%% file: main.tex
\documentclass[12pt,a4paper]{report}
\usepackage[T1]{fontenc}
\usepackage[french]{babel}
\usepackage{hyperref}
\usepackage{mathrsfs} 
\usepackage{numprint}
\usepackage{amsthm}
\usepackage{amsmath}      
\usepackage{mathrsfs}       
\usepackage{amssymb}        
\usepackage{amsfonts}      
\usepackage{multicol}
\usepackage{pdfpages}
\usepackage{a4wide}

\usepackage{pythontex}
\usepackage{graphicx} 
\usepackage{wrapfig} 
\usepackage{stmaryrd}
\usepackage{enumitem}   
\usepackage{comment}

\input{commandes_latex}

\usepackage{multirow} 
\usepackage{multicol}
\usepackage{scrextend} 
\usepackage{tikz}
\usepackage{graphicx}
\usepackage{colortbl}
\usepackage{array}
\usepackage[absolute]{textpos} 
\usepackage{color} 
\definecolor{Prune}{RGB}{99,0,60}
\usepackage{mdframed}
\usepackage{fancyhdr}
\usepackage{geometry}
\usepackage{ragged2e}

\begin{document}

\input{Chapitres/Page_garde_nouvelle} 

\normalsize
\thispagestyle{empty}
\newgeometry{top=3.25 cm, bottom=3.25cm, left=2.6cm, right=2.6cm}

\Ifthispageodd{}{\newpage \thispagestyle{empty} \null\newpage}

\input{Chapitres/Remerciements}

\newpage 
\Ifthispageodd{}{\newpage \null\newpage}
\addcontentsline{toc}{chapter}{Table des matières}

\tableofcontents

\newpage
\Ifthispageodd{}{\newpage \null\newpage}
\input{Chapitres/0_Introduction}

\newpage
\Ifthispageodd{}{\newpage \null\newpage}
\input{Chapitres/1_Outils}
\newpage
\Ifthispageodd{}{\newpage \null\newpage}
\input{Chapitres/2_Outils_bis}

\newpage
\Ifthispageodd{}{\newpage \null\newpage}
\input{Chapitres/3_Exposants_diophantiens_somme}
\newpage
\Ifthispageodd{}{\newpage \null\newpage}
\input{Chapitres/4_Approximation_droites}
\newpage
\Ifthispageodd{}{\newpage \null\newpage}
\input{Chapitres/5_Approximation_1er_angle}
\newpage
\Ifthispageodd{}{\newpage \null\newpage}
\input{Chapitres/6_Approx}
\newpage
\Ifthispageodd{}{\newpage \null\newpage}
\input{Chapitres/8_Independance_exposants}

\newpage
\Ifthispageodd{}{\newpage \null\newpage}
\input{Chapitres/7_Dernier_angle_autre_exemple}

\newpage
\Ifthispageodd{}{\newpage \null\newpage}

\phantomsection

\addcontentsline{toc}{chapter}{Bibliographie}
\bibliographystyle{smfplain.bst}
\bibliography{bibliographie.bib}

\end{document}

%% file: commandes_latex.tex
\newcommand{\R}[0]{\mathbb{R}}
\newcommand{\RQ}[0]{\mathbb{R} \setminus \mathbb{Q}}
\newcommand{\N}[0]{\mathbb{N}}
\newcommand{\Z}[0]{\mathbb{Z}}
\newcommand{\Q}[0]{\mathbb{Q}}
\renewcommand{\P}[0]{\mathbb{P}}
\newcommand{\Nx}[0]{\mathbb{N^*}}
\newcommand{\Rx}[0]{\mathbb{R^*}}

\newcommand{\RR}[0]{\mathcal{R}}

\newcommand{\Mat}[0]{\text{Mat}}
\newcommand{\Ind}[0]{\text{Ind}}
\newcommand{\MM}[0]{\mathcal{M}}
\newcommand{\FF}[0]{\mathcal{F}}
\renewcommand{\AA}[0]{\mathcal{A}}

\newcommand{\CC}[0]{\mathcal{C}}
\newcommand{\DD}[0]{\mathcal{D}}
\newcommand{\EE}[0]{\mathcal{E}}
\newcommand{\HH}[0]{\mathcal{H}}
\newcommand{\GG}[0]{\mathcal{G}}
\newcommand{\PP}[0]{\mathcal{P}}

\newcommand{\GL}[0]{\text{GL}}
\newcommand{\udots}[0]{\reflectbox{$\ddots$}}
\renewcommand{\mod}[0]{\text{ mod }}
\newcommand{\partfrac}[1]{\left\{ #1 \right\} }

\renewcommand{\SS}[0]{\mathcal{S}}
\newcommand{\II}[0]{\mathcal{I}}

\newcommand{\floor}[1]{\left\lfloor #1 \right\rfloor}

\let\oldforall\forall
\renewcommand{\forall}{\oldforall \, }

\newcommand{\Ztir}[0]{\mathbb{Z}\text{-}}
\newcommand{\Zbase}[0]{\mathbb{Z}\text{-base}}

\newcommand{\tir}[0]{\text{-}}

\newcommand{\pgcd}[0]{\text{pgcd}}
\newcommand{\vol}[0]{\text{vol}}
\newcommand{\covol}[0]{\text{covol}}
\newcommand{\Vect}[0]{\text{Vect}}
\newcommand{\rang}[0]{\text{rang}}

\newenvironment{preuve}[1][]{\paragraph{Preuve#1.}}{\hfill$\blacksquare$} 
\newenvironment{theor}[1][]{\paragraph{Théorème #1.}}{\hfill}
\newenvironment{propr}[1][]{\paragraph{Propriété #1.}}{\hfill}
\newenvironment{lemr}[1][]{\paragraph{Lemme #1.}}{\hfill}

\theoremstyle{definition}
\newtheorem{theo}{Théorème}[chapter]
\newtheorem{req}[theo]{Remarque}
\newtheorem{prob}[theo]{Problème}
\newtheorem{cor}[theo]{Corollaire}
\newtheorem{prop}[theo]{Propriété}
\newtheorem{lem}[theo]{Lemme}
\newtheorem{defi}[theo]{Définition}
\newtheorem{ex}[theo]{Exemple}
\newtheorem{conjecture}[theo]{Conjecure}

\newcounter{constante}[chapter]
\newcommand{\cons}[0]{\refstepcounter{constante}c_{\theconstante}}

\usepackage{dsfont}

%% file: Chapitres/Page_garde_nouvelle.tex
\begin{titlepage}

\newgeometry{left=6cm,bottom=2cm, top=1cm, right=1cm}

\tikz[remember picture,overlay] \node[opacity=1,inner sep=0pt] at (-7mm,-135mm){\includegraphics{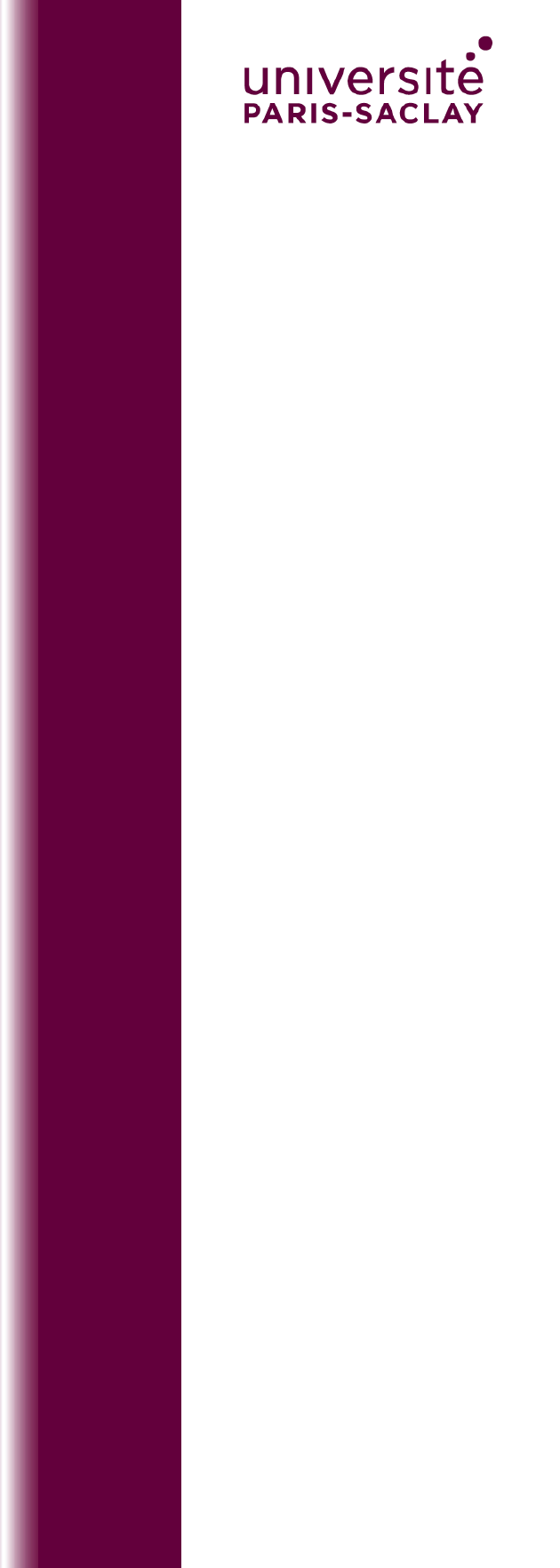}};

\fontfamily{fvs}\fontseries{m}\selectfont


\color{white}

\begin{picture}(0,0)
\put(-142,-743){\rotatebox{90}{\Large \textsc{THESE DE DOCTORAT}}} \\
\put(-110,-743){\rotatebox{90}{NNT : 2024UPASM009 }}
\end{picture}
 


\flushright
\vspace{10mm} 
\color{Prune}

\fontsize{22}{26}\selectfont
  \Huge Approximation de sous-espaces vectoriels de $\R^n$ par des sous-espaces rationnels
\\

\normalsize
\color{black}
\Large{\textit{Approximation of linear subspaces by \\rational linear subspaces}} \\

\fontsize{8}{12}\selectfont

\vspace{1.5cm}

\normalsize
\textbf{Thèse de doctorat de l'université Paris-Saclay}\\

\vspace{6mm}

\small École doctorale n$^{\circ}$ 574 : mathématiques Hadamard (EDMH) \\
\small Spécialité de doctorat: Mathématiques Fondamentales \\
\small Graduate School : Mathématiques\\ Référent : Faculté des sciences d’Orsay \\
\vspace{6mm}

\footnotesize Thèse préparée dans le Laboratoire de math\'ematiques d'Orsay (Université Paris-Saclay, CNRS), sous la direction de \textbf{Stéphane FISCHLER},\\ Ma\^itre de conférences \\
\vspace{15mm}

\textbf{Thèse soutenue à Paris-Saclay, le 05 Juin 2024, par}\\
\bigskip
\Large {\color{Prune} \textbf{Gaétan GUILLOT}} 

\vspace{\fill} 

\bigskip

\flushleft
\small {\color{Prune} \textbf{Composition du jury}}\\
{\color{Prune} \scriptsize {Membres du jury avec voix délibérative}} \\
\vspace{2mm}
\scriptsize
\begin{tabular}{|p{7cm}l}
\arrayrulecolor{Prune}
\textbf{Frédéric PAULIN} &  Président \\ 
Professeur, Université Paris-Saclay    &   \\ 
\textbf{Yann BUGEAUD} &  Rapporteur \& Examinateur \\ 
Professeur, Université de Strasbourg   &   \\ 
\textbf{Nikolay MOSHCHEVITIN} &  Rapporteur \& Examinateur \\ 
Projektassistent, Technische Universit\"at Wien  &   \\ 
\textbf{Jean-Benoît BOST} &  Examinateur \\ 
Professeur, Université Paris-Saclay   &   \\ 
\textbf{Nicolas DE SAXCE} &  Examinateur \\ 
Chargé de recherche, Université Paris-Nord   &   \\

\end{tabular} 

\end{titlepage}


\thispagestyle{empty}
\newgeometry{top=1.5cm, bottom=1.25cm, left=2cm, right=2cm}

\noindent 
\includegraphics[height=2.45cm]{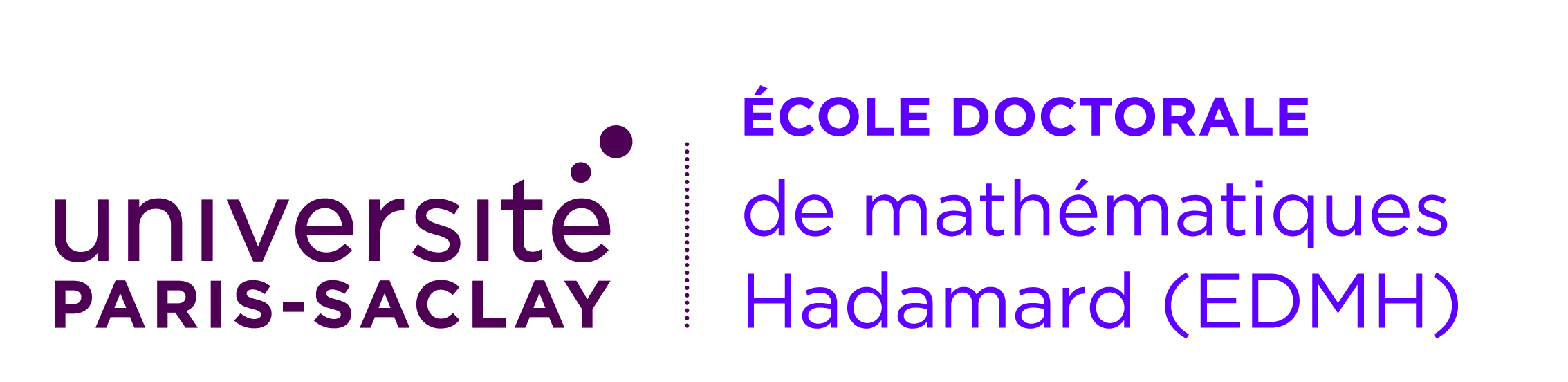}
\vspace{1cm}

\small
\setcounter{page}{0}
\begin{mdframed}[linecolor=Prune,linewidth=1]

\textbf{Titre:} 
Approximation de sous-espaces vectoriels de $\R^n$ par des sous-espaces rationnels

\noindent \textbf{Mots clés:} Arithmétique, Approximation diophantienne, Géométrie des nombres

\vspace{-.5cm}
\begin{multicols}{2}
\noindent \textbf{Résumé:} Pour $A$ un sous-espace vectoriel de $\R^n$ de dimension $d$ et   $B$ un sous-espace rationnel de dimension $e$, on définit $\min(d,e)$ angles $\psi_j(A,B)$ pour $j \in \llbracket 1, \min(d,e) \rrbracket $ qui rendent compte de la proximité entre $A$ et $B$. On étudie alors l'exposant diophantien $\mu_n(A|e)_j$ défini comme la borne supérieure des $\mu > 0$ tels qu'il existe une infinité d'espaces $B$ de dimension $e$ pour lesquels $\psi_j(A,B)$ est inférieur à $H(B)^{-\mu}$, où $H(B)$ est la hauteur de l'espace rationnel $B$. 

On montre d'abord une formule permettant de calculer, sous certaines hypothèses, $\mu_n(A|e)_j$ dans le cas où $A = \bigoplus\limits_{\ell=1}^d A_\ell$ avec $A_\ell$ des droites  de $\R^n$.
Ensuite on construit plusieurs sous-espaces vectoriels de $\R^n$ de dimension $d$ dont on peut prescrire les valeurs prises par $\mu_n(A|e)_j$ pour divers choix de $(e,j)$. 
De ces constructions, on déduit enfin des résultats sur l'indépendance algébrique de familles de fonctions de la forme $\mu_n(\cdot|e)_j$.

\end{multicols}

\end{mdframed}

\vspace{8mm} 

\begin{mdframed}[linecolor=Prune,linewidth=1]

\textbf{Title:} Approximation of linear subspaces by rational linear subspaces


\vspace{-.5cm}
\begin{multicols}{2}
\noindent \textbf{Abstract:} 
For $A$ a vector subspace of $\mathbb{R}^n$ with dimension $d$ and $B$ a rational subspace with dimension $e$, we define $\min(d,e)$ angles $\psi_j(A,B)$ for $j \in \llbracket 1, \min(d,e) \rrbracket$ that capture the proximity between $A$ and $B$. We then study the Diophantine exponent $\mu_n(A|e)_j$, defined as the supremum of $\mu > 0$ such that there exist infinitely many spaces $B$ of dimension $e$ for which $\psi_j(A,B)$ is less than $H(B)^{-\mu}$, where $H(B)$ is the height of the rational space $B$.

We first present a formula allowing the computation, under certain assumptions, of $\mu_n(A|e)j$ when $A = \bigoplus\limits_{\ell=1}^d A_\ell$ with $A_\ell$  a line in $\mathbb{R}^n$.
Then, we construct several vector subspaces of $\mathbb{R}^n$ with dimension $d$ for which we can prescribe the values taken by $\mu_n(A|e)_j$ for various choices of $(e,j)$.
From these constructions, we finally obtain results on the algebraic independence of families of functions of the form $\mu_n(\cdot|e)_j$.
\end{multicols}

\end{mdframed}

\newpage

\hbox{\includegraphics[width=8cm]{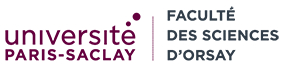}}\vspace*{2cm}

\hbox{\includegraphics[width=5cm]{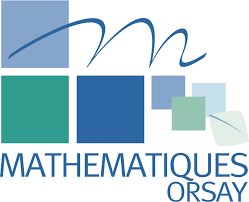}}\vspace*{2cm}

\hbox{\includegraphics[width=4cm]{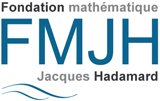}}\vspace*{0.5cm}

%% file: Chapitres/Remerciements.tex
\chapter*{Remerciements}

Tout d'abord, je veux adresser mes plus sincères remerciements à mon directeur de thèse Stéphane Fischler.  Je n'aurais pas pu être mieux encadré durant ces trois années et j'en mesure la chance. Sa disponibilité, ses conseils  et son accompagnement ont été précieux. Merci pour le temps passé à toutes les relectures de ce manuscrit. 

\bigskip

Merci aux rapporteurs de cette thèse, Yann Bugeaud et Nikolay Moshchevitin, pour leur lecture et l'intérêt porté à mon travail. Merci pour les retours et remarques constructives. Je remercie Jean-Benoît Bost, Frédéric Paulin et Nicolas de Saxcé d'avoir accepté de faire partie du jury, je suis honoré de leur présence. 

\bigskip 

J'adresse mes remerciements à tous les enseignants, chercheurs et personnels que j'ai croisés et qui m'ont aidé durant ces années. Un merci particulier à mes "grands frères" de thèse, Elio Joseph et Anthony Poëls, qui ont toujours répondu à mes questions du M2 jusqu'à cette soutenance. 

\bigskip

J'ai eu la chance d'avoir pu compter sur mes amis pendant cette thèse. Dresser une liste est délicat, mais je veux remercier ici certaines personnes dans un ordre plus ou moins chronologique. Tout d'abord, les amis de toujours : Stéphen et Léo, vous accompagner à la batterie depuis tout ce temps est un plaisir sans cesse renouvelé, bientôt les sitcoms raconteront nos vies. Avec Clément G., Alexis, Thomas le groupe formé depuis la maternelle s'est agrandi au fil des ans, entre autres par  Clément M., Lisa, Arthur, merci à vous pour nos  régulières retrouvailles. Je complète cette liste presque exclusivement mayennaise avec les amies du lycée Amélie, Lucie et Jade ; toujours présentes. Mon aventure mathématique a  commencé en prépa (ou au PMU) aux cotés de  Laurine,  Romain, Alexis G., Alexis B. ; des Bretons au grand coeur qui m'ont accueilli (presque) comme un des leurs. Bien entendu, je dois des remerciements aux colocs d'Orsay, Charles et Anthony ; devenus compagnons lointains de thèse, nos soirées passées sur un certain canapé évoqueront toujours en moi une douce nostalgie. Merci Anthony pour les corrections d'orthographe. Depuis la licence, j'ai  eu la chance de partager la compagnie des (plus ou moins) matheux : Philippe, Jordan, Keivan, Marion, Jeanne, Donatien, Pauline, Nicolas, Mathilde. Merci enfin à Thibault d'avoir réussi à me faire sortir du bureau du labo.
Présents depuis l'école maternelle ou depuis le doctorat, votre compagnie est précieuse. Merci à vous pour tous les bons moments passés. Je veux remercier en particulier Philippe pour sa grande influence sur mon parcours (au moins mathématique), pas merci pour les tests de mémoire pourris, merci pour les goûters. 

\newpage
Merci enfin à ma famille, tout particulièrement mon neveu préféré, mes frères et soeurs (j'espère qu'Elise a su réfréner ses envies de déguisements). Merci enfin à mes parents, sans qui je n'aurais sans doute pas développé le goût des sciences, pour leur présence et leur soutien durant ces longues années d'études.

%% file: Chapitres/0_Introduction.tex
\chapter{Introduction}
Cette thèse a pour thème l'approximation diophantienne. Après avoir rappelé le cas classique on introduit le problème étudié ici, l'approximation rationnelle de sous-espaces vectoriels de $\R^n$. On énonce ensuite les résultats connus pour ce problème. Enfin, on présente les principaux résultats obtenus et le plan de la thèse.

\section{Approximation diophantienne classique}\label{section_classique}
 L'approximation diophantienne est l'étude de l'approche des nombres réels par des nombres rationnels. 

Pour un réel $\xi$ donné, on peut rendre la quantité $\left|\xi - \frac{p}{q}\right|$, avec $p \in \Z$ et $q \in \Nx$, aussi petite que l'on veut par densité de $\Q$ dans $\R$. Cependant cela peut se faire au prix d'un dénominateur $q$ possiblement très grand. On cherche donc à étudier qualitativement cette approximation en cherchant des nombres rationnels $\frac{p}{q}$ approchant bien $\xi$, les moins \og compliqués\fg{} possibles.

Ici, la complexité d'un rationnel est la taille de son dénominateur. On introduit alors l'exposant d'irrationalité d'un réel $\xi$ qui est la borne supérieure de l'ensemble des réels strictement positifs $\mu$ tels qu'il existe une infinité de rationnels $\frac{p}{q}$ vérifiant :
\begin{align}
 \left|\xi - \frac{p}{q}\right| \leq \frac{1}{q^\mu}.
\end{align}
On note cette borne $\mu(\xi) \in [1,+\infty]$. On cite quelques résultats sur cet exposant d’irrationalité, on peut les retrouver dans \cite{Hindry_Silverman}, partie D ou dans \cite{Hardy_Wright}, Chapitre XI.
\\ Soit $\xi$ un nombre rationnel de la forme $\frac{a}{b}$ et $\frac{p}{q} \neq \xi$ on remarque que :
\begin{align*}
 \left|\xi - \frac{p}{q}\right| = \left| \frac{aq-bp}{bq}\right| \geq \frac{1}{bq}
\end{align*}
et on en déduit donc que $\mu(\xi) = 1 $ pour tout $\xi \in \Q$.

 Pour tout $\xi \in \RQ$, on a $\mu(\xi) \geq 2 $. 
En effet,le principe des tiroirs de Dirichlet implique que pour tout $Q > 1$ il existe un couple $(p,q) \in \Z \times \Nx$ tel que :
\begin{align*}
 1 \leq q \leq Q \text{ et } \left|q\xi - p\right| \leq \frac{1}{Q}.
\end{align*}

\bigskip

Un théorème plus profond de Roth (\cite{Roth}, \cite{Hindry_Silverman}) implique que pour tout nombre irrationnel algébrique $\xi$, on a $\mu(\xi) = 2$. La réciproque est fausse, on a en effet $\mu(e) = 2$, voir par exemple \cite{Borwein_Borwein} chapitre 11.3.

La théorie des fractions continues, voir \cite{Hardy_Wright} ou \cite{Casselsdio}, est importante en approximation diophantienne classique. Elle permet en effet de fournir les meilleures approximations de nombres réels. Notamment cette théorie montre :
\begin{align}\label{spectre_classique}
 \Big \lbrace \mu(\xi), \xi \in \R \Big\rbrace = \{1\} \cup [2, + \infty]
\end{align}
On peut aussi construire explicitement des réels $\xi$ pour lesquels $\mu(\xi)$ est égal à un réel $\mu \geq 2$ fixé (voir la remarque~\ref{3req_sigma_transcendant} au chapitre \ref{chap3} page \pageref{3req_sigma_transcendant}). Cette thèse s'intéresse à des généralisations de cette approche.
\begin{req}
 Les réels $\xi$ vérifiant $\mu(\xi) = + \infty$ sont appelés nombres de Liouville. Historiquement, ces nombres furent les premiers exemples permettant de montrer l'existence de nombres transcendants \cite{Liouville1844}, \cite{Liouville1851}.
\end{req}


\section{Approximation rationnelle de sous-espaces vectoriels}
En 1967, Schmidt introduit une généralisation de l'approximation diophantienne classique \cite{Schmidt}. Le but est d'approcher des sous-espaces vectoriels de $\R^n$ par des sous-espaces dits \og rationnels\fg{}. 

Soit $n \geq 2$ un entier. Soient $d,e \in \llbracket 1,n-1 \rrbracket$ deux entiers. On cherche ici à approcher $A$, un sous-espace vectoriel de $\R^n$ de dimension $d$ par des sous-espaces vectoriels de $\R^n$ de dimension $e$, rationnels et les moins \og compliqués\fg{} possibles.

Définissons ces notions. On utilise ici des notations provenant de \cite{Schmidt} et \cite{Joseph-these}. Dans la suite, et sauf mention du contraire, toutes les normes $\| \cdot \| $ sont les normes euclidiennes canoniques sur les espaces considérés.

\subsection{Espaces rationnels et hauteur}
\begin{defi}
 Soit $B$ un sous-espace vectoriel de $\R^n$. On dit que $B$ est rationnel si $B$ admet une base de vecteurs à coordonnées dans $\Q$. Une telle base est appelée base rationnelle de $B$. 
 \\On note $\RR_n(e)$ l'ensemble des sous-espaces rationnels de $\R^n$ de dimension $e$.
\end{defi}

\begin{req}
 On peut voir la rationnalité d'un espace par l'équivalence suivante $B \in \RR_n(e)$ si et seulement si $B $  est l'ensemble des solutions de  $n-e$ équations indépendantes à coefficients dans $\Q.$
On renvoie aussi à la notion de $\Q\tir$structure définie dans \cite{Bourbaki} (paragraphe $8$). La $\Q\tir$structure $\Q^n$ dans $\R^n$ permet alors de définir les sous-espaces rationnels de $\R^n$.
\end{req}

\bigskip
On définit maintenant la complexité d'un tel espace, utilisée par Schmidt dans \cite{Schmidt}. A cet effet on rappelle la construction des coordonnées grassmaniennes associées à un espace. On peut retrouver cette construction et les preuves des résultats dans \cite{caldero_germoni}.

\bigskip
 Pour $n$ et $r$ deux entiers, on note $\PP(r,n)$ l'ensemble des parties à $r$ éléments de $\llbracket 1, n \rrbracket$. Pour $n, m, r$ trois entiers, $I \in \PP(r,n) $, $J \in \PP(r,m)$, et $M \in \MM_{n,m}(\R)$ on note $\Delta_{I,J}(M)$ le mineur de $M$ associé aux lignes indexées par $I$ et aux colonnes indexées par $J$. Enfin on pose la matrice :
 \begin{align*}
 \Delta^r(M) = (\Delta_{I,J}(M))_{I \in \PP(r,n), J \in \PP(r,m)}.
 \end{align*}
 On remarque que si $m = r \leq n $ alors $\Delta^r(M) \in \R^{N}$ avec $N= \binom{n}{r}$ en indexant $\PP(r,n)$ par $\llbracket 1, N \rrbracket$ avec l'ordre lexicographique.

 Ces quantités permettent d'énoncer le théorème suivant qui généralise le calcul d'un déterminant par développement par rapport à une colonne ou une ligne. 

 \begin{theo}[Développement de Laplace] \label{theo_laplace}
 Soit $M \in \MM_n(\R)$, $r \in \llbracket 1, n-1 \rrbracket$ et $J \in \PP(r,n)$. Alors :
 \begin{align*}
 \det(M ) =  \sum\limits_{I \in \PP(r,n)} (-1)^{\ell(I) + \ell(J)} \Delta_{I,J}(M) \Delta_{\Bar{I},\Bar{J}}(M)
 \end{align*}
 où $ \Bar{I} = \llbracket 1, n \rrbracket \smallsetminus I$ et $\ell(I) = \# \lbrace (i,j) \in I \times \Bar{I}, i > j \rbrace$.
 \end{theo}

 On peut maintenant énoncer le théorème principal qui permet la définition des coordonnées grassmaniennes.
\begin{theo}
 Soit $B$ un sous-espace vectoriel de $\R^n$ de dimension $e$. Soit $X_1, \ldots, X_e$ une base de $B$. On note $M_B \in \MM_{n,e}(\R)$ la matrice dont les colonnes sont les $X_i$. 
 \\Alors, en notant $N = \binom{n}{e}$, la droite engendrée par $\Delta^e(M_B) \in \R^N$ ne dépend que de $B$.
\end{theo}

 On appelle alors coordonnées grassmaniennes (ou coordonnées de Plücker) de $B$ le point de l'espace projectif $\P^N(\R)$ associé à cette droite. Enfin on appelle représentant des coordonnées grassmaniennes associé à la base $X_1, \ldots, X_e$ le vecteur $\Delta^e(M_B)$ de $\R^N$.

\begin{req}
 L'ensemble des points de $\P^N(\R)$ correspondant à des coordonnées grassmaniennes est défini par un système d'équations appelé relations de Plücker.
\end{req}

\begin{defi}\label{def_hauteur}
 Soit $B$ un sous-espace rationnel de $\R^n$ de dimension $e$.\\On note $(\eta_1, \ldots, \eta_N) \in \Q^N$ un représentant des coordonnées grassmaniennes associé à une base rationnelle de $B$ avec $N = \binom{n}{e}$. On pose $\alpha = \eta_1 \Z + \ldots + \eta_N \Z$ l'idéal fractionnaire de $\Z$ engendré par ces coordonnées. 
 \\ On appelle hauteur de $B$ et on la note $H(B)$ la quantité : 
 \begin{align*}
 H(B) = \frac{1}{N(\alpha)} \| (\eta_1, \ldots, \eta_N) \|.
 \end{align*}
 où $N(\alpha)$ est la norme de l'idéal fractionnaire $\alpha$.
\end{defi}

\begin{req}
 Comme $\alpha$ est de la forme $q\Z$ avec $q \in \Q$, on a $N(\alpha) = |q|$. En outre, on peut choisir $(\eta_1, \ldots, \eta_N)$ de sorte que les coordonnées $\eta_i$ soient entières et premières entre elles dans leur ensemble. Dans ce cas $\alpha =\Z$ et donc $H(B) = \| (\eta_1, \ldots, \eta_N) \|.$
\end{req}

\subsection{Proximité entre deux espaces}
Pour pouvoir parler d'approximation d'un espace par un autre, il reste à définir la notion de proximité entre deux sous-espaces vectoriels. Pour cela on s'intéresse d'abord à l'angle entre deux vecteurs. 

\begin{defi}
 Soit $X,Y \in \R^n \smallsetminus \{ 0 \}$. On note : 
 \begin{align*}
 \omega(X,Y) = \frac{\| X \wedge Y \|}{ \|X \| \cdot\|Y \| }
 \end{align*}
où $\wedge$ est le produit extérieur sur $\R^n$. 
\\Géométriquement $\omega(X,Y)$ est le sinus de l'angle géométrique entre les vecteurs $X$ et $Y$.
\end{defi}

Soit $A$ et $B$ deux sous-espaces vectoriels de $\R^n$ de dimensions respectives $d$ et $e$. On pose $t = \min(d,e)$. On définit maintenant la quantité : 
\begin{align}
 \omega_1(A,B) = \min\limits_{\substack{X \in A \smallsetminus \{ 0 \} \\Y \in B \smallsetminus \{ 0 \}}} \omega(X,Y).
\end{align}
On pose $X_1 \in A$ et $Y_1 \in B$ des vecteurs non nuls réalisant ce minimum. On peut maintenant définir les quantités $\omega_2(A,B), \ldots, \omega_t(A,B)$ par récurrence.
\\On suppose que l'on a construit $\omega_1(A,B), \ldots, \omega_j(A,B)$ pour $j \in \llbracket 1, t-1\rrbracket$ ainsi que les couples $(X_1, Y_1), \ldots, (X_j, Y_j)$ qui les réalisent. 
\\Soit $A_j$ et $B_j$ respectivement l'orthogonal de $\Vect(X_1, \ldots, X_j)$ dans $A$ et de $\Vect(Y_1, \ldots, Y_j)$ dans $B$. On pose maintenant :
\begin{align}
 \omega_{j+1}(A,B) = \min\limits_{\substack{X \in A_j \smallsetminus \{ 0 \} \\Y \in B_j \smallsetminus \{ 0 \}}} \omega(X,Y).
\end{align}
et $(X_{j+1}, Y_{j+1})\in A \times B $ réalisant ce minimum. \\
Une fois la récurrence terminée, on a donc construit $t$ quantités $\omega_i(A,B)$ et des familles orthogonales $(X_1, \ldots, X_t)$ et $(Y_1, \ldots, Y_t)$ vérifiant :
\begin{align*}
 \omega_1(A,B) \leq \ldots \leq \omega_{t}(A,B), \\
 \forall i \in \llbracket 1,t \rrbracket, \quad \omega_i(A,B) = \omega(X_i, Y_i).
\end{align*}

Remarquons maintenant que si $d +e > n $ alors $A \cap B \neq \{ 0 \}$ et donc en particulier $\omega_1(A,B) = \ldots = \omega_{d+e-n}(A,B) = 0$. 
Pour cette raison Schmidt \cite{Schmidt} et Joseph \cite{Joseph-these} se placent dans le cadre $d +e \leq n$. 

\bigskip
 Dans cette thèse on se place dans le cas général avec $d, e \in \llbracket 1,n-1 \rrbracket $ quelconques. 
Les premiers angles nuls ne sont alors pas intéressants à étudier. On introduit pour cela la notation suivante : pour $e,d, n$ trois entiers, on pose 
:
\begin{align*}
 g(d,e,n) = \max(0, d+e -n).
\end{align*} 
\\Les quantités intéressantes à étudier sont alors les $\psi_j(A,B)$ définis par : 
\begin{align*}
 \forall j \in \llbracket 1, t - g (d,e,n) \rrbracket, \quad \psi_j(A,B) = \omega_{j+g(d,e,n)}(A,B). 
\end{align*}
On dit que $A$ et $B$ sont proches au $j\tir$ème angle si $\psi_j(A,B)$ est petit.

\subsection{Exposants diophantiens}\label{1section_exposants_diopha}
On généralise alors l'approche de l'approximation diophantienne classique dans le cadre des sous-espaces vectoriels.
\\ Il reste néanmoins à définir l'irrationalité d'un sous-espace vectoriel $A$ que l'on veut approcher. 

\begin{defi}
 Soient $d \in \llbracket 1, n-1 \rrbracket$ et $A$ un sous-espace vectoriel de dimension $d$ de $\R^n$. Pour $e \in \llbracket 1, n-1 \rrbracket$ et $j \in \llbracket 1, \min(d,e) - g(d,e,n) \rrbracket $, on dit que $A$ est $(e,j)\tir$irrationnel si pour tout espace rationnel $B$ de dimension $e$ on a :
 \begin{align*}
 \dim(A \cap B) < j + g(d,e,n).
 \end{align*}
On note $\II_n(d,e)_j$ l'ensemble des espaces $(e,j)\tir$irrationnels de dimension $d$.
\end{defi}

\begin{req}
 Pour $A$ un sous-espace vectoriel de dimension $d $, on a équivalence entre les assertions suivantes :
\begin{enumerate}[label = (\roman*)]
 \item $A \in \II_n(d,e)_j$
 \item $ \forall B \in \RR_n(e), \psi_j(A,B) \neq 0$
 \item $ \forall B \in \RR_n(e), \omega_{j+g(d,e,n)}(A,B) \neq 0$.
\end{enumerate}
\end{req}

\begin{req}
La notion de sous-espace $(e,j)\tir$irrationnel a été définie par Joseph dans \cite{Joseph-these} seulement pour le cas $g(d,e,n) = 0 $. Schmidt \cite{Schmidt} évoque simplement la condition $\dim(A \cap B) < j$ dans ses corollaires ; il se place aussi dans le cas $d +e \leq n$ c'est-à-dire $g(d,e,n) = 0$.
\end{req}
On a maintenant tous les outils pour définir les exposants diophantiens d'un sous-espace $A \in \II_n(d,e)_j$, pour $d,e \in \llbracket 1, n-1 \rrbracket$ et $ j \in \llbracket 1, \min(d,e) - g(d,e,n) \rrbracket $. 
\\On note $\mu_n(A|e)_j$ la borne supérieure (éventuellement infinie) de l'ensemble des $ \mu > 0$ tels qu'il existe une infinité de $B \in \RR_n(e)$ vérifiant :
\begin{align}\label{ineg_def_mu}
 \psi_j(A,B) \leq \frac{1}{H(B)^\mu}.
\end{align}

\begin{req}
 On pourrait remplacer dans la définition, l'inégalité (\ref{ineg_def_mu}) par $$\psi_j(A,B) \leq \frac{c}{H(B)^\mu}$$ où $c > 0 $ est une constante indépendante de $B$. 
\end{req}

\begin{req}\label{1req_intro_generalisation}
 Ces exposants sont bien une généralisation de l'exposant d'irrationalité en approximation diophantienne classique. En effet on a : 
 \begin{align}
 \forall \xi \in \RQ, \quad \mu(\xi) = \mu_2\bigg(\Vect \begin{pmatrix}
 1 \\ \xi 
 \end{pmatrix}, 1\bigg)_1. 
 \end{align}
 \end{req}
\newpage
Il émerge alors plusieurs problèmes : 

\begin{prob}
 Déterminer $\mu_n(A|e)_j$ pour $A \in \II_n(d,e)_j$ donné.
\end{prob} 

\begin{prob}\label{prob_spec}
 Déterminer le spectre $\SS(n,d,e,j) = \Big\lbrace \mu_n(A|e)_j, A \in \II_n(d,e)_j \Big\rbrace $ en fonction de $(n,d,e,j)$.
\end{prob}

\begin{prob}\label{prob_murond}
 Déterminer la quantité $\inf\limits_{A \in \II_n(d,e)_j} \mu_n(A|e)_j $.
\end{prob}

\begin{req}
En particulier le problème \ref{prob_spec} de détermination du spectre dans le cas $n = 2, d= e=j=1$ est résolu en (\ref{spectre_classique}); on a $\SS(2,1,1,1)= \left[2, +\infty \right]$.
\end{req}

Le problème $\ref{prob_murond}$ n'est pas l'objet d'étude de cette thèse. Le lecteur pourra consulter les travaux de Schmidt \cite{Schmidt}, Moshchevitin \cite{Moshchevitin}, Joseph \cite{Joseph-these}, \cite{joseph_exposants}, \cite{joseph_spectre} et De Saxcé \cite{Saxce_HDR},\cite{Saxce} à ce sujet.

\bigskip
Enfin, dans l'inégalité (\ref{ineg_def_mu}), la proximité de deux espaces est comparée avec une puissance de la hauteur. Dans l'esprit de la conjecture de Duffin-Schaeffer \cite{Duffin_Schaeffer} (démontrée par Koukoulopoulos-Maynard \cite{Koukoulopoulos_Maynard}), on peut s'intéresser à des inégalités du type :
\begin{align*}
 \psi_j(A,B) \leq \phi(t)
\end{align*}
où $\phi$ est une fonction strictement décroissante quelconque avec $H(B) \leq t$.
\\On peut mentionner les récents travaux de Chebotarenko \cite{Chebotarenko} à ce sujet dans le cas $$n = 4, \quad d= e= 2, \quad j=1.$$


\section{Résultats connus sur le spectre des exposants}

\subsection{Valeurs prises par un exposant fixé}
Laurent trouve le spectre dans le cas où l'on approche une droite dans \cite{Laurent_omega_un}.
\begin{theo}[Laurent 2009]
Soit $n \geq 2$ alors pour tout $e \in \llbracket 1,n-1 \rrbracket$ on a : 
\begin{align*}
 \SS(n,1,e,1) = \left[ \frac{n}{n-e}, + \infty \right].
\end{align*}
\end{theo}

Joseph \cite{Joseph-these} prouve le résultat suivant dans le cas $e = d$ et pour le dernier angle.

\begin{theo}[Joseph 2021]\label{Joseph_spec} Soit $n \geq 2 $ et $d \in \llbracket 1, \lfloor n/2 \rfloor \rrbracket$ alors :
\begin{align*}
 \left[ 1 + \frac{1}{2d} + \sqrt{1+\frac{1}{4d^2}}, + \infty \right] \subset \SS(n,d,d,d).
\end{align*}
 
\end{theo}

\begin{req}
 La preuve du théorème~\ref{Joseph_spec} fournit une construction explicite d'un espace de dimension $d$ d'exposant prescrit. On reprend cette approche dans les chapitres~\ref{chap5},~\ref{chap6},~\ref{chap7}~et~\ref{chap8}.
\end{req}
 
Dans \cite{Saxce}, De Saxcé obtient comme corollaire le résultat suivant qui donne le spectre de l'exposant correspondant au dernier angle pour toutes les dimensions $d$ et $e$, précisant alors ceux de Joseph et Laurent tout en les généralisant. 

\begin{theo}[De Saxcé 2022]
 Soit $n \geq 2$ alors pour tous $e,d \in \llbracket 1, n-1 \rrbracket$ on a :
 \begin{align*}
 \SS(n,d,e,\min(d,e) - g(d,e,n) ) = \left[ \frac{n}{\min(d,e)(n-\max(d,e))}, + \infty \right].
 \end{align*}
\end{theo}

Il est à noter que ce théorème prend en compte le cas $d + e > n$.

\subsection{Spectre joint}
On appelle spectre joint, le spectre d'une famille d'exposants. 
Dans \cite{Laurent_omega_un}, Laurent pose la question de la description de l'ensemble des valeurs des $(n-1)\tir$uplets :
\begin{align*}
 (\mu_n(A | 1)_1, \ldots, \mu_n(A | n-1)_1)
\end{align*}
 quand $A$ décrit $\II_n(1,n-1)_1$.
Dans ce sens, il explicite les inégalités prouvées par Schmidt \cite{Schmidt} (théorèmes $9$ et $10$). 

\begin{theo}[Schmidt 1967, Laurent 2009] \label{Laurent_Schmidt} Soit $n \geq 2$. 
\\Pour tout $A \in \II_n(1,n-1)_1$ on a $\mu_n(A|1)_1 \geq \frac{n}{n-1} $ et :
\begin{align*}
 \forall e \in \llbracket 2, n-1 \rrbracket, \quad \frac{e \mu_n(A|e)_1}{\mu_n(A|e)_1 + e -1 } \leq \mu_n(A|e-1)_1 \leq \frac{(n-e)\mu_n(A|e)_1}{n-e+1}.
\end{align*}
 
\end{theo}

Dans \cite{Roy}, Roy donne le spectre joint complet dans le cas où l'on approche des droites, en montrant que les inégalités du théorème~\ref{Laurent_Schmidt} sont optimales.

\begin{theo}[Roy 2016] \label{1theo_roy_spec_joint}Soit $n \geq 2$. Soit $\mu_1, \ldots, \mu_{n-1} \in [1, + \infty ]$ satisfaisant $\mu_1~\geq~\frac{n}{n-1} $ et :
\begin{align*}
 \forall e \in \llbracket 2, n-1 \rrbracket, \quad \frac{e \mu_e}{\mu_e + e -1 } \leq \mu_{e-1} \leq \frac{(n-e)\mu_e}{n-e+1}.
\end{align*}
 Alors il existe $A \in \II_n(1,n-1)_1$ tel que :
 \begin{align*}
 \forall e \in \llbracket 1, n-1 \rrbracket, \quad \mu_n(A| e)_1 = \mu_e. 
 \end{align*}
\end{theo}

Comme l'ensemble $(\mu_1, \ldots, \mu_{n-1})$ vérifiant les conditions du theorème \ref{1theo_roy_spec_joint} contient un ouvert non vide, on en déduit le corollaire suivant. 

\begin{cor}\label{1cor_ind_alg_exposant_dim1}
 Soit $n \geq 2$. 
 \\ L'image de $\left(\mu_n(\cdot|1)_1, \ldots , \mu_n(\cdot|n-1)_1 \right) : \II_{n}(1,n-1)_1 \to \R^{n-1}$ contient un ouvert non vide de $\R^{n-1}$ et la famille de fonctions $(\mu_n(\cdot|e)_{1} )_{e \in \llbracket 1, n-1 \rrbracket}$ est algébriquement indépendante sur $\R$.
\end{cor}

\subsection{Exposants uniformes}

Cette thèse est consacrée à l'étude des exposants diophantiens définis dans la section~\ref{1section_exposants_diopha}. Il existe aussi d'autres exposants reliés à l'inégalité (\ref{ineg_def_mu}) : les exposants uniformes. Ceux-ci, ainsi que leurs relations avec les exposants $\mu_n(\cdot|e)_j$, ont notamment été étudiés par Khintchine \cite{Khinchine}, Jarník \cite{Jarnik}, Bugeaud et Laurent \cite{Bugeaud_Laurent_2005}, \cite{Bugeaud_Laurent_2010} et Marnat \cite{Marnat_Jarnik}, \cite{Marnat_spec}.
\bigskip
\\ On définit ici les exposants diophantiens uniformes d'une droite ($d = 1$). 
\\Pour $A$ une droite de $\R^n$ et $e \in \llbracket 1 ,n-1 \rrbracket$ tels que $A \in \II_n(1,e)_1$, on note $\widehat{\mu}_n(A|e)_1$ la borne supérieure (éventuellement infinie) de l'ensemble des $ \mu > 0$ tels que pour tout $H > 1$ assez grand, il existe $B \in \RR_n(e)$ vérifiant :
\begin{align*}
 \psi_1(A,B) \leq \frac{1}{H^{\mu}} \text{ et } H(B) \leq H.
\end{align*}

Dans le cas $n = 3 $, Laurent \cite{Laurent_dim_2} décrit le spectre des valeurs prises par $$(\mu_3(A|1)_1, \mu_3(A|2)_1, \widehat{\mu}_3(A|1)_1, \widehat{\mu}_3(A|2)_1 )$$ quand $A$ décrit $\II_{3}(1,2)_1$. Une des égalités pour décrire ce spectre est celle de Jarník \cite{Jarnik} :
\begin{align*}
 \widehat{\mu}_3(A|1)_1 + \widehat{\mu}_3(A|2)_1^{-1} = 1. 
\end{align*}
 En particulier, si $n = 3$ les exposants $\widehat{\mu}_n(\cdot|1)_1$ et $ \widehat{\mu}_n(\cdot|2)_1$ vérifient une relation de dépendance algébrique. Cette relation de dépendance n'existe plus pour les dimensions supérieures. Marnat \cite{Marnat_spec} prouve en effet le théorème suivant. 

 \begin{theo}[Marnat, 2018]
 Soit $n \geq 4$. 
 \\Les fonctions définies sur $\II_n(1,n-1)_1$ $$\mu_n(\cdot|1)_1, \ldots, \mu_n(\cdot|n-1)_1, \widehat{\mu}_n(\cdot|1)_1, \ldots, \widehat{\mu}_n(\cdot|n-1)_1 $$ forment une famille algébriquement indépendante sur $\R$.
 \end{theo}

Dans cette thèse, on établit des résultats similaires sur les familles d'exposants $\mu_n(\cdot|e)_j$ pour des espaces $A$ de dimension $d \geq 1$.

 \section{Principaux résultats de cette thèse}
 Dans ce paragraphe, on énonce les principaux résultats démontrés dans cette thèse. Ici, la numérotation des théorèmes, lemmes et propriétés correspond à celle qui est employée dans les chapitres $2$ à $9$.
 \subsection{Résultats sur l'image du spectre joint}\label{1section_inde_algebrique}
Soit $n \in \Nx$ et $d \in \llbracket 1, n-1 \rrbracket$. On étudie dans cette thèse des familles de fonctions de la forme $\mu_n(\cdot|e)_j$ définie sur $\II_{n}(d,e)_{j-g(d,e,n)}$ avec $ e \in \llbracket 1, n-1 \rrbracket$ et $j \in \llbracket 1+ g(d,e,n), \min(d,e) \rrbracket$. On émet la conjecture suivante qui indiquerait alors qu'il n'y a pas de relation algébrique entre ces fonctions.

\begin{conjecture} La famille de fonctions 
 $\left(\mu_n(\cdot|e)_{j}\right)_{ e \in \llbracket 1, n-1 \rrbracket, j \in \llbracket 1+ g(d,e,n), \min(d,e) \rrbracket }$ est algébriquement indépendante sur $\R$.
 \end{conjecture}
\bigskip
Dans cette direction, on exhibe des sous-familles de $\left(\mu_n(\cdot|e)_{j}\right)_{ e \in \llbracket 1, n-1 \rrbracket, j \in \llbracket 1+ g(d,e,n), \min(d,e) \rrbracket }$ qui sont algébriquement indépendante sur $\R$.
Dans le cadre du premier angle $(j = 1)$, on montre notamment le théorème suivant. Le  théorème~\ref{9prop_premier_angle} étend le résultat du corollaire~\ref{1cor_ind_alg_exposant_dim1}. On fournit ainsi une nouvelle preuve du corollaire~\ref{1cor_ind_alg_exposant_dim1} qui correspond alors au cas particulier $d = 1$ du théorème~\ref{9prop_premier_angle}. 

\begin{theor}[\ref{9prop_premier_angle}]
L'image de $(\mu_n(\cdot|1)_1, \ldots, \mu_n(\cdot|n-d)_1)$ contient un ouvert non vide de $\R^{n-d}$ et la famille de fonctions $(\mu_n(\cdot|e)_{1} )_{e \in \llbracket 1, n-d \rrbracket}$ est algébriquement indépendante sur $\R$.
\end{theor}

\bigskip
Un autre angle naturel à étudier est le dernier angle $(j = \min(d,e) )$ ; en effet $\omega_{\min(d,e)}$ \og ressemble \fg{} à une distance car 
\begin{align*}
 \omega_{\min(d,e)}(A,B) = 0 \Longleftrightarrow A \subset B \text{ ou } B \subset A. 
\end{align*}
En particulier, si $d =e $ alors $\omega_d$ est une distance (voir les lemmes \ref{2lem_ineg_triang} et \ref{lem_inclusion_croissance} pour l'inégalité triangulaire) sur la grassmanienne des sous-espaces vectoriels de dimension $d $ de $\R^n$. 
\\Pour l'étude dans le cadre du dernier angle, on montre les deux résultats suivants, dans lesquels on approxime $A$ par des sous-espaces de dimensions $e \in \llbracket 1, n-d \rrbracket$ (théorème~\ref{9lem_sousfam_pour_theo2}) ou bien $e \in \llbracket d, n-1 \rrbracket $ (théorème~\ref{9lem_sousfam_pour_theo3}). Le  théorème~\ref{9lem_sousfam_pour_theo2} est assez naturel, puisque l'on étudie les valeurs de $e$ pour lesquelles on a $e + d \leq n$.

\begin{theor}[\ref{9lem_sousfam_pour_theo2}]
On suppose que $d$ divise $n$. 
\\Alors l'image de $(\mu_n(\cdot|1)_{\min(d,1)}, \ldots, \mu_n(\cdot|n-d)_{\min(d,n-d)})$ contient un ouvert non vide de $\R^{n-d}$ et la famille $(\mu_n(\cdot|e)_{\min(d,e)} )_{e \in \llbracket 1, n-d \rrbracket}$ est algébriquement indépendante sur $\R$. 
 
\end{theor}

\begin{theor}[\ref{9lem_sousfam_pour_theo3}]
On suppose que $d$ divise $n$. 
\\Alors l'image de $(\mu_n(\cdot|d)_{d}, \ldots, \mu_n(\cdot|n-1)_{d})$ contient un ouvert non vide de $\R^{n-d}$ et la famille $(\mu_n(\cdot|e)_{d} )_{e \in \llbracket d, n-1 \rrbracket}$ est algébriquement indépendante sur $\R$. 
\end{theor}

\bigskip
Ces trois théorèmes se démontrent grâce aux constructions d'espaces avec exposants prescrits qui sont réalisées dans cette thèse et dont les principaux résultats sont énoncés dans la section~\ref{1intro_section_expo_prescrits}.

\subsection{Nouveaux outils}
On développe dans cette thèse de nouveaux outils qui vont permettre de réaliser les différentes constructions d'espaces. 

\bigskip
On établit d'abord qu'un espace $A$ est $(e, j)\tir$irrationnel si et seulement si son orthogonal est $(n - e, j)\tir$irrationnel. Il en découle que si un espace de dimension $d$ est $(n-d, j)\tir$irrationnel alors il est $(e, j)\tir$irrationnel pour tout $e \in \llbracket j, n-j \rrbracket$.

On montre ensuite la propriété suivante, reliant les exposants diophantiens d'un espace $A$ et ceux de son orthogonal $A^\perp$, ce qui permet en particulier d'avoir un premier résultat dans le cas où $d +e > n $. 
\begin{propr}[\ref{3prop_egalit_mun_orhto_mun_normal}]
 Soit $d, e \in \llbracket 1, n-1 \rrbracket$ et $j \in \llbracket 1, \min(d,e)- g(d,e,n) \rrbracket $. Soit $A\in \II_{n}(d,e)_j $ alors
\begin{align*}
 \mu_n(A|e)_j = \mu_n(A^\perp| n -e )_j
\end{align*}
 où $A^\perp$ désigne l'orthogonal de $A$ dans $\R^n$.
\end{propr}

\bigskip

Le lemme suivant permet d'exhiber des nombres transcendants de manière explicite ; on obtient une famille qui est de plus algébriquement indépendante sur $\Q$. Celle-ci sera utilisée comme coordonnées des vecteurs de base des espaces considérés dans cette thèse. Cela permet de faire apparaître des propriétés d'irrationalité de ces espaces.

\begin{lemr}[\ref{3lem_sigma_alg_indep}]
 Soit $e \in \Nx$, $\theta \in \N \smallsetminus \{0,1\}$ et $\alpha = (\alpha_k)_{k \in \N} $ une suite à valeurs dans $\R^{*}_+$ vérifiant $\exists c_{\ref{constante_suite_alpha}} > 1, \quad \forall k \in \N, \quad \alpha_{k+1} > c_{\ref{constante_suite_alpha}}\alpha_{k}.$ Soient $J \subset \Nx $ de cardinal au moins $2$  et $\FF \subset \R$ un ensemble fini de réels.
 \\ Soit $\phi : \N \longrightarrow \llbracket 0,e \rrbracket $ une application telle que :
 \begin{align*}
 \forall i \in \llbracket 0,e \rrbracket, \quad \#\phi^{-1}(\{ i \} ) = \infty.
 \end{align*}
Alors il existe un choix de $e +1 $ suites $(u_k^{i})_{ i \in \llbracket 0,e \rrbracket, k \in \N}$ telles que : 
\begin{align*}
 \text{pour tous } i \in \llbracket 0,e \rrbracket \text{ et } k \in \N, \quad u_k^{i} &\left\{ \begin{array}{cl}
 \in J &\text{ si } \phi(k) = i \\
 = 0 &\text{ sinon,  }
 \end{array} \right. 
\end{align*}
et telles que la famille $\left(\sum\limits_{k =0}^{ + \infty} \frac{u_k^i}{\theta^{\floor{\alpha_k}}}\right)_{ i \in \llbracket 0,e \rrbracket}$ soit algébriquement indépendante sur $\Q(\FF)$. 
\end{lemr}

\bigskip
\bigskip

Ensuite on prouve le lemme suivant qui permet de montrer qu'un certain type d'espaces est $(e,j)\tir$irrationnel. On y trouve l'hypothèse d'indépendance algébrique \ref{1intro_hypothese2i} ; pour certains espaces construits dans cette thèse, celle-ci est notamment vérifiée grâce aux nombres évoqués dans le lemme~\ref{3lem_sigma_alg_indep}. 

\begin{lemr}[\ref{3lem_1_irrat}]
 Soit $1 \leq d \leq n -1$. Soit $M = \begin{pmatrix} G \\ \Sigma \end{pmatrix} \in \mathcal{M}_{n,d}(\R)$ vérifiant :
 \begin{enumerate}[label=(\roman*)]
 \item $G \in \GL_{d}(\R)$ et $\Sigma \in \mathcal{M}_{n-d,d}(\R )$. 
 \item Les coefficients de $\Sigma$ forment une famille algébriquement indépendante sur $\Q(\FF)$ où $\FF$ est la famille des coefficients de $G$. \label{1intro_hypothese2i}
 \end{enumerate}
Alors l'espace engendré par $M$ est $(e,1)$-irrationnel pour tout $e \in \llbracket 1, n-1 \rrbracket$.
\end{lemr}

\bigskip 
Enfin on établit un théorème qui permet de calculer certains exposants diophantiens d'un espace $A$ qui est une somme de droites contenues dans des espaces rationnels supplémentaires. 

\begin{theor}[\ref{theo_somme_sev}]
Soit $d \in \left\llbracket 1, \floor{\frac{n}{2} }\right \rrbracket$. On suppose que l'on a $ \bigoplus\limits_{j=1}^d R_j \subset \R^n $ avec $R_j$ des espaces rationnels de dimension $r_j$.\\ Soit $A = \bigoplus\limits_{j=1}^d A_j$ avec $A_j \subset R_j$ et $\dim(A_j) = 1$. 
Pour $J \subset \llbracket 1, d \rrbracket $, on pose $A_J = \bigoplus\limits_{j \in J} A_j$. 
\\ Soit $e \in \llbracket 1, n-1 \rrbracket $ et $k \in \llbracket 1,\min(d,e) - g(A,e) \rrbracket$. On a alors l'équivalence entre les assertions suivantes :
\begin{enumerate}[label = (\roman*)]
\item $A$ est $(e,k)\tir$irrationel.
\item$\forall J \in \PP(k + g(A,e), d), \quad A_J$ est $(e,k+g(A,e) -g(A_J,e))\tir$irrationel.
\end{enumerate}
En outre, dans ce cas on a :
\begin{align*}
\mu_n(A|e)_k =\max\limits_{J \in \PP(k + g(A,e), d)} \mu_n({A}_{J}| e)_{k+g(A,e) -g(A_J,e)}.
\end{align*}
\end{theor}

\subsection{Construction d'espaces avec exposants prescrits}\label{1intro_section_expo_prescrits}

Une grande partie de cette thèse est consacrée à la construction d'espaces dont on peut calculer les exposants. On réalise d'abord la construction suivante.

\begin{theor}[\ref{5theo_droite_approx}]
 Soit $(\gamma_1, \ldots, \gamma_{n-1}) \in \R^{n-1}$ vérifiant $\gamma_1 \geq 2 +\frac{ \sqrt{5}-1}{2} $ et 
 \begin{align*}
 \forall i \in \llbracket 2, n-1 \rrbracket,& \quad \gamma_i \geq \left(2 +\frac{ \sqrt{5}-1}{2}\right)\gamma_{i-1},\\
 \forall (i,j) \in \llbracket 1, n-2 \rrbracket^2,& \quad i+j \leq n-1 \Longrightarrow \gamma_{i+j} \leq \gamma_i \gamma_j.
 \end{align*}
 On peut alors construire explicitement une droite $A$ de $\R^n$ vérifiant $A \in \II_n(1,n-1)_1$ et :
 \begin{align*}
 \forall e \in \llbracket 1, n-1 \rrbracket, \quad \mu_n(A|e)_1 = \gamma_e.
 \end{align*}
\end{theor}

On construit ainsi une droite dont on peut prescrire les exposants diophantiens pour tout $e \in \llbracket 1, n-1 \rrbracket$. Les valeurs prescrites ne sont cependant pas optimales dans le sens où on connait toutes les valeurs prises par le spectre joint $(\mu_n(A|1)_1, \ldots, \mu_n(A|n-1)_1) $ par le théorème~\ref{1theo_roy_spec_joint}. \\
Cette construction va surtout servir de \og brique de base \fg{} pour les suivantes. En effet, tous les espaces considérés dans la suite seront une somme directe de droites ainsi construites.
\\On montre ensuite le résultat suivant, plus général, qui concerne encore le premier angle $(j =1)$.

\begin{theor}[\ref{6theo_principal_avec_gamma}]
Soit $d \in \llbracket 1 , n-1 \rrbracket$. Il existe une constante $C_d$ dépendant de $d$ et $n$ telle que pour $(\gamma_1, \ldots, \gamma_{n-d}) \in \R^{n-d}$ vérifiant $\gamma_1 \geq C_d $ ainsi que 
 \begin{align*}
 \forall i \in \llbracket 2, n-d \rrbracket,& \quad \gamma_i \geq C_d \gamma_{i-1}, \\
 \forall (i,j) \in \llbracket 1, n-d-1 \rrbracket^2,& \quad i+j \leq n-d \Longrightarrow \gamma_{i+j} \leq \gamma_i \gamma_j.
 \end{align*}
 On peut alors construire explicitement $A \in \II_n(1,n-d)_1$ tel que :
 \begin{align*}
 \forall e \in \llbracket 1, n-d \rrbracket, \quad \mu_n(A|e)_1 = \gamma_e.
 \end{align*}
\end{theor}

\bigskip
En supposant maintenant que $d$ divise $n$, on obtient le théorème suivant. On construit ici un sous-espace $A$ de $\R^n$ de dimension $d$ dont on peut calculer un grand nombre d'exposants diophantiens. On connait en particulier des exposants diophantiens de cet espace correspondant à la fois à $e $ vérifiant $g(d,e,n) = 0 $ et à $e$ vérifiant $g(d,e,n) > 0$. 

 \begin{theor}[\ref{7theo_principale}]
 Soit $n =d (m +1)$, $c_{\ref{7cons_petite_hyp_theoc2c1}} =\left(1+ \frac{1}{m} \right)^{\frac{1}{d}}$ et $1 < c_{\ref{7cons_petite_hyp_theoc2}} < c_{\ref{7cons_petite_hyp_theoc2c1}} $.
\\ Soit $(\beta_{1,1}, \ldots, \beta_{1, m} ) \in \R^m $ tels que :
\begin{align*}
 \min\limits_{\ell \in \llbracket 1, m \rrbracket}(\beta_{1,\ell}) > {(3d)^{\frac{c_{\ref{7cons_petite_hyp_theoc2}}}{c_{\ref{7cons_petite_hyp_theoc2}}-1}}} \text{ et } \min\limits_{\ell \in \llbracket 1, m \rrbracket}(\beta_{1,\ell})^{c_{\ref{7cons_petite_hyp_theoc2c1}}} > \max\limits_{\ell \in \llbracket 1, m \rrbracket}(\beta_{1,\ell})^{c_{\ref{7cons_petite_hyp_theoc2}}}.
\end{align*}
Pour $i \in \llbracket 2,d \rrbracket$, soit $(\beta_{i,1}, \ldots, \beta_{i, m}) \in \R^m$ vérifiant pour tout $i \in \llbracket 1,d-1 \rrbracket$ : 
\begin{align*}
 \min\limits_{\ell \in \llbracket 1, m \rrbracket}(\beta_{i,\ell})^{c_{\ref{7cons_petite_hyp_theoc2c1}}} > \max\limits_{\ell \in \llbracket 1, m \rrbracket}(\beta_{i+1,\ell})
 \\
 \text{ et } \min\limits_{\ell \in \llbracket 1, m \rrbracket}(\beta_{i+1,\ell}) > \max\limits_{\ell \in \llbracket 1, m \rrbracket}(\beta_{i,\ell})^{c_{\ref{7cons_petite_hyp_theoc2}}} .
\end{align*}
Il existe un espace $A$ de dimension $d$ dans $\R^{n}$ tel que pour tous $e \in \llbracket 1, n-1\rrbracket $ et $k \in \llbracket 1 + g(d,e,n), \min(d,e) \rrbracket$ vérifiant $e < k (m+1) $ on a $A \in \II_n(d,e)_{k-g(d,e,n)}$ et : 
\begin{align*}
 \mu_n(A|e)_{k-g(d,e,n)}
 &= \bigg( \sum\limits_{q = 1 + \max(0, e-mk)}^k \frac{1}{ \max\limits_{\ell \in \llbracket 0, m-v_q \rrbracket} (\beta_{q+d-k, \ell +1} \ldots \beta_{q+d-k, \ell + v_q}) }\bigg)^{-1}
\end{align*}
où $v_1, \ldots, v_k$ sont définis en posant $u$ et $v$ tels que $e =k v + u$ soit la division euclidienne de $e$ par $k$ et :
\begin{align*}
 v_q &= \left\{
 \begin{array}{lllc}
 v + 1 &\text{ si } q \in \llbracket 1, u \rrbracket \\
 v &\text{ si } q \in \llbracket u+1, k \rrbracket.
 \end{array}
 \right. 
\end{align*}
\end{theor}

Les théorèmes~\ref{6theo_principal} et~\ref{7theo_principale} impliquent notamment les résultats énoncés dans la section~\ref{1section_inde_algebrique}.

\bigskip
On construit enfin un dernier exemple de sous-espace $A$ de dimension $d$ de $\R^n$ dont on peut calculer les exposants diophantiens au dernier angle $(j= \min(d,e) ) $. Il complète dans ce cas-là le théorème~\ref{7theo_principale}, en permettant d'avoir d'autres restrictions sur les valeurs prises par les exposants.

\begin{theor}[\ref{8theo_dernier_angle}]
On suppose que $d$ divise $n$.
\\Pour tout $\alpha \geq 2d(d+4)$, il existe un sous-espace vectoriel $A$ de $\R^n$ tel que 
\begin{align*}
\forall e \in \llbracket d, n-d \rrbracket, \quad &\mu_n(A| e)_d = \frac{\alpha^{q_e +1 } }{r_e + (d-r_e)\alpha} \\
 \forall e \in \llbracket 1, d -1\rrbracket, \quad &\mu_n(A| e)_e = \frac{\alpha }{r_e } = \frac{\alpha}{e}
\end{align*}
où $q_e$ et $ r_e$ sont le quotient et le reste de la division euclidienne de $e $ par $d$.
\end{theor}

\begin{req}
    Dans cette thèse, les constructions des chapitres $\ref{chap5}, \ref{chap6}, \ref{chap7}$ et $\ref{chap8}$ font appel aux nombres construits dans le lemme \ref{3lem_sigma_alg_indep}. En considérant des nombres de Liouville à la place de ceux-ci, on peut aboutir à des constructions d'espaces dont les exposants diophantiens sont infinis. On peut consulter la fin de la démonstration du théorème \ref{Joseph_spec} dans \cite{Joseph-these}, chapitre 6,  à ce sujet.
\end{req}

\section{Plan de la thèse}

Dans le \textbf{chapitre \ref{chap2}}, on énonce des résultats connus, qui vont nous servir dans cette thèse. Ils ont principalement été énoncés et prouvés par Schmidt \cite{Schmidt} et Joseph \cite{Joseph-these}, \cite{joseph_exposants}, \cite{joseph_spectre}.
\bigskip
\\On introduit de nouveaux outils dans le \textbf{chapitre \ref{chap3}}. On démontre des résultats sur le passage à l'orthogonal (dont la propriété~\ref{3prop_egalit_mun_orhto_mun_normal}). Ensuite, on prouve une propriété permettant de calculer la hauteur d'un sous-espace rationnel $B$ en utilisant certaines projections orthogonales de $\R^n$. On construit explicitement (dans la section~\ref{3section_constr_espac_irr}) des nombres transcendants qui vont servir dans toutes les constructions d'espaces effectuées dans le reste de la thèse. La dernière partie du chapitre est la preuve du lemme~\ref{3lem_1_irrat}, énoncé ci-dessus.
\bigskip
\\Le \textbf{chapitre \ref{chap4}} est entièrement consacré à la preuve du  théorème~\ref{theo_somme_sev} qui permet de calculer les exposants diophantiens de certains sous-espaces vectoriels. 
\bigskip
\\Le  théorème~\ref{5theo_droite_approx} est prouvé dans le \textbf{chapitre \ref{chap5}} ; sa preuve occupe l'intégralité du chapitre. On y construit une droite grâce notamment aux outils développés dans la section~\ref{3section_constr_espac_irr} du chapitre \ref{chap3}. On calcule ensuite les exposants diophantiens associés à cette droite.
\bigskip
\\Le \textbf{chapitre \ref{chap6}} est consacré à la preuve du  théorème~\ref{6theo_principal}. Pour cela on fait une récurrence sur la dimension $d$ de l'espace à construire, en s'appuyant sur le cas $d = 1$, étudié dans le chapitre \ref{chap5}.
\bigskip 
\\Le \textbf{chapitre \ref{chap7}} est dédié à la preuve du  théorème~\ref{7theo_principale}. On y construit un espace $A$ défini comme une somme directe orthogonale de droites en utilisant la construction du chapitre \ref{chap5}. On calcule ensuite les exposants diophantiens en jeu notamment grâce au  théorème~\ref{theo_somme_sev}.
\bigskip
\\On étudie certains spectres joints dans le \textbf{chapitre \ref{chap9}}. Grâce aux résultats des chapitres précédents, on démontre des conditions suffisantes pour que ces spectres contiennent un ouvert non vide. Les théorèmes \ref{9prop_premier_angle}, \ref{9lem_sousfam_pour_theo2} et \ref{9lem_sousfam_pour_theo3} sont en particulier démontrés ici.
\bigskip
\\Enfin le \textbf{chapitre \ref{chap8}} est consacré à la preuve du  théorème~\ref{8theo_dernier_angle}. On utilise encore une fois les droites explicitées au chapitre \ref{chap5} pour construire un espace $A$ de dimension $d$. On calcule ensuite les exposants diophantiens $\mu_n(A|e)_{\min(d,e)}$ pour $e \in \llbracket 1, n-d \rrbracket$, en utilisant notamment les résultats de géométrie des nombres exposés au chapitre \ref{chap2}.

%% file: Chapitres/1_Outils.tex
\chapter{Outils}\label{chap2}

Dans ce chapitre, on rappelle quelques résultats connus qui seront utiles dans le développement de cette thèse. Pour cela on introduit la notion de réseau euclidien et on reprend des propriétés notamment données par Schmidt \cite{Schmidt}, \cite{Schmidt_book} et Joseph \cite{Joseph-these}.

\bigskip
Dans cette thèse, toutes les constantes de la forme $c_{i}$ avec $i \in \N$ sont supposées strictement positives sauf mention du contraire.

\section{Géométrie des nombres }
\subsection{D'autres formules pour la hauteur}
On reprend ici la construction de la hauteur, pour donner une formule que l'on utilisera en pratique. On rappelle que $\| \cdot \| $ désigne la norme euclidienne canonique sur $\R^n$. On identifie de plus $\bigwedge^e \R^n$ à $\R ^{ \binom{n}{e}}$ par l'isomorphisme canonique.
\begin{defi}
 On appelle réseau euclidien de $\R^n$ tout $\Ztir$module libre $\Gamma$ discret de $\R^n$. 
\\Pour $\Gamma$ un réseau euclidien, on appelle rang de $\Gamma$ le rang du $\Ztir$module associé.
\end{defi}

\begin{req}
 Ici on ne demande pas qu'un réseau euclidien de $\R^n$ soit de rang $n$.
\end{req}
\begin{defi}
 Soit $\Gamma$ un réseau euclidien de $\R^n$ de rang $r$ et $X_1, \ldots, X_r$ une $\Zbase$ de $\Gamma$. 
 \\On appelle parallélotope fondamental de $\Gamma$ le pavé semi-ouvert $$P = \Big\lbrace \sum\limits_{i = 1}^r t_iX_i, 0 \leq t_i < 1 \Big\rbrace \subset \Vect(\Gamma).$$ 
 \\ On appelle covolume de $\Gamma$ et on note $\covol(\Gamma)$ la quantité :
 \begin{align*}
 \covol(\Gamma) = \vol(P)
 \end{align*}
 où l'on a muni $\Vect(\Gamma)$ de la mesure de Lebesgue.
\end{defi}

Schmidt \cite{Schmidt} donne alors une autre définition de la hauteur.
\begin{prop}\label{2prop_haut_covol}
 Soit $B$ un sous-espace rationnel de $\R^n$. On a : 
 \begin{align*}
 H(B) = \covol(B \cap \Z^n).
 \end{align*}
\end{prop}

Cette propriété donne un point de vue plus géométrique sur la notion de hauteur que dans la définition~\ref{def_hauteur}.

En pratique pour calculer la hauteur d'un espace $B$, on exhibe souvent une $\Zbase$ de $B \cap \Z^n$. A cet effet, on énonce un lemme montré dans \cite{Cassels_geom} (Corollaire 3 du Théorème 1, page 14).

\bigskip 
\begin{lem}\label{lem_Cassels_compl}
 Soit $\Gamma$ un réseau euclidien de rang $f$ de $\R^n$. Pour $e < f$, soit $X_1,\ldots, X_e$ des vecteurs linéairement indépendants de $\Gamma$. 
\\ Alors il existe des vecteurs $X_{e+1}, \ldots, X_f \in \Gamma$ tels que $X_1, \ldots, X_f$ forme une $\Zbase$ de $\Gamma $ si et seulement si :
\begin{align*}
 \forall \, (u_1, \ldots, u_e) \in \R^e ,\quad u_1 X_1 + \ldots + u_e X_e \in \Gamma \Longrightarrow (u_1, \ldots, u_e) \in \Z^e.
\end{align*}
\end{lem}

Autrement dit, on peut compléter en une $\Zbase$ de $\Gamma$, toute famille $X_1, \ldots, X_e $ qui est une $\Zbase$ du module $\Vect(X_1, \ldots, X_e) \cap \Gamma$. 

\begin{req}\label{2req_sous_Zbase}
 La réciproque est vraie dans le sens où si $X_1, \ldots, X_f$ est une $\Zbase$ de $\Gamma$ alors pour tout $e < f$ la famille $X_1, \ldots, X_e$ est une $\Zbase $ de $\Vect(X_1, \ldots, X_e) \cap \Z^n$.
\end{req}

\bigskip

On peut alors montrer une autre formule pour la hauteur d'un sous-espace rationnel :

\begin{prop}\label{prop_haut_prod_ext}
 Soit $B \in \RR_n(e)$. Soit $X_1, \ldots, X_e $ une $\Zbase$ de $B \cap \Z^n$. Alors :
 \begin{align*}
 H(B) = \| X_1 \wedge \ldots \wedge X_e \|
 \end{align*}
 où $X_1 \wedge \ldots \wedge X_e \in \Lambda^e (\R^n)$ est le $e\tir$vecteur produit extérieur des $X_1, \ldots, X_e$. 
\end{prop}

\begin{req}\label{2req_maj_hauteur_base_entiere}
 De manière générale on a $H(B) = \frac{ \| X_1 \wedge \ldots \wedge X_e \|}{N(\alpha(X))}$ pour toute base rationnelle $X_1, \ldots, X_e $ de $B$, en notant $\alpha(X)$ l'idéal fractionnaire engendré par les coordonnées grassmaniennes associées à cette base et $N(\alpha(X))$ la norme de cet idéal. 
 On déduit donc de cette remarque que pour toute base $Z_1, \ldots, Z_e$ de $B$ formée de vecteurs de $\Z^n$, on a 
 $$H(B) \leq \|Z_1 \wedge \ldots \wedge Z_e \|. $$
 
\end{req}

\begin{preuve}
 Les coordonnées grassmaniennes sont les mineurs maximaux de la matrice dont les colonnes sont les $X_i$ pour $i \in \llbracket 1,e \rrbracket$. De même on peut voir $ X_1 \wedge \ldots \wedge X_e$ comme un vecteur de $\R^{N}$, avec $N = \dbinom{n}{e}$, dont les coordonnées sont précisement (au signe près) ces mineurs. On a alors :
 \begin{align*}
 \| X_1 \wedge \ldots \wedge X_e \| = \| (\eta_1, \ldots, \eta_N) \|
 \end{align*}
 où $\eta_1, \ldots, \eta_N$ sont les coordonnées grassmaniennes associées aux $X_i$.
 \\ D'après la définition~\ref{def_hauteur}, en posant $\alpha = \eta_1 \Z + \ldots + \eta_N \Z$ il reste alors seulement à montrer que $N(\alpha) = 1 $. Or $N(\alpha) = \pgcd(\eta_1, \ldots, \eta_N).$ Le lemme suivant montre que cette quantité est bien égale à $1$.
 
\end{preuve}

\begin{lem}
Soit $n \geq 2$ et $B \in \RR_n(e)$. Alors pour toute $\Zbase$ de $B\cap \Z^n$, en notant $\eta_1, \ldots, \eta_N$ les coordonnées grassmaniennes associées, on a :
 \begin{align*}
 \pgcd(\eta_1, \ldots, \eta_N) = 1.
 \end{align*}
\end{lem}

Ce lemme est prouvé dans \cite{Schmidt_book}, lemme $5H$ page $18$.

\subsection{Théorème d'équidistribution de Weyl}

On énonce un théorème d'équidistribution attribué à Weyl sur les parties fractionnaires d'une famille de nombres ; il est utilisé au chapitre \ref{chap7} dans la démonstration du corollaire~\ref{7cor_min_expos}. On trouve une preuve de ce résultat dans \cite{Kuipers_Niederreiter}, chapitre $1.6$.

\bigskip
Pour $x = (x_1, \ldots, x_n) \in \R^n$ et $y = (y_1, \ldots, y_n) \in \R^n$ vérifiant $0 \leq x_j < y_j \leq 1$ pour tout $j \in \llbracket 1,n \rrbracket$, on pose :
\begin{align*}
 [x,y[ \:= [x_1, y_1[ \times \ldots \times [ x_n, y_n[. 
\end{align*}
Pour $(a_1, \ldots, a_n) \in \R^n$ et $N \in \Nx$ on pose :
\begin{align*}
 A(N,[x,y[) = \# \left\{ 1 \leq k \leq N, (\partfrac{ ka_1}, \ldots, \partfrac{ ka_n} ) \in [x,y[ \right\} 
\end{align*}
 où $\partfrac{\cdot}$ est la partie fractionnaire.
 
\begin{theo}
 Si la famille $1,a_1, \ldots, a_n$ est linéairement indépendante sur $\Q$ alors :
 \begin{align*}
 \underset{N \to + \infty}{\lim} \frac{A(N,[x,y[)}{N} = \prod\limits_{i = 1}^n (y_i-x_i)
 \end{align*}
 pour tout $ [x,y[ \: \subset [0,1[^n$.
\end{theo}

\begin{cor}\label{2cor_kronecker}
 Soit $(a_1, \ldots, a_n) \in \R^n$ tel que la famille $1,a_1, \ldots, a_n$ soit linéairement indépendante sur $\Q$. Alors pour tout $(b_1, \ldots, b_n) \in \R^n$ l'ensemble :
 \begin{align*}
 \left\{ (\partfrac{a_1N +b_1}, \ldots, \partfrac{ a_nN +b_n}), \quad N \in \N \right\}
 \end{align*}
 est dense dans $[0,1[^n$.
\end{cor}


\subsection{Théorème de Minkowski}\label{Theo_Mink}
On énonce ici un théorème fondamental de la géométrie des nombres que l'on utilise dans le chapitre \ref{chap8}.. 

\begin{theo}[Minkowski]\label{2th_Minko}
Soit $n \in \Nx$ et $\Gamma$ un réseau euclidien de rang $n$ de $\R^n$. Soit $\CC$ une partie convexe symétrique de $\R^n$ vérifiant :
\begin{align*}
 \vol(\CC) > 2^n \covol(\Gamma).
\end{align*}
Alors $\CC \cap \Gamma \neq \{0 \}$.
\end{theo}

On énonce également un corollaire dans le cas où $B$ est rationnel et $\Gamma = B \cap \Z^n$, ce qui permet de faire le lien avec la hauteur de l'espace $B$ car alors $\covol(\Gamma) = H(B)$. On peut le retrouver dans \cite{Schmidt_book}, lemme $4B$ page $10$.

\begin{cor} \label{cor_Mink2}
 Soit $B$ un sous-espace rationnel de dimension $e$ de $\R^n$. Il existe alors une constante $\cons \label{cons_minko2} >0 $ ne dépendant que de $n$, et un vecteur $X \in B \cap \Z^n \smallsetminus \{ 0 \}$ vérifiant :
 \begin{align*}
 \|X \| \leq c_{\ref{cons_minko2}} H(B)^{{1}/{e}}.
 \end{align*}

\end{cor}


\subsection{Appartenance à un sous-espace rationnel}
Dans cette section, on montre un lemme qui est utilisé tout au long de cette thèse pour montrer qu'un vecteur entier appartient à un espace rationnel donné.

\begin{lem}\label{2lem_X_in_B}
 Soit $Y \in \Z^n $ et $B$ un sous-espace rationnel de dimension $e$. On note $ (X_1, \ldots, X_e )$ une base de $B$ avec $X_i \in \Z^n $ pour tout $i \in \llbracket 1, e \rrbracket$. On suppose que : 
 \begin{align*}
 \| Y \wedge X_1 \wedge \ldots \wedge X_e \| <1.
 \end{align*}
 Alors $Y \in B$.
\end{lem}

\begin{preuve}
 En notant $\| \cdot \|_{\infty}$ la norme infinie on a :
 \begin{align*}
 \| Y \wedge X_1 \wedge \ldots \wedge X_e \|_{\infty} \leq \| Y \wedge X_1 \wedge \ldots \wedge X_e \| <1.
 \end{align*}
Or $ Y \wedge X_1 \wedge \ldots \wedge X_e $ est à coordonnées entières car les vecteurs considérés le sont. Sa norme infinie est donc nulle et donc :
\begin{align*}
 Y \wedge X_1 \wedge \ldots \wedge X_e = 0.
\end{align*}
Il existe donc une relation de dépendance linéaire entre les vecteurs $Y, X_1, \ldots, X_e$ et comme la famille des $X_i$ est une base de $B$ on a bien :
\begin{align*}
 Y \in \Vect(X_1, \ldots, X_e) = B.
\end{align*}
\end{preuve}

\section{Calculs des angles}
On donne ici quelques outils pour calculer les quantités $\omega(X,Y)$ ou $\psi_j(A,B)$. On peut retrouver ces propriétés démontrées de manière exhaustive dans \cite{Joseph-these} et \cite{Schmidt}.

\begin{lem}[Inégalité triangulaire]\label{2lem_ineg_triang}
 Soit $X,Y, Z\in \R^n \smallsetminus \{0\}$ alors :
 \begin{align*}
 \omega(X,Z) \leq \omega(X,Y) +\omega(Y,Z).
 \end{align*}
\end{lem}

\begin{lem}\label{2lem_omega_XY_diff}
 Soit $X,Y \in \R^n \smallsetminus \{0\}$ alors :
 \begin{align*}
 \omega(X,Y) \leq \frac{\| X - Y \| }{\| X\| }.
 \end{align*}
\end{lem}

\begin{lem}\label{lem_inclusion_croissance}
 Soit $A,B, A',B'$ quatre sous-espaces de $\R^n$ vérifiant $A' \subset A$ et $B' \subset B$. Alors :
 \begin{align*}
 \forall k \in \llbracket 1, \min(\dim(A'), \dim(B')) \rrbracket, \quad\omega_k(A',B') \geq \omega_k(A,B).
 \end{align*}
D'autre part si $\dim(A) \leq \dim(B) $ :
\begin{align*}
 \forall X \in A \smallsetminus \{0 \}, \quad \omega(X, B) = \omega_1(\Vect(X), B) \leq \omega_{\dim(A)}(A,B).
\end{align*}
\end{lem}

On en déduit le résultat suivant sur les exposants.
\begin{cor}\label{2cor_croissance_exposants_inclusion}
 Soit $A' \subset A$ deux sous-espaces de $\R^n$ de dimensions respectives $d'$ et $d $. Alors pour $e$ et $j$ tels que $ A \in \II_n(d,e)_j$ on a :
 \begin{align*}
 A' \in \II_n(d',e)_{j+g(d,e,n)-g(d',e,n)} \text{ et } \mu_n(A|e)_j \geq \mu_n(A'|e)_{j+g(d,e,n)-g(d',e,n)}.
 \end{align*}
\end{cor}

\begin{preuve}
 Soit $e$ et $j$ tels que $ A \in \II_n(d,e)_j$. Alors pour tout sous-espace rationnel $B$ de dimension $e $ on a 
 \begin{align*}
 \dim(A' \cap B) \leq \dim(A \cap B) < j + g(d,e,n) \leq j+g(d,e,n)-g(d',e,n) + g(d',e,n).
 \end{align*}
 On en conclut donc que $A' \in \II_n(d',e)_{j+g(d,e,n)-g(d',e,n)}$. 
 \\ Soit $\varepsilon > 0$ ; par définition de l'exposant diophantien, il existe une infinité d'espaces $B$ rationnels de dimension $e $ tels que :
 \begin{align*}
 \omega_{j + g(d,e,n)}(A',B) = \psi_{j+g(d,e,n)-g(d',e,n)}(A',B) \leq H(B)^{-\mu_n(A'|e)_{j+g(d,e,n)-g(d',e,n)} + \varepsilon }.
 \end{align*}
En appliquant le lemme~\ref{lem_inclusion_croissance} avec $k = j + g(d,e,n)$ on a 
\begin{align*}
 \omega_{j + g(d,e,n)}(A,B) \leq \omega_{j + g(d,e,n)}(A',B) \leq H(B)^{-\mu_n(A'|e)_{j+g(d,e,n)-g(d',e,n)} + \varepsilon }.
\end{align*}
Or $ \omega_{j + g(d,e,n)}(A,B) = \psi_j(A,B)$, on en déduit donc que 
\begin{align*}
 \forall \varepsilon >0, \quad \mu_n(A|e)_j \geq \mu_n(A'|e)_{j+g(d,e,n)-g(d',e,n)} - \varepsilon
\end{align*}
ce qui conclut la preuve en faisant tendre $\varepsilon$ vers $0$.

\end{preuve}

\begin{lem}\label{angle_Xortho_sur_X}
 Soit $X \in \R^n \smallsetminus \{ 0 \}$ et $C$ un sous-espace de $\R^n$. Alors :
 \begin{align*}
 \omega(X, C) = \omega_1(\Vect(X), C) = \omega(X, p_C(X) ) = \frac{\| X - p_C(X) \| }{ \| X \| }
 \end{align*}
 où $p_C$ est la projection orthogonale sur $C$.
\end{lem}

\begin{defi}
 Soit $d,e \in \llbracket 1, n -1 \rrbracket$ avec $d +e \leq n$. Soit $A$ et $B$ des sous-espaces vectoriels de $\R^n$ de dimensions respectives $d $ et $e$. On définit la quantité : 
 \begin{align*}
 \varphi(A,B) = \prod\limits_{i = 1}^t \psi_i(A,B)
 \end{align*}
 avec $t = \min(d,e)$.
\end{defi}

\begin{lem} \label{lem_phi_wedge}
 Soit $d,e \in \llbracket 1, n -1 \rrbracket$ avec $d +e \leq n$. Soit $A$ et $B$ des sous-espaces vectoriels de $\R^n$ de dimensions respectives $d $ et $e$. 
 \\Alors, pour toutes bases $X_1, \ldots, X_d $ de $A$ et $Y_1, \ldots, Y_e$ de $B$, on a :
 \begin{align*}
 \varphi(A,B) = \frac{ \| X_1 \wedge \ldots \wedge X_d \wedge Y_1 \wedge \ldots \wedge Y_e \| }{\| X_1 \wedge \ldots \wedge X_d \| \cdot \| Y_1 \wedge \ldots \wedge Y_e \|}.
 \end{align*}
\end{lem}

On utilise souvent dans cette thèse, la conséquence suivante de ce lemme :

\begin{lem}\label{lem_phi_dim1}
 Soit $X \in \R^n$ et $C$ un sous-espace rationnel de dimension $e$ de $\R^n$. Soit $Z_1, \ldots, Z_e$ une $\Zbase$ de $C \cap \Z^n$. 
 \\Alors :
 \begin{align}\label{ega_phi_dim1}
 \| X \wedge Z_1 \wedge \ldots \wedge Z_e \| = \omega_1(\Vect(X),C) \| X \| H(C). 
 \end{align}
\end{lem}

\begin{preuve}
Si $X \in C$ alors les deux membres de l'égalité (\ref{ega_phi_dim1}) sont nuls.
\\En effet, d'une part on a $\omega_1(\Vect(X),C) = \omega(X, p_C(X) ) = \omega(X,X) = 0$. D'autre part, on a une relation de dépendance linéaire entre $X, Z_1, \ldots, Z_e $ et donc $X \wedge Z_1 \wedge \ldots \wedge Z_e = 0$.
\bigskip
\\ Sinon $X \notin C$ et on peut alors appliquer le lemme~\ref{lem_phi_wedge} avec $A = \Vect(X)$ et \linebreak $C = \Vect(Z_1, \ldots, Z_e)$ : 
 \begin{align*}
 \omega_1(\Vect(X),C) = \psi_1(\Vect(X), C) = \frac{ \| X \wedge Z_1 \wedge \ldots \wedge Z_e \| }{\| X\| \cdot \| Z_1 \wedge \ldots \wedge Z_e \|}.
 \end{align*}
De plus en appliquant le lemme~\ref{prop_haut_prod_ext} on a $H(C) = \| Z_1 \wedge \ldots \wedge Z_e \|$ et donc l'égalité (\ref{ega_phi_dim1}).
 
\end{preuve}

\bigskip
Dans \cite{joseph_exposants}, Joseph montre le résultat suivant sur le comportement de la proximité par l’opération de somme directe :

\begin{prop}\label{2prop_4.5Elio}
 Soient $n \geq 2$, et $F_1, \ldots, F_\ell, B_1, \ldots, B_\ell,2 \ell $
sous-espaces vectoriels de $\R^n$
tels que pour tout $i \in \{1, \ldots\ell \}$, $\dim F_i = \dim B_i = d_i \geq 1$. Supposons que les $F_i$
engendrent un sous-espace de dimension $k = d_1 + \ldots + d_\ell$ et de même pour les $B_i$.
Posons
\begin{align*}
 F = \bigoplus\limits_{i = 1}^\ell F_i \text{ et } B = \bigoplus\limits_{i = 1}^\ell B_i.
\end{align*}
Alors on a 
\begin{align*}
 \omega_k(F,B) \leq c_{F,n} \sum\limits_{i = 1}^\ell \omega_{d_i}(F_i,B_i) 
\end{align*}
où $c_{F,n} > 0$ est une constante qui dépend uniquement de $F_1, \ldots, F_\ell$ et de $n$.
\end{prop}

\begin{req}
Comme dans tous ses travaux, Joseph considère seulement le cas où $2k \leq n$. En reprenant la preuve on remarque que cette hypothèse n'est pas nécessaire. En effet Joseph montre d'abord l'existence de droites $D_{i,1}, \ldots, D_{i,d_i}$ de $F_i$ et $E_{i,1},\ldots,E_{i,d_i}$ de $B_i$ vérifiant :
\begin{align*}
  \sum\limits_{j = 1}^{d_i} \psi_1(E_{i,j}, D_{i,j}) \leq d_i \psi_{d_i}(F_i, B_i) \leq n\psi_{d_i}(F_i, B_i).
\end{align*}
Or le raisonnement tient si on remplace les $\psi_\ell$ par des $\omega_\ell$.
\end{req}

\section{Croissance des exposants diophantiens}
Le problème de déterminer le spectre joint d'une famille d'exposants, pose la question des relations possibles entre ces exposants. Les questions dites de transfert apparaissent alors naturellement. On énonce d'abord une propriété (dite de Going-up) donnée par Schmidt \cite{Schmidt}. Il est à noter qu'il existe une propriété similiaire de Going-down que l'on n'utilise pas dans cette thèse. Cette dernière a été énoncée et montrée dans \cite{Schmidt} (la preuve dans le cas complexe a été corrigée par Poëls dans \cite{Poels_GD}). 

\begin{prop}
 Soient $e \in \llbracket 1, n-1 \rrbracket$, $f \in \llbracket e+ 1, n \rrbracket $ et $B \in \RR_n(e)$. 
\\Alors il existe $\cons \label{cons_going_up} > 0 $ ne dépendant que de $n$ et un espace $C \in \RR_n(f)$ vérifiant :
\begin{align*}
 B \subset C \text{ et } H(C) \leq c_{ \ref{cons_going_up}} H(B) ^{\frac{n-f}{n-e}}.
\end{align*}
\end{prop}

On peut alors montrer une propriété de transfert sur les exposants.

\begin{prop}\label{prop_croissance_exposants}
 Soit $(d,e,f) \in \llbracket 1, n-1 \rrbracket^3$ tels que $e \leq f $ et $k \in \llbracket g(d,f,n)+1,\min(d,e) \rrbracket$.
\\Pour $A \in \II_n(d,f)_{k-g(d,f,n)}$ on a :
\begin{align*}
 A \in \II_n(d,e)_{k-g(d,e,n)} \text{ et } \mu_n(A|e)_{k - g(d,e,n)} \leq \frac{n-f}{n-e} \mu_n(A|f)_{k - g(d,f,n)}.
\end{align*}
\end{prop}

\begin{preuve}
On montre d'abord que $A \in \II_n(d,e)_{k-g(d,e,n)}$. 
\\Soit $B \in \RR_n(e)$, il existe $B' \in \RR_n(f)$ tel que $B \subset B'$. Alors :
\begin{align*}
 \dim(A \cap B) \leq \dim(A \cap B') < k 
\end{align*}
car $A \in \II_n(d,f)_{k-g(d,f,n)}$. On en conclut que $A \in \II_n(d,e)_{k-g(d,e,n)}$.
\\ On montre maintenant la deuxième partie de la propriété. On note $\mu_e = \mu_n(A|e)_{k - g(d,e,n)} $ et $ \mu_f = \mu_n(A|f)_{k - g(d,f,n)} $. \\On suppose d'abord que l'on a $ \mu_e < + \infty $.
\\Soit $\varepsilon > 0 $, il existe une infinité d'espaces $B \in \RR_n(e)$ tel que :
\begin{align}\label{2croissane_expo_1}
 \psi_{k - g(d,e,n)}(A,B) \leq H(B)^{-\mu_e + \varepsilon}.
\end{align}
Pour un tel espace $B$ il existe un espace $C \in \RR_n(f)$ tel que $B \subset C$ et :
\begin{align*}
 H(C) \leq c_{ \ref{cons_going_up}} H(B) ^{\frac{n-f}{n-e}} 
\end{align*}
avec $c_{ \ref{cons_going_up}} $ ne dépendant que de $n$. 
\\On rappelle que $\psi_{k - g(d,e,n)}(A,B) = \omega_{k}(A,B)$ et en utilisant le lemme~\ref{lem_inclusion_croissance} on a 
\begin{align}\label{2croissane_expo_2}
 \omega_{k}(A,C) \leq 
 \omega_{k}(A,B) \leq H(B)^{-\mu_e + \varepsilon} \leq c_{ \ref{cons_croissance_expo}} H(C)^{-\frac{n-e}{n-f}\mu_e + \frac{n-e}{n-f}\varepsilon}
\end{align}
avec $\cons \label{cons_croissance_expo} = c_{\ref{cons_going_up}}^{(\mu_e - \varepsilon)\frac{n-e}{n-f}}$ indépendante de $C$. 

Comme $\omega_{k}(A,C) = \psi_{k - g(d,f,n)}(A,C) $ on a donc l'existence d'une infinité d'espaces $C$ rationnels vérifiant :
\begin{align}\label{2croissane_expo_3}
 \psi_{k - g(d,f,n)}(A,C) \leq c_{ \ref{cons_croissance_expo}} H(C)^{-\frac{n-e}{n-f}\mu_e + \frac{n-e}{n-f}\varepsilon}.
\end{align}
Cette inégalité donne alors $\mu_f \geq \frac{n-e}{n-f}\mu_e - \frac{n-e}{n-f}\varepsilon$ pour tout $\varepsilon >0$ et donc $\mu_f \geq \frac{n-e}{n-f}\mu_e$.
\bigskip \\
Dans le cas où $\mu_e = + \infty $, les inégalités $(\ref{2croissane_expo_1})$, $(\ref{2croissane_expo_2})$ et $(\ref{2croissane_expo_3})$ sont vraies pour tout $\mu_e > 0$. \\On en déduit alors $\mu_f \geq \frac{n-e}{n-f}\mu_e$ pour tout $\mu_e > 0$ et donc $\mu_f = + \infty$.

\end{preuve}

On utilise plus souvent le corollaire suivant, en remarquant que $\frac{n-e}{n-f} \geq 1$ : 

\begin{cor}
 \label{croissance_exposants}
 Soit $(d,e,f) \in \llbracket 1,n-1 \rrbracket^3$ tel que $e \leq f $ et $k \in \llbracket g(d,f,n)+1,\min(d,e) \rrbracket$.
\\Pour $A \in \II_n(d,f)_{k-g(d,f,n)}$ on a :
\begin{align*}
 A \in \II_n(d,e)_{k-g(d,e,n)} \text{ et } \mu_n(A|e)_{k - g(d,e,n)} \leq \mu_n(A|f)_{k - g(d,f,n)}.
\end{align*}
\end{cor}

\section{Inclusion dans un sous-espace
vectoriel rationnel
}
Joseph montre dans \cite{Joseph-these}, \cite{joseph_spectre} le résultat suivant qui permet de calculer les exposants diophantiens d'un espace plongé dans un espace ambiant plus grand.

\begin{prop}\label{2prop_elio_plongement}
 Soient $n \geq 2$ et $k \in \llbracket 2, n\rrbracket$. 
 \\Soient $d, e \in \llbracket 1 + g(d,e,k) - g(d,e,n), k - 1\rrbracket $ et $j \in \llbracket 1, \min(d, e) -g(d,e,n)\rrbracket$. Soit $A$ un sous-espace vectoriel de $\R^n$ de
dimension $d$ tel qu’il existe un sous-espace vectoriel rationnel $F \in \mathcal{R}_n(k)$ vérifiant $A \subset F$.
\\Notons $\Phi$ un isomorphisme rationnel de $F$ dans $\R^k$ et $\widetilde{A} = \Phi(A)$, qui est un sous-espace vectoriel de dimension $d$ de $\R^k$.
Supposons que pour tout sous-espace rationnel $B'$ de $F$ de dimension $e$, on a
\begin{align*}
 \dim(A \cap B') < j +g(d,e,n).
\end{align*} 
Alors $A \in \mathcal{I}_n(d, e)_j
, \widetilde{A} \in \mathcal{I}_k(d, e)_{j+ g(d,e,n)-g(d,e,k)}$ et $
\mu_n(A|e)_j = \mu_k(\widetilde{A}|e)_{j+ g(d,e,n)-g(d,e,k)}$.
\end{prop}

 \begin{req}
Encore une fois, Joseph considère seulement le cas où $d +e \leq k$. En reprenant la preuve on remarque que cette hypothèse n'est pas nécessaire.
\end{req}

Ce théorème est notamment utilisé au chapitre \ref{chap4} dans le cadre du corollaire suivant : 
\begin{cor}\label{cor_incl_Elio}
 Soit $n \in \Nx$ et $\Phi$ un isomorphisme rationnel de $\R^n$ dans lui-même. 
 Alors pour tous $A \in \II_n(d,e)_j$, $e \in \llbracket 1,n-1 \rrbracket $ et $j \in \llbracket 1, \min(d, e) -g(d,e,n)\rrbracket$ on a :
 \begin{align*}
 \Phi(A) \in \II_n(d,e)_j \text{ et } \mu_n(\Phi(A)|e)_j = \mu_n(A|e)_j.
 \end{align*}
 \end{cor}

%% file: Chapitres/2_Outils_bis.tex
\chapter{Nouveaux outils } \label{chap3}

L'objectif de ce chapitre est de donner divers nouveaux outils, qui seront utilisés tout au long de cette thèse. 
On commence d'abord par étudier les liens entre les exposants diophantiens d'un sous-espace $A$ et ceux de son orthogonal $A^\perp$. 
Ensuite, on prouve une propriété permettant de calculer la hauteur d'un sous-espace rationnel $B$ en utilisant certaines projections orthogonales de $\R^n$.
Enfin, on réalise les premières constructions de suites de nombres transcendants de la forme $$  \sum\limits_{k =0}^{ + \infty} \frac{u_k}{\theta^{\floor{\alpha_k}}}$$ avec $(u_k) \in (\Nx)^\N, \theta \in \Nx$ et $(\alpha_k) \in (\Rx)^{\N}$.

\section{Comportement par passage à l'orthogonal}
On développe dans cette section, quelques résultats sur le comportement du caractère irrationnel et des exposants diophantiens d'un sous-espace lorsque l'on considère son orthogonal. Ces idées manipulent simplement les définitions exposées en introduction mais le fait d'inclure le cas $d + e > n $ permet des nouveautés.

\begin{prop}\label{prop_ortho_j_irr}
 Soient $A$ un sous-espace de $\R^n$ de dimension $d$ et $e \in \llbracket 1,n-1 \rrbracket$. On a alors pour tout $ j \in \llbracket 1, \min(d,e) -g(d,e,n) \rrbracket $ : 
 \begin{align*}
 A \text{ est } (e,j)\tir\text{irrationnel} \Longleftrightarrow A^\perp \text{ est } (n-e,j)\tir\text{irrationnel}
 \end{align*}
en notant $A^\perp$ l'orthogonal de $A$ dans $\R^n$.
\end{prop}
\begin{preuve}
 Le lemme~\ref{lem_egal_angle_ortho} de Schmidt énoncé ci-dessous donne la propriété, mais on propose ici une preuve plus directe de ce résultat. 
 \bigskip \\
 Soit $B'$ un sous-espace rationnel de $\R^n$ de dimension $n-e$. Alors en notant $B = B^{'\perp} $, on a $\dim(B) = e$ et $B' = B^\perp$. 
 \\ On utilise maintenant la formule de Grassmann et on a :
 \begin{align*}
 \dim(A^\perp \cap B') &= \dim(A^\perp) + \dim(B') - \dim(A^\perp + B') \\
 &= (n -d) + (n-e) - \dim(A^\perp + B^\perp) \\
 &= (n -d) + (n-e) - \dim((A\cap B)^\perp) \\
 &= (n -d) + (n-e) - n + \dim((A\cap B)) \\
 &= n- d-e + \dim (A\cap B).
 \end{align*}
On rappelle que $g(d,e,n) = \max(0, d + e - n )$ et on remarque que $g(n-d,n-e,n) = g(d,e,n) + n -d -e $. Cela donne donc : 
\begin{align*}
 \dim(A\cap B) < j + g(d,e,n) &\Longleftrightarrow \dim(A^\perp \cap B') < j + g(d,e,n) + n - d -e \\
 & \Longleftrightarrow \dim(A^\perp \cap B') < j + g(n - d,n- e,n).
\end{align*}
La propriété~\ref{prop_ortho_j_irr} découle directement de cette équivalence.

\end{preuve}

\begin{prop}\label{prop_j_irr_descend_monte}
 Soit $A$ un sous-espace vectoriel de $\R^n$ de dimension $d$. 
 \\ Alors pour tout $j \in \llbracket 1, \min(d,n-d) \rrbracket$ :
 \begin{align*}
 A \in \II_n(d,n-d)_j \quad \Longrightarrow \quad \forall e \in \llbracket j,n-j \rrbracket, A \in \II_n(d,e)_j.
 \end{align*}
\end{prop}

\begin{preuve}
\textbullet \, \underline{Premier cas :} Soit $e \in \llbracket j,n-d \rrbracket$. 
 \\Soit $B$ sous-espace rationnel de dimension $e$, il existe alors $C$ un sous-espace rationnel de dimension $n-d $ tel que $B \subset C$. 
 On a alors :
 \begin{align*}
 \dim(A \cap B) \leq \dim(A \cap C) < j + g(d,n-d,n) = j + g(d,e,n).
 \end{align*}
 On a donc montré que si $A \in \II_n(d,n-d)_j$ alors $A \in \II_n(d,e)_j$ pour tout $ e \in \llbracket j,n-d \rrbracket$.
 
\bigskip
\textbullet \, \underline{Second cas :} Soit $e \in \llbracket n-d + 1,n-j \rrbracket$. \\D'après la propriété~\ref{prop_ortho_j_irr}, $A^\perp$ est $(n-d,j)\tir$irrationnel. 
 \\On utilise alors la première partie de la preuve et on a $A^\perp \in \II_n(n-d,f)_j$ pour tout $f \in \llbracket j, d \rrbracket $. En particulier comme $n-e \in \llbracket j, d \rrbracket $ on a :
 \begin{align*}
 A^\perp \in \II_n(n-d,n-e)_j.
 \end{align*}
 En appliquant encore une fois la propriété~\ref{prop_ortho_j_irr} on a bien : 
 \begin{align*}
 A \in \II_n(d,e)_j.
 \end{align*}
\end{preuve}

On utilise dans cette thèse, le cas $j = 1$ énoncé dans le corollaire suivant.

\begin{cor}\label{cor_j_irr_descend_monte}
 Soit $A$ un sous-espace vectoriel de $\R^n$ de dimension $d$. Alors
 \begin{align*}
 A \in \II_n(d,n-d)_1 \quad \Longrightarrow \quad \forall e \in \llbracket 1,n-1 \rrbracket, A \in \II_n(d,e)_1.
 \end{align*}
\end{cor}

\bigskip
De plus Schmidt \cite{Schmidt} donne les lemmes suivants pour le passage à l'orthogonal dans la hauteur et les angles.

\begin{lem}\label{hauteur_ortho}
 Soit $B$ un sous-espace vectoriel rationnel de $\R^n$. Alors :
 \begin{align*}
 H(B^\perp) = H(B).
 \end{align*}
\end{lem}

\begin{lem}\label{lem_egal_angle_ortho}
 Soit $A,B$ deux sous-espaces vectoriels de $\R^n$. 
 \\Alors $\min(d,e)- g(d,e,n) = \min(n-d,n-e)- g(n-d,n-e,n) $ et :
\begin{center}
 $\forall j \in \llbracket 1, \min(d,e)- g(d,e,n) \rrbracket, \quad \psi_j(A,B) = \psi_j(A^\perp, B^\perp)$.
\end{center}
\end{lem}

Ces deux lemmes et la définition des exposants diophantiens permettent donc d'avoir la propriété suivante.

\begin{prop}\label{3prop_egalit_mun_orhto_mun_normal}
 Soit $d, e \in \llbracket 1, n-1 \rrbracket$ et $j \in \llbracket 1, \min(d,e)- g(d,e,n) \rrbracket $. Soit $A\in \II_{n}(d,e)_j $ alors
\begin{align*}
 \mu_n(A|e)_j = \mu_n(A^\perp| n -e )_j. 
\end{align*}
\end{prop}

On en déduit directement le corollaire suivant. On rappelle que le spectre $\SS(n,d,e,j)$ défini dans le problème~\ref{prob_spec}, est l'ensemble des valeurs prises par $\mu_n(A|e)_j$ quand $A$ décrit $\II_n(d,e)_j$.

\begin{cor}
Soit $d, e \in \llbracket 1, n-1 \rrbracket$ et $j \in \llbracket 1, \min(d,e)- g(d,e,n) \rrbracket $. Alors :
\begin{align*}
 \SS(n,d,e,j) = \SS(n,n-d,n-e,j).
\end{align*}
 \end{cor}

\section{Comportement de la hauteur par projection orthogonale} 

Dans cette section, on montre un résultat sur la hauteur d'un espace en la \og décomposant\fg{} sur l'image et le noyau d'une projection orthogonale particulière. Cette propriété est utilisée au chapitre \ref{chap4} dans la démonstration du lemme~\ref{Rec_expo}.

\begin{prop}\label{prop_haut_appli}
 Soit $n \in \N \smallsetminus\{0 \}$ et $e \in \llbracket 1, n \rrbracket$. Soit $ p : \R^n \longrightarrow \R^n$ une projection orthogonale vérifiant $p(\Z^n) \subset \Z^n$. On a alors pour tout $B \in \RR_n(e)$ : 
 \begin{align*}
 H(B) = H(\ker(p) \cap B) \: H(p(B)).
 \end{align*}
\end{prop}

\begin{req}
 La proposition \ref{prop_haut_appli} n'est pas vraie en général pour une projection rationnelle.
 \\ En effet, considérons $C = \R^2$ et $p : \R^2 \longrightarrow \R^2$ la projection orthogonale sur $V = \Vect\begin{pmatrix}
 1 & k 
 \end{pmatrix}^\intercal $ pour $k \in \N \smallsetminus \{ 0 \}$.
On a alors $\ker(p) = \Vect\begin{pmatrix}
 -k & 1 
 \end{pmatrix}^\intercal $ et $p(C) = V$, d'où $H(C) = 1 $ et $H(\ker(p)\cap C)H(p(C)) = 1 + k^2$.
\end{req}

\begin{req}
 La condition $p(\Z^n) \subset \Z^n$ est en fait très restrictive. \\On note $(e_i)_{i \in \llbracket 1, n \rrbracket}$ la base canonique de $\R^n$. Comme $p$ est une projection, la norme subordonnée de $p$ vérifie $\| p \| \leq 1$. 
 \\On a alors pour tout $ i \in \llbracket 1,n \rrbracket $:
 \begin{align*}
 \|p(e_i)\| \leq \|e_i\| =1.
 \end{align*}
 Or $p(e_i) \in \Z^n$ et donc $p(e_i) \in \{0, \pm e_1, \ldots,\pm e_n \}$.
 Les seules projections orthogonales possibles sont donc celles sur les espaces de la forme $\underset{ i \in I}{\Vect}(e_i) $ avec $I \subset \llbracket 1, n \rrbracket$.
 \end{req}

\begin{preuve}
 On travaille ici dans l'espace $B$ ; on considère la restriction de $p$ à $B$ en la notant toujours $p$. 
 \\ On note $t = \dim(p(B))$. L'application $p$ vérifie $p(\Z^n) \subset \Z^n$ donc $p(B)$ est rationnel. Comme $\dim(B) = e$ on a $\dim(\ker(p)) = e -t $.
 \\ Soit $X_1, \ldots, X_{e-t}$ une $\Zbase$ de $\ker(p) \cap \Z^n$. Comme $\ker(p) \cap \Z^n \subset B \cap \Z^n$, d'après le lemme~\ref{lem_Cassels_compl} il existe des vecteurs $X_{e-t+1}, \ldots, X_e$ de sorte que :
 \begin{align*}
 X_1, \ldots, X_e \text{ forme une $\Zbase$ de $B \cap \Z^n$}.
 \end{align*}
 On va montrer que l'on peut remplacer $X_i$ par $p(X_i)$ pour $i \in \llbracket e-t +1, e \rrbracket$. 
 \\En effet on a:
 \begin{align} \label{tilde_X_noyau}
 \forall i \in \llbracket e-t +1, e \rrbracket, \, \widetilde{X_i} = X_i - p(X_i) \in \ker(p) \cap \Z^n
 \end{align}
 car $p$ est une projection et $X_i$ et $p(X_i)$ sont des vecteurs entiers. Pour $i \in \llbracket e-t +1, e \rrbracket $ on écrit donc :
 \begin{align}\label{tilde_X_decom}
 \widetilde{X_i} =  \sum\limits_{k = 1}^{e-t} u_{i,k} X_k
 \end{align}
 avec $u_{i,k} \in \Z$ pour $k \in \llbracket 1, e-t \rrbracket$.
 \\ Tous les $\widetilde{X_i}$ sont donc des combinaisons entières de $X_1, \ldots, X_{e-t}$. En particulier la matrice de passage de $(X_1, \ldots, X_d)$ à $(X_1, \ldots, X_{e-t}, p(X_{e-t+1}), \ldots, p(X_{e}))$ est de la forme 
 \begin{align*}
 \begin{pmatrix}
 I_{e-t} & * \\
 0 & I_{t}
 \end{pmatrix} \in M_e(\Z)
 \end{align*}
d'après $(\ref{tilde_X_noyau})$ et $(\ref{tilde_X_decom})$. Le déterminant de cette matrice étant $1$ on en déduit que la famille $$X_1, \ldots, X_{e-t}, p(X_{e-t+1}), \ldots, p(X_{e})$$ est aussi une $\Zbase$ de $B \cap \Z^n.$
\\De plus on remarque que la famille $p(X_{e-t+1}), \ldots, p(X_{e})$ est une base de $p(B)$ donc c'est aussi une $\Zbase$ de $p(B) \cap \Z^n$.

On peut maintenant calculer les hauteurs : par la propriété~\ref{prop_haut_prod_ext} on a 
\begin{align*}
 H(B) &= \| X_1 \wedge \ldots \wedge X_{e-t} \wedge p(X_{e-t+1}) \wedge \ldots \wedge p(X_{e}) \| \\
 &= \| X_1 \wedge \ldots \wedge X_{e-t} \| \cdot \| p(X_{e-t+1}) \wedge \ldots \wedge p(X_{e}) \|,
\end{align*}
la deuxième inégalité provenant du fait que $\ker(p)$ et $p(B)$ sont des espaces orthogonaux. 
\\On conclut en appliquant encore une fois la propriété~\ref{prop_haut_prod_ext} pour les espaces $\ker(p)$ et $p(B)$. On a $\| X_1 \wedge \ldots \wedge X_{e-t} \| = H(\ker(p)) $ et $\| p(X_{e-t+1}) \wedge \ldots \wedge p(X_{e}) \| = H(p(B))$.
Cela donne donc :
\begin{align*}
 H(B) = H(\ker(p))H(p(B)).
\end{align*}

\end{preuve}

On énonce enfin un corollaire de cette propriété, qui pourrait se montrer par ailleurs en exhibant des $\Z\tir$bases des espaces considérés. Il est utilisé au chapitre \ref{chap7} dans la preuve du lemme~\ref{7lem_haut_CNM}.

\begin{cor}\label{2cor_haut_espace_ortho_somme}
 Soit $k \in \llbracket 1, n \rrbracket$. Soient $B \subset \R^k \times \{0 \}^{n-k}$ et $C \subset \{0\}^k \times \R^{n-k} $ des sous-espaces rationnels de $\R^n$. Alors :
 \begin{align*}
 H(B \oplus C) = H(B)H(C).
 \end{align*}
\end{cor}

\begin{preuve}
 On pose $p$ la projection orthogonale sur $\R^k \times \{0 \}^{n-k}$ alors $\ker(p) \cap (B \oplus C) = C $ et $p(B \oplus C) = B$. On applique la propriété~\ref{prop_haut_appli} pour conlcure.
 
\end{preuve}

\section{Construction de nombres transcendants}\label{3section_constr_espac_irr}
On énonce dans cette section, un lemme qui est utilisé dans beaucoup de constructions que l'on réalise dans cette thèse. On construit une famille de nomnbres transcendants qui est algébriquement indépendante sur $\Q$. En définissant des sous-espaces vectoriels grâce à ces nombres, on montre dans les chapitres suivants (\ref{chap5}, \ref{chap6}, \ref{chap7} et \ref{chap8}) que ces espaces sont $(e,j)\tir$irrationnel pour certains $e$ et $j$.

\begin{defi}
Soit $\theta \in \N \smallsetminus \{0,1 \}$, $u = (u_k)_{k \in \N} \in (\Nx)^\N$ une suite prenant un nombre fini de valeurs et $\alpha = (\alpha_k)_{k \in \N} $ une suite à valeurs dans $\R^{*}_+$ vérifiant 
\begin{align}\label{3condition_suite_alpha}
 \forall k \in \N, \quad \alpha_{k+1} > c_{\ref{constante_suite_alpha} } \alpha_k
\end{align}
avec $\cons \label{constante_suite_alpha} >  1$ indépendante de $k$.
On définit la quantité :
 \begin{align*}
 \sigma(\theta, u, \alpha) = \sum\limits_{k = 0}^{+ \infty} \frac{u_k}{\theta^{\floor{\alpha_k}}} < \infty.
 \end{align*} 
\end{defi}

\begin{req}
Dans la suite, quand la suite $u$ est claire, on notera souvent $\sigma(\theta, \alpha) $ pour $\sigma(\theta, u, \alpha)$. De plus, $u$ n'est jamais nulle, cela implique que $\sigma(\theta, \alpha) \notin \Q$.
\end{req}

\begin{req}\label{3req_sigma_transcendant}
Bugeaud  \cite{Bugeaud_approx_Cantor} prouve dans le cas où $c_{\ref{constante_suite_alpha}} 
 > 2$ :
\begin{align}\label{3sigma_valeur_exposant}
 \mu(\sigma(\theta, u)) = \limsup\limits_{k \to + \infty} \frac{\alpha_{k+1}}{\alpha_{k}} \geq c_{\ref{constante_suite_alpha} }
\end{align}
en s'appuyant sur les fractions continues \cite{Shallit}, \cite{Poorten_Shallit}. En particulier $\mu(\sigma(\theta, u)) > 2$ et donc $\sigma(\theta, u)$ est transcendant par le théorème de Roth (section~\ref{section_classique}). Dans le cas $c_{\ref{constante_suite_alpha}} 
 > 1$, le théorème de Ridout \cite{Ridout} donne le caractère transcendant de $\sigma(\theta, u)$ (voir \cite{Bugeaud_trans} pour plus de détails). 
\bigskip
\\Si $c_{\ref{constante_suite_alpha} } \geq \frac{ 3+ \sqrt{5}}{2}$, on peut calculer $ \mu(\sigma(\theta, u)) $ \og à la main \fg{} en considérant les sommes considérant les sommes tronquées $$\sigma^N(\theta,\alpha) = \sum\limits_{k = 0}^{N} \frac{u_k}{\theta^{\floor{\alpha_k}}} $$ qui sont alors les meilleurs approximations de $\sigma(\theta, \alpha)$, voir aussi \cite{Levesley_Salp_Velani} section $8$ à ce sujet. Dans cette thèse, on généralise ce genre de preuve ; cela explique pourquoi tous les exposants qui sont calculés dans les chapitres suivants sont grands. 

Si $c_{\ref{constante_suite_alpha} } \geq \frac{ 3+ \sqrt{5}}{2}$, on retrouve d'ailleurs le résultat $(\ref{3sigma_valeur_exposant})$ dans le cas particulier de la construction du chapitre~\ref{chap5} pour $n = 2$ (voir la remarque~\ref{5req_n=2}). 
 
\end{req}

\bigskip

\begin{lem}\label{3lem_sigma_alg_indep}
 Soit $e \in \Nx$, $\alpha = (\alpha_k)_{k \in \N} $ une suite à valeurs dans $\R^{*}_+$ vérifiant (\ref{3condition_suite_alpha}) et $\theta \in \N \smallsetminus \{0,1\}$. Soit $J \subset \Nx $ de cardinal au moins $2$ et $\FF \subset \R$ un ensemble fini de réels.
 \\ Soit $\phi : \N \longrightarrow \llbracket 0,e \rrbracket $ une application telle que :
 \begin{align*}
 \forall i \in \llbracket 0,e \rrbracket, \quad \#\phi^{-1}(\{ i \} ) = \infty.
 \end{align*}
Alors il existe un choix de $e +1 $ suites $(u_k^{i})_{ i \in \llbracket 0,e \rrbracket, k \in \N}$ telles que : 
\begin{align}\label{3condition_suite_u}
 \text{pour tous } i \in \llbracket 0,e \rrbracket \text{ et } k \in \N, \quad u_k^{i} &\left\{ \begin{array}{cl}
 \in J &\text{ si } \phi(k) = i \\
 = 0 &\text{ sinon }
 \end{array} \right. , 
\end{align}
et telles que la famille $(\sigma(\theta, (u_k^{i})_{k\in \N}, \alpha))_{ i \in \llbracket 0,e \rrbracket}$ soit algébriquement indépendante sur $\Q(\FF)$.
\end{lem}

\begin{preuve}
On note les quantités $(\sigma(\theta, (u_k^{i})_{k\in \N}, \alpha))_{ i \in \llbracket 0,e \rrbracket }$ par $\sigma_0, \ldots, \sigma_{e}$ et on raisonne par récurrence sur $ t \in \llbracket 0,e \rrbracket$. 
\\ L'ensemble des réels algébriques sur $\Q(\FF)$ est dénombrable car $\#\FF < + \infty$ et l'ensemble des suites $(u_k^0)$ vérifiant (\ref{3condition_suite_u}) est non dénombrable car $\#J \geq 2$. On choisit donc une suite telle que $\sigma_0$ soit transcendant sur $\Q(\FF)$.

\bigskip On suppose maintenant que l'on a construit $\sigma_0, \ldots, \sigma_t$ une famille algébriquement indépendante sur $\Q$ avec $t \in\llbracket 0, e-1\rrbracket$. L'ensemble des réels algébriques sur $\Q(\RR,\sigma_0, \ldots, \sigma_t)$ est dénombrable mais l'ensemble des suites $(u_k^{t+1})_{k \in \N}$ vérifiant $(\ref{3condition_suite_u})$ est non dénombrable car $\#J \geq 2$. 
\\On peut donc choisir une suite telle que $\sigma_{t+1} $ soit transcendant sur $\Q(\RR,\sigma_0, \ldots, \sigma_t)$ ce qui conclut la récurrence.

\end{preuve}

\section{Exemple d'espaces \texorpdfstring{$(e,1)\tir$}{}irrationnels}
Le lemme suivant est utile pour montrer qu'un sous-espace vectoriel est $(e,1)\tir$irrationnel. On reprend ici des idées de la preuve du lemme~$6.3$ de \cite{Joseph-these}, celui-ci étant un cas particulier du lemme~\ref{3lem_1_irrat} dans le cas $n = 2d$. Combiné avec les constructions du lemme~\ref{3lem_sigma_alg_indep}, il permet de construire des espaces que l'on sait $(e,1)\tir$irrationnel. Ce lemme contient notamment le cas $d +e >n $, grâce au corollaire~\ref{cor_j_irr_descend_monte}.

\begin{lem}\label{3lem_1_irrat}
 Soit $1 \leq d \leq n -1$. Soit $M = \begin{pmatrix} G \\ \Sigma \end{pmatrix} \in \mathcal{M}_{n,d}(\R)$ vérifiant :
 \begin{enumerate}[label=(\roman*)]
 \item $G \in \GL_{d}(\R)$ et $\Sigma \in \mathcal{M}_{n-d,d}(\R )$. \label{condition_1}
 \item Les coefficients de $\Sigma$ forment une famille algébriquement indépendante sur $\Q(\FF)$ où $\FF$ est la famille des coefficients de $G$. \label{condition_2}
 \end{enumerate}
Alors pour tout $e \in \llbracket 1, n-1 \rrbracket$, le sous-espace vectoriel engendré par les colonnes de $M$ est $(e,1)$-irrationnel.
\end{lem}

\begin{preuve}
On note $Y_1, \ldots, Y_d$ les colonnes de la matrice $M $ et $A = \Vect(Y_1, \ldots, Y_d)$ l'espace engendré par ces colonnes. 

\bigskip
On montre seulement que $A$ est $(n-d,1)$-irrationnel. Par le corollaire~\ref{cor_j_irr_descend_monte} on a alors que $A$ est $(e,1)$-irrationnel pour tout $e \in \llbracket 1, n-1 \rrbracket$.

 Soit $B$ un espace rationnel de dimension $n-d$. On suppose par l'absurde que $A \cap B \neq \{ 0 \} $. En notant $Z_1, \ldots, Z_{n-d}$ une base rationnelle de $B$, on a 
 \begin{align*}
 Y_1 \wedge \ldots Y_d \wedge Z_1 \wedge \ldots \wedge Z_{n-d} = 0.
 \end{align*}
 Dans la suite, on note $Q = (Z_1 | \cdots | Z_{n-d})$ la matrice dont les colonnes sont les $Z_i$. 
 \\L'égalité $ Y_1 \wedge \ldots Y_d \wedge Z_1 \wedge \ldots \wedge Z_{n-d} = 0 $ donne l'annulation du déterminant de la matrice suivante :
 \begin{align*}
 \begin{pmatrix}
 \begin{matrix}G \\ 
 \Sigma 
 \end{matrix} & \begin{matrix}
 Z_1 & \cdots & Z_{n-d}
 \end{matrix}
 \end{pmatrix} \in \mathcal{M}_n(\R).
 \end{align*}
Or ce déterminant est un polynôme à coefficients rationnels en les coefficients de $M$, car les $Z_i$ sont rationnels.
\\On peut aussi voir ce déterminant comme un polynôme de $\Q(\FF)[X_1, \ldots, X_{d(n-d)}]$ évalué en les $d(n-d)$ coefficients de $\Sigma$. Comme par \ref{condition_2} ces coefficients forment une famille algébriquement indépendante sur $\Q(\FF)$ ce polynôme est identiquement nul. 
\\ On peut alors remplacer les coefficients de $\Sigma$ par n'importe quelle famille de réels et le déterminant sera nul. On a donc :
\begin{align}\label{det_nul}
 \forall \Sigma \in \mathcal{M}_{n-d,d}(\R ), \quad \det \begin{pmatrix}
 \begin{matrix}G \\ 
 \Sigma 
 \end{matrix} & \begin{matrix}
 Z_1 & \cdots & Z_{n-d}
 \end{matrix}
 \end{pmatrix} = 0.
\end{align}

Pour $\Delta$ un mineur de taille $n-d$ de $Q$ on note $\Ind(\Delta)$ l'ensemble des indices $$1 \leq i_1 < \ldots < i_{n-d} \leq n $$ des lignes dont est extrait $\Delta$. On note aussi $\Mat({\Delta}) $ la sous-matrice de $Q$ de taille $n-d$ dont $\Delta$ est le déterminant.
\\ On montre par récurrence forte sur $r \in \llbracket 0, \min(d, n-d) \rrbracket$ la propriété : \\
\og Pour tout mineur $\Delta$ de taille $n-d$ de $Q$ tel que $\#(\Ind(\Delta) \cap \llbracket 1,d \rrbracket )= r $ on a $\Delta = 0$.\fg{} 
\\ \textbullet \, Si $r = 0 $ alors $\Delta$ est le mineur extrait des $n - d $ dernières lignes de $Q$ et donc en prenant $\Sigma = (0) $ dans $(\ref{det_nul})$ on a :
\begin{align*}
 0= \det \begin{pmatrix}
 \begin{matrix}G \\ 
 0
 \end{matrix} & \begin{matrix}
 * \\
 \Mat({\Delta})
 \end{matrix}
 \end{pmatrix} = \det(G) \Delta
\end{align*}
et donc $\Delta = 0 $ car $G$ est inversible.
\\ \textbullet \, Soit $r \in \llbracket 1,\min(d,n-d) \rrbracket$ et on suppose la propriété vraie pour tout $0 \leq k < r$. Soit $\Delta$ un mineur de taille $n-d$ de $Q$ tel que $ \#(\Ind(\Delta) \cap \llbracket 1,d \rrbracket )= r $.
\\Sans perte de généralité, on suppose que $ \Ind(\Delta) \cap \llbracket 1,d \rrbracket=\llbracket 1, r\rrbracket$ et $ \Ind(\Delta) \cap \llbracket d+1,n \rrbracket = \llbracket d+1,n-r \rrbracket$.
\\Pour $i$ un entier et $U$ une matrice on note $U_i$ la $i-$ème ligne de $U$.
\\On a donc $\Mat({\Delta}) = \begin{pmatrix}
 B_{1} \\ \vdots \\ B_r \\ B_{d+1} \\ \vdots \\ B_{{n-r}}
\end{pmatrix}$ avec $B_i \in \MM_{1,n-d}(\Q)$.
\\On définit une matrice $\Sigma = \begin{pmatrix}
 \Sigma_1\\ \vdots \\ \Sigma_{n-d}
\end{pmatrix} \in \mathcal{M}_{n-d,d}(\R )$ en posant :
\begin{align*}
 \Sigma_i = \left\{ \begin{array}{cl}
 G_{i-n+d+r} &\text{ si } i \in \llbracket n-d-r+1,n- d\rrbracket \\
 0 &\text{ sinon }
 \end{array} \right. .
\end{align*}
On a défini $\Sigma$ de sorte que les $r$ dernières lignes de $\Sigma$ soient égales aux $r$ premières lignes de $G$. \\
On utilise maintenant (\ref{det_nul}) : 
\begin{align*}
 0= \det \begin{pmatrix}
 \begin{matrix}G \\ 
 \Sigma
 \end{matrix} & B
 \end{pmatrix} = \det \begin{pmatrix}
 \begin{matrix}G_1\\ \vdots \\ G_r \\ G_{r+1} \\ \vdots \\ G_d \\ 0 \\ \vdots \\ 0 \\ G_1 \\ \vdots \\ G_r 
 \end{matrix} & 
 \begin{matrix} B_1\\ \vdots \\ B_r \\ B_{r+1}\\ \vdots \\ B_d \\ B_{d+1} \\ \vdots \\ B_{n-r} \\ B_{n-r+1} \\ \vdots \\ B_n 
 \end{matrix} 
 \end{pmatrix} = \det \begin{pmatrix}
 \begin{matrix}G_1\\ \vdots \\ G_r \\ G_{r+1} \\ \vdots \\ G_d \\ G_1 \\ \vdots \\ G_r \\ 0 \\ \vdots \\ 0 
 \end{matrix} & 
 \begin{matrix} B_{n-r+1} \\ \vdots \\ B_n \\ B_{r+1} \\ \vdots \\ B_{d} \\ B_1\\ \vdots \\ B_r \\ B_{d+1} \\ \vdots \\ B_{n-r} 
 \end{matrix} 
 \end{pmatrix} 
\end{align*} 
en échangeant des lignes.
\\ On note $M = \begin{pmatrix}
 \begin{matrix}G_1\\ \vdots \\ G_r \\ G_{r+1} \\ \vdots \\ G_d \\ G_1 \\ \vdots \\ G_r \\ 0 \\ \vdots \\ 0 
 \end{matrix} & 
 \begin{matrix} B_{n-r+1} \\ \vdots \\ B_n \\ B_{r+1} \\ \vdots \\ B_{d} \\ B_1\\ \vdots \\ B_r \\ B_{d+1} \\ \vdots \\ B_{n-r} 
 \end{matrix} 
 \end{pmatrix} \in \MM_n(\R) $.
 \\En faisant un développement de Laplace (théorème~\ref{theo_laplace}) avec $J = \llbracket 1,d \rrbracket$ on a: 
\begin{align}\label{Laplace_GB}
 0 =  \sum\limits_{I \in \PP(d,n)} (-1)^{\ell(I) + \ell(J)} \Delta_{I,J}(M) \Delta_{\Bar{I},\Bar{J}}(M).
\end{align}
On étudie alors les mineurs de $M$ selon le choix de $I \in \PP(d,n)$.
\\ Si $I \cap \llbracket d + r +1, n \rrbracket \neq \emptyset $, on a $\Delta_{I,J}(M) = 0 $ car il contient une ligne nulle. 
\\ Sinon $I \cap \llbracket d + r +1, n \rrbracket = \emptyset $ alors on a $\#(I \cap \llbracket r +1, r +d \rrbracket) \geq d- r $. 
\begin{itemize}[label = $\bullet$]
 \item Si $\#(I \cap \llbracket r +1, r +d \rrbracket) = d- r $ alors $I\cap \llbracket 1, r \rrbracket = \llbracket 1, r \rrbracket $. Dans ce cas, si $I \neq \llbracket 1, d\rrbracket$ alors $\Delta_{I,J}(M)=0$ . En effet, on a deux lignes égales dans ce mineur.
 \item Si $\#(I \cap \llbracket r +1, r +d \rrbracket) > d- r $ alors $\# (\Bar{I} \cap \llbracket r +1, r +d \rrbracket) < r $. Dans ce cas, $\Delta_{\Bar{I},\Bar{J}}(M)$ est un mineur de $B$ tel que $\#(\Ind(\Delta_{\Bar{I},\Bar{J}}(M)) \cap \llbracket 1, d \rrbracket ) <r$. Par hypothèse de récurrence, ce mineur est nul.
\end{itemize}
On reprend (\ref{Laplace_GB}) où le seul terme possiblement non nul est alors celui correspondant à $I = \llbracket 1,d \rrbracket$ et on a :
\begin{align*}
 0 = \pm \det(G) \det \begin{pmatrix} B_1\\ \vdots \\ B_r \\ B_{d+1} \\ \vdots \\ B_{n-r} 
 \end{pmatrix} = \pm \det(G) \Delta.
\end{align*}
Or $\det(G) \neq 0 $ par l'hypothèse \ref{condition_1}, donc $\Delta = 0$.

On a donc montré que tout mineur de taille $n-d$ de $Q$ était nul. En particulier $$\rang(B) < n-d $$ ce qui est contradictoire puisque $Z_1, \ldots, Z_{n-d}$ forment une base de $B$.

 \end{preuve}

%% file: Chapitres/3_Exposants_diophantiens_somme.tex
\chapter{Exposants diophantiens d'une somme }\label{chap4}

On montre dans ce chapitre, un résultat permettant de calculer les exposants diophantiens de certaines sommes directes de droites de $\R^n$. Le théorème ainsi montré permet notamment des calculs d'exposants diophantiens dits \og intermédiaires \fg{} (c'est-à-dire $\mu_n(A|e)_j$ avec $ j \neq 1, \min(d,e) )$. 

\section{Résultat principal}

Dans tout ce chapitre, comme $n$ ne varie pas on note $g(A,e) $ la quantité :
\begin{align*}
 g(A,e) = g(\dim(A), e,n) = \max(0, \dim(A) +e -n) .
\end{align*}

\begin{theo}\label{theo_somme_sev}
Soit $d \in \left\llbracket 1, \floor{\frac{n}{2} } \right\rrbracket$. On suppose que l'on a $ \bigoplus\limits_{j=1}^d R_j \subset \R^n $ avec $R_j$ des espaces rationnels de dimension $r_j$.\\ Soit $A = \bigoplus\limits_{j=1}^d A_j$ avec $A_j \subset R_j$ et $\dim(A_j) = 1$. Pour $J \subset \llbracket 1,d \rrbracket $, on pose $A_J = \bigoplus\limits_{j \in J} A_j$. 
\\ Soit $e \in \llbracket 1, n-1 \rrbracket $ et $k \in \llbracket 1,\min(d,e) - g(A,e) \rrbracket$. On a alors l'équivalence entre les assertions suivantes :
\begin{enumerate}[label = (\roman*)]
 \item $A$ est $(e,k)\tir$irrationnel.
 \item $\forall J \in \PP(k + g(A,e), d), \quad A_J $ est $(e,k+g(A,e) -g(A_J,e))\tir$irrationnel.
\end{enumerate}
En outre, dans ce cas on a :
\begin{align*}
\mu_n(A|e)_k = \max\limits_{J \in \PP(k + g(A,e), d)} \mu_n({A}_{J}| e)_{k+g(A,e) -g(A_J,e)}.
\end{align*}
\end{theo}

On rappelle que $\PP(k + g(A,e), d)$ désigne l'ensemble des parties à $k + g(A,e)$ éléments de $\llbracket 1, d\rrbracket$.

\begin{req}
 On remarque que $ \sum\limits_{j = 1}^d r_j \leq n$ et que pour avoir des droites $A_j$ non rationnelles il faut que $r_j \geq 2$. Cela nous restreint alors à considérer des \og petits\fg{} espaces dans le sens où $2d \leq n$.
\end{req}

\begin{ex}
 Soit $n = r_1 + r_2 $ avec $r_1, r_2 \geq 2$, on pose $$R_1 = \R^{r_1} \times \{0\}^{r_2} \text{ et } R_2 = \{0 \}^{r_1} \times \R^{r_2}.$$ Pour $i \in \{1, 2 \}$, soit $A_i$ une droite $(r_i-1, 1)\tir$irrationnel contenue dans $R_i$. En posant $A = A_1 \oplus A_2$, on a alors 
 \begin{align*}
 \forall e \in \llbracket 1, \min(r_1, r_2) -1 \rrbracket, \quad \mu_n(A|e)_1 = \max (\mu_n(A_1|e)_1, \mu_n(A_2|e)_1).
 \end{align*}
\end{ex}
\bigskip
La preuve de ce théorème se fait par récurrence sur $d$ : pour démontrer le théorème avec $A$, on l'applique aux $A_{\llbracket 1,d \rrbracket \smallsetminus \{j\} }$ (voir les détails ci-dessous). Le point clef est la proposition suivante, dont la preuve occupe l'essentiel de ce chapitre : 

\begin{prop}\label{prop_princi} Soient $d, e, k, R_j, A_j$ comme dans le théorème~\ref{theo_somme_sev}. 
\\Soit $J \subset \llbracket 1,d \rrbracket $, tel que $\#J \geq k+ g(A_J,e) + 1$. On a alors équivalence entre les assertions :
\begin{enumerate}[label = (\roman*)]
 \item $A_J$ est $(e,k)\tir$irrationnel. \label{4implic1}
 \item $\forall j \in J $, $A_{J \smallsetminus \{j\}}$ est $(e,k+g(A_J,e) -g(A_{J \smallsetminus \{j\}},e)) \tir$irrationnel. \label{4implic2}
\end{enumerate}

En outre, dans ce cas on a :
\begin{align*}
\mu_n(A_J|e)_k = \max\limits_{j \in J} \mu_n({A}_{J \smallsetminus \{j\}}| e)_{k+g(A_J,e) -g(A_{J \smallsetminus \{j\}},e)} .
\end{align*}
\end{prop}

\bigskip

\begin{preuve}[ (Théorème~\ref{theo_somme_sev})]
Si $d = 1 $ alors $g(A,e) = 0$ et $k = 1$ est la seule valeur possible pour $k$. Les conclusions du théorème sont alors triviales.

 Soit $d \geq 2$, on suppose que pour tout $1 \leq d' <d$ le théorème est vérifié. Par hypothèse on a $k \leq d -g(A,e)$.
 \\ \textbullet \, Si $k = d - g(A,e)$ alors le théorème est évident car $A_J = A $ pour $J = \llbracket 1,d \rrbracket$. 
 \bigskip
 \\ \textbullet \, Sinon $k + g(A,e) + 1 \leq d $ et on peut appliquer la propriété~\ref{prop_princi} avec $ J = \llbracket 1,d \rrbracket$ et donc $A_J = A$. On a alors l'équivalence entre :
\begin{enumerate}[label = (\roman*)]
 \item $A_J$ est $(e,k)\tir$irrationnel.
 \item $\forall j \in \llbracket 1,d \rrbracket $, $A_{J \smallsetminus \{j\}}$ est $(e,k+g(A,e) -g(A_{J\smallsetminus \{j\}},e)) \tir$irrationnel. \label{irr_ii}
\end{enumerate}
Or $A_{J \smallsetminus \{j\}}$ est de dimension $d -1 < d$ et on peut donc appliquer l'hypothèse de récurrence à l'espace $A_{J\smallsetminus \{j\}}$ avec $k' =k+g(A,e) -g(A_{J\smallsetminus \{j\}} $. On a donc l'équivalence de \ref{irr_ii} avec : 
\begin{align*}
 \forall j \in \llbracket 1,d \rrbracket, \forall I \subset J \smallsetminus \{j\}, \#I = k + g(A,e) \Longrightarrow A_{I} \text{ est } (e,k+g(A,e) -g(A_I,e))\text{-irrationnel}. 
\end{align*}
Or cette assertion se reformule immédiatement en :
\begin{align*}
 \forall I \in \PP(k + g(A,e),d), \quad A_{I} \text{ est } (e,k+g(A,e) -g(A_I,e))\text{-irrationnel}. 
\end{align*}
Cela démontre alors l'équivalence cherchée. \\
De plus, dans le cas où ces conditions sont vérifiées on a, par la propriété~\ref{prop_princi} : 
\begin{align*}
 \mu_n(A_J|e)_k &= \max\limits_{j \in J} \mu_n({A}_{J \smallsetminus \{j\}}| e)_{k+g(A_J,e) -g(A_{J \smallsetminus \{j\}},e)} 
 \end{align*}
 et donc par l'hypothèse de récurrence
 \begin{align*}
\mu_n(A_J|e)_k &= \max\limits_{j \in J} \max\limits_{\underset{\# I =k+g(A_J,e) }{I \subset J \smallsetminus \{j \} } } \mu_n({A}_I| e)_{k+g(A_J,e) - g(A_I,e)} \\
 &= \max\limits_{\underset{\# I =k+g(A_J,e) }{I \subset \llbracket 1, d \rrbracket} } \mu_n({A}_I| e)_{k+g(A_J,e) - g(A_I,e)} 
\end{align*}
et donc la récurrence est prouvée car $A_J = A$.

\end{preuve}

Il reste alors à montrer la propriété~\ref{prop_princi}. Le reste du chapitre est consacré à sa preuve. On montre d'abord l'implication $\ref{4implic1} \Longrightarrow \ref{4implic2}$, la plus simple, ainsi qu'une première inégalité sur les exposants dans le lemme~\ref{premiere_implication}.
Puis dans le lemme~\ref{Rec_angles}, on explique comment on peut \og passer \fg{} de l'espace $A_J$ à $A_{J \smallsetminus \{j \}}$ en exhibant une certaine projection orthogonale. Enfin on montre l'implication inverse $\ref{4implic2} \Longrightarrow \ref{4implic1}$ et l'égalité sur les exposants en jeu dans la propriété~\ref{prop_princi} grâce au lemme~\ref{Rec_expo}. Ce dernier fait notamment appel à la propriété~\ref{prop_haut_appli} démontrée au chapitre~\ref{chap3}.


\section{Preuve de la propriété~\ref{prop_princi} }
Dans la suite, on se donne $d, e,k, R_j$ et $A_j$ comme dans le théorème~\ref{theo_somme_sev} et on considère $J \subset \llbracket 1,d \rrbracket $ tel que $\#J \geq k+ g(A_J,e) + 1$.


\subsection{Première implication}

Une implication de la propriété~\ref{prop_princi} est plus simple que l'autre. On la montre ici en utilisant les définitions de l'irrationalité d'un sous-espace et des exposants diophantiens.

\begin{lem}\label{premiere_implication}
Si $A_J$ est $(e,k)\tir$irrationnel alors pour tout $j \in J$ : 
\begin{align*}
 A_{J \smallsetminus \{j \}} \text{ est } (e,k+g(A_J,e) - g(A_{J \smallsetminus \{j \} },e)\tir\text{irrationnel} 
\end{align*}
et de plus 
\begin{align*}
 \mu_n(A_J|e)_k \geq \max\limits_{j \in J} \mu_n({A}_{J \smallsetminus \{j\}}|e)_{k+g(A_J,e) - g(A_{J \smallsetminus \{j \}},e )}.
\end{align*}
\end{lem}

\begin{preuve}
 On suppose que l'espace $A_J$ est $(e,k)\tir$irrationnel et soit $j \in J$. On rappelle que $g(\dim(A_J),e,n) = g(A_J,e)$ et $g(\dim(A_{J \smallsetminus \{j \} },e,n) =g(\dim(A_{J \smallsetminus \{j \} },e)$.

Comme $A_{J \smallsetminus \{j\}} \subset A_J$ le corollaire~\ref{2cor_croissance_exposants_inclusion} donne $A_{J \smallsetminus \{j \}} \in \II_n(A_{J \smallsetminus \{j \}},e)_{k+g(A_{J},e) - g((A_{J \smallsetminus \{j \} },e)} $ ce qui prouve la première partie du lemme. De plus, par ce même corollaire on a 
\begin{align*}
 \mu_n(A_J|e)_k \geq \mu_n({A}_{J \smallsetminus \{j\}}| e)_{k+g(A_J,e) - g(A_{J \smallsetminus \{j \}}| e)}.
 \end{align*}
 Comme cela est valable pour tout $j \in J$ on a donc: 
 \begin{align*}
 \mu_n(A_J|e)_k \geq \max\limits_{j \in J}\mu_n({A}_{J \smallsetminus \{j\}}| e)_{k+g(A_J,e) - g(A_{J \smallsetminus \{j \}}| e)}.
 \end{align*}
 
\end{preuve}

\subsection{Passage à des sous-espaces orthogonaux}
Pour prouver la suite, on se ramène au cas où les $R_j$ sont engendrés par des vecteurs de la base canonique. 
On introduit les espaces $R'_j$ définis comme suit :
\begin{align*}
 &R'_j = \{0\}^{r_1 } \times \ldots \times \{0 \}^{r_{j-1}} \times {\R^{r_j}} \times \{0 \}^{r_{j+1}} \ldots \times \{0 \}^{r_d} \subset \R^n \text{ pour } j \in \llbracket 1, d\rrbracket.
\end{align*}
Pour $j \in \llbracket 1, d\rrbracket$, soit $\varphi_j $ un isomorphisme rationnel de $R_j$ dans $R'_j$. On pose alors $\varphi : \R^n \to \R^n$ l'application définie par :
\begin{align*}
 \forall j \in \llbracket 1,d \rrbracket, \quad \varphi_{|R_j} = \varphi_j.
\end{align*}
 Alors $\varphi$ est un isomorphisme rationnel de $\R^n$ vers $\R^n$. On a donc d'après le corollaire~\ref{cor_incl_Elio}:
\begin{align*}
 \mu_n(\varphi(A)|e)_k = \mu_n(A|e)_k 
 \text{ et } \mu_n(\varphi(A_J)|e)_k = \mu_n(A_J|e)_k
\end{align*}
pour tous $J \subset \llbracket 1, d \rrbracket$, $e \in \llbracket 1, n-1 \rrbracket$ et $k$ tels que les espaces considérés soient $(e,k)\tir$irrationnels.

\bigskip 
On peut donc dorénavant supposer que $R_j =R'_j $ pour tout $j \in \llbracket 1,d \rrbracket$.
\\Pour $J \subset \llbracket1,d \rrbracket$, on note $R_J = \bigoplus\limits_{j \in J} R_j$ et cette somme directe est orthogonale. On note $p_j$ la projection orthogonale sur $R_j$ et $\widehat{p_{j}} = id - p_j$ la projection orthogonale sur $R_j^\perp$.
\begin{req}\label{4norme_widehat_p_j}
 On remarque que l'on a pour tout $j \in \llbracket 1,d \rrbracket$, $\widehat{p_j} =  \sum\limits_{s \neq j } p_s$. Pour tous $X \in \R^n$ et $j \neq s$ on a en particulier :
 \begin{align*}
 \| p_s(X) \| \leq \| \widehat{p_j}(X) \|.
 \end{align*}
\end{req}

\subsection{Etude des angles entre espaces projetés }

Le but de cette partie est démontrer le lemme suivant qui permet \og diminuer\fg{} la dimension de l'espace $A_J$.

\begin{lem}\label{Rec_angles}
Soit $J \subset \llbracket 1,d \rrbracket $ tel que $\#J \geq k + g(A_J,e) + 1 $ et $C$ un sous-espace vectoriel de $\R^n$ de dimension $e$. Alors il existe $j = j(C) \in J$ tel qu'en notant $C_j = \widehat{p_j}(C)$ on ait : 
\begin{align*}
 &\dim(C_j) \geq k +g(A_J,e) \\
 \text{ et } &\omega_{k+g(A_J,e) }(A_{J\smallsetminus \{j \}}, C_j) \leq c_{\ref{cons_lemme_angle1}}\omega_{k+g(A_J,e) }({A}_{J}, C) 
\end{align*}
avec $\cons \label{cons_lemme_angle1} >0 $ une constante ne dépendant que de $n$.
\end{lem}

Pour démontrer le lemme~\ref{Rec_angles}, on constate que pour tout sous-espace vectoriel $C$ de dimension $e$, il existe $C' $ un sous-espace de $C$ de dimension $k+g(A_J,e)$ tel que 
\begin{align*}
 \omega_{k+g(A_J,e) }({A}_{J}, C) = \omega_{k+g(A_J,e)}({A}_{J}, C').
\end{align*} 
On trouve le $j \in J$ du lemme~\ref{Rec_angles} grâce au lemme suivant :

\begin{lem}\label{petit_vecteurs}
Soit $J$ une partie non vide de $\llbracket 1, d \rrbracket $ et $ C'$ un sous-espace vectoriel de $\R^n$ tel que $\dim(C') < \#J $. Posons $ \cons \label{cons_lemme_angle}= \frac{ 1}{\sqrt{n^2 +1 }}$. Alors il existe $j \in J$ tel que, pour tout $X \in C' $ : 
\[ \| \widehat{p_j}(X) \| \geq c_{\ref{cons_lemme_angle}} \|X\|.\]
\end{lem}

\bigskip

\begin{req}
 L'hypothèse $\dim(C') < \#J$ est optimale dans le sens où si $\dim(C') \geq \#J$, on peut avoir $\bigoplus\limits_{j \in J} \Delta_j \subset C'$ avec $\Delta_j \subset R_j$. En particulier, pour $X \in \Delta_j$, on a $\widehat{p_j}(X) = 0$ ce qui contredit le résultat du lemme.
\end{req}

\begin{req}
 Ce lemme est intéressant du fait que la constante $c_{\ref{cons_lemme_angle}}$ ne dépend pas de $C'$. En effet si l'on oublie cette contrainte on peut facilement conclure : comme $\dim(C') <~\#J$, il existe $j \in J$ tel que $R_j \cap C' = \{0 \} $. 
 \\ En particulier $\ker(\widehat{p_j}) \cap C' = \{0 \}$ et donc $\widehat{p_j}$ est injective sur $C'$. On peut alors prendre $$c_{\ref{cons_lemme_angle}} = \min\limits_{X \in C', \|X \| = 1} \|\widehat{p_j} (X) \|$$
 mais cette constante dépend alors de $C'$.
\end{req}

\begin{preuve}
On suppose le contraire par l'absurde.

On a donc l'existence d'une famille $(X_j)_{j \in J}$ de vecteurs de $C'$ telle que pour tout $j \in J $:
\begin{align*}
 \|\widehat{p_j}(X_j) \| < c_{\ref{cons_lemme_angle}} \|X_j\|.
\end{align*}
Les vecteurs $X_j$ sont deux à deux distincts. En effet, si on a $X_j = X_k$ pour $j \neq k$, en utilisant la remarque~\ref{4norme_widehat_p_j} on trouve :
\begin{align*}
 \| X_j \|^2 = \|p_j(X_j)\|^2 + \|\widehat{p_j}(X_j) \|^2 \leq \|\widehat{p_k}(X_j)\|^2 + \|\widehat{p_j}(X_j) \|^2 < 2c_{\ref{cons_lemme_angle}}^2 \|X_j\|^2 
\end{align*}
ce qui donne alors $c_{\ref{cons_lemme_angle}} > \frac{1}{\sqrt{2}}$. Or $c_{\ref{cons_lemme_angle}} = \frac{ 1}{\sqrt{n^2 +1 }} \leq \frac{1}{\sqrt{2}}$, ce qui soulève une contradiction.
\bigskip
\\ Pour tout $j \in J$ on a :
\begin{align*}
 \| \widehat{p_j}(X_j) \|^2 &< c_{\ref{cons_lemme_angle}}^2 \|X_j\|^2 \\
 &= c_{\ref{cons_lemme_angle}}^2(\|p_j(X_j)\|^2 +\|\widehat{p_j}(X_j)\|^2) 
\end{align*}
et donc $c_{\ref{cons_lemme_angle}}^2\| p_j(X_j) \|^2 > (1-c_{\ref{cons_lemme_angle}}^2)\|\widehat{p_j}(X_j)\|^2$. On rappelle que $\frac{\| \cdot \|^2_1}{n} \leq \| \cdot \|^2 \leq n \| \cdot \|^2_\infty $. Cela donne :
\begin{align} \label{inegal_pj_pjchapeau}
 \| p_j(X_j) \|_\infty^2 &>\frac{1-c_{\ref{cons_lemme_angle}}^2}{n^2c_{\ref{cons_lemme_angle}}^2} \| \widehat{p_j}(X_j) \|_1^2 = \| \widehat{p_j}(X_j) \|_1^2
\end{align}
par définition de $c_{\ref{cons_lemme_angle}}$.
\newpage
Notons $X_j = \begin{pmatrix} x_{1,j} & \cdots & x_{n,j} \end{pmatrix}^\intercal$ et $\ell_j$ un indice tel que $ \| p_j(X_j) \|_\infty = | x_{\ell_j,j} |$. On remarque que les indices $\ell_j$ sont distincts car les $R_i$ sont en somme directe, engendrés par des vecteurs de la base canonique et car $p_j(X_j) \in R_j$.

On étudie maintenant la famille $(X_j)$. On note $M \in \MM_{n,\#J}(\R) $ la matrice dont les colonnes sont ces $X_j$ pour $j \in J $. Enfin on note $M' $ la matrice carrée de taille $\#J$ extraite de $M$ dont les lignes correspondent aux lignes de $M$ indexées par les $\ell_j$. On a :
\begin{align*}
 M' = \begin{pmatrix}
 \begin{matrix}
 x_{\ell_1,j_1} \\ \vdots \\ x_{\ell_{\#J},j_1}
 \end{matrix}
 & \cdots &
 \begin{matrix}
 x_{\ell_1,j_{\#J}} \\ \vdots \\ x_{\ell_{\#J},j_{\#J}}
 \end{matrix}
 \end{pmatrix}
\end{align*}
en écrivant $J = \{j_1 < \ldots < j_{\#J} \}$. On va montrer que $M'$ est à diagonale strictement dominante donc inversible. Joseph utilise aussi cette propriété dans le chapitre~$6$ de \cite{Joseph-these} pour montrer qu'une certaine famille est une base.
\\ D'après l'inégalité (\ref{inegal_pj_pjchapeau}), pour tout $i \in \llbracket 1, \# J \rrbracket$ on a :
\begin{align*}
 | x_{\ell_{i},j_i} | = \| p_{j_i}(X_{j_i}) \|_\infty > \| \widehat{p_{j_i}}(X_{j_i}) \|_1 = \bigg\| \sum\limits_{s \neq j_i }p_s(X_{j_i}) \bigg\|_1 \geq \sum\limits_{s \neq {j_i} } |x_{\ell_s,j_i} |.
\end{align*}
Alors $M'$ est une matrice à diagonale strictement dominante donc elle est inversible. En particulier $M$ est de rang $\#J$ et les vecteurs $X_j$ forment une famille libre de $C'$ et donc 
\begin{align*}
 \dim(C') \geq \dim (\Vect_{j \in J} (X_j)) = \#J
\end{align*}
ce qui fournit une contradiction.

\end{preuve}

 Avant de passer à la preuve du lemme~\ref{Rec_angles} on montre le lemme technique suivant.
 
 Dans la suite, pour un vecteur $x \in \R^n$ et un espace $V \subset \R^n$ donné on note $x^V$ le projeté orthogonal de $x$ sur $V$. On note aussi $A_J^\perp$ l'orthogonal de $A_J$ dans $\R^n$.

\begin{lem}\label{lem_minoration_U_j}
 Soit $J \subset \llbracket 1,d \rrbracket$. Pour tous $U \in \R^n$ et $ j \in J $ on a : 
 \begin{align*}
 \|U^{A_J^\perp}\| \geq \bigg\| \widehat{p_j}(U) ^{A_{J \smallsetminus \{j\}} ^\perp}\bigg\|.
 \end{align*}
\end{lem} 

\begin{preuve}
On introduit, pour $A$ et $R$ deux sous-espaces vectoriels de $\R^n$, la notation suivante :
\begin{align*}
 A^{\perp_R} = A^\perp \cap R.
\end{align*}
On a donc les relations suivantes :
\begin{align}
A_J^\perp &= R_{\llbracket 1,d \rrbracket \smallsetminus J} \oplus \overset{\perp}{\bigoplus\limits_{i \in J}} A_i^{\perp_{R_i}}, \label{4ortho_A_J}\\
 A_{J \smallsetminus \{j\}}^\perp &= R_{\llbracket 1,d \rrbracket \smallsetminus {J} \cup \{j\}} \oplus \overset{\perp}{\bigoplus\limits_{i \in {J \smallsetminus \{j\}}}} A_i^{\perp_{R_i}}. \label{4ortho_A_J_j}
\end{align}
 Soit $U\in C$. Comme $p_j(U) \in R_j$ on a d'après $(\ref{4ortho_A_J})$ : 
 \begin{align*}
 p_j(U)^{A_J^\perp} = p_j(U)^{A_j^{\perp_{R_j}}} \in R_j. 
 \end{align*}
En effet $R_j$ est orthogonal à toutes les autres composantes de $A_J^\perp$.
\\ De même, comme $\widehat{p_j}(U) \in R_{\llbracket 1,d \rrbracket \smallsetminus \{j\} }$, on a $\widehat{p_j}(U)^{A_i^{\perp_{R_j}}} =0 $ et donc : 
\begin{align*}
 \widehat{p_j}(U)^{A_J^\perp} =\widehat{p_j}(U)^{R_{\llbracket 1,d \rrbracket \smallsetminus J } } + \sum\limits_{i \in J \smallsetminus \{j\}} \widehat{p_j}(U)^{A_i^{\perp_{R_i}} } \in R_{\llbracket 1,d \rrbracket \smallsetminus \{j\}}. 
\end{align*}
De plus comme $\widehat{p_j}(U)^{R_{\llbracket 1,d \rrbracket \smallsetminus J \cup \{j\}} } = \widehat{p_j}(U)^{R_{\llbracket 1,d \rrbracket \smallsetminus J }} + \widehat{p_j}(U)^{R_j} = \widehat{p_j}(U)^{R_{\llbracket 1,d \rrbracket \smallsetminus J }} $, en utilisant $(\ref{4ortho_A_J_j})$ on a :
\begin{align}\label{4egalité_sans_j}
 \widehat{p_j}(U)^{A_J^\perp} =\widehat{p_j}(U)^{R_{\llbracket 1,d \rrbracket \smallsetminus J \cup \{j \} } } + \sum\limits_{i \in J \smallsetminus \{j\}} \widehat{p_j}(U)^{A_i^{\perp_{R_i}} } = \widehat{p_j}(U) ^{A_{J \smallsetminus \{j\}} ^\perp}.
\end{align}
Les vecteurs $p_j(U) ^{A_J^\perp} \in R_j $ et $ \widehat{p_j}(U) ^{A_J^\perp} \in R_{\llbracket 1,d \rrbracket \smallsetminus \{j\}} = R_j^\perp$ sont donc orthogonaux. On peut donc minorer la norme de $U^{A_J^\perp} $ de la façon suivante :
\begin{align*}
 \|U^{A_J^\perp}\|^2 = \| p_j(U) ^{A_J^\perp}+ \widehat{p_j}(U) ^{A_J^\perp} \|^2 = \| p_j(U) ^{A_J^\perp} \|^2+ \| \widehat{p_j}(U) ^{A_J^\perp} \|^2 \geq \| \widehat{p_j}(U) ^{A_J^\perp} \|^2.
\end{align*}
Le lemme est donc montré car on a déjà vu que $\widehat{p_j}(U) ^{A_J^\perp} = \widehat{p_j}(U) ^{A_{J \smallsetminus \{j\}}^\perp} $ en (\ref{4egalité_sans_j}).

\end{preuve}

On a maintenant tous les outils pour montrer le lemme~\ref{Rec_angles}.
\begin{preuve}[ (lemme~\ref{Rec_angles})]
Soit $J \subset \llbracket 1,d \rrbracket$ tel que $\#J \geq k + g(A_J,e) + 1$. On rappelle que l'on note $C'$ un sous-espace vectoriel de $C$ de dimension $k + g(A_J,e)$ tel que :
\begin{align*}
 \omega_{k + g(A_J,e)}(A_J,C) = \omega_{k + g(A_J,e)}(A_J,C').
\end{align*}
Comme $\dim(C') < \#J$, le lemme~\ref{petit_vecteurs} donne $j \in J$ tel que pour tout $X \in C'$ 
\begin{align*}
 \| \widehat{p_j}(X) \| \geq c_{\ref{cons_lemme_angle}} \|X\|.
\end{align*}
On étudie l'espace $C'_j = \widehat{p_j}(C')$.
\\Comme pour tout $X \in C' \smallsetminus \{ 0 \}, \widehat{p_j}(X) \neq 0$, l'espace $C'_j$ est de dimension $k + g(A_j,e)$ comme $C'$ . Enfin $C'_j \subset C_j = \widehat{p_j}(C)$ et donc $ \dim(C_j) \geq \dim(C'_j) = k+g(A_J,e)$ ce qui donne la première conclusion du lemme. 

\bigskip 
On étudie maintenant l'espace $ A_{J \smallsetminus \{j\}}$. On remarque d'abord que $ A_{J \smallsetminus \{j\}} = \widehat{p_j}(A_J)$. 
\\Soit $U_j \in C'_j \smallsetminus \{0 \}$, il existe $U \in C'$ tel que $\widehat{p_j}(U) = U_j$.
\\ On rappelle que, d'après lemme~\ref{angle_Xortho_sur_X} : 
\begin{align*}
 \omega(U_j, A_{J \smallsetminus \{j\}}) &= \frac{\| U_j ^{A_{J \smallsetminus \{j\}}^\perp} \| }{\| U_j \| }.
\end{align*}
On utilise maintenant le lemme~\ref{lem_minoration_U_j} qui donne $\| U_j ^{A_{J \smallsetminus \{j\}}^\perp} \| \leq \| U ^{A_{J}^\perp} \|$ et on rappelle que $ \| U_j \| \geq c_{\ref{cons_lemme_angle}} \| U \|$. On a donc :
\begin{align*}
 \omega(U_j, A_{J \smallsetminus \{j\}}) &\leq \frac{\| U ^{A_J^\perp} \| }{c_{\ref{cons_lemme_angle}}\| U \| } \\
 &= \frac{1}{c_{\ref{cons_lemme_angle}}}\omega(U,A_J) \\
 &\leq \frac{1}{c_{\ref{cons_lemme_angle}}}\omega_{k + g(A_J,e)}(C',A_J) 
\end{align*}

car $ U \in C'$. La dernière inégalité provient du fait que $$\dim(C') = k + g(A_J,e) \leq \#J = \dim(A_J)$$ et on utilise le lemme~\ref{lem_inclusion_croissance}.
\\On rappelle que $\omega_{k + g(A_J,e)}(C',A_J) 
 = \omega_{k + g(A_J,e)}(C,A_J) $ par construction de $C'$ et donc pour tout $U_j \in C'_j \smallsetminus \{0 \}$ 
 \begin{align*}
 \omega(U_j,A_{J \smallsetminus \{j\}}) \leq \frac{1}{c_{\ref{cons_lemme_angle}}} \omega_{k + g(A_J,e)}(C,A_J).
 \end{align*}
Comme $k + g(A_J,e) = \dim(C'_j) \leq \#J -1 = \dim(A_{J\smallsetminus \{j \} })$, il existe $U_j \in C'_j \smallsetminus \{0 \}$ tel que 
\begin{align*}
 \omega(U_j, A_{J \smallsetminus \{j \} }) = \omega_{k+g(A_J,e)}(C'_j, A_{J\smallsetminus\{j \}} ).
\end{align*}
On a donc :
\begin{align*}
 \omega_{k+g(A_J,e)}(C'_j, A_{J\smallsetminus\{j \}} ) = \omega(U_j,A_{J \smallsetminus \{j\}}) \leq \frac{1}{c_{\ref{cons_lemme_angle}}}\omega_{k+g(A_J,e)}(C,A_J).
\end{align*}
Or $C'_j \subset C_j$ donc $\omega_{k+g(A_J,e)}(C'_j, A_{J\smallsetminus\{j \}} ) \geq \omega_{k+g(A_J,e)}(C_j, A_{J\smallsetminus\{j \}} ) $ d'après le lemme~\ref{lem_inclusion_croissance}. On a donc :
\begin{align*}
 \omega_{k+g(A_J,e)}(C_j, A_{J\smallsetminus\{j \}} ) \leq \frac{1}{c_{\ref{cons_lemme_angle}}}\omega_{k+g(A_J,e)}(C,A_J).
\end{align*}
Le lemme~\ref{Rec_angles} est prouvé avec $c_{\ref{cons_lemme_angle1}} = c_{\ref{cons_lemme_angle}}^{-1}$.

\end{preuve}

\subsection{Seconde implication et conclusion}
On a tous les outils pour montrer la seconde implication et le résultat sur les exposants. 

\begin{lem}\label{lem_irration_descente}
Supposons que pour tout $j \in J$, 
$$A_{J \smallsetminus \{j\}} \text{ est } (e,k+g(A_J,e ) -g(A_{J \smallsetminus \{j\}},e ))\tir\text{irrationnel.} $$ Alors $A_J$ est $(e,k)\tir$irrationnel.
\end{lem}

\begin{preuve}
 On raisonne par contraposition et on suppose que $A_J$ n'est pas $(e,k)\tir$irrationnel. 
 \\ Il existe alors $C$ un espace rationnel de dimension $e$ tel que $\dim(A_J \cap C) \geq k + g(A_J,e)$ ce qui est équivalent à $\omega_{k+g(A_J,e)}(A_J,C) = 0$.
\\Or d'après le lemme~\ref{Rec_angles} il existe $j \in J$ tel que : 
\begin{align*}
 &\dim(C_j) \geq k +g(A_J,e)\\
 \text{ et } &\omega_{k+g(A_J,e)}({A}_{J\smallsetminus \{j \}}, C_j) \leq c_{\ref{cons_lemme_angle1}} \omega_{k+g(A_J,e)}({A}_{J}, C) 
\end{align*}
avec $C_j = \widehat{p_j}(C)$.
L'espace $C_j$ est rationnel et de dimension inférieure ou égale à $e$. Il existe donc $\widetilde{C} $ un espace rationnel de dimension $e$ tel que $C_j \subset \widetilde{C}$. On a alors d'après le lemme~\ref{lem_inclusion_croissance} :
\begin{align*}
 0 \leq c_{\ref{cons_lemme_angle}} \omega_{k+g(A_J,e)}({A}_{J\smallsetminus \{j \}}, \widetilde{C}) \leq c_{\ref{cons_lemme_angle}} \omega_{k+g(A_J,e)}({A}_{J\smallsetminus \{j \}}, C_j) &\leq \omega_{k+g(A_J,e)}({A}_{J}, C) =0
\end{align*}
et donc
\begin{align*}
 \omega_{k+g(A_J,e) -g(A_{J \smallsetminus \{j\}},e) + g(A_{J \smallsetminus \{j\}},e) }({A}_{J\smallsetminus \{j \}}, \widetilde{C}) = \omega_{k+g(A_J,e)}({A}_{J\smallsetminus \{j \}}, \widetilde{C})=0
\end{align*} ce qui signifie que $A_{J\smallsetminus \{j\}}$ n'est pas $(e,k+g(A_J,e) -g(A_{J \smallsetminus \{j\}},e) )\tir$irrationnel.

\end{preuve}

\begin{lem}\label{Rec_expo}
Soit $J \subset \llbracket 1,d \rrbracket$ tel que $\# J \geq k + g(A_J,e) + 1$ alors:
\begin{align*}
 \max\limits_{j \in J} \mu_n({A}_{J \smallsetminus \{j\}}| e)_{k+g(A_J,e) - g(A_{J \smallsetminus \{j \}},e )} \geq \mu_n(A_J|e)_k.
\end{align*}
\end{lem}

\begin{preuve}
Soit $\varepsilon > 0 $. On note $\beta = \mu_n(A_J|e)_k$. Il existe alors une infinité de sous-espaces rationnels $C$ de dimension $e$ tels que :
\begin{align}\label{4min_psik_AJC}
 \psi_k(A_J,C) \leq H(C)^{-\beta + \varepsilon}.
\end{align}
D'après le lemme~\ref{Rec_angles}, pour tout tel espace $C$, il existe $j \in J$ tel que : \begin{align*}
 &\dim(C_j) \geq k +g(A_J,e) \\
 \text{ et } & \omega_{k +g(A_J,e)}({A}_{J\smallsetminus \{j \}}, C_j) \leq c_{\ref{cons_lemme_angle}}\omega_{k +g(A_J,e)}({A}_{J}, C) 
\end{align*}
en notant $C_j = \widehat{p_j}(C)$ et avec $c_{\ref{cons_lemme_angle}} > 0$ ne dépendant que de $n$.
\\Comme les $j \in J$ sont en nombre fini, il existe $j \in J$ tel que pour une infinité d'espaces rationnels $C$ on ait ces inégalités avec $C_j$.\\
 Comme le nombre de dimensions possibles pour les $C_j$ est fini, il existe $t \in \llbracket k + g(A_J,e), e \rrbracket$ tel que pour une infinité d'espaces rationnels $C$ de dimension $e$ vérifiant $(\ref{4min_psik_AJC})$, on ait :
\begin{align}\label{4inega_angle}
 &\dim(C_j) = t \nonumber \\
 \text{ et } & \omega_{k +g(A_J,e)}({A}_{J\smallsetminus \{j \}}, C_j) \leq c_{\ref{cons_lemme_angle}} \omega_{k +g(A_J,e)}({A}_{J}, C).
\end{align}
D'après la propriété~\ref{prop_haut_appli} on a 
\begin{align}\label{4inega_hauteur}
 H(C) = H(C_j) H(\ker(\widehat{p_j}) \cap C) \geq H(C_j).
\end{align} 
On note $\gamma_{e,j} = \mu_n(A_{J\smallsetminus \{j \}}|e)_{k + g(A_J,e)- g(A_{J\smallsetminus \{j \}},e)} $ et $\gamma_{t,j} = \mu_n(A_{J\smallsetminus \{j\}}|t)_{k + g(A_J,e) - g(A_{J\smallsetminus \{j \}},t)} $. 
D'après le corollaire~\ref{croissance_exposants} on a 
\begin{align}\label{4inega_exposant}
 \gamma_{e,j} \geq \gamma_{t,j} .
\end{align}
En reprenant la défintion de l'exposant diophantien $\gamma_{t,j} $, il existe donc une constante $\cons \label{4cons_gamma_t} > 0$ dépendant seulement de $A$ et de $\varepsilon$ telle que pour tout espace $B$ rationnel de dimension $t$ on a : 
\begin{align}\label{4ineg_gamma_t}
 c_{\ref{4cons_gamma_t}}H(B) ^{-\gamma_{t,j} - \varepsilon} \leq 
 \psi_{k + g(A_J,e) - g(A_{J\smallsetminus \{j\}},t)}({A}_{J\smallsetminus \{j\}}, B). 
\end{align}
Donc pour tout espace $C$ de dimension $e$ vérifiant $\dim(C_j) =t $ et $(\ref{4min_psik_AJC})$ on a :
\begin{align}\label{4inega_gamma_combine}
 c_{\ref{4cons_gamma_t}}H(C )^{-\gamma_{e,j} - \varepsilon} \leq c_{\ref{4cons_gamma_t}}H(C_j)^{-\gamma_{e,j} - \varepsilon} \leq c_{\ref{4cons_gamma_t}}H(C_j)^{-\gamma_{t,j} - \varepsilon} \leq 
 \psi_{k + g(A_J,e) - g(A_{J\smallsetminus \{j\}},t)}({A}_{J\smallsetminus \{j\}}, C_j), 
\end{align}
en utilisant les inégalités (\ref{4inega_hauteur}), (\ref{4inega_exposant}) et (\ref{4ineg_gamma_t}). 
\\ On rappelle que $ \psi_{k + g(A_J,e) - g(A_{J\smallsetminus \{j\}},t)}({A}_{J\smallsetminus \{j\}}, C_j) = \omega_{k +g(A_J,e)}({A}_{J\smallsetminus \{j \}}, C_j)$ si $\dim(C_j) = t$. En regroupant alors les inégalités (\ref{4min_psik_AJC}), (\ref{4inega_angle}) et (\ref{4inega_gamma_combine}) on trouve :
\begin{align*}
 c_{\ref{4cons_gamma_t}}c_{\ref{cons_lemme_angle}}^{-1} H(C )^{-\gamma_{e,j} - \varepsilon} \leq \omega_{k +g(A_J,e)}(A_J, C) \leq \psi_{k }(A_J,C) \leq H(C)^{-\beta + \varepsilon}
\end{align*}
et donc en faisant tendre $H(C)$ vers l'infini on obtient $\gamma_{e,j} \geq \beta - 2 \varepsilon$. 

Comme cette inégalité est vraie pour tout $\varepsilon$ on a donc $\gamma_{e,j} \geq \beta $ et donc en particulier : 
\begin{align*}
 \max\limits_{j \in J}\gamma_{e,j} = \max\limits_{j \in J} \mu_n(A_{J\smallsetminus \{j\}}|e)_{k + g(A_J,e) - g(A_{J\smallsetminus \{j \}},e)} \geq \mu_n(A_J|e)_k.
\end{align*}

\end{preuve}

Les lemmes \ref{premiere_implication}, \ref{lem_irration_descente}, \ref{Rec_expo} donnent alors la propriété~\ref{prop_princi}.

%% file: Chapitres/4_Approximation_droites.tex
\chapter{Construction de droites avec exposants diophantiens prescrits}\label{chap5}

Dans ce chapitre, on présente notre première construction d'espace avec exposants prescrits. On y étudie le cas où l'on veut construire une droite de $\R^n$. Ce cas est important car il est utilisé dans toutes les autres constructions de cette thèse (chapitres \ref{chap6}, \ref{chap7} et \ref{chap8}). 
\\Ici on a $d = 1 $ et $e \in \llbracket 1, n-1 \rrbracket$ donc $g(d,e,n) = 0$.

\bigskip

Ce chapitre établit le théorème suivant.

\begin{theo}\label{5theo_droite_approx}
 Soit $(\gamma_1, \ldots, \gamma_{n-1}) \in \R^{n-1}$ vérifiant $\gamma_1 \geq 2 +\frac{ \sqrt{5}-1}{2} $ et 
 \begin{align}
 \forall i \in \llbracket 2, n-1 \rrbracket,& \quad \gamma_i \geq \left(2 +\frac{ \sqrt{5}-1}{2}\right)\gamma_{i-1},\label{5hyp_gamma_th1} \\
 \forall (i,j) \in \llbracket 1, n-2 \rrbracket^2,& \quad i+j \leq n-1 \Longrightarrow \gamma_{i+j} \leq \gamma_i \gamma_j. \label{5hyp_gamma_th2}
 \end{align}
 On peut alors construire explicitement une droite $A$ de $\R^n$ vérifiant $A \in \II_n(1,n-1)_1$ et :
 \begin{align*}
 \forall e \in \llbracket 1, n-1 \rrbracket, \quad \mu_n(A|e)_1 = \gamma_e.
 \end{align*}
\end{theo}

\bigskip

\begin{req}\label{5req_n=2}
 Avec $n = 2 $ et $\gamma_1 \geq 2 +\frac{ \sqrt{5}-1}{2} \approx 2,618\ldots$, on construit une droite $A = \Vect \begin{pmatrix}
 1 \\ \xi
 \end{pmatrix}$ vérifiant $\mu_2(A|1)_1 = \gamma_1$. En utilisant la remarque~$\ref{1req_intro_generalisation}$ cela donne en particulier $\mu(\xi) = \gamma_1$ comme énoncé dans la remarque~\ref{3req_sigma_transcendant}. 
\end{req}

\bigskip
\begin{req}
L'hypothèse $(\ref{5hyp_gamma_th2})$, qui stipule que la suite $(\gamma_i)$ est sous-multiplicative, assure l'égalité $\mu_n(A|e)_1 = \gamma_e$. Sans cette hypothèse, on aurait seulement la minoration $\mu_n(A|e)_1 \geq \gamma_e$ comme cela est détaillé dans la preuve, section~\ref{5section_prop_technique}.
\end{req}

\bigskip
\begin{req}
 La preuve de ce théorème reprend des idées de la preuve du théorème~\ref{Joseph_spec} de Joseph. On mentionne qu'il existe aussi une droite $A$ de $\R^n$ vérifiant :
 \begin{align*}
 \forall e \in \llbracket 1, n-1 \rrbracket, \quad \mu_n(A|e)_1 =+ \infty.
 \end{align*}
 En effet, Joseph construit une droite $A$ vérifiant $\mu_n(A|1)_1 = + \infty$ et la propriété~\ref{croissance_exposants} donne $\mu_n(A|e)_1 = + \infty $ pour tout $e \in \llbracket 1, n-1 \rrbracket$.
\end{req}

\begin{req}
 Les valeurs prescrites des exposants diophantiens dans le théorème~\ref{5theo_droite_approx} ne sont pas optimales dans le sens où l'on connait toutes les valeurs possibles prises par le spectre joint $(\mu_n(A|1)_1, \ldots, \mu_n(A|n-1)_1) $ par le théorème~\ref{1theo_roy_spec_joint}. Ici, les hypothèses $(\ref{5hyp_gamma_th1})$ et $(\ref{5hyp_gamma_th2})$ sont plus restrictives.
 \\ On utilise, dans cette thèse, la preuve du théorème~\ref{5theo_droite_approx} pour réaliser les constructions effectuées dans les chapitres suivants.
\end{req}

\bigskip 
Pour prouver le théorème~\ref{5theo_droite_approx}, on montre un résultat plus général à savoir la propriété~\ref{5prop_technique}. On y construit une droite $A$ de $\R^n$ dont on peut prescrire les exposants diophantiens. Cette construction, réalisée dans la section~\ref{5section_BN}, fait appel aux nombres $\sigma(\theta, u, \alpha)$ que l'on a définis dans le chapitre \ref{chap3} (section~\ref{3section_constr_espac_irr}).
\\On minore l'exposant $\mu_n(A|e)_1$ dans le corollaire~\ref{5cor_min} en exhibant des sous-espaces rationnels $B_{N,e}$, définis dans la section~\ref{5sect_construction_B_N}, qui approchent bien $A$. En montrant ensuite dans le lemme~\ref{5lem_meilleurs_espaces} que ces sous-espaces sont les \og meilleurs \fg{} approximations de $A$, on majore $\mu_n(A|e)_1$ dans le corollaire~\ref{5cor_maj}. Pour la majoration de $\mu_n(A|e)_1$, on a besoin de minorer $\psi_1(A, B_{N,e})$, ce qui se fait en considérant, dans le lemme \ref{5min_psi_A_BN}, les mineurs de taille $2$ de certaines matrices et en montrant que l'un d'eux est forcément assez grand

\section{Propriété technique}\label{5section_prop_technique}
On montre le théorème~\ref{5theo_droite_approx} grâce à la propriété technique suivante dont la preuve occupe la plus grande partie de chapitre.

\begin{prop}\label{5prop_technique}
 Soit $(\beta_N)_{N \in \Nx} \in [ 2 +\frac{ \sqrt{5}-1}{2}, + \infty [ ^{\Nx}$ une suite périodique de période $T \in \Nx$.\\ Il existe une droite $A$ de $\R^n$ vérifiant $A \in \II_n(1,n-1)_1$ et :
 \begin{align*}
 \forall e \in \llbracket 1, n-1 \rrbracket, \quad \mu_n(A|e)_1 = \max\limits_{i \in \llbracket 0, T-1\rrbracket} \beta_{i+1}\ldots\beta_{i+e}.
 \end{align*}
 De plus la construction de $A$ est explicite.
\end{prop}

\begin{cor}\label{5cor_technique}
 Soit $(\beta_1, \ldots, \beta_{n-1}) \in [ 2 +\frac{ \sqrt{5}-1}{2}, + \infty [ ^{n-1}$.\\ Il existe une droite $A$ de $\R^n$ vérifiant $A \in \II_n(1,n-1)_1$ et :
 \begin{align*}
 \forall e \in \llbracket 1, n-1 \rrbracket, \quad \mu_n(A|e)_1 = \max\limits_{i \in \llbracket 0, n-1-e\rrbracket} \beta_{i+1}\ldots\beta_{i+e}.
 \end{align*}
\end{cor}

\begin{preuve}
 Ce corollaire se déduit directement de la propriété~\ref{5prop_technique} en introduisant la suite $(\beta_N)_{N \in \Nx}$ qui complète $\beta_1, \ldots, \beta_{n-1}$ définie par :
 \begin{align*}
 \forall j \in \llbracket n, 2n-2 \rrbracket, \quad \beta_j = 2 +\frac{ \sqrt{5}-1}{2}
 \end{align*}
 et 
 \begin{align*}
 \forall i \in \llbracket 1, 2n-2 \rrbracket, \forall k \in \N, \quad \beta_{i + k(2n-2)} = \beta_i.
\end{align*}
Cette suite $(\beta_N)_{N \in \Nx}$ ainsi construite est donc $(2n-2)\tir$périodique et à valeurs dans $[ 2 +\frac{ \sqrt{5}-1}{2}, + \infty [$ . \\
 D'après la propriété~\ref{5prop_technique}, il existe une droite $A$ de $\R^n$ vérifiant $A \in \II_n(1,n-1)_1$ et :
 \begin{align*}
 \forall e \in \llbracket 1, n-1 \rrbracket, \quad \mu_n(A|e)_1 = \max\limits_{i \in \llbracket 0, 2n-3\rrbracket} \beta_{i+1}\ldots\beta_{i+e}.
 \end{align*}
Il reste alors à montrer que 
\begin{align*}
 \max\limits_{i \in \llbracket 0, 2n-3\rrbracket} \beta_{i+1}\ldots\beta_{i+e} \leq \max\limits_{i \in \llbracket 0, n-1-e\rrbracket} \beta_{i+1}\ldots\beta_{i+e},
 \end{align*}
 l'autre inégalité étant triviale. \\
Soit $k \in \llbracket 0, 2n-3\rrbracket$ tel que $\beta_{k+1}\ldots\beta_{k+e} = \max\limits_{i \in \llbracket 0, 2n-3\rrbracket} \beta_{i+1}\ldots\beta_{i+e} $. 
On distingue $4$ cas. 
\\ \textbullet \, Si $k \in \llbracket 0, n-e -1 \rrbracket$, alors par définition
\begin{align*}
 \beta_{k+1}\ldots\beta_{k+e} \leq \max\limits_{i \in \llbracket 0, n-1-e\rrbracket} \beta_{i+1}\ldots\beta_{i+e}.
\end{align*}
\\ \textbullet \, Si $k \in \llbracket n-e, n-2 \rrbracket$, on a :
\begin{align*}
 \beta_{k+1}\ldots\beta_{k+e} &= \prod\limits_{\ell = 1}^{n-1-k} \beta_{k + \ell } \times \prod\limits_{\ell = n-k}^{e } \beta_{k + \ell }.
\end{align*}
Comme pour tout $ 
\ell \in \llbracket n-k, e \rrbracket$, $\beta_{k + \ell} = 2 +\frac{ \sqrt{5}-1}{2} \leq \beta_{k+\ell -e }$, on a
\begin{align*}
 \beta_{k+1}\ldots\beta_{k+e} &\leq \prod\limits_{\ell = 1}^{n-1-k} \beta_{k + \ell } \times \prod\limits_{\ell = n-k}^{e } \beta_{k + \ell -e } \\
 &\leq \prod\limits_{\ell = n-k-e}^{n-1-k} \beta_{k + \ell } \\
 &\leq \beta_{i + 1} \ldots \beta_{i +e }
\end{align*}
avec $i = n-1-e$ et donc $\beta_{k+1}\ldots\beta_{k+e} \leq \max\limits_{i \in \llbracket 0, n-1-e\rrbracket} \beta_{i+1}\ldots\beta_{i+e}.$
\\ \textbullet \, Si $k \in \llbracket n-1, 2n-2-e \rrbracket$, on a $\beta_{k +\ell } = 2 +\frac{ \sqrt{5}-1}{2} $ pour tout $\ell \in \llbracket 1, e\rrbracket$ et donc :
\begin{align*}
 \beta_{k+1}\ldots\beta_{k+e} &\leq \max\limits_{i \in \llbracket 0, n-1-e\rrbracket} \beta_{i+1}\ldots\beta_{i+e}.
\end{align*}
\\ \textbullet \, Si $k \in \llbracket 2n-1-e, 2n-3 \rrbracket$, on a :
\begin{align*}
 \beta_{k+1}\ldots\beta_{k+e} &= \prod\limits_{\ell = 1}^{2n-2-k} \beta_{k + \ell } \times \prod\limits_{\ell = 2n-1-k}^{e } \beta_{k + \ell } 
 \\
 &= \prod\limits_{\ell = 1}^{2n-2-k} \beta_{k + \ell } \times \prod\limits_{\ell = 1-k}^{e -2n +2} \beta_{k + \ell } 
\end{align*}
en utilisant la périodicité de $\beta$. 
\\Comme pour tout $ 
\ell \in \llbracket 1, 2n-2-k \rrbracket$, $\beta_{k + \ell} = 2 +\frac{ \sqrt{5}-1}{2} \leq \beta_{k+ \ell + e -2n +2 }$, on a
\begin{align*}
 \beta_{k+1}\ldots\beta_{k+e} &\leq \prod\limits_{\ell = 1}^{2n-2-k} \beta_{k+ \ell + e -2n +2 } \times \prod\limits_{\ell = 1-k}^{e -2n +2} \beta_{k + \ell } \\
 &\leq \prod\limits_{\ell = e -2n + 3 }^{e-k} \beta_{k+\ell } \times \prod\limits_{\ell = 1-k}^{e -2n +2} \beta_{k + \ell } \\
 &\leq \beta_{i + 1} \ldots \beta_{i +e }
\end{align*}
avec $i = 1$ et donc $\beta_{k+1}\ldots\beta_{k+e} \leq \max\limits_{i \in \llbracket 0, n-1-e\rrbracket} \beta_{i+1}\ldots\beta_{i+e}.$

\end{preuve}

On peut alors déduire le théorème~\ref{5theo_droite_approx} du corollaire~\ref{5cor_technique}.

\begin{preuve}[ (théorème~\ref{5theo_droite_approx})]\label{5preuve_du_theo}
On construit $(\beta_1, \ldots, \beta_{n-1}) \in [ 2 +\frac{ \sqrt{5}-1}{2}, + \infty [ ^{n-1} $ de sorte à appliquer le corollaire~\ref{5cor_technique}. 
\\ On pose $\gamma_0 =1 $ et 
$$\forall j \in \llbracket 1,n-1 \rrbracket, \quad \beta_j = \gamma_j \gamma_{j-1}^{-1}.$$ 
On a $\beta_1 = \gamma_1 \geq 2 +\frac{ \sqrt{5}-1}{2} $ par hypothèse. D'après l'hypothèse $(\ref{5hyp_gamma_th1})$, on a $\gamma_{j} \geq \left(2 +\frac{ \sqrt{5}-1}{2}\right)\gamma_{j-1}$ pour tout $j \in \llbracket 2, n-1 \rrbracket$. On en déduit alors :
\begin{align*}
 \forall j \in \llbracket 1, n-1 \rrbracket, \quad \beta_j \geq 2 +\frac{ \sqrt{5}-1}{2}.
\end{align*}
Le corollaire $\ref{5cor_technique}$ donne l'existence d'une droite
$A$ de $\R^n$ vérifiant $A \in \II_n(1,n-1)_1$ et :
 \begin{align*}
 \forall e \in \llbracket 1, n-1 \rrbracket, \quad \mu_n(A|e)_1 = \max\limits_{i \in \llbracket 0, n-1-e\rrbracket} \beta_{i+1}\ldots\beta_{i+e}.
 \end{align*}
Il reste alors à montrer que $\max\limits_{i \in \llbracket 0, n-1-e\rrbracket} \beta_{i+1}\ldots\beta_{i+e}= \gamma_e$ pour tout $ e \in \llbracket 1, n-1 \rrbracket$. 
On remarque que pour $i \in \llbracket 0, n-1-e\rrbracket$ :
\begin{align*}
 \beta_{i+1}\ldots\beta_{i+e} = \frac{\gamma_{i+e}}{\gamma_{i}}.
\end{align*}
En prenant $i = 0$ on a donc $\max\limits_{i \in \llbracket 0, n-1-e\rrbracket} \beta_{i+1}\ldots\beta_{i+e} \geq \gamma_e$. 
\\ D'autre part pour $i \in \llbracket 1, n-1-e \rrbracket$, l'hypothèse $(\ref{5hyp_gamma_th2})$ donne $\gamma_{i+e} \leq \gamma_{i}\gamma_{e}$ et donc en particulier 
\begin{align*}
 \beta_{i+1}\ldots\beta_{i+e} \leq \gamma_e.
\end{align*}
Alors $\max\limits_{i \in \llbracket 0, n-1-e\rrbracket} \beta_{i+1}\ldots\beta_{i+e} \leq \gamma_e$.

\end{preuve}

Le reste du chapitre est consacré à la preuve de la propriété $\ref{5prop_technique}$. 

On se donne $(\beta_N)_{N \in \Nx} \in [ 2 +\frac{ \sqrt{5}-1}{2}, + \infty [ ^{\Nx}$ une suite périodique de période $T \in \Nx$.

\section{Construction de la droite}\label{5section_BN}

Dans cette section, on construit la droite $A$ de la propriété~\ref{5prop_technique} ainsi que des vecteurs $X_N$ pour $N \in \N$, qui sont ceux réalisant les meilleures approximations de $A$. 

On introduit la suite $\alpha = (\alpha_k)_{k \in \N}$ définie par :
\begin{align*}
 \alpha_0 &= 1, \\
 \forall k \in \N, \quad \alpha_{k+1} &= \beta_{k+1}\alpha_k.
\end{align*}

\begin{req}\label{5rq_ecart_alpha}
 Pour tout $k \in \N$ on a $\floor{\alpha_{k+1}} - \floor{\alpha_k} \geq 1$. 
\\En effet, le cas $k = 0 $ est trivial car $\beta_1 > 2$ et pour $k \in \Nx$ :
 \begin{align*}
 \floor{\alpha_{k+1 }} - \floor{\alpha_k} \geq \alpha_{k+1 } -1 - \alpha_k \geq (\beta_{k+1}-1)\alpha_k -1 \geq 1
 \end{align*}
 car $\beta_{k+1} \geq 2$ et $\alpha_k \geq 2$.
\end{req}

Soit $\theta$ un nombre premier supérieur ou égal à $5$ et $\phi : \N \to \llbracket 0, n-2 \rrbracket $ définie par :
\begin{align*}
 \phi(k) =(k \mod (n-1)) \in \llbracket 0, n-2 \rrbracket
\end{align*}
où $k \mod (n-1)$ est le reste de la division euclidienne de $k$ par $n-1$.
\\On remarque que pour tout sous-ensemble de $\N$ de la forme $\mathcal{X} = \{x, x+1, \ldots, x + n -2\}$ avec $x \in \N$ on a $\phi_{|\mathcal{X}}$ bijective.\\
D'après le lemme~\ref{3lem_sigma_alg_indep}, comme $\beta_N \geq 2 + \frac{\sqrt{5}-1}{2} > 1 $, il existe $n-1$ suites $u^0,\ldots, u^{n-2}$ vérifiant : 
\begin{align}\label{5cons_suite_u}
 \forall j \in \llbracket 0, n-2 \rrbracket, \forall k \in \N,\quad u^j_k \left\{ \begin{array}{ll}
 \in \{1,2\} &\text{ si } \phi(k) = j \\
    = 0   &\text{ sinon }
 \end{array} \right.
\end{align}
d'où $u^j_k \neq 0 $ si et seulement si $j = k \mod (n-1)$ ; 
et telles que la famille $(\sigma_0, \ldots, \sigma_{n-2})$ soit algébriquement indépendante sur $\Q$ avec :
\begin{align*}
 \forall j \in \llbracket 0, n-2 \rrbracket, \quad \sigma_j = \sigma(\theta, u^j, \alpha) =  \sum\limits_{k = 0}^{+ \infty} \frac{u^j_k}{\theta^{\floor{\alpha_k}}}.
\end{align*}
On pose alors $A = \Vect(Y)$ où :
\begin{align*}
 Y = \begin{pmatrix}
 1 \\
 \sigma_0 \\
 \vdots \\
 \sigma_{n-2}
 \end{pmatrix}. 
\end{align*}
D'après le lemme~\ref{3lem_1_irrat}, l'espace $A$ est $(e,1)\tir$irrationnel pour tout $e \in \llbracket 1, n-1 \rrbracket$. 
Dans le reste de ce chapitre, on va montrer que $\mu_n(A|e)_1 = \max\limits_{i \in \llbracket 0, T-1\rrbracket} \beta_{i+1}\ldots\beta_{i+e}$ pour $e \in \llbracket 1, n-1 \rrbracket$.

Pour cela on pose, pour $N \in \N$, le vecteur tronqué :
\begin{align}\label{5def_XN_ chap5}
 X_N = \theta^{\floor{\alpha_N}} \begin{pmatrix}
 1 \\
 \sigma_{0,N} \\
 \vdots \\
 \sigma_{n-2,N}
 \end{pmatrix} = \theta^{\floor{\alpha_N}} \begin{pmatrix}
 1 \\
 \sum\limits_{k = 0}^{N} \frac{u^0_k}{\theta^{\floor{\alpha_k}}} \\
 \vdots \\
 \sum\limits_{k = 0}^{N} \frac{u^{n-2}_k}{\theta^{\floor{\alpha_k}}}
 \end{pmatrix} \in \Z^n
\end{align}
où l'on a posé $\sigma_{j,N} = \sum\limits_{k = 0}^{N} \frac{u^{j}_k}{\theta^{\floor{\alpha_k}}} $ pour $j \in \llbracket 0, n-2 \rrbracket$. On a en particulier : 
\begin{align*}
 \left| \sigma_j - \sigma_{j,N} \right| = \left| \sum\limits_{k = N+1}^{+ \infty} \frac{u^{j}_k}{\theta^{\floor{\alpha_k}}} \right| 
 \leq 2 \sum\limits_{k = \floor{\alpha_{N+1}} }^{+ \infty} \frac{1}{\theta^{k} } \leq \frac{2 \theta}{\theta - 1 } \frac{ 1} {\theta^{\floor{\alpha_{N+1}}}} 
 \end{align*}
car $\floor{\alpha_{k+1}} - \floor{\alpha_{k}} \geq 1 $ d'après la remarque~\ref{5rq_ecart_alpha}. On a donc :
\begin{align}\label{5sigma_moins_sigmaN}
 \left| \sigma_j - \sigma_{j,N} \right| \leq\frac{ c_{\ref{5cons_sigma_moins_sigmaN}}} {\theta^{{\alpha_{N+1}}}} 
\end{align}
avec $\cons = \frac{2\theta^2}{\theta -1} \label{5cons_sigma_moins_sigmaN}$.\\
En particulier on a :
\begin{align}\label{5convergence_XN_Y}
 \theta^{- \floor{\alpha_N}} X_N\underset{N\to+\infty}{\longrightarrow} Y .
\end{align}
La preuve de la propriété~\ref{5prop_technique} se déroule en montrant que les \og meilleurs\fg \, espaces approchant $A$ sont ceux engendrés par ces vecteurs $X_N$. On pose pour cela pour $N \in \N$ et $e \in \llbracket 1, n-1 \rrbracket$ : 
\begin{align*}
 B_{N,e} = \Vect(X_N, \ldots, X_{N+e-1}).
\end{align*}

Dans la suite on fixe $e \in \llbracket 1, n-1 \rrbracket$ et on montre que $\mu_n(A|e)_1 = \max\limits_{i \in \llbracket 0, T-1 \rrbracket} \beta_{i+1}\ldots\beta_{i+e}$.

\section{L'espace \texorpdfstring{$B_{N,e}$}{}}\label{5sect_construction_B_N}
On étudie ici plus précisement l'espace $ B_{N,e} = \Vect(X_N, \ldots, X_{N+e-1})$ en donnant en particulier une $\Zbase$ de $B_{N,e} \cap \Z^n$ et la hauteur $H(B_{N,e}).$

\bigskip
Pour $N \in \N$ on remarque que :
\begin{align}\label{5rec_X_N}
 X_{N+1} = \theta^{\floor{\alpha_{N+1}} -\floor{\alpha_N}} X_N + w_{N+1}
\end{align}
avec $w_{N+1} = \begin{pmatrix}
 0 & u^0_{N+1} & \cdots & u^{n-2}_{N+1}
\end{pmatrix}^\intercal \in \Z^n$.

D'après la construction des suites $u^j$ en $(\ref{5cons_suite_u})$, le vecteur $w_{N+1}$ comporte exactement une coordonnée non nulle qui est $u^j_{N+1}$ pour $j = \phi(N+1)$. On pose alors le vecteur suivant : 
\begin{align*}
 v_{N+1} = \frac{1}{u^{\phi(N+1)}_{N+1}}w_{N+1}.
\end{align*}
On a défini des vecteurs $v_{N+1}, \ldots, v_{N+e-1}$ et on constate que ce sont des vecteurs deux à deux distincts de la base canonique, et que le premier vecteur $(1,0, \ldots, 0)^\intercal$ de cette base n'en fait pas partie. De plus on a :
\begin{align}\label{5B_N_base_v}
 B_{N,e} = \Vect(X_N, v_{N+1}, \ldots, v_{N+e-1})
\end{align}
en utilisant la relation (\ref{5rec_X_N}) et la définition de $B_{N,e}$. On en déduit que $\dim(B_{N,e}) = e$ car la première coordonnée de $X_N$ est non nulle.

\begin{lem}\label{5Z_base_B_N}
 Soit $N \in \N$. Alors 
 \begin{align*}
 X_N, v_{N+1}, \ldots, v_{N+e-1} \text{ est une } \Zbase \text{ de } B_{N,e} \cap \Z^n.
 \end{align*}
\end{lem}

\begin{preuve}
 Les vecteurs $X_N, v_{N+1}, \ldots, v_{N+e-1}$ sont entiers et d'après $(\ref{5B_N_base_v})$ ils engendrent $B_{N,e}$. \\
 Soit $(a, b_1, \ldots,b_{e-1}) \in \R^e$ tels que :
 \begin{align*}
 X = aX_N + b_1v_{N+1} + \ldots + b_{e-1}v_{N+e-1} \in B_{N,e} \cap \Z^n.
 \end{align*}
 Pour prouver le lemme, il reste à montrer que $(a, b_1, \ldots,b_{e-1}) \in \Z^e$.
\\ La première coordonnée de $X$ est $a\theta^{\floor{\alpha_N}}$, on a donc :
\begin{align*}
 a\theta^{\floor{\alpha_N}} \in \Z.
\end{align*}
On pose $i = \phi(N)$, on a donc $u^i_N \in \{ 1,2\} $ et pour tout $j \in \llbracket 1, n-2 \rrbracket $, $u^i_{N+j} = 0$.
\\On étudie la $(i+1)\tir$ième coordonnée de $X$ et on a :
\begin{align}\label{5a_theta_sigma}
 a\theta^{\floor{\alpha_N}} \sigma_{i,N} + \sum\limits_{j = 1}^{e-1} b_j \frac{u^i_{N+j}}{u^{\phi(N+j)}_{N+j} } = a\theta^{\floor{\alpha_N}} \sigma_{i,N} \in \Z.
\end{align}
Or 
\begin{align}\label{5a_theta_sigma2}
 \sigma_{i,N} = \sum\limits_{k = 0}^N \frac{u^i_k}{\theta^{\floor{\alpha_k}}} = \frac{1}{\theta^{\floor{\alpha_N}}}\sum\limits_{k = 0}^N \theta^{\floor{\alpha_N} - \floor{\alpha_k}}u^i_k = \frac{1}{\theta^{\floor{\alpha_N}}}\left(\theta U_N + u^i_{N} \right)
\end{align}
avec $U_N = \sum\limits_{k = 0}^{N-1} \theta^{\floor{\alpha_N} - \floor{\alpha_k}-1}u^i_k \in \Z$. 
\\On note $V_N = \theta U_N + u^i_{N} $. Les relations $(\ref{5a_theta_sigma})$ et $(\ref{5a_theta_sigma2})$ donnent alors $aV_N \in \Z$.
\bigskip 
\\D'autre part, comme $u^i_{N} \in \{ 1,2\} $ et $\theta > 2$, on a $\pgcd(\theta, V_N) = 1 $ en particulier :
\begin{align*}
 \pgcd(\theta^{\floor{\alpha_N}}, V_N) = 1.
\end{align*}
Par le théorème de Bézout, il existe $ p_1$ et $p_2$ deux entiers tels que $p_1 \theta^{\floor{\alpha_N}} + p_2 V_N = 1$. On a alors 
\begin{align*}
 a = a (p_1 \theta^{\floor{\alpha_N}} + p_2 V_N) = p_1 a \theta^{\floor{\alpha_N }}+ p_2 a V_N \in \Z
\end{align*}
car $a\theta^{\floor{\alpha_N}} \in \Z$ et $aV_N \in \Z$. 

En particulier $X - a X_N = b_1v_{N+1} + \ldots + b_{e-1}v_{N+e-1} \in \Z$. Or les $v_{N+j}$ sont des vecteurs distincts de la base canonique, donc $b_j \in \Z$ pour tout $j \in \llbracket 1,e-1 \rrbracket$.

On conclut que $(a, b_1, \ldots,b_{e-1}) \in \Z^e$ et donc le lemme est prouvé.

\end{preuve}

On peut maintenant calculer la hauteur de l'espace $B_{N,e}$. Pour cela il faut estimer la norme des vecteurs $X_N$.

\begin{lem}\label{5lem_norme_XN}
 Il existe des constantes $\cons \label{5cons_minor_norme_XN} > 0 $ et $\cons \label{5cons_major_norme_XN} > 0 $ indépendantes de $N$ telles que :
 \begin{align*}
 \forall N \in \N, \quad c_{\ref{5cons_minor_norme_XN}} \theta^{{\alpha_N}} \leq \| X_N \| \leq c_{\ref{5cons_major_norme_XN}} \theta^{{\alpha_N}}.
 \end{align*}
\end{lem}

\begin{preuve}
 D'après $(\ref{5convergence_XN_Y}) $ on a $ \theta^{-\floor{\alpha_N}}X_N \underset{N\to+\infty}{\longrightarrow} Y $ et donc 
\begin{align*}
 \theta^{-\floor{\alpha_N}}\|X_N\| \underset{N\to+\infty}{\longrightarrow} \|Y \|.
\end{align*}
Comme $\frac{\theta^{{\alpha_N}} }{\theta} \leq \theta^{\floor{\alpha_N}} \leq \theta^{\alpha_N}$, cette relation asymptotique donne l'existence de constantes $c_{\ref{5cons_minor_norme_XN}} > 0 $ et $c_{\ref{5cons_major_norme_XN}} > 0 $ telles que :
 \begin{align*}
 \forall N \in \N, \quad c_{\ref{5cons_minor_norme_XN}} \theta^{{\alpha_N}} \leq \| X_N \| \leq c_{\ref{5cons_major_norme_XN}} \theta^{{\alpha_N}}.
 \end{align*}
 
\end{preuve}

\begin{lem}\label{5lem_haut_BN}
 Il existe des constantes $\cons \label{5cons_minor_haut_BN} > 0 $ et $\cons \label{5cons_major_haut_BN} > 0 $ indépendantes de $N$ telles que :
 \begin{align*}
 \forall N \in \N, \quad c_{\ref{5cons_minor_haut_BN}} \theta^{{\alpha_N}} \leq H(B_{N,e}) \leq c_{\ref{5cons_major_haut_BN}} \theta^{{\alpha_N}}.
 \end{align*}
\end{lem}

\begin{preuve}
 D'après le lemme~\ref{5B_N_base_v}, les vecteurs $X_N, v_{N+1}, \ldots, v_{N+e-1}$ forment une $\Zbase$ de $B_{N,e}\cap \Z^n$. On a donc par le lemme~\ref{prop_haut_prod_ext}:
 \begin{align*}
 H(B_{N,e}) = \| X_N \wedge v_{N+1} \wedge \ldots \wedge v_{N+e-1} \|.
 \end{align*}
On en déduit :
\begin{align*}
 H(B_{N,e}) \leq \| X_N \| \| v_{N+1} \| \ldots \| v_{N+e-1} \| \leq \| X_N \|
\end{align*}
et on a la majoration du lemme en prenant $c_{\ref{5cons_major_haut_BN}} =c_{\ref{5cons_major_norme_XN}} $ la constante du lemme~\ref{5lem_norme_XN}.

D'autre part on sait que les $v_{N+j}$ pour $j \in \llbracket 1, e-1 \rrbracket$ forment une famille libre. De plus $Y \notin \Vect(v_{N+1}, \ldots, v_{N+e-1}) $ car la première coordonnée de $Y $ est $1$ et celle des $v_{N+j}$ est nulle, donc $\|Y \wedge v_{N+1} \wedge \ldots \wedge v_{N+e-1} \| > 0$. En outre, le $(e-1)\tir$uplet $(v_{N+1}, \ldots, v_{N+e-1})$ ne peut prendre qu'un nombre fini de valeurs quand $N$ varie. Il existe donc une constante $\cons \label{5cons_v_non_col} > 0 $ indépendante de $N$ telle que :
\begin{align*}
 \forall N \in \N, \quad \|Y \wedge v_{N+1} \wedge \ldots \wedge v_{N+e-1} \| \geq c_{ \ref{5cons_v_non_col}}.
\end{align*}
Comme $\theta^{-\floor{\alpha_N}}X_N \underset{N \to \infty}{\longrightarrow} Y$ d'après $(\ref{5convergence_XN_Y})$, il existe $N_0 \in N$ tel que :
\begin{align*}
 \forall N \geq N_0, \quad \| \theta^{-\floor{\alpha_N}}X_N - Y \| \leq \frac{c_{ \ref{5cons_v_non_col}}}{2}. 
\end{align*}
On minore alors $H(B_{N,e})$ pour $N \geq N_0$ :
\begin{align*}
 H(B_{N,e}) &= \| X_N \wedge v_{N+1} \wedge \ldots \wedge v_{N+e-1} \| \\
 &= \theta^{\floor{\alpha_N}} \| (\theta^{-\floor{\alpha_N}}X_N - Y + Y )\wedge v_{N+1} \wedge \ldots \wedge v_{N+e-1} \| \\
 &\geq \theta^{\floor{\alpha_N}} \left(\| Y \wedge v_{N+1} \wedge \ldots \wedge v_{N+e-1} \| - \|(\theta^{-\floor{\alpha_N}}X_N - Y) \wedge v_{N+1} \wedge \ldots \wedge v_{N+e-1} \| \right) \\
 &\geq \frac{c_{ \ref{5cons_v_non_col}}}{2\theta}\theta^{{\alpha_N}} 
\end{align*}
d'après les inégalités précédentes. 
\\On pose $\cons \label{5cons_min_N<N0} = \min\limits_{0 \leq N <N_0} \theta^{-{\alpha_N}}H(B_{N,e}) > 0$ et $c_{\ref{5cons_minor_haut_BN}} = \min(c_{\ref{5cons_min_N<N0}}, \frac{c_{ \ref{5cons_v_non_col}}}{2\theta}) > 0$. 
\\On a alors pour tout $N \in \N$
\begin{align*}
 H(B_{N,e}) \geq c_{\ref{5cons_minor_haut_BN}} \theta^{{\alpha_N}} 
\end{align*}
et la minoration du lemme est prouvée.

\end{preuve}

\section{Minoration de l'exposant}\label{5Sect_Minoration_exposant}
On montre dans cette section, que les espaces $B_{N,e}$ (définis dans la section~\ref{5sect_construction_B_N}) réalisent une bonne approximation de $A$ au premier angle. Cela permet alors de minorer l'exposant $\mu_n(A|e)_1$.

\begin{lem}\label{5lem_omega_XN_Y}
 Il existe une constante $\cons \label{5cons_omega_XN_Y} > 0$ indépendante de $N$ telle que pour tout $N \in \N$ 
 \begin{align*}
 \omega(Y,X_N) \leq c_{\ref{5cons_omega_XN_Y}}\theta^{-{\alpha_{N+1}}}.
 \end{align*}
\end{lem}

\begin{preuve}
 Soit $N \in \N$. D'après le lemme~\ref{2lem_omega_XY_diff} on a 
 \begin{align*}
 \omega(Y,X_N) = \omega(Y,\theta^{-\floor{\alpha_{N}}}X_N) \leq \frac{ \| Y- \theta^{-\floor{\alpha_{N}}}X_N \| }{\|Y \|}.
 \end{align*}
 Or $Y- \theta^{-\floor{\alpha_{N}}}X_N = \begin{pmatrix}
 0 &
 \sum\limits_{k = N+1}^{+ \infty} \frac{u^0_k}{\theta^{\floor{\alpha_k}}}
 &
 \cdots &
 \sum\limits_{k = N+1}^{+ \infty} \frac{u^{n-2}_k}{\theta^{\floor{\alpha_k}}}
 \end{pmatrix}^\intercal$ et pour tout $j \in \llbracket 0, n-2 \rrbracket $ :
 \begin{align*}
 \left| \sum\limits_{k = N+1}^{+ \infty} \frac{u^{j}_k}{\theta^{\floor{\alpha_k}}} \right| 
 = \left| \sigma_j - \sigma_{j,N} \right| \leq\frac{ c_{\ref{5cons_sigma_moins_sigmaN}}}{\theta^{{\alpha_{N+1}}}}
 \end{align*}
d'après $(\ref{5sigma_moins_sigmaN})$. On a donc :
\begin{align*}
 \| Y- \theta^{-\floor{\alpha_{N}}}X_N \| \leq \frac{ c_{\ref{5cons_sigma_moins_sigmaN}}}{\theta^{{\alpha_{N+1}}}} \sqrt{n-1}.
\end{align*}
On a donc prouvé le lemme avec $c_{\ref{5cons_omega_XN_Y}} = \frac{ \sqrt{n-1}c_{\ref{5cons_sigma_moins_sigmaN}} }{\|Y \| }$.

\end{preuve}

\begin{lem}\label{5maj_psi_A_BN}
 Il existe une constante $\cons \label{5cons_major_prox_A_BN} > 0$ indépendante de $N$ telle que pour tout $N \in \N$ :
 \begin{align*}
 \psi_1(A,B_{N,e}) \leq 
 c_{\ref{5cons_major_prox_A_BN}} H(B_{N,e})^{-\beta_{N+1}\ldots\beta_{N+e}}.
 \end{align*}
\end{lem}

\begin{preuve}
Soit $N \in \N$. Comme $X_{N+e-1} \in B_{N,e}$ et $Y \in A $ on a :
\begin{align}\label{5Y_X_n+e-1}
 \psi_1(A,B_{N,e}) \leq \omega(Y, X_{N+e-1} )\leq c_{\ref{5cons_omega_XN_Y}} \theta^{-{\alpha_{N+e}}} = c_{\ref{5cons_omega_XN_Y}}\theta^{-{\alpha_{N}}\beta_{N+1} \ldots \beta_{N+e}} 
\end{align}
 d'après le lemme~\ref{5lem_omega_XN_Y}.

On rappelle que $H(B_{N,e}) \leq c_{\ref{5cons_major_haut_BN}} \theta^{{\alpha_N}}$ d'après le lemmme \ref{5lem_haut_BN}. On a donc :
\begin{align*}
 \psi_1(A,B_{N,e}) \leq (c_{\ref{5cons_omega_XN_Y}} c_{\ref{5cons_major_haut_BN}}^{\beta_{N+1} \ldots \beta_{N+e}})H(B_{N,e})^{-\beta_{N+1} \ldots \beta_{N+e}}.
\end{align*}
On pose $c_{\ref{5cons_major_prox_A_BN}} = c_{\ref{5cons_omega_XN_Y}} c_{\ref{5cons_major_haut_BN}}^{(\max\limits_{i \in \llbracket 1, T \rrbracket} \beta_{i})^e} \geq c_{\ref{5cons_omega_XN_Y}} c_{\ref{5cons_major_haut_BN}}^{\beta_{N+1} \ldots \beta_{N+e}}$ et le lemme est prouvé.

\end{preuve}

On rappelle que les $\beta_{N+1},\ldots, \beta_{N+e}$ sont à valeurs dans $\{ \beta_1, \ldots, \beta_{T}\}$. 
On en déduit alors le corollaire suivant.
\begin{cor}\label{5cor_min}
 On a :
 \begin{align*}
 \mu_n(A|e)_1 \geq \max\limits_{i \in \llbracket 0, T-1 \rrbracket } \beta_{i+1}\ldots\beta_{i+e}.
 \end{align*}
\end{cor}

\begin{preuve}
 Comme la suite des $(\beta_k)$ est $T\tir$périodique, on a :
 \begin{align*}
 \forall i \in \llbracket 0, T \rrbracket, \forall N \in \N, \quad N \equiv i \mod (T -1) \Longrightarrow \beta_{N+1}\ldots\beta_{N+e} = \beta_{i+1}\ldots\beta_{i+e}.
 \end{align*}
 Le lemme~\ref{5maj_psi_A_BN} donne donc, pour tout $i \in \llbracket 0, T \rrbracket $, une infinité d'espaces $B$ rationnels de dimension $e$ tels que :
 \begin{align*}
 \psi_1(A,B) \leq c_{\ref{5cons_major_prox_A_BN}} H(B)^{- \beta_{i+1}\ldots\beta_{i+e}}.
 \end{align*}
 On a donc en prenant $i$ qui réalise le maximum des $\beta_{i+1}\ldots\beta_{i+e}$ :
 \begin{align*}
 \mu_n(A|e)_1 \geq \max\limits_{i \in \llbracket 0, T-1 \rrbracket } \beta_{i+1}\ldots\beta_{i+e}.
 \end{align*}

\end{preuve}

\section{Majoration de l'exposant}

On montre dans cette section, que les espaces $B_{N,e}$ (définis dans la section~\ref{5sect_construction_B_N}) réalisent en fait la meilleure approximation de $A$ au premier angle. Cela permet alors de majorer l'exposant $\mu_n(A|e)_1$.

\begin{lem}\label{5min_psi_A_BN}
 Il existe une constante $\cons \label{5cons_minor_prox_A_BN} > 0$ indépendante de $N$ telle que pour tout $N \in \N$ :
 \begin{align*}
 \psi_1(A,B_{N,e}) \geq 
 c_{\ref{5cons_minor_prox_A_BN}} \theta^{-{\alpha_{N+e}}}.
 \end{align*}
 
\end{lem}

\begin{preuve}
 On a montré dans la section~\ref{5sect_construction_B_N} que $B_{N,e} = \Vect(X_N, v_{N+1}, \ldots, v_{N+e-1})$. Comme $X_{N+e-1} - X_N \in \Vect(v_{N+1}, \ldots, v_{N+e-1})$, d'après $(\ref{5rec_X_N})$, on a : 
 \begin{align*}
 B_{N,e} = \Vect(X_{N+e-1}, v_{N+1}, \ldots, v_{N+e-1}) .
 \end{align*}
Soit $X \in B_{N,e} \smallsetminus \{0\}$. On va montrer que 
\begin{align*}
 \omega(X,Y) \geq c_{\ref{5cons_minor_prox_A_BN}} \frac{1}{\theta^{\floor{\alpha_{N+e}} }}
\end{align*} et donc en particulier $\omega(X,Y) \geq c_{\ref{5cons_minor_prox_A_BN}} \frac{1}{\theta^{{\alpha_{N+e}} }}$. Comme $ \psi_1(A,B_{N,e}) = \min\limits_{X \in B_{N,e} \smallsetminus \{0\}} \omega(X,Y)$ cela prouve le lemme.
\\On écrit $X$ dans la base $X_{N+e-1}, v_{N+1}, \ldots, v_{N+e-1}$ :
\begin{align*}
 X = a Z_{N+e-1} + b_1 v_{N+1} + \ldots + b_{e-1} v_{N+e-1}.
\end{align*}
avec $Z_{N+e-1} = \theta^{-\floor{\alpha_{N+e-1}}} X_{N+e-1}$. 
\\Quitte à renormaliser $X$, on peut supposer que :
\begin{align*}
 a^2 + \sum\limits_{j = 1}^{e-1} b_j^2 = 1.
\end{align*}
Cette hypothèse donne en particulier $\|X\| \leq \cons \label{5cons_maj_X}$ avec $c_{\ref{5cons_maj_X}}$ indépendante de $N$ car les vecteurs $v_{N+j}$ sont de norme $1$ et $Z_{N+e-1} \to Y$ d'après $(\ref{5convergence_XN_Y})$.
\\On cherche alors à minorer la quantité $\| X \wedge Y \|$. On rappelle que l'on a :
\begin{align*}
 X = \begin{pmatrix}
 a \\
 a\sigma_{0,N+e-1} +  \sum\limits_{j = 1}^{e-1} b_ju^0_{N+j} \\
 \vdots \\
 a\sigma_{n-2,N+e-1} +  \sum\limits_{j = 1}^{e-1} b_ju^{n-2}_{N+j} 
 \end{pmatrix}
 \text{ et } Y = \begin{pmatrix}
 1 \\
 \sigma_{0} \\
 \vdots \\
 \sigma_{n-2}
 \end{pmatrix}
\end{align*}
et que l'on peut minorer $\| X \wedge Y \|$ par la valeur absolue de tout mineur de taille $2$ de la matrice $(X|Y)$. On a donc :
\begin{align*}
 \forall i \in \llbracket 0, n-2 \rrbracket, \quad \| X \wedge Y \| &\geq \left| \det\begin{pmatrix} a & 1 \\
a\sigma_{i,N+e-1} +  \sum\limits_{j = 1}^{e-1} b_ju^i_{N+j} & \sigma_i \end{pmatrix} \right|\\ &= \left| a (\sigma_i - \sigma_{i,N+e-1}) -  \sum\limits_{j = 1}^{e-1} b_ju^i_{N+j} \right|.
\end{align*}
Comme $ \sigma_i - \sigma_{i,N+e-1} = \sum\limits_{k = N+e}^{+\infty} \frac{u^i_k}{\theta^{\floor{\alpha_k}}}$ on obtient finalement :
\begin{align}\label{5minor_prod_ext_X_Y}
 \forall i \in \llbracket 0, n-2 \rrbracket, \quad \| X \wedge Y \| \geq \left| a  \sum\limits_{k = N+e}^{+\infty} \frac{u^i_k}{\theta^{\floor{\alpha_k}}} -  \sum\limits_{j = 1}^{e-1} b_ju^i_{N+j} \right|.
\end{align}
On distingue maintenant deux cas.
\\ \textbullet \: \underline{Premier cas :} pour tout $j \in \llbracket 1, e-1 \rrbracket$, on a $ |b_j| < \frac{1}{\theta^{\floor{\alpha_{N+e}}}} $. On a alors :
\begin{align*}
 a^2 = 1 - \sum\limits_{j = 1}^{e-1} b_i^2 \geq 1 - \frac{e-1}{\theta^{2e}}
\end{align*}
car $\floor{\alpha_{N+e}} \geq e$ d'après la remarque~\ref{5rq_ecart_alpha}. On pose $\cons \label{5cons_cas<} = \sqrt{1 - \frac{e-1}{\theta^{2e}}} > 0 $ et donc $|a| \geq c_{\ref{5cons_cas<}}$. 
\\Soit $i = \phi(N+e) \in \llbracket 0,n-2 \rrbracket$. On a alors 
\begin{align*}
 u^i_{N+e} \geq 1 \text{ et } u^i_{N+1} = \ldots = u^i_{N+e-1} =0.
\end{align*} On utilise $(\ref{5minor_prod_ext_X_Y})$ avec ce $i$ et on trouve :
\begin{align*}
 \| X \wedge Y \| \geq \left| a  \sum\limits_{k = N+e}^{+\infty} \frac{u^i_k}{\theta^{\floor{\alpha_k}}} \right| \geq |a| \left| \frac{u^i_{N+e}}{\theta^{\floor{\alpha_{N+e}}}}\right| \geq c_{\ref{5cons_cas<}} \frac{1}{\theta^{\floor{\alpha_{N+e}}}}.
\end{align*}
\\ \textbullet \: \underline{Second cas :} il existe $j_0 \in \llbracket 1,e -1 \rrbracket$ tel que $|b_{j_0}| \geq \frac{1}{\theta^{\floor{\alpha_{N+e}}}} $.
\\ Soit $i = \phi(N+j_0)$. On a alors 
\begin{align*}
 u^{i}_{N+j_0} \geq 1 \text{ et } u^i_{N+1} = \ldots =u^i_{N+j_0-1} =u^i_{N+j_0+1} = \ldots = u^i_{N+e-1} = u^i_{N+e} = 0.
\end{align*}
En appliquant $(\ref{5minor_prod_ext_X_Y})$ avec $i$ on trouve : 
\begin{align*}
 \| X \wedge Y \| &\geq \left| a  \sum\limits_{k = N+e+1}^{+\infty} \frac{u^i_k}{\theta^{\floor{\alpha_k}}} - b_{j_0} u^i_{N+j_0}\right| \\
 &\geq |b_{j_0}| - |a|  \sum\limits_{k = N+e+1}^{+\infty} \frac{u^i_k}{\theta^{\floor{\alpha_k}}} 
 \\
 &\geq \frac{1}{\theta^{\floor{\alpha_{N+e}}}} - \frac{1}{\theta^{\floor{\alpha_{N+e+1}}}} \frac{2\theta}{\theta-1}
\end{align*}
car $|a| \leq 1 $ par l'hypothèse $a^2 + \sum\limits_{j = 1}^{e-1} b_j^2 = 1$. Or, par la remarque~\ref{5rq_ecart_alpha}
 \begin{align*}
 \frac{1}{\theta^{\floor{\alpha_{N+e+1}}}} \frac{2\theta}{\theta-1} \leq \frac{1}{\theta^{\floor{\alpha_{N+e}} +1}} \frac{2\theta}{\theta-1} \leq \frac{1}{\theta^{\floor{\alpha_{N+e}} }} \frac{2}{\theta-1} \leq \frac{1}{2\theta^{\floor{\alpha_{N+e}} }}.
 \end{align*}
Cela donne alors 
\begin{align*}
 \| X \wedge Y \| &\geq \frac{1}{\theta^{\floor{\alpha_{N+e}}}} - \frac{1}{2\theta^{\floor{\alpha_{N+e}} }} = \frac{1}{2\theta^{\floor{\alpha_{N+e}} }}. 
\end{align*}
Soit maintenant $c_{\ref{5cons_minor_prox_A_BN}} = \frac{\min(c_{\ref{5cons_cas<}}, \frac{1}{2})}{\|Y\|c_{\ref{5cons_maj_X}}}$. Comme $\|X\| \leq c_{\ref{5cons_maj_X}}$ et d'après les deux cas que l'on a étudiés on a :
\begin{align*}
 \omega(X,Y) = \frac{\| X \wedge Y\|}{\| X\| \cdot\|Y\|} \geq \frac{\min(c_{\ref{5cons_cas<}}, \frac{1}{2})}{c_{\ref{5cons_maj_X}}\cdot \|Y\|}\frac{1}{\theta^{\floor{\alpha_{N+e}} }} = c_{\ref{5cons_minor_prox_A_BN}} \frac{1}{\theta^{\floor{\alpha_{N+e}} }}
\end{align*}
ce qui termine la preuve du lemme.

\end{preuve}

\begin{req}
 Les lemmes \ref{5lem_omega_XN_Y} et \ref{5min_psi_A_BN} donnent $$\cons \theta^{-\alpha_{N+e}} \leq \psi_1(A,B_{N,e}) \leq \cons \theta^{-\alpha_{N+e}}. $$
\end{req}
\bigskip

On montre maintenant que les espaces $B_{N,e}$ sont les \og meilleurs\fg \, approchant $A$ dans le sens du lemme suivant. On note $ K_e = \max\limits_{i \in \llbracket 0, T-1 \rrbracket } \beta_{i+1}\ldots\beta_{i+e}$.

\begin{lem}\label{5lem_meilleurs_espaces}
 Soit $\varepsilon > 0 $ et $B$ un sous-espace rationnel de dimension $e$ tel que :
 \begin{align}\label{5hyp_lem_meilleurs_espaces}
 \psi_1(A,B) \leq H(B)^{-K_e - \varepsilon}.
 \end{align}
 Alors si $H(B)$ est assez grand en fonction de $\varepsilon$ et $A$, il existe $N \in \N$ tel que $B =B_{N,e}$.
\end{lem}

\begin{preuve}
Soit $N \in \N$ l'entier vérifiant :
\begin{align}\label{5choix_du_N}
 \theta^{{\alpha_{N+e-1}}} \leq H(B)^{K_e + \frac{\varepsilon}{2} -1 } < \theta^{{\alpha_{N+e}}}.
\end{align}
On montre alors que ce $N$ convient si $H(B)$ est assez grand. 
\bigskip

Soit $Z_1, \ldots, Z_e$ une $\Zbase$ de $B \cap \Z^n$. On va montrer que si $H(B)$ est assez grand, alors pour tout $i \in \llbracket 0, e-1 \rrbracket$ la quantité suivante :
 \begin{align*}
 \DD_{N+i} = \| X_{N+i} \wedge Z_1 \ldots \wedge Z_e \| 
 \end{align*}
 s'annule. D'après le lemme~\ref{lem_phi_dim1} on a : 
 \begin{align*}
 \DD_{N+i} = \omega_1(X_{N+i}, B) \| X_{N+i} \| H(B).
 \end{align*}
 Or $ \omega_1(X_{N+i}, B) \leq \omega(X_{N+i}, Y) + \omega_1(A, B)$ par l'inégalité triangulaire (lemme~\ref{2lem_ineg_triang}). De plus $\| X_{N+i} \| \leq c_{\ref{5cons_major_norme_XN}} \theta^{{\alpha_{N+i}}} $ d'après le lemme~\ref{5lem_norme_XN} et $ \omega_1(A,B) = \psi_1(A, B)$. Ainsi :
 \begin{align*}
 \DD_{N+i} \leq c_{\ref{5cons_major_norme_XN}} \theta^{{\alpha_{N+i}}}\left(\omega(X_{N+i}, Y) + \psi_1(A, B)\right) H(B).
 \end{align*}
 En utlisant le lemme~\ref{5lem_omega_XN_Y} qui donne $ \omega(X_{N+i}, Y) \leq c_{\ref{5cons_omega_XN_Y}}\theta^{-{\alpha_{N+i+1}}} $ et l'hypothèse (\ref{5hyp_lem_meilleurs_espaces}) on a 
 \begin{align*}
 \DD_{N+i} &\leq c_{\ref{5cons_major_norme_XN}}\theta^{{\alpha_{N+i}}}\left(c_{\ref{5cons_omega_XN_Y}} \theta^{-{\alpha_{N+i+1}}}+ H(B)^{-K_e- \varepsilon}\right) H(B) \\
 &\leq c_{\ref{cons_maj_DNi}}\left(\theta^{-{\alpha_{N+i+1}}+{\alpha_{N+i}}} H(B)+ \theta^{{\alpha_{N+i}}}H(B)^{-K_e - \varepsilon +1 }\right)
 \end{align*}
avec $\cons >0 \label{cons_maj_DNi}$ indépendante de $B$.
\\ D'après le choix de $N$ en $(\ref{5choix_du_N})$ on a $\theta^{{\alpha_{N+i}}} \leq \theta^{{\alpha_{N+e-1}}} \leq H(B)^{K_e+ \frac{\varepsilon}{2} -1 }$ et $ \theta \geq H(B)^{\frac{K_e+ \frac{\varepsilon}{2} -1 }{{\alpha_{N+e}}} } $. Alors pour tout $i \in \llbracket 0,e-1 \rrbracket$ : 
\begin{align}\label{5inegal_D_N_i}
 \DD_{N+i} \leq c_{\ref{cons_maj_DNi}}\left(H(B)^{\gamma_i } + H(B)^{- \frac{\varepsilon}{2} }\right)
\end{align}
en posant $\gamma_i = 1 +\frac{(-{\alpha_{N+i+1}}+{\alpha_{N+i}})(K_e + \frac{\varepsilon}{2} -1) }{{\alpha_{N+e}}}$.
\\ On remarque tout d'abord que $-{\alpha_{N+i+1}}+{\alpha_{N+i}}$ est maximal pour $i = 0$. En effet si $i \geq 1$ :
\begin{align*}
 {\alpha_{N+1}}+{\alpha_{N+i}} \leq 2\alpha_{N+i} \leq \alpha_{N+i+1} \leq {\alpha_{N+i+1}} +{\alpha_{N}}
\end{align*}
en utilisant le fait que $\beta_{N+i+1} \geq 2 $ et $\alpha_N \geq 0$.

\bigskip 
 On a donc $\gamma_i \leq \gamma_0$ pour tout $ i \in \llbracket 0, e-1 \rrbracket$ et on va montrer que $\gamma_0 \leq -c_{\ref{cons_epsilon_negatif}} \varepsilon$ avec une certaine constante $c_{\ref{cons_epsilon_negatif}} > 0 $ indépendante de $B$ qu'il reste à définir.\\
On a :
\begin{align*}
 \gamma_0 = 1 +\frac{(-{\alpha_{N+1}}+{\alpha_{N}})(K_e + \frac{\varepsilon}{2} -1) }{{\alpha_{N+e}}} &= 1 +\frac{(-\beta_{N+1} + 1)(K_e + \frac{\varepsilon}{2} -1) }{{\beta_{N+1}\ldots\beta_{N+e}}} ;
\end{align*}
enfin :
\begin{align}\label{5maj_separe_Ke_eps}
 \gamma_0 &= \frac{\beta_{N+1}\ldots\beta_{N+e} + (-\beta_{N+1} + 1)(K_e + \frac{\varepsilon}{2} -1) }{{\beta_{N+1}\ldots\beta_{N+e}}} \nonumber \\ 
 &\leq \frac{ (-\beta_{N+1}+ 1) (K_e - 1 ) + \beta_{N+1}\ldots\beta_{N+e} }{\beta_{N+1}\ldots\beta_{N+e}} - c_{\ref{cons_epsilon_negatif}} \varepsilon.
\end{align}
avec $\cons \label{cons_epsilon_negatif} = \frac{1+\sqrt{5}}{4K_e} 
 \leq \frac{(\beta_{N+1} - 1)}{2\beta_{N+1}\ldots\beta_{N+e}}$ indépendante de $B$. On majore maintenant l'autre terme :
 \begin{align*}
 (-\beta_{N+1}+ 1) (K_e - 1 ) + \beta_{N+1}\ldots\beta_{N+e} \leq (-\beta_{N+1}+ 1) (K_e - 1 ) + K_e \\
 = -K_e(\beta_{N+1} -2) + \beta_{N+1} -1 .
 \end{align*}
 Or $K_e = \max\limits_{i \in \llbracket 0,T-1 \rrbracket } \beta_{i+1}\ldots\beta_{i+e}\geq \beta_{N+1} > 2$ et donc 
 \begin{align*}
 -K_e(\beta_{N+1} -2) + \beta_{N+1} -1 &\leq - \beta_{N+1}(\beta_{N+1} -2) + \beta_{N+1} -1 \\
 &\leq -\beta_{N+1}^2 + 3\beta_{N+1} -1 \\
 &\leq 0
 \end{align*}
car $\beta_{N+1} \geq 2 +\frac{\sqrt{5}-1}{2} = \frac{3+\sqrt{5}}{2}$. 
\\On a donc en reprenant l'inégalité $(\ref{5maj_separe_Ke_eps})$, pour tout $i \in \llbracket 0,e-1 \rrbracket$ :
\begin{align*}
 \gamma_i \leq \gamma_0 \leq -c_{\ref{cons_epsilon_negatif}} \varepsilon
\end{align*}
avec $c_{\ref{cons_epsilon_negatif}} > 0$ indépendante de $B$. 
\bigskip
\\Revenons à la majoration $(\ref{5inegal_D_N_i})$. Pour tout $i \in \llbracket 0,e-1 \rrbracket$ on a :
\begin{align*}
 \DD_{N+i} &\leq c_{\ref{cons_maj_DNi}}\left(H(B)^{\gamma_i } + H(B)^{- \frac{\varepsilon}{2} }\right) 
 \leq c_{\ref{cons_maj_DNi}}\left(H(B)^{-c_{\ref{cons_epsilon_negatif}} \varepsilon} + H(B)^{- \frac{\varepsilon}{2} }\right).
\end{align*}
Le terme de droite tend vers $0$ quand $H(B) \to + \infty$. En particulier, si $H(B)$ est assez grand en fonction de $c_{\ref{cons_maj_DNi}},c_{\ref{cons_epsilon_negatif}} $ et $\varepsilon$ on a :
\begin{align*}
 \forall i \in \llbracket 0,e-1 \rrbracket, \quad \| X_{N+i} \wedge Z_1 \ldots \wedge Z_e \| = \DD_{N+i} < 1.
\end{align*}
Si tel est le cas, on a, d'après le lemme~\ref{2lem_X_in_B} : 
\begin{align*}
 \forall i \in \llbracket 0,e-1 \rrbracket, \quad X_{N+i} \in B.
\end{align*}
On rappelle que $B_{N,e} = \Vect(X_N, \ldots, X_{N+e-1})$. On a alors montré que, si $H(B)$ est assez grand, $B_{N,e} \subset B$ pour $N$ vérifiant $ \theta^{{\alpha_{N+e-1}}} \leq H(B)^{K_e + \frac{\varepsilon}{2} -1 } < \theta^{{\alpha_{N+e}}}$.
\\ Par égalité des dimensions, on a alors $B_{N,e} =B$ et le lemme est prouvé.

\end{preuve}

Le corollaire suivant termine la preuve de la propriété~\ref{5prop_technique} compte tenu de la minoration déjà obtenue au corollaire~\ref{5cor_min}.

\begin{cor}\label{5cor_maj}
 On a :
 \begin{align*}
 \mu_n(A|e)_1 \leq K_e = \max\limits_{i \in \llbracket 0,T-1 \rrbracket } \beta_{i+1}\ldots\beta_{i+e}.
 \end{align*}
\end{cor}

\begin{preuve}
 On suppose par l'absurde que $\mu_n(A|e)_1 > K_e.$ Il existe alors $\varepsilon > 0 $ tel que :
 \begin{align*}
 \mu_n(A|e)_1 \geq K_e + 2\varepsilon.
 \end{align*}
Par définition de l'exposant diophantien, il existe une infinité d'espaces rationnels $B$ de dimension $e$ vérifiant :
\begin{align}\label{5inegl_psi_1_A_B}
 0 < \psi_1(A,B) \leq H(B)^{-\mu_n(A|e)_1 + \varepsilon} \leq H(B)^{- K_e - \varepsilon}.
\end{align}
D'après le lemme~\ref{5lem_meilleurs_espaces}, si $H(B)$ est assez grand et vérifie $(\ref{5inegl_psi_1_A_B})$ alors $B = B_{N,e}$ avec $N \in \N$. Il existe alors une infinité d'entiers $N \in \N$ tels que :
\begin{align}\label{5infinite_N_psi}
 0 < \psi_1(A,B_{N,e}) \leq H(B_{N,e})^{- K_e - \varepsilon}.
\end{align}
D'autre part d'après le lemme~\ref{5min_psi_A_BN} on a :
\begin{align}\label{5min_psi_1_A_BN_new}
 \forall N \in \N, \quad \psi_1(A,B_{N,e}) \geq c_{ \ref{5cons_minor_prox_A_BN}} \theta^{- {\alpha_{N+e}}} = c_{ \ref{5cons_minor_prox_A_BN}}\theta^{- {\alpha_{N}}\beta_{N+1} \ldots \beta_{N+e}} \geq c_{ \ref{5cons_minor_prox_A_BN}}\theta^{- {\alpha_{N}}K_e}.
\end{align}
De plus d'après le lemme~\ref{5lem_haut_BN} on a $H(B_{N,e}) \geq c_{\ref{5cons_minor_haut_BN}} \theta^{{\alpha_N}}.$ En regroupant alors les inégalités $(\ref{5infinite_N_psi}) $ et $(\ref{5min_psi_1_A_BN_new})$ on trouve :
\begin{align*}
c_{\ref{5cons_minor_haut_BN}}^{K_e}H(B_{N,e}) ^{- K_e} \leq c_{\ref{5cons_minor_prox_A_BN}}^{-1}H(B_{N,e})^{- K_e - \varepsilon}.
\end{align*}
On rappelle que cette inégalité est valable pour une infinité d'entiers $N$. On a donc :
\begin{align*}
c_{\ref{5cons_minor_haut_BN}}^{K_e} c_{\ref{5cons_minor_prox_A_BN}} \leq H(B_{N,e})^{- \varepsilon}
\end{align*}
et en faisant tendre $N$ vers $+ \infty$, $H(B_{N,e}) \to + \infty$ et alors $c_{\ref{5cons_minor_haut_BN}}^{K_e} c_{\ref{5cons_minor_prox_A_BN}} = 0$. 
\\ Or $c_{\ref{5cons_minor_haut_BN}} > 0 $ et $c_{\ref{5cons_minor_prox_A_BN}} >0$ ce qui soulève une contradiction et prouve le corollaire. 

\end{preuve}

%% file: Chapitres/5_Approximation_1er_angle.tex
\chapter{Construction d'espaces avec exposants prescrits au premier angle}\label{chap6}

Dans ce chapitre, on construit des espaces avec exposants prescrits dans le prolongement du chapitre~\ref{chap5}. Ici, les espaces $A$ construits sont de dimension $d$ quelconque et on connait les exposants diophantiens pour le premier angle dans le cadre où on approxime $A$ par des espaces rationnels de dimension $e$ avec $d +e \leq n$. On fixe $n \in \Nx$.

\bigskip
On pose $C_{1} = 2 + \frac{\sqrt{5}-1}{2} $ et pour tout entier $ d \in \llbracket 2, n-1 \rrbracket$, on pose 
\begin{align}\label{6definition_Cd}
 C_{d} = 5n^2 (C_{d-1})^{2n}.
\end{align}
On va montrer dans ce chapitre, le résultat suivant.

\begin{theo}\label{6theo_principal_avec_gamma}
Soient $d \in \llbracket 1 , n-1 \rrbracket$ et $(\gamma_1, \ldots, \gamma_{n-d}) \in \R^{n-d}$ vérifiant $\gamma_1 \geq C_d $ ainsi que 
 \begin{align}
 \forall i \in \llbracket 2, n-d \rrbracket,& \quad \gamma_i \geq C_d \gamma_{i-1},\label{6hyp_gamma_th1} \\
 \forall (i,j) \in \llbracket 1, n-d-1 \rrbracket^2,& \quad i+j \leq n-d \Longrightarrow \gamma_{i+j} \leq \gamma_i \gamma_j. \label{6hyp_gamma_th2}
 \end{align}
 On peut alors construire explicitement $A \in \II_n(1,n-d)_1$ tel que :
 \begin{align*}
 \forall e \in \llbracket 1, n-d \rrbracket, \quad \mu_n(A|e)_1 = \gamma_e.
 \end{align*}
\end{theo}

\bigskip
\begin{req}
 Ce théorème se limite au cas du premier angle $\psi_1$, ici égal à $\omega_1$. Cela s'explique notamment par le fait que l'on dispose des relations suivantes :
 \begin{align*}
 \varphi(A,B) = \psi_1(A,B) \ldots \psi_{\min(\dim(A), \dim(B))}(A, B) \geq \psi_1(A,B)^{\min(\dim(A), \dim(B))}
 \end{align*}
 et \begin{align*}
 \| X \wedge Z_1 \wedge \ldots \wedge Z_e \| = \omega_1(\Vect(X),C) \| X \| H(C)
 \end{align*}
 avec $Z_1, \ldots, Z_e $ une $\Zbase$ de $C \cap \Z^n$.
 \\Ces relations particulières permettent alors les calculs de ce chapitre.
\end{req}

Le  théorème~\ref{6theo_principal_avec_gamma} se déduit directement du théorème suivant. On peut reprendre la preuve du  théorème~\ref{5theo_droite_approx} page \pageref{5preuve_du_theo} pour les détails.

\begin{theo}\label{6theo_principal}
 Soit $d \in \llbracket 1, n-1 \rrbracket $. 
 \\ Il existe une constante explicite $C_d >0$, telle que pour tout $ (\beta_1, \ldots, \beta_{n-d}) \in [C_{d}, +\infty[^{n-d}$, il existe $A \in \II_n(d,n-d)_1$ vérifiant
 \begin{align*}
 \forall e \in \llbracket 1, n-d \rrbracket, \quad \mu_n(A|e)_1 = \max\limits_{i \in \llbracket 0, n-d-e \rrbracket } \beta_{i+1}\ldots\beta_{i+e}.
 \end{align*}
 De plus la construction de $A$ est explicite.
\end{theo}

\bigskip
La preuve de ce  théorème~\ref{6theo_principal} se fait par récurrence sur $d \in \llbracket 1, n-1 \rrbracket$, le cas $d = 1$ étant en fait traité dans le chapitre~\ref{chap5}. Pour $d > 1$, on construit alors, dans la section~\ref{6section_construction_espace_A} un espace $A$ de dimension $d$ comme somme directe d'une droite $\Vect(Y_1)$ (construite comme dans le chapitre~\ref{chap5}) et d'un espace de dimension $d -1$ dont on connaît les exposants diophantiens par hypothèse de récurrence.
\\On calcule ensuite $\mu_n(A|e)_1$, pour $e \in \llbracket 1, n-d \rrbracket$. Pour cela on montre en fait que les meilleurs espaces approximant $A$ sont ceux qui sont proches de $\Vect(Y_1)$ et donc $ \mu_n(A|e)_1 = \mu_n(\Vect(Y_1)|e)_1$. 
\\En effet, $\Vect(Y_1)$ est très bien approché (avec les $\mu_n(\Vect(Y_1)|e)_1$ souhaités) tandis que l'espace de dimension $d-1$ obtenu par hypothèse de récurrence est beaucoup moins bien approché. La minoration $\mu_n(A|e)_1 \geq \mu_n(\Vect(Y_1)|e)_1$ est relativement triviale et est traitée en fin de section~\ref{6section_construction_espace_A}. Pour la majoration, on retrouve la structure de la preuve présentée dans la chapitre~\ref{chap5} : on exhibe d'abord des espaces rationnels $B_{N,e}$ approchant bien $\Vect(Y_1)$ (section~\ref{6sectionBNE}) et on montre ensuite que ce sont les \og meilleurs \fg{} par le lemme~\ref{6lem_meilleurs_espaces}. En estimant alors précisement $\psi_1(A,B_{N,e})$ dans le lemme~\ref{6min_psi_A_BN} (en considérant notamment les mineurs de taille $2$ de certaines matrices) on obtient la majoration dans le corollaire~\ref{6cor_maj_expos}. Du fait de la dimension $d > 1$, les preuves diffèrent cependant du chapitre~\ref{chap5} dans leur technicité. 

\section{Construction de l'espace \texorpdfstring{$A$}{}}
\subsection{Hypothèse de récurrence et initialisation}
Pour $d \in \llbracket 1, n-1 \rrbracket$ on pose l'hypothèse de récurrence $\HH(d)$ suivante : 
\\ \textbf{Pour tout $(\beta_1, \ldots, \beta_{n-d}) \in [C_{d}, +\infty[^{n-d}$ il existe un espace $A$ de dimension $d$ engendré par des vecteurs de la forme :}
\begin{align*}
 Y_1 = \begin{pmatrix}
 1 \\
 0 \\ 
 \vdots \\
 0 \\
 \tau_{n-d,1} \\
 \vdots \\ 
 \tau_{1,1}
 \end{pmatrix}, \quad Y_2 = \begin{pmatrix}
 1\\
 0\\ 
 \vdots \\
 \tau_{n-d +1,2} \\
 \tau_{n-d,2} \\
 \vdots \\ 
 \tau_{1,2}
 \end{pmatrix}, \ldots, \quad Y_d = \begin{pmatrix}
 1\\
 \tau_{n-1,d} \\
 \tau_{n-2,d} \\
 \vdots \\
 \\
 \vdots \\ 
 \tau_{1,d}
 \end{pmatrix}
\end{align*}
\textbf{avec $\tau_{i,j} >0 $ tels que la famille $(\tau_{i,j})$ soit algébriquement indépendante sur $\Q$, qui vérifie tout $e \in \llbracket 1, n-d \rrbracket $ : }
\begin{align*}
 \mu_n(A|e)_1 = \max\limits_{i \in \llbracket 0, n-d-e \rrbracket } \beta_{i+1}\ldots\beta_{i+e}.
\end{align*}

Le cas $d = 1$ est traité dans le chapitre~\ref{chap5}. En effet, soit $(\beta_1, \ldots,\beta_{n-1}) \in [C_{1}, +\infty [ ^{n-1} $ avec $C_1 = 2 + \frac{\sqrt{5}-1}{2}$. Le corollaire~\ref{5cor_technique} (démontré au chapitre~\ref{chap5}) fournit un vecteur $Y_1 = \begin{pmatrix}1 &\sigma_0 & \cdots & \sigma_{n-2} \end{pmatrix}^\intercal $ avec $\sigma_0, \ldots, \sigma_{n-2}$ algébriquement indépendants sur $\Q$ tel que
\begin{align*}
 \forall e \in \llbracket 1, n-1 \rrbracket, \quad \mu_n(\Vect(Y_1)|e)_1 = \max\limits_{i \in \llbracket 0, n-1-e \rrbracket } \beta_{i+1}\ldots\beta_{i+e}.
\end{align*}

\subsection{Espace \texorpdfstring{$A$}{} et minoration de l'exposant}\label{6section_construction_espace_A}
Dans cette section, on construit les vecteurs $Y_1, \ldots, Y_d$ et $A = \Vect(Y_1, \ldots, Y_d)$ en utilisant l'hypothèse de récurrence.
\bigskip 

Le cas $d =1$ étant traité, soit $d \in \llbracket 2, n -1 \rrbracket$. On suppose que $\HH(d-1)$ est vraie. 
\\On l'applique avec $\beta'_1 = \ldots = \beta'_{n-d+1} = C_{d-1}$. Il existe alors des vecteurs $Y_2, \ldots, Y_d \in \R^n$ de la forme 
\begin{align*}
 Y_2 = \begin{pmatrix}
 1\\
 0\\ 
 \vdots \\
 \tau_{n-d +1,2} \\
 \vdots \\ 
 \tau_{1,2}
 \end{pmatrix}, \quad Y_3 = \begin{pmatrix}
 1 \\
 \vdots \\
 0 \\
 \tau_{n-d+2,3} \\
 \tau_{n-d+1,3} \\
 \vdots \\ 
 \tau_{1,3}
 \end{pmatrix}, \ldots, \quad
 Y_d = \begin{pmatrix}
 1\\
 \tau_{n-1,d} \\
 \tau_{n-2,d} \\
 \vdots \\
 \\
 \vdots \\ 
 \tau_{1,d}
 \end{pmatrix}
\end{align*}
avec $\tau_{i,j} >0 $ et tels que la famille $(\tau_{i,j})$ soit algébriquement indépendante sur $\Q$. L'espace $A' = \Vect(Y_2, \ldots, Y_d)$ vérifie :
\begin{align}\label{6hyp_rec}
 \forall e \in \llbracket 1, n-d+1 \rrbracket \quad \mu_n(A'|e)_1 = (C_{d-1})^e.
\end{align}
Dans la suite, on note $\lambda_e = (C_{d-1})^e $ pour $e \in \llbracket 1, n-d+1 \rrbracket$.

Soit $ (\beta_1, \ldots, \beta_{n-d}) \in [C_{d}, +\infty[^{n-d}$. On construit maintenant $Y_1 \in \R^n$ comme dans le chapitre~\ref{chap5}.
\\On étend la suite $\beta = (\beta_k)_{k \in \Nx}$ par :
\begin{align*}
 \forall i \in \llbracket n-d, 2n-2d \rrbracket, \quad \beta_i &= C_d,
 \end{align*}
puis par périodicité :
\begin{align*}
 \forall i \in \llbracket 1, 2n-2d \rrbracket, \forall k \in \N, \quad \beta_{i + k(2n-2d)} = \beta_i.
\end{align*}
La suite $(\beta_k)_{k \in \Nx}$ ne prend alors ses valeurs que dans $\{ \beta_1, \ldots, \beta_{n-d}, C_d \}$. \\

On introduit aussi la suite $\alpha = (\alpha_k)_{k \in \N}$ définie par :
\begin{align*}
 \alpha_0 &= 1, \\
 \forall k \in \N, \alpha_{k+1} &= \beta_{k+1}\alpha_k.
\end{align*}

 Soit $\theta$ un nombre premier supérieur ou égal à $5$ et $\phi : \N \to \llbracket 0, n-d-1 \rrbracket $ définie par :
\begin{align*}
 \phi(k) =(k \mod (n-d)) \in \llbracket 0, n-d-1 \rrbracket.
\end{align*}

D'après le lemme~\ref{3lem_sigma_alg_indep} il existe $n-d$ suites $u^0,\ldots, u^{n-d-1}$ vérifiant : 
\begin{align}\label{6cons_suite_u}
 \forall j \in \llbracket 0, n-d-1 \rrbracket, \forall k \in \N, \quad u^j_k \left\{ \begin{array}{ll}
 \in \{1, 2 \} &\text{ si } \phi(k) = j \\
 = 0 &\text{ sinon }
 \end{array} \right. 
\end{align}
d'où $u^j_k \neq 0 $ si et seulement $j = k \mod (n-d-1)$ ; et telles que la famille $(\sigma_0, \ldots, \sigma_{n-d-1})$ soit algébriquement indépendante sur $\Q(\FF )$ où $\FF =(\tau_{i,j})$ et où on a noté :
\begin{align*}
 \forall j \in \llbracket 0, n-d-1 \rrbracket, \quad \sigma_j = \sigma(\theta, u^j, \alpha) =  \sum\limits_{k = 0}^{+ \infty} \frac{u^j_k}{\theta^{\floor{\alpha_k}}}.
\end{align*}
On pose alors $Y_1 = \begin{pmatrix}
 1 & 0 & \cdots & 0 &
 \sigma_0 &
 \cdots &
 \sigma_{n-d-1}
 \end{pmatrix}^\intercal $.
 La preuve du corollaire~\ref{5cor_technique} et notamment la construction de la section~\ref{5section_BN} donnent alors 
 \begin{align}\label{6exposant_Y1}
 \forall e \in \llbracket 1,n-d \rrbracket, \quad \mu_n(\Vect(Y_1)|e)_1 = \max\limits_{i \in \llbracket 0, n-d-e \rrbracket } \beta_{i+1}\ldots\beta_{i+e}
 \end{align}
 en plongeant $Y_1$ dans $\R \times \{0\}^{d-1} \times \R^{n-d} $ et en utilisant la propriété~\ref{2prop_4.5Elio} avec $k = n-d +1 $. 

 \bigskip
 
 On pose alors $A = \Vect(Y_1, \ldots, Y_d)$ et on note $K_e $ les quantités 
 \begin{align*}
 K_e = \max\limits_{i \in \llbracket 0, n-d-e \rrbracket } \beta_{i+1}\ldots\beta_{i+e}
 \end{align*} pour $e \in \llbracket 1, n-d \rrbracket$. 

\begin{lem}
 L'espace $A$ est $(n-d,1)\tir$irrationnel. 
\end{lem}

\begin{preuve}
 On pose $\AA = \begin{pmatrix}
 Y_1 & \cdots & Y_d
 \end{pmatrix}$ la matrice dont les colonnes sont les $Y_j$. On peut écrire $\AA = \begin{pmatrix}
 G \\ \Sigma
 \end{pmatrix}$ avec :
 \begin{align*}
 G = \begin{pmatrix}
 1 & 1 & \cdots & \cdots & 1 \\
 0 & & & 0& \tau_{n-1,d} \\
 \vdots & \cdots & 0 & \tau_{n-2,d-1} & \tau_{n-2,d} \\
 \vdots & \udots & & & \vdots \\
 0 & \tau_{n-d+1, 2} & \cdots & \cdots & \tau_{n-d+1,d}
 \end{pmatrix}
 \text{ et } \Sigma = \begin{pmatrix}
 \sigma_0 & \tau_{n-d,2} & \cdots & \tau_{n-d,d} \\
 & &\vdots & \\
 \sigma_{n-d-1} & \tau_{1,2} & \cdots & \tau_{1,d}
 \end{pmatrix}.
 \end{align*}
 On a $\Sigma \in \mathcal{M}_{n-d,d}(\R )$ et $G \in \GL_{d}(\R)$ car $\det(G) = \pm  \prod\limits_{i=1}^{d-1} \tau_{n-d+i, i+1} \neq 0 $ en développant par rapport à la première colonne. 
 \\De plus, par construction, la famille $\{ \sigma_k \} \cup \{ \tau_{i,j}\}$ est algébriquement indépendante sur $\Q$. En particulier on remarque que les coefficients de $\Sigma$ forment une famille algébriquement indépendante sur $\Q(\FF)$ où $\FF$ est la famille des coefficients de $G$.
 \\Alors l'espace $A$ engendré par les vecteurs $Y_1, \ldots, Y_d$ est $(n-d,1)\tir$irrationnel par le lemme~\ref{3lem_1_irrat}.
 
\end{preuve}

\bigskip

On va maintenant calculer $\mu_n(A|e)_1$ pour $e \in \llbracket 1, n-d \rrbracket$. On fixe $e \in \llbracket 1, n-d \rrbracket$.


En utilisant $(\ref{6exposant_Y1})$ on a 
\begin{align*}
 \mu_n(\Vect(Y_1)|e)_1 = K_e.
\end{align*}
On remarque que $\Vect(Y_1) \subset A$ et on applique alors le corollaire $\ref{2cor_croissance_exposants_inclusion}$ en remarquant que $g(1,e,n) = g(d,e,n) = 0$, et on a alors
\begin{align}\label{6minoration_exposant_par_Ke}
 \mu_n(A|e)_1 \geq \mu_n(\Vect(Y_1)|e)_1 = K_e.
\end{align}

Il reste alors à majorer l'exposant $\mu_n(A|e)_1$ par $K_e$ et c'est l'objet du reste du chapitre.
\section{L'espace \texorpdfstring{$B_{N,e}$}{}}\label{6sectionBNE}
Dans ce chapitre, on montre que le vecteur de $A$ réalisant les meilleures approximations par des espaces rationnels est le vecteur $Y_1$. On se ramène alors en quelque sorte au cas $d= 1$ traité dans le chapitre~\ref{chap5}, le vecteur $Y_1$ étant simplement le vecteur $Y$ étudié au chapitre~\ref{chap5} auquel on a ajouté des coordonnées nulles. 
\\Les \frquote{meilleurs} espaces sont donc similaires à ceux étudiés dans la section~\ref{5sect_construction_B_N} du chapitre~\ref{chap5}, c'est pourquoi on les note de la même façon. On rappelle ici leur construction et les propriétés associés.

\bigskip
On pose, pour $N \in \N$, le vecteur \frquote{tronqué} :
\begin{align*}
 X_N = \theta^{\floor{\alpha_N}} \begin{pmatrix}
 1 \\
 0 \\
 \vdots \\
 0 \\
 \sigma_{0,N} \\
 \vdots \\
 \sigma_{n-d-1,N}
 \end{pmatrix} \in \Z^n
\end{align*}
où l'on a posé $\sigma_{j,N} = \sum\limits_{k = 0}^{N} \frac{u^{j}_k}{\theta^{\floor{\alpha_k}}} \in \frac{1}{ \theta^{\floor{\alpha_N}}} \Z $ pour $j \in \llbracket 0, n-d-1 \rrbracket$. 
\\Par la construction des $X_N$ on a :
\begin{align*}
 \theta^{-\floor{\alpha_N}}X_N \underset{N \to + \infty}{\longrightarrow} Y_1
\end{align*}
et donc 
\begin{align}\label{6norme_X_N}
 c_{\ref{6con_min_norme_XN}} \theta^{{\alpha_N}} \leq \| X_N \| \leq c_{\ref{6con_maj_norme_XN}}\theta^{{\alpha_N}}
\end{align} 
pour tout $N \in \N$ avec $ \cons \label{6con_min_norme_XN} >0 $ et $\cons \label{6con_maj_norme_XN}>0$ indépendantes de $N$, comme énoncé dans le lemme~\ref{5lem_norme_XN}. Enfin on a 
\begin{align}\label{6X_N_diff_tronque}
 \| Y_1 - \theta^{-\floor{\alpha_N}}X_N \| \leq c_{\ref{6cons_diff_tronque}} \theta^{-\alpha_{N+1}}
\end{align}
avec $\cons > 0 \label{6cons_diff_tronque}$ indépendante de $N$.

\bigskip
Soit alors, pour $N \in \N$, $B_{N,e }= \Vect(X_N, X_{N+1}, \ldots, X_{N+e-1})$.
\begin{req}\label{6req_section_BN}
 On a $B_{N,e }= \Vect(X_N, v_{N+1}, \ldots, v_{N+e-1})$ où les $v_j$ sont des vecteurs distincts de la base canonique (voir section~\ref{5sect_construction_B_N}).
De plus on a montré dans le lemme~\ref{5Z_base_B_N} que ces vecteurs forment une $\Zbase$ de $B_{N,e} \cap \Z^n$ et 
\begin{align}\label{6hauteur_B_N}
 c_{\ref{6cons_minor_haut_BN}} \theta^{{\alpha_N}} \leq H(B_{N,e}) \leq c_{\ref{6cons_major_haut_BN}} \theta^{{\alpha_N}}
\end{align}
grâce au lemme~\ref{5lem_haut_BN}, avec $\cons \label{6cons_minor_haut_BN} $ et $\cons \label{6cons_major_haut_BN}$ indépendante de $N$. 

\end{req}

\section{Majoration de l'exposant}
C'est dans cette section que la preuve diffère du chapitre~\ref{chap5} et où on utilise l'hypothèse de récurrence.

\subsection{Espaces de meilleures approximations}
On montre ici que les espaces $B_{N,e}$ sont ceux qui réalisent les meilleurs approximations de $A$ au premier angle. On rappelle que l'entier $e \in \llbracket 1, n-d \rrbracket$ est fixé.

\begin{lem}\label{6lem_meilleurs_espaces}
 Soit $\varepsilon > 0 $ et $B$ un sous-espace rationnel de dimension $e$ tel que :
 \begin{align}\label{6hyp_lem_meilleurs_espaces}
 \psi_1(A,B) \leq H(B)^{-K_e- \varepsilon}.
 \end{align}
 Alors si $H(B)$ est assez grand en fonction de $\varepsilon$ et $A$, il existe $N \in \N$ tel que $B =B_{N,e}$.
\end{lem}

Soit $Z_1, \ldots, Z_e$ une $\Zbase$ de $B \cap \Z^n$. On reprend le schéma de la preuve du lemme~\ref{5lem_meilleurs_espaces} en montrant que si $H(B)$ est assez grand, alors pour tout $i \in \llbracket 0, e-1 \rrbracket$ la quantité suivante :
 \begin{align*}
 \DD_{N+i} = \| X_{N+i} \wedge Z_1 \ldots \wedge Z_e \| 
 \end{align*}
 s'annule pour un certain $N \in \N $ à déterminer.

On montre dans ce but, deux lemmes préalables, en rappellant la notation $\lambda_e = (C_{d-1})^e$.

\begin{lem}\label{6lem_maj_Y1_Y_d_Z1_Ze}
 Soit $\varepsilon >0$. On suppose que $B$ est un espace rationnel de dimension $e$ tel que $\psi_1(A,B) \leq H(B)^{-\lambda_e -\varepsilon}$. 
 \\Alors si $H(B)$ est assez grand en fonction de $A$ et de $\varepsilon$, on a 
 \begin{align*}
 \forall \delta \in \left]0, \frac{\varepsilon}{2}\right[, \quad \| Y_1 \wedge \ldots \wedge Y_d \wedge Z_1 \wedge \ldots \wedge Z_e \| 
 &\leq c_{\ref{6cons_major_Y_Z}} \psi_1(A,B) H(B)^{1+\lambda_e + \delta} 
 \end{align*}
 avec $\cons \label{6cons_major_Y_Z} >0 $ indépendante de $B$ mais pouvant dépendre de $\delta$.
\end{lem}

\begin{req}
 Le lemme~\ref{lem_phi_wedge} donne $$\| Y_1 \wedge \ldots \wedge Y_d \wedge Z_1 \wedge \ldots \wedge Z_e \| = \varphi(A,B) H(B)\| Y_1 \wedge \ldots \wedge Y_d\|. $$
 Ici l'intérêt est de faire apparaître le premier angle $\psi_1(A,B)$.
\end{req}
\begin{preuve}
 Soit $Y =  \sum\limits_{j=1}^d a_j Y_j \in A$ de norme $1$ tel que $\psi_1(A,B) = \omega(Y,p_B(Y))$. \\
On utilise la relation suivante 
\begin{align}\label{6relation_a_1_X}
 \omega  \bigg( \sum\limits_{j = 2}^d a_j Y_j, p_B(Y) ) - \omega(Y, p_B(Y) \bigg) \leq \omega \bigg(Y,  \sum\limits_{j = 2}^d a_j Y_j \bigg)
\end{align}
 provenant de l'inégalité triangulaire (lemme~\ref{2lem_ineg_triang}). Or 
\begin{align}\label{6maj_angle_YajYj}
 \omega \bigg(Y,  \sum\limits_{j = 2}^d a_j Y_j \bigg) = \frac{\left\|  \sum\limits_{j = 1}^d a_j Y_j \wedge  \sum\limits_{j = 2}^d a_j Y_j \right\| }{\left\|  \sum\limits_{j = 2}^d a_j Y_j \right\| } = \frac{\left\| a_1 Y_1 \wedge  \sum\limits_{j = 2}^d a_j Y_j \right\| }{\left\|  \sum\limits_{j = 2}^d a_j Y_j \right\| } \leq \left\| a_1 Y_1 \right\|.
\end{align}
D'autre part, soit $\delta \in \left] 0, \frac{\varepsilon}{2} \right[$. Comme $ \sum\limits_{j = 2}^d a_j Y_j \in A' = \Vect(Y_2, \ldots, Y_d)$ et à condition que $ \sum\limits_{j = 2}^d a_j Y_j \neq 0 $, l'hypothèse de récurrence $(\ref{6hyp_rec})$ donne : 
\begin{align*}
 \omega \bigg( \sum\limits_{j = 2}^d a_j Y_j, p_B(Y) \bigg) \geq \psi_1(A',B ) \geq c_{\ref{6cons_Aprime}} H(B)^{-\mu_n(A'|e)_1-\delta} = c_{\ref{6cons_Aprime}} H(B)^{-\lambda_e-\delta}
\end{align*}
avec $\cons > 0 \label{6cons_Aprime}$ ne dépendant que de $A'$ et $\delta$.
\\ En reprenant $(\ref{6relation_a_1_X})$ et $(\ref{6maj_angle_YajYj})$ on a donc :
\begin{align*}
 c_{\ref{6cons_Aprime}} H(B)^{-\lambda_e-\delta} - \psi_1(A,B) \leq |a_1| \|Y_1\|. 
\end{align*}
On rappelle que $\delta < \frac{\varepsilon}{2}$ et que par hypothèse $\psi_1(A,B) \leq H(B)^{-\lambda_e-\varepsilon}$. 
\\Comme $H(B)$ est assez grand on a $H(B)^{-\frac{\varepsilon}{2}} \leq \frac{c_{\ref{6cons_Aprime}}}{2}$ d'où :
\begin{align*}
 \frac{c_{\ref{6cons_Aprime}}}{2} H(B)^{-\lambda_e-\delta} &\leq H(B)^{-\lambda_e-\delta}\left(c_{\ref{6cons_Aprime}} - H(B)^{-\frac{\varepsilon}{2}}\right) \\
 &\leq c_{\ref{6cons_Aprime}} H(B)^{-\lambda_e-\delta} - H(B)^{ -\lambda_e-\varepsilon} \\
 &\leq |a_1| \|Y_1\|
\end{align*}
et donc $ c_{\ref{6cons_minor_a_1}}H(B)^{-\lambda_e-\delta} \leq |a_1|$ avec $\cons \label{6cons_minor_a_1} >0 $ indépendante de $B$ mais dépendante de $\delta$. De plus, cette inégalité est toujours vraie si $ \sum\limits_{j = 2}^d a_j Y_j = 0 $, car dans ce cas $a_1 = \|Y_1 \|^{-1}$. \\
\bigskip 

En particulier $a_1 \neq 0$. On pose $D_{Y,Z} = \| Y_1 \wedge \ldots \wedge Y_d \wedge Z_1 \wedge \ldots \wedge Z_e \|$ et on calcule alors :
\begin{align*}
D_{Y,Z} &= \frac{1}{|a_1|}\| a_1Y_1 \wedge Y_2 \wedge \ldots \wedge Y_d \wedge Z_1 \wedge \ldots \wedge Z_e \| \\
 &= \frac{1}{|a_1|}\bigg\| \bigg(a_1Y_1 +  \sum\limits_{j = 2}^d a_j Y_j - p_B(X)\bigg) \wedge Y_2 \wedge \ldots \wedge Y_d \wedge Z_1 \wedge \ldots \wedge Z_e \bigg\| \\
 &= \frac{1}{|a_1|}\| \left(Y - p_B(Y)\right) \wedge Y_2 \wedge \ldots \wedge Y_d \wedge Z_1 \wedge \ldots \wedge Z_e \|
\end{align*}
car $ \sum\limits_{j = 2}^d a_j Y_j - p_B(X) \in \Vect(Y_2, \ldots, Y_d, Z_1, \ldots, Z_e)$. On peut alors majorer : 
\begin{align*}
 D_{Y,Z} &\leq \frac{\| Y - p_B(Y)\| \cdot \| Y_2 \wedge \ldots \wedge Y_d \| \cdot \| Z_1 \wedge \ldots \wedge Z_e \| }{|a_1|} \\
 &\leq \frac{\psi_1(A,B) \| Y_2 \wedge \ldots \wedge Y_d \| H(B) }{ c_{\ref{6cons_minor_a_1}}H(B)^{-\lambda_e-\delta}} 
\end{align*}
car $Z_1, \ldots, Z_e$ est une $\Zbase$ de $B \cap \Z^n$ (propriété~\ref{prop_haut_prod_ext}).
\\On a donc :
\begin{align*}
 \forall \delta \in \left] 0, \frac{\varepsilon}{2} \right[, \quad \| Y_1 \wedge \ldots \wedge Y_d \wedge Z_1 \wedge \ldots \wedge Z_e \| 
 &\leq c_{\ref{6cons_major_Y_Z}} \psi_1(A,B) H(B)^{1+\lambda_e+\delta} 
\end{align*}avec $c_{ \ref{6cons_major_Y_Z}} = \| Y_2 \wedge \ldots \wedge Y_d \| c_{\ref{6cons_minor_a_1}}^{-1}$.

\end{preuve}

\begin{lem}\label{6lem_maj_DN}
 Soit $\varepsilon >0$. On suppose que $B$ est un espace rationnel de dimension $e$ tel que $\psi_1(A,B) \leq H(B)^{-\lambda_e - \varepsilon}$. On rappelle que $Z_1, \ldots, Z_e$ est une $\Zbase$ de $B \cap \Z^n$.
 \\Pour $N$ un entier, on pose $\DD_N =\| X_{N} \wedge Z_1 \ldots \wedge Z_e \| $.
 \\Alors si $H(B)$ est assez grand en fonction de $A$ et de $\varepsilon$, on a pour tout $N \in \N$ :
 \begin{align*}
 \forall \delta \in \left]0, \frac{\varepsilon}{2}\right[, \quad \DD_N \leq c_{\ref{6cons_maj_D_N_lem}}\theta^{{\alpha_{N}} (1 + t\lambda_{e+1}+t\delta) } H(B)^{1+ t\lambda_{e+1}+t\delta} \left(\psi_1(A,B) H(B)^{ \lambda_e+\delta } + \frac{1 }{\theta^{{\alpha_{N+1}}}} \right)
 \end{align*}
 en posant $t = \min(d-1,e+1)$, et avec $\cons \label{6cons_maj_D_N_lem} >0$ indépendante de $N$ et de $B$ mais dépendante de $\delta$ .
\end{lem}

\begin{preuve}
On pose $\EE_N = \| X_N \wedge Y_2 \wedge \ldots \wedge Y_d \wedge Z_1 \ldots \wedge Z_e \|$. 
\\On utilise le lemme~\ref{lem_phi_wedge} et on a :
\begin{align}\label{6lien_phi_EEN}
 \EE_N &= \| Y_2 \wedge \ldots \wedge Y_d \wedge X_N \wedge Z_1 \ldots \wedge Z_e \| \nonumber \\
 &=\varphi(A', C_N ) \| Y_2 \wedge \ldots \wedge Y_d \| \cdot \| X_N \wedge Z_1 \ldots \wedge Z_e \|
\end{align}
avec $C_N = \Vect(X_N,Z_1 \ldots, Z_e)$ et $A' = \Vect(Y_2, \ldots, Y_d)$.
\bigskip
\\On rappelle que $\varphi(A', C_N ) = \psi_1(A',C_N) \ldots \psi_u(A',C_N) \geq \psi_1(A',C_N) ^u $ avec $u = \min(d-1, \dim(C_N))$. Comme $\dim(C_N) \in \{e,e+1\}$ on a $u \leq \min(d-1,e+1) = t $ et donc 
\begin{align*}
 \varphi(A', C_N ) \geq \psi_1(A',C_N)^t
\end{align*}
car $\psi_1(A',C_N)\leq 1$. L'inéquation $(\ref{6lien_phi_EEN})$ donne alors :
\begin{align}\label{6maj_DN_1}
 \DD_N = \| X_N \wedge Z_1 \ldots \wedge Z_e \| \leq \frac{c_{\ref{6cons_ext_Y2_Yd}}\EE_N}{\psi_1(A',C_N)^t}
\end{align}
avec $\cons = \| Y_2 \wedge \ldots \wedge Y_d \|^{-1} \label{6cons_ext_Y2_Yd}$ une constante indépendante de $N$. Enfin comme $C_N$ est un espace rationnel de dimension $f \in \{e, e+1\}$, par l'hypothèse de récurrence $(\ref{6hyp_rec})$ on a :
\begin{align*}
 \forall \delta >0, \quad \psi_1(A',C_N) \geq c_{\ref{6cons_minor_angle_A'_CN}}H(C_N)^{-\lambda_f - \delta} \geq c_{\ref{6cons_minor_angle_A'_CN}}H(C_N)^{-\lambda_{e+1} - \delta}
\end{align*}
avec $\cons \label{6cons_minor_angle_A'_CN} >0$ indépendante de $N$ et dépendante de $\delta$. Dans la suite on fixe $\delta \in \left] 0, \frac{\varepsilon}{2} \right[$.
\\Comme $X_N, Z_1, \ldots, Z_e$ sont des vecteurs entiers on a d'après la remarque~\ref{2req_maj_hauteur_base_entiere} :
\begin{align*}
 H(C_N) \leq\|X_N\| \cdot \| Z_1 \wedge \ldots \wedge Z_e \| \leq c_{ \ref{6con_maj_norme_XN}} \theta^{\alpha_N} H(B)
\end{align*}
avec $c_{ \ref{6con_maj_norme_XN}} $ indépendante de $N$ provenant de $(\ref{6norme_X_N})$. L'inégalité $(\ref{6maj_DN_1})$ devient alors :
\begin{align}\label{6maj_DN_2}
 \DD_N \leq c_{\ref{6cons_ext_Y2_Yd}} c_{\ref{6cons_minor_angle_A'_CN}}^{-t} c_{ \ref{6con_maj_norme_XN}}^{t\lambda_{e+1} + t\delta }\theta^{\alpha_N(t\lambda_{e+1} +t\delta) } H(B)^{t\lambda_{e+1} + t\delta}\EE_N.
\end{align}
D'autre part, on majore $\EE_N$. Pour $N \in \N$, on pose $Z_N = \theta^{-\floor{\alpha_N}} X_N$ et $W_N = Y_1 - Z_N$. On a alors 
\begin{align*}
 \EE_N&= \theta^{\floor{\alpha_N}}\| Z_N\wedge Y_2 \wedge \ldots \wedge Y_d \wedge Z_1 \ldots \wedge Z_e \| \\
 &= \theta^{\floor{\alpha_N}} \| (Y_1- W_N) \wedge Y_2 \wedge \ldots \wedge Y_d \wedge Z_1 \ldots \wedge Z_e \| \\
 &\leq\theta^{{\alpha_N}}\left(\| Y_1 \wedge Y_2 \wedge \ldots \wedge Y_d \wedge Z_1 \ldots \wedge Z_e \| + \| W_N \wedge Y_2 \wedge \ldots \wedge Y_d \wedge Z_1 \ldots \wedge Z_e \|\right).
\end{align*}
On étudie les deux termes séparement. Premièrement :
\begin{align*}
 \| Y_1 \wedge Y_2 \wedge \ldots \wedge Y_d \wedge Z_1 \ldots \wedge Z_e \| \leq c_{\ref{6cons_major_Y_Z}} \psi_1(A,B) H(B)^{1+\lambda_e + \delta} 
\end{align*}
d'après le lemme~\ref{6lem_maj_Y1_Y_d_Z1_Ze}. Deuxièmement, d'après $(\ref{6X_N_diff_tronque})$ on a $\|W_N\| = \| Y_1 -Z_N \| \leq c_{\ref{6cons_diff_tronque}} \theta^{-\alpha_{N+1}}$ donc
\begin{align*}
 \| W_N \wedge Y_2 \wedge \ldots \wedge Y_d \wedge Z_1 \ldots \wedge Z_e \|&\leq \| W_N \| \cdot \| Y_2 \wedge \ldots \wedge Y_d \| \cdot \| Z_1 \ldots \wedge Z_e \| 
 \\ &\leq c_{\ref{6cons_diff_tronque}} \theta^{-\alpha_{N+1}} c_{\ref{6cons_ext_Y2_Yd2}}H(B)
\end{align*}
avec $ \cons \label{6cons_ext_Y2_Yd2} =\| Y_2 \wedge \ldots \wedge Y_d \| $. Ces deux inégalités permettent de majorer $\EE_N$ :
\begin{align*}
 \EE_N \leq c_{\ref{6cons_pour_pas_que_ca_deborde}}\left(\psi_1(A,B) H(B)^{ \lambda_e + \delta } + \frac{1 }{\theta^{{\alpha_{N+1}}}} \right)
\end{align*}
avec $\cons = \max(c_{\ref{6cons_major_Y_Z}},c_{\ref{6cons_diff_tronque}} c_{\ref{6cons_ext_Y2_Yd2}}) \label{6cons_pour_pas_que_ca_deborde}$
et donc en reprenant $(\ref{6maj_DN_2})$ :
\begin{align*}
 \DD_N &\leq c_{\ref{6cons_ext_Y2_Yd}} c_{\ref{6cons_minor_angle_A'_CN}}^{-t} c_{ \ref{6con_maj_norme_XN}}^{t\lambda_{e+1} + t\delta }\theta^{\alpha_N(t\lambda_{e+1} +t\delta) } H(B)^{t\lambda_{e+1} + t\delta}\theta^{{\alpha_N}} c_{\ref{6cons_pour_pas_que_ca_deborde}}\left(\psi_1(A,B) H(B)^{ \lambda_e + \delta } + \frac{1 }{\theta^{{\alpha_{N+1}}}} \right) \\
 &\leq c_{\ref{6cons_maj_D_N_lem}} \theta^{{\alpha_{N}} (1 + t\lambda_{e+1}+t\delta) } H(B)^{1+ t\lambda_{e+1}+t\delta} \left(\psi_1(A,B) H(B)^{ \lambda_e + \delta } + \frac{1 }{\theta^{{\alpha_{N+1}}}} \right)
\end{align*}
avec $c_{\ref{6cons_maj_D_N_lem}} = \ c_{\ref{6cons_ext_Y2_Yd}} c_{\ref{6cons_minor_angle_A'_CN}}^{-t} c_{ \ref{6con_maj_norme_XN}}^{t\lambda_{e+1} } c_{ \ref{6con_maj_norme_XN}}^{t \delta}c_{\ref{6cons_pour_pas_que_ca_deborde}}$. 

\end{preuve}

 On a maintenant tous les outils pour prouver le lemme~\ref{6lem_meilleurs_espaces}. 

 \begin{preuve}[ (lemme~\ref{6lem_meilleurs_espaces})] 
 On pose 
 \begin{align}\label{6choix_delta}
 \delta = \min\left(\frac{\varepsilon}{4(t+1)}, \frac{n(C_{d-1})^n -1 }{t}\right).
 \end{align}
 Soit $N \in \N $ l'entier verifiant :
 \begin{align}\label{6choix_N}
 \theta^{\alpha_{N+e-1}(1+t \lambda_{e+1} + t\delta)} \leq H(B)^{K_e - 1-t\lambda_{e+1} - \lambda_e + \frac{\varepsilon}{2}} < \theta^{\alpha_{N+e}(1+t \lambda_{e+1} + t\delta)}
 \end{align}
 où $t = \min(d-1,e+1)$. Ce choix a un sens ; en effet $K_e - 1-t\lambda_{e+1} - \lambda_e + \frac{\varepsilon}{2} > 0$ car $K_e \geq C_d = 5n^2(C_{d-1})^{2n} \geq 1 + n(C_{d-1})^{e+1} + (C_{d-1})^{e}$.
 
 \bigskip
 
 On applique maintenant le lemme~\ref{6lem_maj_DN}. 
 \\ On a bien $\psi_1(A,B) \leq H(B)^{-\lambda_e - \varepsilon}$. En effet, par hypothèse
 \begin{align*}
 \psi_1(A,B) \leq H(B)^{-K_e- \varepsilon}
 \end{align*}et $K_e \geq (C_d)^e \geq (C_{d-1})^e \geq \lambda_e$. 
 \\Si $H(B)$ est assez grand, le lemme~\ref{6lem_maj_DN} donne alors pour tout $i \in \llbracket 0, e-1 \rrbracket$ :
 \begin{align*}
 D_{N+i} &\leq c_{\ref{6cons_maj_D_N_lem}}\theta^{{\alpha_{N+i}} (1 + t\lambda_{e+1} + t\delta) } H(B)^{1+ t\lambda_{e+1} + t\delta} \left(\psi_1(A,B) H(B)^{ \lambda_e + \delta } + \frac{1 }{\theta^{{\alpha_{N+i+1}}}} \right) \\
 &\leq c_{\ref{6cons_maj_D_N_lem}}\theta^{{\alpha_{N+i}} (1 + t\lambda_{e+1} + t\delta) } H(B)^{1+ t\lambda_{e+1}+ t\delta} \left(H(B)^{-K_e- \varepsilon} H(B)^{ \lambda_e + \delta} + \frac{1 }{\theta^{{\alpha_{N+i+1}}}} \right). 
 \end{align*}
Par le choix de $N$ en $(\ref{6choix_N})$ on a pour tout $i \in \llbracket 0, e-1 \rrbracket$ :
\begin{align*}
 \theta^{\alpha_{N+i}(1+t \lambda_{e+1} + t\delta)} \leq H(B)^{K_e - 1-t\lambda_{e+1} - \lambda_e + \frac{\varepsilon}{2} }
\end{align*}
 par croissance de la suite $(\alpha_N)$. 
 De plus, par choix de $\delta$ en $(\ref{6choix_delta})$, on a $\frac{-\varepsilon}{2} + t\delta + \delta \leq \frac{-\varepsilon}{4}$.\\ Alors pour tout $i \in \llbracket 0, e-1 \rrbracket$:
\begin{align}\label{6maj_DNi_23}
 D_{N+i} &\leq c_{\ref{6cons_maj_D_N_lem}} \left(H(B)^{\frac{-\varepsilon}{4} } + \theta^{{\alpha_{N+i}} (1 + t\lambda_{e+1} + t\delta) - \alpha_{N+i+1} } H(B)^{1+ t\lambda_{e+1}+ t\delta}\right).
\end{align}
On s'intéresse dorénavant au second terme $\GG_{N+i} = \theta^{{\alpha_{N+i}} (1 + t\lambda_{e+1} + t\delta) - \alpha_{N+i+1} } H(B)^{1+ t\lambda_{e+1}+ t\delta}$. 
\\ On note $ \eta_e = 1+t \lambda_{e+1} + t\delta \leq 2n(C_{d-1})^{n}$ par choix de $\delta $ en $(\ref{6choix_delta})$. On a 
\begin{align*}
 \alpha_{N+i} \eta_e - \alpha_{N+i+1} = \alpha_N (\beta_{N+1} \ldots \beta_{N+i})(\eta_e - \beta_{N+i+1}) \leq 0
\end{align*}
car d'après $(\ref{6definition_Cd})$ :
\begin{align}\label{6diff_Cd_etae}
 C_d = 5n^2(C_{d-1})^{2n} \geq (2n(C_{d-1})^{n} +1 )2n(C_{d-1})^{n} \geq \eta_e^2 + \eta_e . 
\end{align}
Le choix de $N$ donne une minoration de $\theta$ :
\begin{align*}
 \theta \geq H(B) ^{ \frac{K_e -\eta_e + \frac{\varepsilon}{2}}{\alpha_{N+e}\eta_e}}
\end{align*}et on peut donc majorer :
\begin{align*}
 \theta^{{\alpha_{N+i}} \eta_e - \alpha_{N+i+1} } \leq H(B)^{ \frac{(K_e - \eta_e+ \frac{\varepsilon}{2})(\alpha_{N+i} \eta_e - \alpha_{N+i+1})}{\alpha_{N+e}\eta_e}}
\end{align*}
et donc : 
\begin{align*}
 \GG_{N+i} \leq H(B)^{ \frac{(K_e -\eta_e - \lambda_e + \frac{\varepsilon}{2} )(\alpha_{N+i} \eta_e - \alpha_{N+i+1}) } {\alpha_{N+e}\eta_e} + \eta_e } .
\end{align*}
On étudie alors l'exposant $\gamma_i = \frac{(K_e -\eta_e - \lambda_e + \frac{\varepsilon}{2} )(\alpha_{N+i} \eta_e - \alpha_{N+i+1}) } {\alpha_{N+e}\eta_e} + \eta_e $ pour $i \in \llbracket0,e-1 \rrbracket$.
\bigskip
On pose $\cons \label{6cons_epsilon_negatif} = \frac{1}{2K_e}$ une constante indépendante de $N$. D'après $(\ref{6diff_Cd_etae})$ on a $$c_{ \ref{6cons_epsilon_negatif}} \leq \frac{\eta_e^2}{2K_e\eta_e} \leq \frac{ C_d - \eta_e}{2K_e \eta_e} \leq \frac{(\beta_{N+i+1} - \eta_e)(2\beta_{N+1} \ldots \beta_{N+i})}{\beta_{N+1} \ldots \beta_{N+e}\eta_e}.$$ 
On majore alors $\gamma_i$, en remarquant que $\beta_{N+1} \ldots \beta_{N+e} \leq K_e$ :
\begin{align}\label{6maj_gamma0}
 \gamma_i &= \frac{(K_e -\eta_e - \lambda_e + \frac{\varepsilon}{2} )(\alpha_{N+i} \eta_e - \alpha_{N+i+1}) } {\alpha_{N+e}\eta_e} + \eta_e \nonumber \\
 &= \frac{(K_e -\eta_e - \lambda_e + \frac{\varepsilon}{2} )(\eta_e - \beta_{N+1})(\beta_{N+1} \ldots \beta_{N+i+1}) } {\beta_{N+1} \ldots \beta_{N+e}\eta_e} + \eta_e \nonumber \\
 &\leq \frac{(K_e -\eta_e - \lambda_e )(\eta_e - \beta_{N+i+1})(\beta_{N+1} \ldots \beta_{N+i}) + \beta_{N+1} \ldots \beta_{N+e}\eta_e^2 } {\beta_{N+1} \ldots \beta_{N+e}\eta_e} - c_{ \ref{6cons_epsilon_negatif}} \varepsilon \nonumber \\
 &\leq \frac{(K_e -\eta_e - \lambda_e )(\eta_e - \beta_{N+i+1})(\beta_{N+1} \ldots \beta_{N+i}) + K_e\eta_e^2 } {\beta_{N+1} \ldots \beta_{N+e}\eta_e} - c_{ \ref{6cons_epsilon_negatif}} \varepsilon.
\end{align}
On majore maintenant le premier terme en remarquant que celui-ci est maximal pour $i= 0$ car $\beta_{N+1} \ldots \beta_{N+i} \geq 0$, 
\begin{align*}
 (K_e -\eta_e - \lambda_e )(\eta_e - \beta_{N+1}) +K_e\eta_e^2 &\leq (K_e -\eta_e - \lambda_e )(\eta_e - C_d) + K_e\eta_e^2 \\
 &\leq (K_e -\eta_e - \lambda_e )\eta_e^2 + K_e\eta_e^2 \\
 &\leq -(\eta_e + \lambda_e )\eta_e \\
 &\leq 0
\end{align*}
en utilisant $(\ref{6diff_Cd_etae})$. \\
En reprenant $(\ref{6maj_gamma0})$ on a donc pour tout $i \in \llbracket 0,e-1 \rrbracket $: 
\begin{align*}
 \gamma_i \leq - c_{ \ref{6cons_epsilon_negatif}} \varepsilon.
\end{align*}
Enfin l'inégalité $(\ref{6maj_DNi_23})$ donne pour tout $i \in \llbracket 0, e-1 \rrbracket$:
\begin{align*}
 D_{N+i} &\leq c_{\ref{6cons_maj_D_N_lem}} \left(H(B)^{\frac{-\varepsilon}{4}} + H(B)^{ -c_{ \ref{6cons_epsilon_negatif}} \varepsilon}\right).
\end{align*}

 En particulier, si $H(B)$ est assez grand en fonction de $c_{\ref{6cons_maj_D_N_lem}},c_{\ref{6cons_epsilon_negatif}} $ et $\varepsilon$, on a :
\begin{align*}
 \forall i \in \llbracket 0,e-1 \rrbracket, \quad\| X_{N+i} \wedge Z_1 \ldots \wedge Z_e \| = \DD_{N+i} < 1.
\end{align*}
On a donc, d'après le lemme~\ref{2lem_X_in_B} : 
\begin{align*}
 \forall i \in \llbracket 0,e-1 \rrbracket, \quad X_{N+i} \in B.
\end{align*}
On rappelle que $B_{N,e} = \Vect(X_N, \ldots, X_{N+e-1})$. On a alors montré que, si $H(B)$ est assez grand, $B_{N,e} \subset B$ pour $N$ vérifiant $(\ref{6choix_N})$.
\\ Par égalité des dimensions, on a alors $B_{N,e} =B$ et le lemme~\ref{6lem_meilleurs_espaces} est prouvé.

 \end{preuve}

\subsection{Minoration du premier angle et conclusion}

Dans cette section, on minore l'angle $\psi_1(A, B_{N,e})$ ; cela nous permet de majorer $\mu_n(A|e)_1$ et de conclure la preuve du théorème~\ref{6theo_principal}.

\begin{lem}\label{6min_psi_A_BN}
 Il existe une constante $\cons \label{6cons_min_angle_1} >0 $ indépendante de $N$ telle que pour tout $N \in \N$: 
 \begin{align*}
 \psi_1(A, B_{N,e}) \geq c_{\ref{6cons_min_angle_1}} \theta^{ - \alpha_{N+e}}.
 \end{align*}
\end{lem}

\begin{preuve}
 Soit $X \in B_{N,e}$ et $Y \in A$ non nuls tels que :
 \begin{align*}
 \omega(X,Y) = \psi_1(A,B_{N,e}).
 \end{align*}
On rappelle (voir section~\ref{6sectionBNE}) que $B_{N,e }= \Vect(X_N, v_{N+1}, \ldots, v_{N+e-1})$ \\où $v_j = \begin{pmatrix}
 0 &
 \cdots & 0 & \rho_{0, j} & \cdots & \rho_{n-d-1, j}
\end{pmatrix}^\intercal $ avec 
\begin{align*}
 \forall i \in \llbracket 0, n-d-1 \rrbracket, \quad \rho_{i, j} = \frac{1}{u^{\phi(j)}_{j}}u^i_{j} \in \{0,1 \}.
\end{align*}
D'après la remarque~\ref{6req_section_BN}, les $v_j$ sont des vecteurs distincts de la base canonique. 
\\On pose alors $a, a_1, \ldots, a_{e-1} \in \R$ tels que :
\begin{align*}
 X = a\theta^{-\floor{\alpha_N}}X_N + \sum\limits_{i =1}^{e-1} a_iv_{N+i} = \begin{pmatrix}
a \\
0 \\
\vdots \\
0 \\
a\sigma_{0,N} +  \sum\limits_{i = 1 }^{e-1} a_i \rho_{0,N+i} \\
\vdots \\ 
a\sigma_{n-d-1,N} + \sum\limits_{i = 1 }^{e-1} a_i \rho_{n-d-1,N+i}
\end{pmatrix}.
\end{align*}
De même, on pose $b_1, \ldots, b_d \in \R$ tels que :
\begin{align*}
 Y = \sum\limits_{j = 1}^d b_jY_j = \begin{pmatrix}
  \sum\limits_{j=1}^d b_j \\
b_d \tau_{n-1,d} \\
\vdots \\
 \sum\limits_{j = 2}^d b_j \tau_{n-d+1,j} \\
b_1\sigma_0 +  \sum\limits_{j =2 }^d b_j \tau_{n-d,j} \\
\vdots \\ 
b_1\sigma_{n-d-1} +  \sum\limits_{j =2 }^d b_j\tau_{1,j}
\end{pmatrix}.
\end{align*}
Quitte à renormaliser les vecteurs $X$ et $Y$, ce qui ne change pas $\omega(X,Y)$, on suppose que l'on a :
\begin{align*}
 a^2 + \sum\limits_{i=1}^{e-1} a_i^2 = \sum\limits_{j=1}^d b_j^2 = 1.
\end{align*}
En particulier, cette hypothèse donne :
\begin{align}\label{6maj_prodXY}
 \| X \| \cdot \|Y\| \leq c_{\ref{6cons_maj_prodXY}} 
\end{align}
avec $\cons \label{6cons_maj_prodXY} >0$ une constante indépendante de $N$. 
\\ On cherche alors à minorer $\|X \wedge Y \|$. On utilise ici, le fait que les coordonnées du vecteur $X \wedge Y$ sont les mineurs de taille $2$ de la matrice $\begin{pmatrix}
 X & Y
\end{pmatrix} \in \MM_{n,2}(\R)$. En particulier, $\| X \wedge Y \| $ est minoré par la valeur absolue de chacun de ces mineurs. 
On pose les quantités suivantes : 
\begin{align*}
 \sigma &= \max(1, \max\limits_{i \in \llbracket 0, n-d-1\rrbracket} (\sigma_i )), \\
 \tau &= \min\limits_{ j \in \llbracket 2, d\rrbracket} \min\limits_{i \in \llbracket 1, n-1-d+j\rrbracket}(\tau_{i, j} ), \\
 T &= \max\limits_{ j \in \llbracket 2, d\rrbracket} \max\limits_{i \in \llbracket 1, n-1-d+j\rrbracket}(\tau_{i, j} ),\\
 s &= \frac{1}{\theta^{\floor{\alpha_{n-d-1}}}}.
\end{align*}
On suppose dorénavant $N \geq n-d-1$, alors pour tout $i \in \llbracket 0, n-d-1 \rrbracket $ il existe $ k \in \llbracket 0, N \rrbracket$ tel que $\phi(k) = i$, d'où $u^i_k \in \{1,2\}$ et 
\begin{align}\label{6s<sigma}
 s \leq \sigma_{i, N}. \end{align}
Enfin on définit la quantité suivante :
\begin{align}\label{6defM}
 M = \frac{\tau}{4(d-1)(T+\sigma)(1 + \frac{T}{\tau})^{d-2} } >0.
\end{align}
On va montrer
\begin{align}\label{6minor_voulue_XYwedge}
 \| X \wedge Y \| \geq c_{\ref{6cons_minor_XwedgeY}} \theta^{ - \alpha_{N+e}}
\end{align}
avec $\cons \label{6cons_minor_XwedgeY} > 0 $ indépendante de $N$, en faisant une disjonction de cas suivant les valeurs prises par les $b_j$.
\\ \textbullet \, \underline{Premier cas :} Si il existe $k \in \llbracket 2, d \rrbracket$ tel que $\left| 
 \sum\limits_{j = k}^d b_j \tau_{n-d+k-1,j} \right| \geq M \theta^{- \alpha_{N+e}}$, on étudie le mineur de la matrice $\begin{pmatrix}
 X & Y
\end{pmatrix} $ correspondant aux lignes $d- k +2 $ et $d+1 +\ell $ avec $\ell \in \llbracket 0, n-d-1 \rrbracket$. On a alors :
\begin{align}\label{6minor_cas_1_bminor}
 \forall \ell \in \llbracket 0, n-d-1 \rrbracket, \quad \| X \wedge Y \| &\geq \left| \det \begin{pmatrix}
 0 & \sum\limits_{j = k}^d b_j \tau_{n-d+k-1,j} \\
 a\sigma_{\ell,N} +  \sum\limits_{i = 1 }^{e-1} a_i \rho_{\ell,N+ i} & b_1\sigma_\ell +  \sum\limits_{j =2 }^d b_j \tau_{n-d- \ell,j}
 \end{pmatrix} \right| \nonumber \\
 &= \left| \sum\limits_{j = k}^d b_j \tau_{n-d+k-1,j} \right| \cdot \left| a\sigma_{\ell,N} +  \sum\limits_{i = 1 }^{e-1} a_i \rho_{\ell,N+i} \right| \nonumber \\
 &\geq M \theta^{- \alpha_{N+e}}\left| a\sigma_{\ell,N} +  \sum\limits_{i = 1 }^{e-1} a_i \rho_{\ell,N+i} \right|.
\end{align}
On distingue alors deux cas : $|a| \geq \frac{1}{2\sigma\sqrt{e}}$ et $|a| < \frac{1}{2\sigma\sqrt{e}}$. 
\\ Dans le premier cas, on choisit $\ell = \phi(N+e)$. On a en particulier :
\begin{align*}
 \rho_{\ell,N+1} = \ldots = \rho_{\ell,N+e-1} = 0
\end{align*}
et donc d'après $(\ref{6minor_cas_1_bminor})$ et en utilisant $(\ref{6s<sigma})$ :
\begin{align*}
 \| X \wedge Y \| \geq M \theta^{- \alpha_{N+e}}\left| a\sigma_{\ell,N} \right| \geq \frac{M\sigma_{\ell,N}}{2\sigma\sqrt{e}} \theta^{- \alpha_{N+e}} \geq \frac{Ms}{2\sigma\sqrt{e}} \theta^{- \alpha_{N+e}}
\end{align*}
ce qui donne $(\ref{6minor_voulue_XYwedge})$.
Dans le deuxième cas on a $|a| < \frac{1}{2\sigma\sqrt{e}} $ et donc :
\begin{align*}
  \sum\limits_{i = 1}^{e-1}a_i^2 = 1 - a^2 \geq 1- \frac{1}{4\sigma^2 e} \geq \frac{e-1}{e}
\end{align*}
puisque $\sigma \geq 1$. \\
En particulier il existe $i \in \llbracket 1, e-1 \rrbracket$ tel que $|a_i| \geq \frac{1}{\sqrt{e}}$.
\\On choisit alors $\ell = \phi(N+i)$. En particulier :
\begin{align*}
 \rho_{\ell, N+i } = 1 \text{ et } \rho_{\ell,N+1} = \ldots =\rho_{\ell,N+i-1} = \rho_{\ell,N+i+1}= \ldots = \rho_{\ell,N+e-1} = 0.
\end{align*}
D'après $(\ref{6minor_cas_1_bminor})$ on a :
\begin{align*}
 \| X \wedge Y \| &\geq M \theta^{- \alpha_{N+e}}\left| a\sigma_{\ell,N} + a_i \right| \\
 &\geq M \theta^{- \alpha_{N+e}} \left(|a_i| - |a|\sigma \right) \\
 &\geq M \theta^{- \alpha_{N+e}} \left(\frac{1}{\sqrt{e}} - \frac{1}{2\sqrt{e}}\right) \\
 &\geq \frac{M}{2\sqrt{e}} \theta^{- \alpha_{N+e}}.
\end{align*}
\\ \textbullet \, \underline{Second cas :} Si l'hypothèse $ \exists \: k \in \llbracket 2, d \rrbracket, \left| 
 \sum\limits_{j = k}^d b_j \tau_{n-d+k-1,j} \right| \geq M \theta^{- \alpha_{N+e}} $ du premier cas n'est pas respectée on a alors :
\begin{align*}
 \forall k \in \llbracket 2, d \rrbracket, \quad \left| 
 \sum\limits_{j = k}^d b_j \tau_{n-d+k-1,j} \right| < M \theta^{- \alpha_{N+e}}.
\end{align*}
On montrera plus loin (dans le lemme~\ref{6lem_maj_b}) qu'on a alors : 
\begin{align}\label{6resultat_sur_familleb_admis}
 \forall j \in \llbracket 2, d \rrbracket, \quad |b_j| \leq \frac{M(1 + \frac{T}{\tau})^{d-2} \theta^{-\alpha_{N+e}}}{\tau}, \quad \left| \sum\limits_{j=1}^d b_j \right| \geq \frac{1}{2} \text{ et } |b_1| \geq \frac{3}{4}
\end{align}
car $M = \frac{\tau}{4(d-1)(T+\sigma)(1 + \frac{T}{\tau})^{d-2} }\leq \frac{\tau}{4(d-1)(1 + \frac{T}{\tau})^{d-2} \theta^{-\alpha_{N+e}}} $. On admet pour l'instant ce résultat. 
\\On considère le mineur de la matrice $\begin{pmatrix}
 X & Y
\end{pmatrix} $ correspondant aux lignes $1$ et $d+1 +\ell $ avec $\ell \in \llbracket 0, n-d-1 \rrbracket$. On a alors pour tout $\ell \in \llbracket 0, n-d-1 \rrbracket,$
\begin{align}\label{6minor_cas_2_bminor}
 \| X \wedge Y \| &\geq \left| \det \begin{pmatrix}
 a & \sum\limits_{j = 1}^d b_j \\
 a\sigma_{\ell,N} +  \sum\limits_{i = 1 }^{e-1} a_i \rho_{\ell,N+i} & b_1\sigma_\ell +  \sum\limits_{j =2 }^d b_j \tau_{n-d- \ell,j}
 \end{pmatrix} \right|. 
\end{align}
On distingue alors deux cas : $|a| \geq K$ et $|a| < K $ avec $K = \frac{1}{4(\sigma + (T + \sigma)(d-1))\sqrt{e}}$. 
\\ Dans le premier cas, on choisit $ \ell = \phi(N+e)$. On a en particulier :
\begin{align*}
 \rho_{\ell,N+e} = 1 \text{ et } \rho_{\ell,N+1} = \ldots = \rho_{\ell,N+e-1} = 0
\end{align*}
et $\sigma_{\ell,N} = \sum\limits_{k = 0}^{N} \frac{u^{\ell}_k}{\theta^{\floor{\alpha_k}}} = \sum\limits_{k = 0}^{N+e-1 } \frac{u^{\ell}_k}{\theta^{\floor{\alpha_k}}} $ car alors $u^{\ell}_{N+1} = \ldots = u^{\ell}_{N+e-1} = 0 $. On en déduit que
\begin{align*}
 | \sigma_\ell - \sigma_{\ell, N} | = \sum\limits_{k = N+e}^{+ \infty } \frac{u^{\ell}_k}{\theta^{\floor{\alpha_k}}} \geq\frac{u^{\ell}_{N+e}}{\theta^{\floor{\alpha_{N+e}}}} \geq \frac{1}{\theta^{{\alpha_{N+e}}}}.
\end{align*}
L'inégalité $(\ref{6minor_cas_2_bminor})$ donne alors :
\begin{align*}
 \| X \wedge Y \| &\geq \left| a\left(b_1\sigma_\ell +  \sum\limits_{j =2 }^d b_j \tau_{n-d- \ell,j} \right) - a\sigma_{\ell,N} \sum\limits_{j = 1}^d b_j\right| \\
 &\geq |a| \left| b_1 (\sigma_\ell - \sigma_{\ell, N}) -  \sum\limits_{j =2 }^d b_j (\sigma_{\ell, N} - \tau_{n-d- \ell,j}) \right| \\
\end{align*}
et donc en rappellant que $\sigma$ et $T$ majorent respectivement les $\sigma_{\ell} $ et les $\tau_{\ell,j}$ on a 
\begin{align*}
 \| X \wedge Y\| &\geq |a| \left(|b_1|\frac{1}{\theta^{{\alpha_{N+e}}}} - (\sigma +T)  \sum\limits_{j =2 }^d |b_j| \right) \\
 &\geq |a|\left(\frac{3}{4\theta^{{\alpha_{N+e}}}} - (\sigma +T)(d-1) \frac{M(1 + \frac{T}{\tau})^{d-2} }{\tau \theta^{\alpha_{N+e}}} \right) \\
 &\geq \frac{|a|}{\theta^{\alpha_{N+e}}} \left(\frac{3}{4} - \frac{1}{4} 
 \right)\\
 &\geq \frac{K}{2\theta^{\alpha_{N+e}}}
\end{align*}
en utilisant $(\ref{6resultat_sur_familleb_admis})$ et par définition de $M$ en $(\ref{6defM})$.

 Dans le deuxième cas, on a $|a| < K = \frac{1}{4(\sigma + (T + \sigma)(d-1))\sqrt{e}}$ et donc 
\begin{align*}
  \sum\limits_{i = 1}^{e-1}a_i^2 = 1 - a^2 \geq 1- \frac{1}{16(\sigma + (T + \sigma)(d-1))^2 e} \geq \frac{e-1}{e}. 
\end{align*}
En particulier il existe $i \in \llbracket 1, e-1 \rrbracket$ tel que $|a_i| \geq \frac{1}{\sqrt{e}}$.
\\On choisit alors $\ell = \phi(N+i)$. En particulier :
\begin{align*}
 \rho_{\ell, N+i } = 1 \text{ et } \rho_{\ell,N+1} = \ldots \rho_{\ell,N+i-1} = \rho_{\ell,N+i+1}= \ldots = \rho_{\ell,N+e-1} = 0.
\end{align*}
D'après $(\ref{6minor_cas_2_bminor})$ et $(\ref{6resultat_sur_familleb_admis})$ on a :
\begin{align*}
 \| X \wedge Y \| &\geq \left| a\left(b_1\sigma_\ell +  \sum\limits_{j =2 }^d b_j \tau_{n-d- \ell,j} \right) - (a\sigma_{\ell,N} + a_i ) \sum\limits_{j = 1}^d b_j\right| \\
 &\geq \left| a_i \sum\limits_{j = 1}^d b_j\right| - \left|ab_1(\sigma_\ell - \sigma_{\ell,N}) \right| - \left| a\sum\limits_{j = 2}^d b_j(\tau_{n-d- \ell,j} - \sigma_{\ell, N}) \right| \\ 
 &\geq \frac{|a_i|}{2} - |a| (\sigma + (T + \sigma)\sum\limits_{j = 2}^d |b_j|) \\
 &\geq \frac{1}{2\sqrt{e}} - |a| (\sigma + (T + \sigma)(d-1)) \\
 &\geq \frac{1}{4\sqrt{e}}
\end{align*}
en majorant grossièremment $\sum\limits_{j = 2}^d |b_j|$ par $ d-1$ car $\sum\limits_{j=1}^d b_j^2 = 1$.
En particulier $$\|X \wedge Y \| \geq \theta^{ - \alpha_{N+e}}.$$

 Dans tous les cas on trouve, pour tout $N \geq n-d -1$ :
\begin{align*}
 \| X \wedge Y \| \geq c_{\ref{6cons_minor_XwedgeY}} \theta^{ - \alpha_{N+e}}
\end{align*}
avec $c_{\ref{6cons_minor_XwedgeY}} > 0 $ indépendante de $N$. Quitte à diminuer $c_{\ref{6cons_minor_XwedgeY}}$, comme on a $\psi_1(A,B_{N,e}) > 0 $ par $(e,j)\tir$irrationalité de $A$, on peut supposer cette inégalité vraie pour tout $N \in \N$.
\\ En utilisant $(\ref{6maj_prodXY})$ on trouve donc :
\begin{align*}
 \psi_1(A, B_{N,e}) = \omega(X,Y) = \frac{\| X \wedge Y \| }{\| X \| \cdot \|Y\|} \geq c_{\ref{6cons_minor_XwedgeY}} c_{\ref{6cons_maj_prodXY}} ^{-1} \theta^{ - \alpha_{N+e}}.
\end{align*}
Le lemme est donc prouvé avec $c_{\ref{6cons_min_angle_1}} = c_{\ref{6cons_minor_XwedgeY}} c_{\ref{6cons_maj_prodXY}} ^{-1}$.

\end{preuve}

On va maintenant montrer le résultat utilisé dans la preuve précédente. 

\begin{lem}\label{6lem_maj_b}
Soit $b_1, \ldots, b_d $ vérifiant $ \sum\limits_{j=1}^d b_j^2 = 1$. \\ On suppose qu'il existe $0 <M \leq \frac{\tau}{4(d-1)(1 + \frac{T}{\tau})^{d-2} \theta^{-\alpha_{N+e}}} $ tel que :
 \begin{align}\label{6lem_maj_b_hyp}
 \forall k \in \llbracket 2, d \rrbracket, \quad \left| 
 \sum\limits_{j = k}^d b_j \tau_{n-d+k-1,j} \right| < M \theta^{- \alpha_{N+e}}.
\end{align}
Alors 
\begin{align*}
 |b_1| \geq \frac{3}{4}, \quad \forall j \in \llbracket 2, d \rrbracket, \quad |b_j| \leq \frac{M(1 + \frac{T}{\tau})^{d-2} \theta^{-\alpha_{N+e}}}{\tau} \text{ et } \left| \sum\limits_{j=1}^d b_j \right| \geq \frac{1}{2} 
\end{align*}
où $\tau$ et $T$ désignent respectivement le minimum et le maximum de la famille des $\{ \tau_{i,j}\}$.
\end{lem}

\begin{preuve}
 On raisonne par récurrence descendante sur $i \in \llbracket 2, d\rrbracket$. On montre un résultat plus fin : 
 \begin{align}\label{6lem_maj_b_hyp_rec}
 \forall i \in \llbracket 2, d\rrbracket, \quad |b_i| \leq \frac{M(1 + \frac{T}{\tau})^{d-i} \theta^{-\alpha_{N+e}}}{\tau}.
 \end{align}
Si $i = d $ alors l'hypothèse $(\ref{6lem_maj_b_hyp})$ avec $ k =d$ donne :
\begin{align*}
 |b_d| \leq \frac{ M \theta^{- \alpha_{N+e}}}{ \tau_{n-1,d}} \leq \frac{M(1 + \frac{T}{\tau})^{d-d} \theta^{-\alpha_{N+e}}}{\tau}.
\end{align*}
Soit $i\in \llbracket 2, d-1 \rrbracket $, on suppose que pour tout $i' > i $ l'inégalité $(\ref{6lem_maj_b_hyp_rec})$ est vérifiée. On applique alors l'hypothèse $(\ref{6lem_maj_b_hyp})$ avec $ k = i$ :
\begin{align*}
 \left| \sum\limits_{j= i}^d b_j \tau_{n-d+i-1, j} \right| < M \theta^{- \alpha_{N+e}}.
\end{align*}
En particulier on a donc $ | b_i \tau_{n-d+i-1, i}| \leq M \theta^{- \alpha_{N+e}} + \sum\limits_{j= i+1}^d \left| b_j \tau_{n-d+i-1, j} \right| $ et donc en utilisant l'hypothèse de récurrence :
\begin{align*}
 | b_i | &\leq \frac{1}{\tau_{n-d+i-1, i}}\left(M \theta^{- \alpha_{N+e}} + \sum\limits_{j= i+1}^d \frac{M(1 + \frac{T}{\tau})^{d-j} \theta^{-\alpha_{N+e}}}{\tau}T \right) \\
 &\leq \frac{M}{\tau} \theta^{- \alpha_{N+e}} \left(1 + \frac{T}{\tau} \left(\frac{1 -(1 + \frac{T}{\tau})^{d-i} }{-\frac{T}{\tau}}\right)\right) \\
 &= \frac{M}{\tau}\theta^{- \alpha_{N+e}}\left(1 + \frac{T}{\tau}\right)^{d-i}
\end{align*}
ce qui prouve donc $(\ref{6lem_maj_b_hyp_rec})$ pour tout $i \in \llbracket 2,d \rrbracket$ et donc la deuxième inégalité du lemme. 

En particulier on a :
\begin{align*}
 \forall i \in \llbracket 2, d\rrbracket, \quad |b_i| \leq \frac{M(1 + \frac{T}{\tau})^{d-2} \theta^{-\alpha_{N+e}}}{\tau} \leq \frac{1}{4(d-1)}.
\end{align*}
On utilise le fait que $\sum\limits_{i = 1}^d b_i^2 = 1$ pour avoir :
\begin{align*}
 |b_1|^2 = 1 - \sum\limits_{i = 2}^d b_i^2 \geq 1 - (d-1)\frac{1}{16(d-1)^2} \geq \frac{3}{4}.
\end{align*}
De plus :
\begin{align*}
 \left| \sum\limits_{j=1}^d b_j \right| \geq |b_1| - \left| \sum\limits_{j=2}^d b_j \right| \geq \frac{3}{4} - (d-1)\frac{1}{4(d-1)} = \frac{1}{2}.
\end{align*}
\end{preuve}

Le corollaire suivant termine la preuve du théorème~\ref{6theo_principal} compte tenu de la minoration $(\ref{6minoration_exposant_par_Ke})$ obtenue à la fin de la section~\ref{6section_construction_espace_A}.

\begin{cor}\label{6cor_maj_expos}
 On a :
 \begin{align*}
 \mu_n(A|e)_1 \leq K_e = \max\limits_{i \in \llbracket 0, n-d-e \rrbracket } \beta_{i+1}\ldots\beta_{i+e}.
 \end{align*}
\end{cor}

\begin{preuve}
 On suppose par l'absurde que $\mu_n(A|e)_1 > K_e.$ Il existe alors $\varepsilon > 0 $ tel que :
 \begin{align*}
 \mu_n(A|e)_1 \geq K_e + 2\varepsilon.
 \end{align*}
Par définition de l'exposant diophantien, il existe une infinité d'espaces rationnels $B$ de dimension $e$ vérifiant :
\begin{align}\label{6inegl_psi_1_A_B}
 0 < \psi_1(A,B) \leq H(B)^{-\mu_n(A|e)_1 + \varepsilon} \leq H(B)^{- K_e - \varepsilon}.
\end{align}
D'après le lemme~\ref{6lem_meilleurs_espaces}, si $H(B)$ est assez grand on a $B = B_{N,e}$ avec $N \in \N$ pour les espaces vérifiant $(\ref{6inegl_psi_1_A_B})$. Il existe alors une infinité d'entiers $N \in \N$ tels que :
\begin{align}\label{6infinite_N_psi}
 0 < \psi_1(A,B_{N,e}) \leq H(B_{N,e})^{- K_e - \varepsilon}.
\end{align}
D'autre part d'après le lemme~\ref{6min_psi_A_BN} on a :
\begin{align}\label{6min_psi_1_A_BN_new}
 \forall N \in \N, \quad \psi_1(A,B_{N,e}) \geq c_{\ref{6cons_min_angle_1}} \theta^{- {\alpha_{N+e}}} = c_{\ref{6cons_min_angle_1}}\theta^{- {\alpha_{N}}\beta_{N+1} \ldots \beta_{N+e}} \geq c_{\ref{6cons_min_angle_1}}\theta^{- {\alpha_{N}}K_e}.
\end{align}
De plus d'après $(\ref{6hauteur_B_N})$ on a $H(B_{N,e}) \geq c_{\ref{6cons_minor_haut_BN}} \theta^{{\alpha_N}}.$ En regroupant alors les inégalités $(\ref{6infinite_N_psi}) $ et $(\ref{6min_psi_1_A_BN_new})$ on trouve :
\begin{align*}
c_{\ref{6cons_min_angle_1}}c_{\ref{6cons_minor_haut_BN}}^{K_e}H(B_{N,e}) ^{- K_e} \leq H(B_{N,e})^{- K_e - \varepsilon}.
\end{align*}
On rappelle que cette inégalité est valable pour une infinité d'entiers $N$. En faisant tendre $N$ vers $+ \infty$, $H(B_{N,e}) \to + \infty$ et alors $c_{\ref{6cons_min_angle_1}}c_{\ref{6cons_minor_haut_BN}}^{K_e} = 0$. 
\\ Or $c_{\ref{6cons_min_angle_1}} >0$ et $c_{\ref{6cons_minor_haut_BN}} > 0 $ ce qui soulève une contradiction et prouve le corollaire. 

\end{preuve}

%% file: Chapitres/6_Approx.tex
\chapter{Construction d'espaces avec exposants prescrits à plusieurs angles}\label{chap7}

Dans ce chapitre, on construit des espaces dont on peut calculer une famille d'exposants diophantiens correspondant à des angles différents.

On utilise encore une fois les constructions du chapitre~\ref{chap5} et le résultat du chapitre~\ref{chap4} pour calculer ces exposants.

On fixe $n \in \Nx$ et dans tout ce chapitre, comme $n$ ne varie pas on note $g(A,e) $ la quantité :
\begin{align*}
 g(A,e) = g(\dim(A), e,n) = \max(0, \dim(A) +e -n) .
\end{align*}
On introduit aussi, pour $e, \ell $ deux entiers, la quantité :
\begin{align*}
 f(e,\ell) = \max(0, e-\ell).
\end{align*}
On suppose ici que $n = (m+1)d$ avec $m, d \in (\Nx)^2$. On a alors le  théorème~suivant.

\begin{theo}\label{7theo_principale}
 Soit $\cons \label{7cons_petite_hyp_theoc2c1} =\left(1+ \frac{1}{m} \right)^{\frac{1}{d}}$ et $1 < \cons \label{7cons_petite_hyp_theoc2} < c_{\ref{7cons_petite_hyp_theoc2c1}} $.
\\ Soit $(\beta_{1,1}, \ldots, \beta_{1, m} ) \in \R^m $ tels que :
\begin{align}\label{7hypothèse_prop_princi1}
 \min\limits_{\ell \in \llbracket 1, m \rrbracket}(\beta_{1,\ell}) > {(3d)^{\frac{c_{\ref{7cons_petite_hyp_theoc2}}}{c_{\ref{7cons_petite_hyp_theoc2}}-1}}} \text{ et } \min\limits_{\ell \in \llbracket 1, m \rrbracket}(\beta_{1,\ell})^{c_{\ref{7cons_petite_hyp_theoc2c1}}} > \max\limits_{\ell \in \llbracket 1, m \rrbracket}(\beta_{1,\ell})^{c_{\ref{7cons_petite_hyp_theoc2}}}.
\end{align}
Pour $i \in \llbracket 2,d \rrbracket$, soit $(\beta_{i,1}, \ldots, \beta_{i, m}) \in \R^m$ vérifiant pour tout $i \in \llbracket 1,d-1 \rrbracket$ : 
\begin{align}
 \min\limits_{\ell \in \llbracket 1, m \rrbracket}(\beta_{i,\ell})^{c_{\ref{7cons_petite_hyp_theoc2c1}}} > \max\limits_{\ell \in \llbracket 1, m \rrbracket}(\beta_{i+1,\ell}) \label{7hypothèse_prop_princi2}
 \\
 \text{ et } \min\limits_{\ell \in \llbracket 1, m \rrbracket}(\beta_{i+1,\ell}) > \max\limits_{\ell \in \llbracket 1, m \rrbracket}(\beta_{i,\ell})^{c_{\ref{7cons_petite_hyp_theoc2}}} \label{7hypothèse_prop_princi3}.
\end{align}
Il existe un espace $A$ de dimension $d$ dans $\R^{n}$ tel que pour tous $e \in \llbracket 1, n-1\rrbracket $ et $k \in \llbracket 1 + g(A,e), \min(d,e) \rrbracket$ vérifiant $e < k (m+1) $ on a $A \in \II_n(d,e)_{k-g(A,e)}$ et : 
\begin{align*}
 \mu_n(A|e)_{k-g(A,e)}
 &= \frac{1}{\sum\limits_{q = 1 + f(e, mk)}^k \frac{1}{K_{q+d-k,v_q} }} 
\end{align*}
où $v_1, \ldots, v_k$ sont définis en posant $u$ et $v$ tels que $e =k v + u$ soit la division euclidienne de $e$ par $k$ et :
\begin{align}\label{7eq_def_vj}
 v_q &= \left\{
 \begin{array}{ll}
 v + 1 &\text{ si } q \in \llbracket 1, u \rrbracket \\
 v &\text{ si } q \in \llbracket u+1, k \rrbracket
 \end{array}
 \right. 
\end{align}
et enfin 
\begin{align*}
 \forall i \in \llbracket 1, d \rrbracket, \quad \forall v \in \llbracket 1, m \rrbracket, \quad K_{i,v} = \max\limits_{\ell \in \llbracket 0, m-v \rrbracket} \beta_{i, \ell +1} \ldots \beta_{i, \ell + v}.
\end{align*}
\end{theo}

\begin{req}\label{7req_somme_vq_egal_e}
 On a $ \sum\limits_{q = 1 }^k v_q = u(v+1) + (k-u)v = e $ et $v_q \leq m+1$ pour tout $q \in \llbracket 1, k \rrbracket$ car $e < k(m+1)$.
\end{req}

\begin{req}
Si $d =1$, on construit dans le théorème~\ref{7theo_principale} une droite $A$ vérifiant pour tout $e \in \llbracket 1, n-1 \rrbracket $ 
\begin{align*}
 \mu_n(A|e)_1 = K_{1,v_1} = \max\limits_{\ell \in \llbracket 0, n-1-e \rrbracket} \beta_{i, \ell +1} \ldots \beta_{i, \ell + e}.
\end{align*}
On retrouve alors le corollaire~\ref{5cor_technique} avec cependant ici une hypothèse $(\ref{7hypothèse_prop_princi1})$ plus restrictive.
\end{req}

\begin{req}
On peut reformuler le  théorème~en considérant la quantité $\nu_n(A|e)_k = (\mu_n(A|e)_k)^{-1}$. On va construire ici l'espace $A$ comme somme orthogonale de droites $A_q$. Pour $ k \in \llbracket 1 + g(A,e), \min(d,e) \rrbracket$, on note $B_q$ un espace rationnel de dimension $v_q$ réalisant la meilleure approximation de $A_q$ ; par construction des droites $A_q$, les espaces $B_q$ sont aussi deux à deux orthogonaux. Enfin on pose $B = \bigoplus\limits_{q =1}^k B_q$. Par orthogonalité, on a donc $H(B) = H(B_1) \ldots H(B_d) $. Les hypothèses \eqref{7hypothèse_prop_princi1}, \eqref{7hypothèse_prop_princi2} et \eqref{7hypothèse_prop_princi3}, et la construction des $B_q$ donnent alors 
\begin{align*}
    \forall q \in \llbracket 1 , f(e,mk) \rrbracket, \quad H(B_q) = 1
\end{align*}et que les quantités $\omega_1(A_q,B_q) $ pour  $ q \in \llbracket 1+ f(e,mk),k \rrbracket $ se comportent toutes comme  $\omega_k(A,B) $ (asymptotiquement en $H(B)$). 
On va montrer que les espaces $B$ ainsi construits sont les meilleures approximations de $A$ ; cela donne alors que $(\omega_k(A,B)) ^{\nu_n(A|e)_{k-g(A,e)}}$ se comporte comme $ H(B)$ quand $H(B) \to + \infty$. Comme par définition des exposants diophantiens, on a $(\omega_1(A_q,B_q))^{\nu_n(A_{q}|v_q)_1} $ de l'ordre de $ H(B_q)$, et alors
\begin{align*}
 \nu_n(A|e)_{k-g(A,e)} = \sum\limits_{q = 1 + f(e, mk)}^k \nu_n(A_{q+d-k}|v_q)_1
\end{align*}
avec $ \nu_n(A_{q+d-k}|v_q)_1 = K_{q+d-k,v_q}^{-1}. $
\end{req}

On énonce deux corollaires directs de ce théorème. Le premier montre que si on se restreint à des petites valeurs de $e$ (au plus $\frac{n-d}{d} = m$) alors on connait explicitement tous les exposants diophantiens $\mu_n(A|e)_k$ associés ; on a en particulier $g(d,e,n) =0 $ dans cette situation.

\begin{cor}
 Sous les hypothèses $(\ref{7hypothèse_prop_princi1})$, $(\ref{7hypothèse_prop_princi2})$ et $(\ref{7hypothèse_prop_princi3})$, il existe un espace $A$ de dimension $d$ tel que pour tout $e \in \llbracket 1, m \rrbracket$ et tout $k \in \llbracket 1, \min(d,e) \rrbracket$ on a 
 \begin{align*}
 \mu_n(A|e)_{k}
 &= \frac{1}{\sum\limits_{q = 1 }^k \frac{1}{K_{q+d-k,v_q} }}.
 \end{align*}

\end{cor}

Enfin, en prenant $\alpha_i = \beta_{i,1} = \ldots = \beta_{i,m}$ pour tout $i \in \llbracket 1,d \rrbracket$ dans le théorème~\ref{7theo_principale}, on a le corollaire suivant (qui concerne les mêmes exposants que le théorème~\ref{7theo_principale} mais possède moins de paramètres).

\begin{cor}\label{7cor_exposantsprescrits}
Soit $c_{\ref{7cons_petite_hyp_theoc2c1}}$ et $c_{\ref{7cons_petite_hyp_theoc2}}$ les constantes du  théorème~$\ref{7theo_principale}$.
\\
 Soit $\alpha_d \geq \ldots \geq \alpha_1$ tels que $\alpha_{1} > (3d)^{\frac{c_{\ref{7cons_petite_hyp_theoc2}}}{c_{\ref{7cons_petite_hyp_theoc2}}-1}} $ et pour tout $i \in \llbracket 1,d-1 \rrbracket, \alpha_{i}^{c_{\ref{7cons_petite_hyp_theoc2}}} < \alpha_{i+1} < \alpha_{i}^{c_{\ref{7cons_petite_hyp_theoc2c1}}}$.
\\Il existe un espace $A$ de dimension $d$ dans $\R^{n}$ tel que pour tous $e \in \llbracket 1, n-1\rrbracket $ et $k \in~\llbracket 1 +~g(A,e), \min(d,e) \rrbracket$ vérifiant $e < k (m+1) $ on a : 
\begin{align*}
 \mu_n(A|e)_{k -g(A,e)}
 &= \frac{1}{ \sum\limits_{q = 1 + f(e, mk)}^k \frac{1}{\alpha_{q+d-k}^{v_{q}}} } 
\end{align*}
\end{cor}

\bigskip

Pour démontrer le théorème~\ref{7theo_principale}, on construit dans la section~\ref{7section_def_A} un espace $A$ de dimension $d$ comme somme directe orthogonale de droites $\Vect(Y_j)$ pour $j \in \llbracket 1 ,d \rrbracket$. Ces droites sont construites de façon similaire à ce qui a été fait au chapitre~\ref{chap5}, on utilise d'ailleurs beaucoup de résultats provenant de celui-ci. Une grande partie de la démonstration est occupée par la preuve de la propriété~\ref{7prop_int}, dans la section~\ref{7section_prop_inter}. 
\\Cette propriété consiste en le calcul des exposants $\mu_n(A_J|e)_{k - g(A_J,e)}$ avec $A_J \subset A$ un sous-espace de dimension $k$ ; on peut en effet calculer ceux-ci car ils correspondent alors au dernier angle puisque $k = \dim(A_J)$. On raisonne alors par double inégalité, en exhibant des espaces rationnels approchant bien $A_J$ (lemme~\ref{7lem_prox_AJ_CN_H(CN)}) et en montrant ensuite que ce sont les \og meilleurs \fg{} (lemme~\ref{7lem_meilleur_espaces}). 
\\Enfin on conclut la preuve du théorème~\ref{7theo_principale} dans la section~\ref{7section_calcul_expo}, en utilisant le théorème~\ref{theo_somme_sev} qui permet de calculer $\mu_n(A|e)_{k - g(A,e)}$ en fonctions des $\mu_n(A_J|e)_{k - g(A_J,e)}$.

\section{Etude des \texorpdfstring{$\beta_{i,\ell}$}{} }
Les hypothèses $(\ref{7hypothèse_prop_princi1})$, $(\ref{7hypothèse_prop_princi2})$ et $(\ref{7hypothèse_prop_princi3})$ permettent d'obtenir une propriété d'indépendance linéaire sur les $\beta_{i,\ell}$ que l'on développe dans cette section. On définit de plus des suites $(\alpha_{i,N})_{N \in \N}$ pour $i \in \llbracket 1,d \rrbracket$ ; nécessaires à la construction de l'espace $A$ dans la section suivante.

\subsection{Hypothèse d'indépendance linéaire}
Tout d'abord pour tout $i \in \llbracket 1, d \rrbracket$, on introduit un terme $\beta_{i, m+1}$ pour avoir une propriété d'indépendance linéaire sur la suite des $(\beta_{i, \ell})_{i \in \llbracket 1, d \rrbracket, \ell \in \llbracket 1, m+1 \rrbracket}$.

\bigskip
 On choisit $(\beta_{i, m+1})_{i \in \llbracket 1, d \rrbracket}$ de sorte que :
\begin{align}
 &\min\limits_{\ell \in \llbracket 1, m +1 \rrbracket} \beta_{1, \ell} \geq \max((\max\limits_{\ell \in \llbracket 1, m +1 \rrbracket} \beta_{1, \ell})^{\frac{c_{\ref{7cons_petite_hyp_theoc2}}}{c_{\ref{7cons_petite_hyp_theoc2c1}}}},{(3d)^{\frac{c_{\ref{7cons_petite_hyp_theoc2}}}{c_{\ref{7cons_petite_hyp_theoc2}}-1} }}) \label{7_1hypot_m=beta_m+1},
\end{align}
ainsi que pour tout $ i \in \llbracket 1,d-1 \rrbracket, $
\begin{align}
\min\limits_{\ell \in \llbracket 1, m +1 \rrbracket}(\beta_{i,\ell})^{ {c_{\ref{7cons_petite_hyp_theoc2c1}}}} \geq \max\limits_{\ell \in \llbracket 1, m +1\rrbracket}(\beta_{i+1,\ell})\label{7_2ahypot_m=beta_m+1}, \\
 \min\limits_{\ell \in \llbracket 1, m +1 \rrbracket}(\beta_{i+1,\ell})\geq \max\limits_{\ell \in \llbracket 1, m +1 \rrbracket}(\beta_{i,\ell})^{{c_{\ref{7cons_petite_hyp_theoc2}} } } \label{7_2bhypot_m=beta_m+1}, 
\end{align}
et 
\begin{align}
&\text{la famille } \{1\} \cup \left(\frac{\log(E_i)}{\log(E_j)}\right)_{i,j \in \llbracket 1,d \rrbracket^2, i \neq j} \text{ est linéairement indépendante sur $\Q$} \label{7_3hypot_m=beta_m+1}
\end{align}
en notant $E_i = \beta_{i,1}\ldots \beta_{i,m}(\beta_{i,m+1})^m$ pour $i \in \llbracket 1,d \rrbracket$.

\bigskip 
De tels $\beta_{i, m+1}$ existent car les inégalités des hypothèses $(\ref{7hypothèse_prop_princi1})$, $(\ref{7hypothèse_prop_princi2})$ 
 et $(\ref{7hypothèse_prop_princi3})$ sont strictes. L'ensemble des $(\beta_{i, m+1})_{i \in \llbracket 1, d \rrbracket}$ vérifiant $(\ref{7_1hypot_m=beta_m+1})$, $(\ref{7_2ahypot_m=beta_m+1})$ et $(\ref{7_2bhypot_m=beta_m+1})$ contient donc un ouvert non vide, on peut choisir $(\beta_{1,m+1}, \ldots, \beta_{d,m+1})$ vérifiant $(\ref{7_3hypot_m=beta_m+1})$ dans celui-ci car l'ensemble des $d\tir$uplets pour lesquels $(\ref{7_3hypot_m=beta_m+1})$ n'est pas vérifiée est une réunion dénombrable d'hypersurfaces de $\R^d$. 
\bigskip

\begin{req}\label{7req_choix_constante}
 On a choisi $c_{\ref{7cons_petite_hyp_theoc2c1}} = \left(\frac{m+1}{m}\right)^{\frac{1}{d}}$ de sorte à avoir 
 $$ \forall v \in \llbracket 1, m +1 \rrbracket, \quad (v-1)c_{\ref{7cons_petite_hyp_theoc2c1}}^d \leq v.$$
\end{req}

\subsection{Prolongement des \texorpdfstring{$\beta_{i,\ell}$}{} et étude des \texorpdfstring{$K_{j,v}$}{} }\label{7section_def_beta}

 On pose pour $i \in \llbracket 1, d \rrbracket$ et $ \ell \in \llbracket m+2, 2m \rrbracket$, $\beta_{i, \ell} = \beta_{i, m+1}$. \\
Enfin, on étend par périodicité la suite des $\beta_{i, \ell}$ en posant pour tout $i \in \llbracket 1, d \rrbracket $ :
\begin{align*}
 \forall \ell \in \llbracket 1, 2m \rrbracket, \quad \forall p \in \N, \quad \beta_{i, \ell + 2mp } = \beta_{i, \ell}.
\end{align*}

\begin{req}\label{7req_max_2m=m}
On a par périodicité
$$ \forall i \in \llbracket 1, d \rrbracket, \quad \max\limits_{\ell \in \Nx} \beta_{i,\ell+1}\ldots \beta_{i,\ell+v}= \max\limits_{\ell \in \llbracket 0, 2m -1 \rrbracket} \beta_{i,\ell+1}\ldots \beta_{i,\ell+v}.$$
Comme $\beta_{i,m+1} = \beta_{i,m+1 } = \ldots = \beta_{i,2m}$ et pour tout $ \ell \in \llbracket 1, m \rrbracket$, $\beta_{i,m+1} \leq \beta_{i, \ell} $ on a donc
$$ \forall i \in \llbracket 1, d \rrbracket, \quad \forall v \in \llbracket 1, m +1 \rrbracket, \quad \max\limits_{\ell \in \llbracket 0, 2m -1 \rrbracket} \beta_{i,\ell+1}\ldots \beta_{i,\ell+v}= \max\limits_{\ell \in \llbracket 0, m + 1 -v \rrbracket} \beta_{i,\ell+1}\ldots \beta_{i,\ell+v}.$$
En particulier, si $v \leq m $, on remarque que $$\max\limits_{\ell \in \llbracket 0, m + 1 -v \rrbracket} \beta_{i,\ell+1}\ldots \beta_{i,\ell+v} = \max\limits_{\ell \in \llbracket 0, m -v \rrbracket} \beta_{i,\ell+1}\ldots \beta_{i,\ell+v} = K_{i,v}.$$
\end{req}
On étend alors la définition de $K_{i,v}$ à $v=0$ et $v= m+1$ par :
\begin{align*}
 \forall i \in \llbracket 1, d \rrbracket, \quad \forall v \in \llbracket 0, m +1 \rrbracket, \quad K_{i,v} = \max\limits_{\ell \in \llbracket 0, m+1-v \rrbracket} \beta_{i, \ell +1} \ldots \beta_{i, \ell + v} = \max\limits_{\ell \in \Nx} \beta_{i,\ell+1}\ldots \beta_{i,\ell+v} 
\end{align*}
avec la convention qu'un produit vide est égal à $1$ si bien que $K_{i,0} = 1 $.

\begin{prop}\label{7lem_min_KKi}
 Soit $k \in \llbracket 1, d \rrbracket $ et $e \in \llbracket 1, k(m+1)-1 \rrbracket$. Soit $1 \leq j_1 \leq \ldots \leq j_k \leq d $. 
 \\Alors pour tout $q \in \llbracket 1, k \rrbracket $ on a :
 \begin{align*}
 \left(1 - \frac{1}{\min\limits_{\ell \in \llbracket 1, m+1\rrbracket}\beta_{j_q, \ell}} \right)\left( \frac{1}{\sum\limits_{\ell =1+ f(e,mk) }^k \frac{1}{K_{j_{\ell},v_{\ell} }} } -1 \right) - K_{j_q,v_q-1} \geq 0
 \end{align*}
où $v_{1}, \ldots, v_{k}$ sont définis dans le théorème~\ref{7theo_principale}.
\end{prop}

\begin{req}
 Les hypothèses $(\ref{7hypothèse_prop_princi1})$ et $(\ref{7hypothèse_prop_princi2})$ sont plus fortes que celle énoncée dans la propriété $\ref{7lem_min_KKi}$. En fait, la propriété $\ref{7lem_min_KKi}$ est suffisante pour montrer le théorème~\ref{7theo_principale}. 
\end{req}

Avant de prouver l'inégalité de la propriété~$\ref{7lem_min_KKi}$, on va démontrer deux lemmes.

\begin{lem}\label{7lem_beta_minor}
 Pour tout $\beta \geq \min\limits_{\ell \in \llbracket 1, m+1 \rrbracket} \beta_{1, \ell}$ on a :
 \begin{align*}
 \beta - d(1+ 2 \beta^{\frac{1}{c_{\ref{7cons_petite_hyp_theoc2}}}}) \geq 0.
 \end{align*}
\end{lem}

\begin{preuve}
 Comme $\min\limits_{\ell \in \llbracket 1, m+1 \rrbracket} \beta_{1, \ell} \geq (3d)^{\frac{c_{\ref{7cons_petite_hyp_theoc2}}}{c_{\ref{7cons_petite_hyp_theoc2}} - 1 }} \geq 1$ par l'hypothèse $(\ref{7_1hypot_m=beta_m+1})$, on a :
 \begin{align*}
 \beta^{\frac{1}{c_{\ref{7cons_petite_hyp_theoc2}}}}(\beta^{\frac{c_{\ref{7cons_petite_hyp_theoc2}} -1 }{c_{\ref{7cons_petite_hyp_theoc2}}}} -2d ) \geq 1(3d - 2d) \geq d. 
 \end{align*}
 Or $ \beta^{\frac{1}{c_{\ref{7cons_petite_hyp_theoc2}}}}(\beta^{\frac{c_{\ref{7cons_petite_hyp_theoc2}} -1 }{c_{\ref{7cons_petite_hyp_theoc2}}}} -2d )= \beta - 2 d \beta^{\frac{1}{c_{\ref{7cons_petite_hyp_theoc2}}}}$ et le lemme est prouvé.
 
\end{preuve}

\begin{lem}\label{7lem_croissance_Kj}
 Soit $k \in \llbracket 1, d \rrbracket $ et $e \in \llbracket 1, k(m+1)-1 \rrbracket$. Soit $1 \leq j_1 < \ldots < j_k \leq d $. On note $u$ le reste de la division euclidienne de $e$ par $k$. 
 \\Alors pour tout $q \in \llbracket 1, k \rrbracket $ on a :
 \begin{align*}
 (K_{j_q,v_q -1})^{c_{\ref{7cons_petite_hyp_theoc2}}} \leq K_{j_{u+1}, v_{u+1}} \leq K_{j_q, v_q}.
 \end{align*}
\end{lem}

\begin{preuve}
Il suffit de montrer les deux inégalités suivantes pour conclure : 
\begin{align}
 &(\max\limits_{\ell \in \llbracket 1, m+1\rrbracket}\beta_{j_{q}, \ell})^{(v_{q}-1)c_{\ref{7cons_petite_hyp_theoc2}} } \leq (\min\limits_{\ell \in \llbracket 1, m+1\rrbracket}\beta_{j_{u+1}, \ell})^{v_{u+1} } \label{7ineg_aprouver_1}\\
 \text{ et } &(\max\limits_{\ell \in \llbracket 1, m+1\rrbracket}\beta_{j_{u+1}, \ell})^{v_{u+1} } \leq (\min\limits_{\ell \in \llbracket 1, m+1\rrbracket}\beta_{j_{q}, \ell})^{v_{q} }. \label{7ineg_aprouver_2}
\end{align}
En effet, on a $(\min \limits_{\ell \in \llbracket 1, m+1\rrbracket}\beta_{j, \ell})^{v} \leq K_{j,v} \leq (\max\limits_{\ell \in \llbracket 1, m+1\rrbracket}\beta_{j, \ell})^{v} $ pour tous $j \in \llbracket 1, d \rrbracket$ et $ v \in \llbracket 1, m+1 \rrbracket$. 

 \bigskip
 On distingue trois cas selon la valeur de $q$.
 \\ \textbullet \, \underline{Si $q < u+1$} alors par définition des $v_j$ on a :
 \begin{align*}
 v_q = v_{u+1} + 1 \in \llbracket 2, m+1 \rrbracket. 
 \end{align*}
 On utilise alors l'hypothèse $(\ref{7_2bhypot_m=beta_m+1})$, qui appliquée $j_{u+1}-j_q $ fois, donne
 \begin{align*}
 (\max\limits_{\ell \in \llbracket 1, m+1\rrbracket}\beta_{j_{q}, \ell})^{(v_{q}-1)c_{\ref{7cons_petite_hyp_theoc2}}^{j_{u+1}-j_q} } \leq (\min\limits_{\ell \in \llbracket 1, m+1\rrbracket}\beta_{j_{u+1}, \ell})^{v_{q} -1} = (\min\limits_{\ell \in \llbracket 1, m+1\rrbracket}\beta_{j_{u+1}, \ell})^{v_{u+1} };
 \end{align*}
 cela donne $(\ref{7ineg_aprouver_1})$ car $j_{u+1}-j_q \geq 1$. \\
 Pour l'autre inégalité, on utilise cette fois-ci l'hypothèse $(\ref{7_2ahypot_m=beta_m+1})$, qui appliquée $j_{u+1}-j_q $ fois, donne :
 \begin{align*}
 (\max\limits_{\ell \in \llbracket 1, m+1\rrbracket}\beta_{j_{u+1}, \ell})^{v_{u+1} } \leq (\min\limits_{\ell \in \llbracket 1, m+1\rrbracket}\beta_{j_{q}, \ell})^{(v_{q}-1)c_{\ref{7cons_petite_hyp_theoc2c1}}^{j_{u+1}-j_q } }.
 \end{align*}
Or, comme $j_{u+1}-j_q \leq d$ on a $(v_{q}-1)c_{\ref{7cons_petite_hyp_theoc2c1}}^{j_{u+1}-j_q } \leq (v_{q}-1)c_{\ref{7cons_petite_hyp_theoc2c1}}^{d} \leq v_q$ d'après la remarque~\ref{7req_choix_constante} et donc 
 \begin{align*}
 (\max\limits_{\ell \in \llbracket 1, m+1\rrbracket}\beta_{j_{u+1}, \ell})^{v_{u+1} } \leq (\min\limits_{\ell \in \llbracket 1, m+1\rrbracket}\beta_{j_{q}, \ell})^{v_{q}}
 \end{align*}
 ce qui prouve $(\ref{7ineg_aprouver_2})$ dans le cas $q < u+1$.
 \\ \textbullet \, \underline{Si $q = u+1$} l'inégalité $(\ref{7ineg_aprouver_2})$ est triviale. 
 D'autre part, les hypothèses $(\ref{7_2ahypot_m=beta_m+1})$ et $(\ref{7_2bhypot_m=beta_m+1})$ combinées 
 donnent 
 $$ (\min\limits_{\ell \in \llbracket 1, m+1\rrbracket}\beta_{j_{u+1}, \ell})\geq (\max\limits_{\ell \in \llbracket 1, m+1\rrbracket}\beta_{j_{u+1}, \ell})^{\frac{c_{\ref{7cons_petite_hyp_theoc2}}}{c_{\ref{7cons_petite_hyp_theoc2c1}}}}$$
et donc 
\begin{align*}
 (\min\limits_{\ell \in \llbracket 1, m+1\rrbracket}\beta_{j_{u+1}, \ell})^{v_{u+1}} \geq (\max\limits_{\ell \in \llbracket 1, m+1\rrbracket}\beta_{j_{u+1}, \ell})^{v_{u+1} {\frac{c_{\ref{7cons_petite_hyp_theoc2}}}{c_{\ref{7cons_petite_hyp_theoc2c1}}}} } \geq (\max\limits_{\ell \in \llbracket 1, m+1\rrbracket}\beta_{j_{u+1}, \ell})^{(v_{u+1}-1) c_{\ref{7cons_petite_hyp_theoc2}} }
\end{align*}
car $ (v-1) c_{\ref{7cons_petite_hyp_theoc2c1}} \leq (v-1) c_{\ref{7cons_petite_hyp_theoc2c1}}^d \leq v$ pour tout $v \in \llbracket 1, m +1 \rrbracket$ par la remarque~\ref{7req_choix_constante}. Cela donne $(\ref{7ineg_aprouver_2})$ dans le cas $q = u+1$.
\\ \textbullet \, \underline{Si $q > u+1$} alors par définition des $v_j$ on a :
 \begin{align*}
 v_q = v_{u+1} \in \llbracket 1, m \rrbracket. 
 \end{align*}
L'inégalité $(\ref{7ineg_aprouver_2})$ est claire car $\max\limits_{\ell \in \llbracket 1, m+1\rrbracket}\beta_{j_{u+1}, \ell} \leq (\min\limits_{\ell \in \llbracket 1, m+1\rrbracket}\beta_{j_{q}, \ell})^{c_{\ref{7cons_petite_hyp_theoc2}}^{-(j_q-j_{u+1})}} \leq \min\limits_{\ell \in \llbracket 1, m+1\rrbracket}\beta_{j_{q}, \ell}$ en appliquant $j_q-j_{u+1} \geq 1 $ fois l'hypothèse $(\ref{7_2bhypot_m=beta_m+1})$.
\\ L'inégalité $(\ref{7ineg_aprouver_1})$ provient de $(\ref{7_2ahypot_m=beta_m+1})$ appliquéé $j_q-j_{u+1}$ fois : 
\begin{align*}
 (\max\limits_{\ell \in \llbracket 1, m+1\rrbracket}\beta_{j_{q}, \ell})^{(v_{q}-1)c_{\ref{7cons_petite_hyp_theoc2}} } \leq (\min\limits_{\ell \in \llbracket 1, m+1\rrbracket}\beta_{j_{u+1}, \ell})^{(v_{q}-1)c_{\ref{7cons_petite_hyp_theoc2}}c_{\ref{7cons_petite_hyp_theoc2c1}}^{j_q-j_{u+1}}}\leq (\min\limits_{\ell \in \llbracket 1, m+1\rrbracket}\beta_{j_{u+1}, \ell})^{v_{q}}
\end{align*}
car $(v_{q}-1)c_{\ref{7cons_petite_hyp_theoc2}}c_{\ref{7cons_petite_hyp_theoc2c1}}^{j_q-j_{u+1}} \leq (v_{q}-1)c_{\ref{7cons_petite_hyp_theoc2c1}}^{d}$ car $c_{\ref{7cons_petite_hyp_theoc2}}\leq c_{\ref{7cons_petite_hyp_theoc2c1}}$ et $ (v_{q}-1)c_{\ref{7cons_petite_hyp_theoc2c1}}^{d}\leq v_q$ par la remarque~\ref{7req_choix_constante}. Cela donne $(\ref{7ineg_aprouver_1})$ dans le cas $q > u+1$.

 \end{preuve}

On peut alors prouver la propriété~\ref{7lem_min_KKi}. 

\begin{preuve}[ (propriété~\ref{7lem_min_KKi})]
Soit $q \in \llbracket 1, k\rrbracket$. On note $u$ le reste de la division euclienne de $e $ par $k$. \\
En notant $f = f(e,mk) $ on a :
\begin{align*} \frac{1}{\sum\limits_{\ell =1+ f }^k \frac{1}{K_{j_{\ell},v_{\ell} }} } \geq \frac{1}{(k-f) \frac{1}{ \min\limits_{\ell \in \llbracket 1 + f, k \rrbracket} K_{j_{\ell},v_{\ell} }} } \geq \frac{ K_{j_{u+1},v_{u+1}} }{ k-f} \geq \frac{ K_{j_{u+1},v_{u+1}} }{ d}
\end{align*}
car d'après le lemme~\ref{7lem_croissance_Kj}, $K_{j_{u+1},v_{u+1}} = \min\limits_{\ell \in \llbracket 1 + f, k \rrbracket} K_{j_{\ell},v_{\ell} } $.
\\ Enfin comme $\min\limits_{\ell \in \llbracket 1, m+1\rrbracket} \beta_{j_q, \ell} \geq \min\limits_{\ell \in \llbracket 1, m+1\rrbracket} \beta_{1, \ell} \geq 3d $ par les hypothèses $ (\ref{7_1hypot_m=beta_m+1})$ et $ (\ref{7_2bhypot_m=beta_m+1})$, on a 
\begin{align*}
 \left(1 - \frac{1}{\min\limits_{\ell \in \llbracket 1, m+1\rrbracket}\beta_{j_q, \ell}} \right) \geq \frac{1}{2}.
\end{align*}
On peut donc minorer la quantité $ \left(1 - \frac{1}{\min\limits_{\ell \in \llbracket 1, m+1\rrbracket}\beta_{j_q, \ell}} \right)\left( \frac{1}{\sum\limits_{\ell =1+ f(e,mk) }^k \frac{1}{K_{j_{\ell},v_{\ell} }} } -1 \right) - K_{j_q,v_q-1}$ par 
\begin{align*}
 \frac{1}{2}\left(\frac{ K_{j_{u+1},v_{u+1}} }{ d} -1 \right) - K_{j_q,v_q-1} &= \frac{K_{j_{u+1},v_{u+1}} -d(1 +2 K_{j_q,v_q-1}) }{2d}.
\end{align*}
Il reste alors à montrer que $K_{j_{u+1},v_{u+1}} -d(1 +2 K_{j_q,v_q-1}) \geq 0$ pour conclure. \\
Si $v_q = 1$ alors $K_{j_q,v_q-1} = 1$ et $K_{j_{u+1},v_{u+1}} \geq \min\limits_{\ell \in \llbracket 1, m+1\rrbracket} \beta_{1, \ell} \geq 3d$, donc la propriété~\ref{7lem_min_KKi} est prouvée dans ce cas.
\\ Sinon $v_q \geq 2$ et en appliquant la minoration du lemme~\ref{7lem_croissance_Kj} on a 
\begin{align*}
 {K_{j_{u+1},v_{u+1}} -d(1 +2 K_{j_q,v_q-1}) } \geq {K_{j_q,v_q-1}^{c_{\ref{7cons_petite_hyp_theoc2}}} -d(1 +2 K_{j_q,v_q-1}) }.
\end{align*}

On conclut la preuve en appliquant le lemme~\ref{7lem_beta_minor} avec $\beta = K_{j_q,v_q-1}^{c_{\ref{7cons_petite_hyp_theoc2}}} \geq \min\limits_{\ell \in \llbracket 1, m+1 \rrbracket} \beta_{1, \ell} $ par $(\ref{7_2bhypot_m=beta_m+1})$, et on trouve 
\begin{align*}
 K_{j_q,v_q-1}^{c_{\ref{7cons_petite_hyp_theoc2}}} -d(1 +2 K_{j_q,v_q-1}) \geq 0.
\end{align*}

\end{preuve}

\subsection{Définition des suites \texorpdfstring{$\alpha_{i,{N}}$}{} }

On introduit les suites $(\alpha_{i,N})_{N \in \N}$ pour $i \in \llbracket 1,d \rrbracket$ définies par :
\begin{align*}
 \alpha_{i,0} &= 1
 \\
 \forall N \in \N, \quad \alpha_{i, N+1 } &= \beta_{i,N+1} \alpha_{i,N}
\end{align*}

On rappelle que l'on note $E_i= \beta_{i,1}\ldots \beta_{i,m}(\beta_{i,m+1})^m$ pour $i \in \llbracket 1,d \rrbracket$. \\ En fait on a $E_i= \beta_{i,1}\ldots \beta_{i,2m}$ d'après les définitions du début de la section~\ref{7section_def_beta} et donc :
\begin{align}\label{7def_alphaN_formule_EI}
 \forall i \in \llbracket1, d\rrbracket, \quad \forall N \in \N, \quad \alpha_{i, N} = E_i^{ \floor{ \frac{N}{2m}}}\beta_{i,1} \ldots \beta_{i,N \mod 2m }= E_i^{ \floor{ \frac{N}{2m}}} \alpha_{i, N \mod 2m}
\end{align}
avec $N \mod 2m$ l'entier $k \in \llbracket 0, 2m-1 \rrbracket$ tel que $N \equiv k \mod (2m)$.

On remarque alors que :
\begin{align*}
 \forall i \in \llbracket 1, d \rrbracket, \quad \forall v \in \llbracket 1, m+1 \rrbracket, \quad K_{i,v} = \max\limits_{\ell \in \llbracket 0, m + 1-v \rrbracket} \beta_{i,\ell+1}\ldots \beta_{i,\ell+v} = \max\limits_{\ell \in \llbracket 0, m+1-v \rrbracket} \frac{\alpha_{i, \ell +v}}{\alpha_{i, \ell}} .
\end{align*}

\section{Construction de l'espace \texorpdfstring{$A$}{}}\label{7section_def_A}

On utilise les constructions que l'on a faites au chapitre~\ref{chap5} pour définir l'espace $A$ du théorème~\ref{7theo_principale}. On rappelle ensuite des propriétés sur les espaces dont on va montrer dans la suite qu'ils sont les meilleurs approximations de $A$.
\bigskip
 \\Pour tous $i \in \llbracket 1, d\rrbracket $ et $\ell \in \Nx $, on a $\beta_{i, \ell} \geq 3d \geq 2 + \frac{\sqrt{5}-1}{2}$ et la suite $(\beta_{i,N})$ est $(2m)\tir$périodique, on peut donc appliquer la propriété~\ref{5prop_technique}.
 \\Soit $i \in \llbracket 1, d \rrbracket$. D'après la propriété~\ref{5prop_technique}, il existe une droite $\Delta'_i = \Vect(Y'_i) \subset \R^{m+1}$ telle que : 
 \begin{align*}
 \forall e \in \llbracket 1, m \rrbracket, \quad \mu_{m+1}(\Delta'_i|e)_1 = \max\limits_{\ell\in \llbracket 0, 2m-1 \rrbracket} \beta_{i,\ell+1}\ldots \beta_{i,\ell+e}= \max\limits_{\ell \in \llbracket 0, m -e\rrbracket} \beta_{i,\ell+1}\ldots \beta_{i,\ell+e}
 \end{align*}
 en utilisant la remarque~\ref{7req_max_2m=m}. \\
 De plus, d'après la preuve de cette propriété et en fixant $\theta$ un nombre premier supérieur à $5$, $Y'_i$ est de la forme :
\begin{align*}
 Y'_i = \begin{pmatrix}
 1 & \sigma_{1,i} & \cdots & \sigma_{m,i}
 \end{pmatrix}^\intercal
\end{align*}
avec $\sigma_{1,i}, \ldots, \sigma_{m, i} $ algébriquement indépendants sur $\Q$.

Puisque $ n = (m+1)d$ on pose pour $ i \in \llbracket 1, d \rrbracket $ :
\begin{align*}
 Y_i = \begin{pmatrix}
 0 \\
 \vdots \\
 0 \\
 Y'_i \\
 0 \\
 \vdots \\
 0
 \end{pmatrix} \in \{ 0 \}^{(i-1)(m+1)} \times \R^{m+1} \times \{0\}^{n - i(m+1)}
\end{align*}
en plaçant les coordonnées de $Y'_i$ entre les lignes $(i-1)(m+1) + 1 $ et $i(m+1)$ de sorte que les vecteurs $Y_i$ soient deux à deux orthogonaux. 
\\ Enfin, les \og meilleurs\fg{} vecteurs approchant les $Y_i$ sont les $X_{N,i}$ avec :
\begin{align*}
 X_{N,i} = \begin{pmatrix}
 0 \\
 \vdots \\
 0 \\
 X'_{N,i} \\
 0 \\
 \vdots \\
 0
 \end{pmatrix} 
\end{align*}
où les $ X'_{N,i}$ sont définis en $(\ref{5def_XN_ chap5})$ dans le chapitre~\ref{chap5}. 
\\Pour $v \in \llbracket 1, m +1 \rrbracket$ et $N \in \N$, on pose :
\begin{align}
 B^i_{N, v} = \Vect(X_{N, i}, \ldots, X_{N+v-1,i})
\end{align}
Les vecteurs $X_{N,i}$ sont construits de sorte que $\dim(B^i_{N,v}) = v $, voir $(\ref{5B_N_base_v})$.
\\ Les propriétés $\ref{5lem_norme_XN}, \ref{5lem_haut_BN},\ref{5lem_omega_XN_Y}$ et $\ref{5min_psi_A_BN}$ donnent alors :
\begin{prop}\label{7prop_BNV_chap4}
 Soit $i \in \llbracket 1, d\rrbracket$. Pour $v \in \llbracket 1, m \rrbracket$ on a :
 \begin{align*}
 c_{\ref{7cons_norme_XNi_minor}}\theta^{\alpha_{i,N}} \leq 
 \:& \| X_{N,i} \| \leq c_{\ref{7cons_norme_XNi_major}}\theta^{\alpha_{i,N}}, \\
 c_{\ref{7cons_hauteur_BNvi_minor}}\theta^{\alpha_{i,N}} \leq 
 \:& H (B^i_{N, v}) \leq c_{\ref{7cons_hauteur_BNvi_major}}\theta^{\alpha_{i,N}},\\
 c_{\ref{7cons_angle_Yi_XNi_minor}} \theta^{-{\alpha_{i,N+1}}} \leq \:& \omega(Y_i, X_{N,i}) \leq c_{\ref{7cons_angle_Yi_XNi_major}} \theta^{-{\alpha_{i,N+1}}}, \\
 c_{\ref{7cons_psi1_Yi_BNvi_minor}} \theta^{-{\alpha_{i,N+v}}} \leq \:& \psi_1(\Vect(Y_i), \Vect(B^i_{N, v})) \leq c_{\ref{7cons_psi1_Yi_BNvi_major}} \theta^{-{\alpha_{i,N+v}}}.
 \end{align*}
 avec $\cons \label{7cons_norme_XNi_minor}$, $\cons \label{7cons_norme_XNi_major}$, $\cons \label{7cons_hauteur_BNvi_minor}$, $\cons \label{7cons_hauteur_BNvi_major}$, $\cons \label{7cons_angle_Yi_XNi_minor}$, $\cons \label{7cons_angle_Yi_XNi_major}$, $\cons \label{7cons_psi1_Yi_BNvi_minor}$ et $\cons \label{7cons_psi1_Yi_BNvi_major}$ des constantes indépendantes de $N$. 
\end{prop}

\bigskip
On pose enfin $A = \Vect(Y_1, \ldots, Y_d) $ et pour $J \subset \llbracket 1,d \rrbracket$ :
\begin{align}
 A_J = \Vect_{j \in J} (Y_j).
\end{align}
Comme les vecteurs $Y_j$ sont deux à deux orthogonaux, on a $\dim(A) = d$ et $\dim(A_J) = \#J$.

\begin{lem}\label{7lem_AJ_irr}
 Soit $J \subset \llbracket 1,d \rrbracket$ non vide et $e \in \llbracket \#J, \#J(m+1) -1 \rrbracket$.
 \\Alors $A_J$ est $(e,\#J - g(A_J,e)) -$irrationnel.
\end{lem}

\begin{preuve}
On suppose le contraire par l'absurde. Il existe alors $B$ un sous-espace rationnel de dimension $e$ tel que $\dim(A_J \cap B) \geq \#J - g(A_J,e) + g(A_J,e) = \#J$.
\\Commme $\dim(A_J) = \#J $ on a alors $A_J \cap B = A_J$.
\\ En particulier $Y_j\in B$ pour tout $j \in J$. Soit $X_1, \ldots, X_e$ une base rationnelle de $B$. 
\\On a $ Y_j \wedge X_1 \wedge \ldots \wedge X_e = 0$ pour tout $j \in J $, c'est-à-dire l'annulation de tous les mineurs de taille $e+1$ de la matrice $(Y_j \mid X_1 \mid \ldots \mid X_e)$.
\\ Soit $N$ un mineur de taille $e$ extrait de $M= (X_1 \mid \ldots \mid X_e)$. Comme $e \leq\#J(m+1) -1$, on peut trouver $j \in J$ tel que les lignes des coefficients non nuls de $Y_j$ ne soient pas toutes des lignes de la matrice extraite de $M$ correspondant à $N$. \\
On considère le mineur de taille $e+1$ de $(Y_j
\mid X_1 \mid \ldots \mid X_e)$ obtenu en considérant les lignes de $N$ et la ligne d'un coefficient non nul de $Y_j$ dont on vient de donner l'existence. On développe ce mineur par rapport à la première colonne et on trouve :
\begin{align*}
 0 = \tau N + \tau_1 N_1 + \ldots + \tau_{m}N_{m}
\end{align*}
où $\{\tau, \tau_1, \ldots, \tau_{m}\} = \{1, \sigma_{1,i}, \ldots, \sigma_{m,i}\}$ et les $N_j$ sont au signe près, des mineurs de taille~$e$ de $(X_1 \mid \ldots \mid X_e)$.
\\ Comme $\sigma_{1,i}, \ldots, \sigma_{m,i}$ sont algébriquement indépendants, $\tau, \tau_1, \ldots, \tau_{m}$ sont linéairement indépendants sur $\Q$ et on trouve $N = N_1 = \ldots = N_{m} = 0$ et en particulier $N= 0$. 
\\Tout mineur de taille $e$ de la matrice $(X_1 \mid \ldots \mid X_e)$ est nul donc $\dim(B) <e$ ce qui est contradictoire.

\end{preuve}


\section{Propriété intermédiaire}\label{7section_prop_inter}
Dans cette section, on démontre la propriété suivante :
\begin{prop}\label{7prop_int}
Soit $ J \subset \llbracket 1, d \rrbracket$ de cardinal $k \in \llbracket 1, d \rrbracket$. On écrit $ J = \{j_1, \ldots, j_k\}$ avec $j_1 < \ldots< j_k $.
\\ On a alors :
\begin{align*}
 \forall e \in \llbracket k, k(m+1) -1 \rrbracket, \quad \mu_n(A_J| e)_{k - g(A_J,e)} = \frac{1}{\sum\limits_{q = 1+ f(e,mk)}^k \frac{1}{K_{j_q,v_q }} }.
\end{align*}
\end{prop} 
Pour $e \in \llbracket k, k(m+1) -1 \rrbracket$, on cherche donc à calculer l'exposant diophantien $\mu_n(A_J|e)_{k - g(A_J,e)} $, le lemme $\ref{7lem_AJ_irr}$ donnant $A_J \in \II_n(k,e)_{k - g(A_J,e)}$. 

\bigskip

Toute la section~\ref{7section_prop_inter} est consacrée à la preuve de cette propriété.\\
Soit $J \subset \llbracket 1, d \rrbracket$ de cardinal $k$. Sans perte de généralité, on suppose dorénavant que $J = \llbracket 1, k \rrbracket$. D'après les notations de la propriété~\ref{7prop_int} on a $ j_{q} =q $ pour tout $q \in \llbracket 1, k \rrbracket$.
\\On fixe aussi $e \in \llbracket k, k(m+1) -1 \rrbracket $.


 \subsection{Hauteur des meilleurs espaces \texorpdfstring{$C^J_{N}$}{} }

On construit les \og meilleurs\fg{} espaces approchant $A_J$.
Soit $i \in \llbracket 1, d \rrbracket$. Pour $N_i$ un entier, on pose 
\begin{align*}
 C^i_{N_i} = \Vect(X_{N_i,i}, \ldots, X_{N_i +v_i -1, i}) = B^i_{N_i, v_i}
\end{align*}
avec $v_i $ défini par les relations en $(\ref{7eq_def_vj})$.

Comme $v_i \leq m +1$, on a en particulier $\dim(C^i_{N_i}) =v_i$. \\
Pour $N= (N_j)_{ j \in \llbracket 1, k \rrbracket} \in \N^{ k}$ on pose 
$$C^J_{N}= \bigoplus\limits_{j \in \llbracket 1, k \rrbracket} C_{N_j}^j .$$ 
Comme les espaces $C_{N_j}^j$ sont deux à deux orthogonaux, on a $\dim(C^J_{N}) = \sum\limits_{j = 1}^k v_j = e $ d'après la remarque $\ref{7req_somme_vq_egal_e}$. 

 On pose pour simplifier les notations 
\begin{align*}
 f = f(e, mk) = \max(0, e-km) \text{ et } g = g(A_J, e) = \max(0, k + e - n). 
\end{align*}

\bigskip 

 On étudie alors $C^J_{N}$.

\begin{lem}\label{7lem_vj=m+1}
 On note $u$ le reste de la division euclidienne de $e $ par $k$. \\
Si $f= 0 $ alors pour tout $j \in \llbracket 1, k \rrbracket$, on a 
 \begin{align*}
 v_j \leq m.
 \end{align*}
Si $f > 0$ alors $u = f$ et 
 \begin{align*}
 &\forall j \in \llbracket 1, f \rrbracket, \quad v_j = m+1 \\
 &\forall j \in \llbracket f +1, k \rrbracket, \quad v_j =m,
 \end{align*}
donc, si $f > 0$, pour tout $j \in \llbracket 1, f \rrbracket$ on a $$C^j_{N_j} = \{ 0 \}^{(j-1)(m+1)} \times \R^{m+1} \times \{0\}^{n - j(m+1)}. $$
\end{lem}

\begin{req}
 La dernière partie de ce lemme reste vraie si $f = 0$, l'intervalle $\llbracket 1,f \rrbracket$ étant alors vide.
 
\end{req}

\begin{preuve}
On note $f = f(e, mk)$ dans cette preuve. \\
 Si $f = 0$ alors $e -km \leq 0$.
 En écrivant la divison euclidienne $e = kv + u $, on a $v \leq m-1 $ ou $(v=m$ et $ u = 0$). 
 \\ Comme pour tout $j \in \llbracket 1, u\rrbracket $, $v_j = v+1$ et pour tout $j \in \llbracket u+1, k\rrbracket $, $v_j = v$ on a bien $$v_j \leq m$$ pour tout $j \in \llbracket 1, k \rrbracket$.
 
 \bigskip
 Sinon $f = e - km > 0$. 
 \\La division euclidienne $e = kv + u $ de $e $ par $k$ vérifie alors $v \geq m.$ \\ De plus comme $e \leq k(m+1)-1 $ on a $v = m$. Enfin 
 \begin{align*}
 u = e - kv = e- km = f. 
 \end{align*}
 Par définition des $v_j$ on a, pour tout $j \in \llbracket 1, f\rrbracket $, $$v_j = v + 1 = m+1$$ et pour $j \in \llbracket f +1,k\rrbracket $, $$v_j = v= m.$$ 
 
Pour $ j \in \llbracket 1, f \rrbracket $, on a par définition des $X_{N,j} $, $$ C^j_{N_j} \subset \{ 0 \}^{(j-1)(m+1)} \times \R^{m+1} \times \{0\}^{n - j(m+1)}.$$ 
Or $\dim(C^j_{N_j} ) =v_j $ et donc pour $j \in \llbracket 1, f \rrbracket$
 \begin{align*}
 \dim(C^j_{N_j} ) = v_j = m+1. 
 \end{align*}
 Par égalité des dimensions, cela prouve la dernière partie du lemme.
 
\end{preuve}

\begin{lem}\label{7lem_haut_CNM}
 On a 
 $$c_{\ref{7cons_haut_CN_minor}} \theta^{\left(\sum\limits_{j = f+1 }^k \alpha_{j,N_j}\right) } \leq H(C^J_{N}) \leq c_{\ref{7cons_haut_CN_major}} \theta^{\left(\sum\limits_{j = f+1 }^k \alpha_{j,N_j}\right) }$$
 avec $\cons \label{7cons_haut_CN_minor}$ et $\cons \label{7cons_haut_CN_major}$ indépendantes de $N$. 
\end{lem}

\begin{preuve}
 On rappelle que $C^J_{N}= \bigoplus\limits_{j \in \llbracket 1, k \rrbracket} C_{N_j}^j$ et que cette somme est orthogonale. On a donc d'après le corollaire~\ref{2cor_haut_espace_ortho_somme} :
 \begin{align*}
 H(C^J_{N}) = \prod\limits_{j = 1}^k H(C^j_{N_j}).
 \end{align*}
 D'après le lemme~\ref{7lem_vj=m+1}, pour tout $j \in \llbracket 1, f \rrbracket$, on a $H(C^j_{N_j}) = 1$.
 \\Par ailleurs pour $j \in \llbracket f+ 1, k \rrbracket, C^j_{N_j} = B^j_{N_j, v_j}$ avec $ v_j \leq m$. La propriété~\ref{7prop_BNV_chap4} donne donc $ c_{\ref{7cons_hauteur_BNvi_minor}}\theta^{\alpha_{j,N_j}} \leq H(C^j_{N_j}) \leq c_{\ref{7cons_hauteur_BNvi_major}}\theta^{\alpha_{j,N_j}} $ ce qui prouve le lemme.
 
\end{preuve}

\subsection{Estimation de l'angle \texorpdfstring{$\psi_{k-g}(A_J, C^J_{N}) $}{} }

On établit dans cette section, des relations entre l'angle $\psi_{k-g}(A_J, C^J_N)$ et la hauteur $H(C^J_N)$ en fonction des suites $(\alpha_{j,N_j})$ pour $j \in J$.

\begin{lem}\label{7lem_minor_prox_Y_j_CN}
Il existe $\cons \label{7cons_min_Yj_CNM} > 0 $ indépendante de $N = (N_1, \ldots,N_k)$ telle que 
\begin{align*}
 \forall i \in \llbracket f+1, k \rrbracket, \quad \omega_1(\Vect(Y_i), C^J_{N}) \geq c_{ \ref{7cons_min_Yj_CNM}}\theta^{-{\alpha_{i,M_i+1}}} 
\end{align*}
 en posant pour $j \in \llbracket 1,k \rrbracket$ et $N_j \in \N$, $M_j = N_j + v_j -1 $.
\end{lem}

\begin{preuve}
Soit $i \in \llbracket f+ 1, k \rrbracket $. \\
 D'après la propriété~\ref{7prop_BNV_chap4}, il existe une constante $\cons \label{7cons_minoration_prox_YJ_BN}$ telle que 
 \begin{align}\label{7inega_de_la_prop_quonavaitavant}
 \omega_1(Y_i, C^j_{N_i}) = \psi_1(Y_i, B^i_{N_i, v_i}) \geq c_{ \ref{7cons_minoration_prox_YJ_BN}} \theta^{-{\alpha_{i,N_i+ v_i}}} = c_{\ref{7cons_minoration_prox_YJ_BN}} \theta^{-{\alpha_{i,M_i +1}}}.
 \end{align}
On peut choisir $c_{\ref{7cons_minoration_prox_YJ_BN}}$ indépendante de $i$, en effet $i$ ne prend qu'un nombre fini de valeurs.\\
 Soit $X \in C^J_{N} \smallsetminus \{0\} $. On écrit 
 $$ X = \sum\limits_{j=1}^k V_j$$
 avec $V_j \in C^j_{N_j}$ pour tout $j \in \llbracket 1, k \rrbracket$. On va minorer $ \omega(Y_i, X) = \frac{\| Y_i \wedge X\|}{\| Y_i \| \cdot \|X\|}$.
\bigskip
 \\ Tout d'abord, on a pour tout $j \in \llbracket 1, k \rrbracket $, $V_j \in \{ 0 \}^{(j-1)(m+1)} \times \R^{m+1} \times \{0\}^{n - j(m+1)}$. En décomposant chaque $V_j$ dans la base canonique, on voit que cette décomposition fait intervenir des vecteurs distincts pour chaque $j \in \llbracket 1, k \rrbracket$. A fortiori les vecteurs de la base canonique de $\bigwedge^2 \R^n$ intervenant dans la décomposition de $Y_i \wedge V_j $ sont aussi distincts pour chaque $j \in \llbracket 1, k \rrbracket$. On a donc les vecteurs $Y_i \wedge V_j $ deux à deux orthogonaux.
 \\D'après le  théorème~de Pythagore on a alors 
 $$\| Y_i \wedge X\|^2 = \| \sum\limits_{j=1}^k (Y_i \wedge V_j ) \|^2 = \sum\limits_{j=1}^k \|Y_i \wedge V_j \|^2.$$ 
 En particulier pour tout $j \in \llbracket 1, k \rrbracket$, on a $\| Y_i \wedge X\| \geq \|Y_i \wedge V_j \|$. 
 \\Les $V_j$ sont deux à deux orthogonaux et donc $\|X\|^2 = \sum\limits_{j=1}^k \|V_j\|^2$. Il existe alors $j_0 \in \llbracket 1, k \rrbracket$ tel que $\|V_{j_0} \| \geq k^{\frac{-1}{2}}\|X\|$. 
\\ \textbullet \, Si $j_0 = i$ alors :
\begin{align*}
 \omega(Y_i, X) &\geq \frac{\| Y_i \wedge V_i\|}{\| Y_i \| \cdot \|X\|} \\
 &\geq k^{-\frac{1}{2}} \frac{\| Y_i \wedge V_i\|}{\| Y_i \| \cdot \|V_i\|} \\
 &\geq k^{-\frac{1}{2}} \omega(Y_i,V_i) \\
 &\geq k^{-\frac{1}{2}} \omega_1(Y_i, C^j_{N_i}) \\
 &\geq k^{-\frac{1}{2}}c_{\ref{7cons_minoration_prox_YJ_BN}} \theta^{-{\alpha_{i,M_i +1}}}
\end{align*}
en utilisant $(\ref{7inega_de_la_prop_quonavaitavant})$.
 \\ \textbullet \, Si $j_0 \neq i$ alors en étudiant les mineurs de taille $2$ de $(Y_i \mid X)$ où on extrait la ligne correspondant au $1$ de $Y_i$ et une autre ligne correspondant à une coordonnée non nulle de $V_{j_0}$, on a $\| Y_i \wedge X\|^2 \geq  \sum\limits_{v} (1 \times |v|)^2 = \| V_{j_0} \|^2$ où les $v$ sont les coordonnées non nulles de $V_{j_0}$. Cela donne :
\begin{align*}
 \omega(Y_i,X)
 &\geq \frac{\| V_{j_0} \|}{\|Y_i\| \|X\|} \\
 &\geq \frac{k^{\frac{-1}{2}}}{ \|Y_i\|} \\
 &\geq \frac{k^{\frac{-1}{2}}}{ \|Y_i\|} \theta^{- {\alpha_{i,M_{i} +1}} }.
\end{align*}
On pose alors $\cons \label{7cons_eni} = \min(\min\limits_{i = 1} ^k \frac{k^{\frac{-1}{2}}}{ \|Y_i\|}, k^{\frac{1}{2}} c_{\ref{7cons_minoration_prox_YJ_BN}} )$ et on a :
\begin{align*}
 \forall i \in \llbracket f+1, k \rrbracket, \quad \omega_1(\Vect(Y_i), C^J_{N})= \min\limits_{X \in C^J_{N} \smallsetminus \{0 \} } \omega(Y_i,X) \geq c_{\ref{7cons_eni}}\theta^{- {\alpha_{i,M_{i} +1}} }.
\end{align*}

\end{preuve}

\begin{lem}\label{7lem_prox_AJ_CN_theta}
On a 
\begin{align}\label{7trad_approx_Aj_C}
 c_{\ref{7cons_min_angle_AJC}} \theta^{- \min\limits_{j = f+1}^k \alpha_{j,M_j+1}} \leq \psi_{k-g}(A_J, C^J_{N}) \leq c_{\ref{7cons_maj_angle_AJC}}\theta^{- \min\limits_{j = f+1}^k \alpha_{j,M_j+1}}.
\end{align}
avec $\cons \label{7cons_min_angle_AJC}, \cons \label{7cons_maj_angle_AJC} $ indépendantes des $N_j$.
\end{lem}

\begin{preuve}

On rappelle que $$\psi_{k-g}(A_J, C^J_{N}) = \omega_k(A_J, C^J_{N})$$ car $g = g(A_J,e) $ et $\dim(C^J_{N}) = e$.
\\ D'après le lemme~\ref{7lem_vj=m+1}, pour tout $j \in \llbracket 1,f \rrbracket$
$$Y_j \in C^j_{N_j} \subset C^J_{N}$$ car $Y_j \in \{ 0 \}^{(j-1)(m+1)} \times \R^{m+1} \times \{0\}^{n - j(m+1)}$. 
\\ On applique la propriété~\ref{2prop_4.5Elio} aux espaces $$A_J = \bigoplus\limits_{j =1}^k \Vect(Y_j) \text{ et }\bigoplus\limits_{j =1}^f \Vect(Y_j) \oplus \bigoplus\limits_{j =f+1}^k \Vect(X_{M_j,j}) \subset C^J_{N} $$
en rappelant la notation $M_j = N_j + v_j -1$.
On a alors 
\begin{align*}
 \omega_k(A_J, C^J_{N}) &\leq \omega_k \bigg(A_J, \bigoplus\limits_{j =1}^f \Vect(Y_j) \oplus \bigoplus\limits_{j =f+1}^k \Vect(X_{M_j,j}) \bigg) \\
 &\leq c_{\ref{7cons_elio_maj}} \bigg(\sum\limits_{j = 1 }^f \omega(Y_j, Y_j) + \sum\limits_{j = f+ 1 }^k \omega(Y_j, X_{M_j,j}) \bigg) \\
 &= c_{\ref{7cons_elio_maj}} \sum\limits_{j = f+ 1 }^k \omega(Y_j, X_{M_j,j})
\end{align*}
avec $\cons \label{7cons_elio_maj} >0 $ ne dépendant que de $Y_1, \ldots, Y_k$ et $n$. 
\\ Or d'après la propriété~\ref{7prop_BNV_chap4}, pour tout $j \in \llbracket f+1, k \rrbracket$
\begin{align*}
 \omega(Y_j, X_{j,M_j}) \leq c_{ \ref{7cons_angle_Yi_XNi_major}} \theta^{-{\alpha_{j,M_j+1}}}
\end{align*}
avec $c_{\ref{7cons_angle_Yi_XNi_major}} $ une constante indépendante des $N_j$. On a alors 
\begin{align*}
 \omega_k(A_J, C^J_{N}) &\leq c_{\ref{7cons_elio_maj}}c_{ \ref{7cons_angle_Yi_XNi_major}} \sum\limits_{j = f+ 1 }^k \theta^{-{\alpha_{j,M_j+1}}} \leq k c_{\ref{7cons_elio_maj}}c_{ \ref{7cons_angle_Yi_XNi_major}} \theta^{-\min\limits_{j = f+ 1 }^k {\alpha_{j,M_j+1}}}.
\end{align*}
La majoration de $(\ref{7trad_approx_Aj_C})$ est donc prouvée avec $c_{ \ref{7cons_maj_angle_AJC}} = k c_{\ref{7cons_elio_maj}}c_{ \ref{7cons_angle_Yi_XNi_major}} $.

\bigskip 
On montre maintenant la minoration. 
\\ Soit $j \in \llbracket f+1, k \rrbracket$. Comme $Y_j \in A_J$, on a d'après le lemme~\ref{lem_inclusion_croissance} :
\begin{align*}
 \omega_1(\Vect(Y_j), C^J_{N}) \leq \omega_k(A_J, C^J_{N})
\end{align*}
car $ \dim(A_J) = k $ et $\dim(C^J_{N}) =e \geq k$. \\Or d'après le lemme~\ref{7lem_minor_prox_Y_j_CN}, on a $\omega_1(\Vect(Y_j), C^J_{N}) \geq c_{ \ref{7cons_min_Yj_CNM}}\theta^{-{\alpha_{j,M_j+1}}}$ avec $ c_{ \ref{7cons_min_Yj_CNM}}$ indépendante de $(N,M)$ et de $j$. On a donc 
\begin{align*}
 \omega_k(A_J, C^J_{N}) &\geq \max\limits_{j \in \llbracket f+1, k \rrbracket} c_{ \ref{7cons_min_Yj_CNM}}\theta^{-{\alpha_{j,M_j+1}}} \\
 &= c_{ \ref{7cons_min_Yj_CNM}}\theta^{- \min\limits_{j = f+1}^k \alpha_{j,M_j+1}}.
\end{align*}
Comme $ \omega_k(A_J, C^J_{N}) =\psi_{k-g}(A_J, C^J_{N})$, on a la minoration, ce qui conclut la preuve du lemme.

\end{preuve}

\begin{lem}\label{7lem_prox_AJ_CN_H(CN)}
 On a :
 \begin{align*}
 c_{\ref{7cons_psik_AJ_CNJ_minor}} H(C^J_{N}) ^{ \frac{ - \min\limits_{j = f+1}^k \alpha_{j,M_j+1}} {\sum\limits_{j = f+1 }^k \alpha_{j,N_j} } } \leq \psi_{k-g}(A_J, C^J_{N}) \leq c_{\ref{7cons_psik_AJ_CNJ_major}} H(C^J_{N}) ^{ \frac{ - \min\limits_{j = f+1}^k \alpha_{j,M_j+1}} {\sum\limits_{j = f+1 }^k \alpha_{j,N_j} } }
 \end{align*}
avec $\cons \label{7cons_psik_AJ_CNJ_minor}$ et $\cons \label{7cons_psik_AJ_CNJ_major}$ indépendantes de $N$.
\end{lem}

\begin{preuve}
 Les lemmes $\ref{7lem_haut_CNM}$ et $\ref{7lem_prox_AJ_CN_theta}$ donnent :
 \begin{align*}
 c_{\ref{7cons_haut_CN_minor}} \theta^{\left(\sum\limits_{j = f+1 }^k \alpha_{j,N_j}\right) } \leq \: & H(C^J_{N}) \leq c_{\ref{7cons_haut_CN_major}} \theta^{\left(\sum\limits_{j = f+1 }^k \alpha_{j,N_j}\right) } 
 \\\text{ et } c_{\ref{7cons_min_angle_AJC}} \theta^{- \min\limits_{j = f+1}^k \alpha_{j,M_j+1}} \leq \: & \psi_{k-g}(A_J, C^J_{N}) \leq c_{\ref{7cons_maj_angle_AJC}}\theta^{- \min\limits_{j = f+1}^k \alpha_{j,M_j+1}}.
 \end{align*}
 Ces deux estimations regroupées donnent le lemme.
 
\end{preuve}

\subsection{Minoration de l'exposant }
Le lemme~\ref{7lem_prox_AJ_CN_H(CN)} permet alors de minorer l'exposant $\mu_n(A_J| e)_{k - g(A_J,e)} $ en considérant certains $N = (N_1, \ldots, N_k) \in \N^k$ .

\begin{cor}\label{7cor_min_expos}
 On a $$ \mu_n(A_J| e)_{k - g(A_J,e)} \geq \frac{1}{\sum\limits_{i = f+1}^k \frac{1}{K_{i,v_{i} }} }$$
 avec $K_{q,v_q} = \max\limits_{\ell \in \llbracket 0, m-1 \rrbracket} \beta_{q, \ell +1} \ldots \beta_{q, \ell + v_q}$.
\end{cor}

\begin{preuve}
On rappelle que si $i \in \llbracket f+1, k \rrbracket$, d'après le lemme~\ref{7lem_vj=m+1} on a $v_i \in \llbracket 1, m \rrbracket $ et alors :
\begin{align*}
 \max\limits_{\ell \in \llbracket 0, m -1\rrbracket} \beta_{i,\ell+1}\ldots \beta_{i,\ell+v_i} = \max\limits_{\ell \in \llbracket 0, m-1 \rrbracket} \frac{\alpha_{i, \ell +v_i}}{\alpha_{i, \ell}} .
\end{align*}
Pour $i \in \llbracket f+1, k \rrbracket$, on note alors $L_i \in \llbracket 0, m-1 \rrbracket$ un entier tel que $$ K_{i,v_i}= \max\limits_{\ell \in \llbracket 0, m -1\rrbracket} \beta_{i,\ell+1}\ldots \beta_{i,\ell+v_i} = \frac{\alpha_{i, L_i +v_i}}{\alpha_{i, L_i}}.$$

 On rappelle que pour $j \in \llbracket 1,k \rrbracket$ et $N_j \in \N$, on a $M_j = N_j + v_j -1 $ où $v_j \in \llbracket 1, m+1 \rrbracket $ est défini en $(\ref{7eq_def_vj})$.

 Pour $(N_1, \ldots, N_f )\in \N^{f}$ et $ N_{f+1} \in \Nx$ un multiple de $2m$ fixés, on pose pour $i \in \llbracket f+ 2, k \rrbracket $ :
\begin{align}\label{7eq_def_Ni}
 &N_i = 2m \floor{ \frac{N_{f+1}\log(E_{f+1})}{2m\log(E_i)} + \frac{\log(\alpha_{f+1,v_{f+1}-1})}{\log(E_i)} } + L_i.
\end{align}
On sait, par le lemme~\ref{7lem_prox_AJ_CN_H(CN)}, que :
\begin{align}\label{7minor_mu_dans_preuve}
 \mu_n(A_J| e)_{k - g(A_J,e)} &\geq \limsup\limits_{\underset{2m | N_{f+1}}{N_{f+1} \to + \infty}} \frac{\min\limits_{i = f+1}^k \alpha_{i,M_i +1 } }{ \sum\limits_{i = f+1}^k \alpha_{i,N_i} }
\end{align}
où les $N_{f+2}, \ldots, N_k, M_{f+1}, \ldots, M_k$ sont ceux définis ci-dessus en fonction de $N_{f+1}$ qui sera choisi plus tard. 
La suite de la preuve est consacrée à montrer que cette limite supérieure est minorée par $ \frac{1}{\sum\limits_{i = f+1}^k \frac{1}{K_{i,v_{i} }} }$.
\\On fixe alors $(N_1, \ldots, N_f \in \N^{f})$ et $ N_{f+1} \in \Nx$ un multiple de $2m$ et on étudie $$\frac{\min\limits_{i = f+1}^k \alpha_{i,M_i +1 } }{ \sum\limits_{i = f+1}^k \alpha_{i,N_i} }.$$\\
On rappelle que l'on note $E_i= \beta_{i,1}\ldots \beta_{i,m}(\beta_{i,m+1})^m$ pour $i \in \llbracket 1,d \rrbracket$ et que
\begin{align}
 \text{la famille } \{1\} \cup \left(\frac{\log(E_{f+1})}{\log(E_j)}\right)_{j \in \llbracket f+2,k \rrbracket} \text{ est linéairement indépendante sur $\Q$} \label{7indp_line_Ej}
\end{align}
d'après $(\ref{7_3hypot_m=beta_m+1})$, par choix des $\beta_{i,m+1}$. 

Pour $i \in \llbracket f+ 2, k \rrbracket $ on définit :
\begin{align}\label{7eq_def_deltai}
 &\delta_i = \partfrac{ \frac{N_{f+1}\log(E_{f+1})}{2m\log(E_i)} + \frac{\log(\alpha_{f+1,v_{f+1}+ L_{f+1}-1})}{\log(E_i)} } \in [0,1[
\end{align}
où $\partfrac{u} = u - \floor{u}$ représente la partie fractionnaire de $u \in \R$. 
L'écriture $(\ref{7def_alphaN_formule_EI})$ donne pour $i \in \llbracket f+ 1,k \rrbracket $:
\begin{align}\label{7rel_alphaNi_Ei}
 \alpha_{i,N_i} = E_i^{\floor{\frac{N_i}{2m}}} \alpha_{i, N_i \mod 2m} = E_i^{\floor{\frac{N_i}{2m}}} \alpha_{i, L_i}
\end{align}
car $2m$ divise $N_i-L_i$ d'après $(\ref{7eq_def_Ni})$. De même :
\begin{align}\label{7rel_alphaMi_Ei}
 \alpha_{i,M_i+1} = E_i^{\floor{\frac{N_i + v_i }{2m}}} \alpha_{i, N_i + v_i \mod 2m} = E_i^{\floor{\frac{N_i}{2m}}} \alpha_{i, L_i + v_i } 
\end{align}
car $0 \leq L_i + v_i < 2m$ pour $i \in \llbracket f+1, k \rrbracket$ En effet $L_i \in \llbracket 0, m-1 \rrbracket$ et $v_i \in \llbracket 1, m \rrbracket$ pour tout $i \in \llbracket f+1, k \rrbracket$ d'après le lemme~\ref{7lem_vj=m+1}.
\\D'après les dépendances entre $N_i$ et $N_{f+1}$ on a : 
\begin{align} \label{7rel_alphaNf_et_alpha_Ni}
 \alpha_{i,N_i} &= E_i^{\floor{\frac{N_i}{2m}}} \alpha_{i, L_i} \nonumber \\
 &= E_i^{\floor{ \frac{N_{f+1}\log(E_{f+1})}{2m\log(E_i)} + \frac{\log(\alpha_{f+1,v_{f+1}-1})}{\log(E_i)} } }\alpha_{i, L_i}\nonumber \\
 &= E_i^{{ \frac{N_{f+1}\log(E_{f+1})}{2m\log(E_i) } + \frac{\log(\alpha_{f+1,v_{f+1}-1})}{\log(E_i)} } - \delta_i }\alpha_{i, L_i} \nonumber \\
 &= E_{f+1}^{ \frac{N_{f+1}}{2m} }\alpha_{f+1,v_{f+1}-1} E_i^{-\delta_i}\alpha_{i, L_i} \nonumber \\
 &= \frac{\alpha_{f+1, N_{f+1}} \alpha_{f+1,v_{f+1}-1} E_i^{-\delta_i}\alpha_{i, L_i}}{\alpha_{f+1, L_{f+1}}}
\end{align}
car $E_{f+1}^{ \frac{N_{f+1}}{2m} } = E_{f+1}^{\floor{\frac{N_{f+1}}{2m}}} $ comme $2m | N_{f+1}$ et $\alpha_{f+1, N_{f+1}} = E_{f+1}^{\floor{\frac{N_{f+1}}{2m}}} \alpha_{f+1, L_{f+1} } $ d'après $(\ref{7rel_alphaNi_Ei})$.
\\ De même en utilisant $(\ref{7rel_alphaMi_Ei})$ on a pour tout $i \in \llbracket f+1, k \rrbracket$ 
\begin{align}\label{7rel_alphaNf_et_alpha_Mi}
 \alpha_{i,M_i+1} =\frac{\alpha_{f+1, N_{f+1}} \alpha_{f+1,v_{f+1}-1} E_i^{-\delta_i}\alpha_{i, L_i+v_i}}{\alpha_{f+1, L_{f+1}}}.
\end{align}
On étudie donc :
\begin{align*}
\frac{\min\limits_{i = f+1}^k \alpha_{i,M_i +1 } }{ \sum\limits_{i = f+1}^k \alpha_{i,N_i} } &= \frac{\min\limits_{i = f+1}^k \frac{\alpha_{i,M_i +1 }}{\alpha_{f+1, N_{f+1}}} }{ \sum\limits_{i = f+1}^k \frac{\alpha_{i,N_i}}{\alpha_{f+1, N_{f+1}}} } \\
 &= \frac{\min(\frac{\alpha_{f+1, L_{f+1} + v_{f+1} }}{\alpha_{f+1, L_{f+1}}},\min\limits_{i = f+2}^k (\frac{\alpha_{f+1, L_{f+1} + v_{f+1} -1}}{\alpha_{f+1, L_{f+1}}} E_i^{-\delta_i}\alpha_{i, L_i +v_i} ) ) }{ 1+ \sum\limits_{i = f+2}^k \frac{\alpha_{f+1, L_{f+1} + v_{f+1} -1}}{\alpha_{f+1, L_{f+1}}} E_i^{-\delta_i}\alpha_{i, L_i } } 
\end{align*}
en utilisant les relations $(\ref{7rel_alphaNf_et_alpha_Ni})$ et $(\ref{7rel_alphaNf_et_alpha_Mi})$.\\
On rappelle que $K_{f+1, v_{f+1}} = \frac{\alpha_{f+1, L_{f+1} + v_{f+1} }}{\alpha_{f+1, L_{f+1}}}$ et $\frac{\alpha_{f+1, L_{f+1} + v_{f+1} } } {\beta_{f+1, L_{f+1} + v_{f+1} }}= \alpha_{f+1, L_{f+1} + v_{f+1} -1 }$. On a donc 
\begin{align*}
 \frac{\min\limits_{i = f+1}^k \alpha_{i,M_i +1 } }{ \sum\limits_{i = f+1}^k \alpha_{i,N_i} } &= \frac{\min(K_{f+1, v_{f+1} },\min\limits_{i = f+2}^k (\frac{K_{f+1, v_{f+1} } }{\beta_{f+1, L_{f+1} + v_{f+1} }}E_i^{-\delta_i}\alpha_{i, L_i +v_i} ) ) }{ 1+ \sum\limits_{i = f+2}^k \frac{K_{f+1, v_{f+1} } }{\beta_{f+1, L_{f+1} + v_{f+1} }} E_i^{-\delta_i}\alpha_{i, L_i } } \\ 
 &= \frac{\min(1,\min\limits_{i = f+2}^k (\frac{E_i^{-\delta_i}\alpha_{i, L_i +v_i} }{\beta_{f+1, L_{f+1} + v_{f+1} }} ) ) }{ \frac{1}{K_{f+1, v_{f+1}}} + \sum\limits_{i = f+2}^k \frac{ E_i^{-\delta_i}\alpha_{i, L_i } } {\beta_{f+1, L_{f+1} + v_{f+1} }}} .
\end{align*}
Enfin comme $\frac{\alpha_{i, L_i + v_i} }{K_{i,v_i}} = \alpha_{i, L_i }$ pour tout $i \in \llbracket f+2, k \rrbracket$, on a finalement 
\begin{align*}
 \frac{\min\limits_{i = f+1}^k \alpha_{i,M_i +1 } }{ \sum\limits_{i = f+1}^k \alpha_{i,N_i} } &= \frac{\min(1,\min\limits_{i = f+2}^k (\frac{E_i^{-\delta_i}\alpha_{i, L_i +v_i} }{\beta_{f+1, L_{f+1} + v_{f+1} }} ) ) }{ \frac{1}{K_{f+1, v_{f+1}}} + \sum\limits_{i = f+2}^k \frac{ E_i^{-\delta_i}\alpha_{i, L_i + v_i} } {\beta_{f+1, L_{f+1} + v_{f+1} }} \frac{1}{K_{i,v_i}}}.
\end{align*}
On a donc montré 
\begin{align*}
 \limsup\limits_{\underset{2m | N_{f+1}}{N_{f+1} \to + \infty}} \frac{\min\limits_{i = f+1}^k \alpha_{i,M_i +1 } }{ \sum\limits_{i = f+1}^k \alpha_{i,N_i} } = \limsup\limits_{\underset{2m | N_{f+1}}{N_{f+1} \to + \infty}} \frac{\min(1,\min\limits_{i = f+2}^k (\frac{E_i^{-\delta_i}\alpha_{i, L_i +v_i} }{\beta_{f+1, L_{f+1} + v_{f+1} }} ) ) }{ \frac{1}{K_{f+1, v_{f+1}}} + \sum\limits_{i = f+2}^k \frac{ E_i^{-\delta_i}\alpha_{i, L_i + v_i} } {\beta_{f+1, L_{f+1} + v_{f+1} }} \frac{1}{K_{i,v_i}}}.
\end{align*}
\bigskip
\\Le corollaire~\ref{2cor_kronecker}, $(\ref{7indp_line_Ej})$ et la définition des $\delta_i$ en $(\ref{7eq_def_deltai})$ donnent :
\begin{align*}
 \left\{ (\delta_{f+2},\ldots,\delta_{k}), N_{f+1} \in \Nx, 2m|N_{f+1} \right\} \text{ dense dans } [0,1[^{k - f- 1}.
\end{align*} On a donc, en utilisant cette densité : 
\begin{align*}
 \limsup\limits_{\underset{2m | N_{f+1}}{N_{f+1} \to + \infty}} \frac{\min\limits_{i = f+1}^k \alpha_{i,M_i +1 } }{ \sum\limits_{i = f+1}^k \alpha_{i,N_i} } &= \sup\limits_{(\delta_i) \in [0,1[^{k - f-1} } \frac{\min(1,\min\limits_{i = f+2}^k (\frac{E_i^{-\delta_i}\alpha_{i, L_i +v_i} }{\beta_{f+1, L_{f+1} + v_{f+1} }} ) ) }{ \frac{1}{K_{f+1, v_{f+1}}} + \sum\limits_{i = f+2}^k \frac{ E_i^{-\delta_i}\alpha_{i, L_i + v_i} } {\beta_{f+1, L_{f+1} + v_{f+1} }} \frac{1}{K_{i,v_i}}}.
\end{align*}
Pour tout $i \in \llbracket f+2, k \rrbracket$ et $\delta_i \in [0,1[$, on pose $$u_i = \frac{ E_i^{-\delta_i}\alpha_{i, L_i + v_i} } {\beta_{f+1, L_{f+1} + v_{f+1} }} \in \left] \frac{ \alpha_{i, L_i + v_i} } {E_i\beta_{f+1, L_{f+1} + v_{f+1} }}, \frac{ \alpha_{i, L_i + v_i} } {\beta_{f+1, L_{f+1} + v_{f+1} }}\right] $$
et $u_i$ prend toutes les valeurs de l'intervalle $\left] \frac{ \alpha_{i, L_i + v_i} } {E_i\beta_{f+1, L_{f+1} + v_{f+1} }}, \frac{ \alpha_{i, L_i + v_i} } {\beta_{f+1, L_{f+1} + v_{f+1} }}\right]$ quand $\delta_i$ parcourt $[0,1[$. 
\\De plus $1 \in \left] \frac{ \alpha_{i, L_i + v_i} } {E_i\beta_{f+1, L_{f+1} + v_{f+1} }}, \frac{ \alpha_{i, L_i + v_i} } {\beta_{f+1, L_{f+1} + v_{f+1} }}\right]$ pour tout $i$ car $\alpha_{i, L_i + v_i} \leq E_i = \alpha_{i,2m}$ et $1 \leq \beta_{f+1, L_{f+1} + v_{f+1} } \leq \min\limits_{\ell \in \llbracket 1, 2m\rrbracket} \beta_{i, \ell} \leq \alpha_{i, L_i + v_i}$ d'après l'hypothèse $(\ref{7_2bhypot_m=beta_m+1})$. On a donc 
\begin{align*}
 \sup\limits_{(\delta_i) \in [0,1[^{k - f-1} } \frac{\min(1,\min\limits_{i = f+2}^k u_i ) }{ \frac{1}{K_{f+1, v_{f+1}}} + \sum\limits_{i = f+2}^k u_i \frac{1}{K_{i,v_i}}} &\geq \frac{\min(1,\min\limits_{i = f+2}^k 1 ) }{ \frac{1}{K_{f+1, v_{f+1}}} + \sum\limits_{i = f+2}^k \frac{1}{K_{i,v_i}}}= \frac{1}{\sum\limits_{i = f+1}^k \frac{1}{K_{i,v_{i} }} }. 
\end{align*}
On a donc montré que $\mu_n(A_J|e)_{k - g(A_J,e)} \geq \frac{1}{\sum\limits_{i = f+1}^k \frac{1}{K_{i,v_{i} }} }$ en reprenant l'équation $(\ref{7minor_mu_dans_preuve})$ ce qui termine la preuve du lemme.

\end{preuve}

\subsection{Majoration de l'exposant}
Dans cette section, on montre que les espaces $C^J_{N}$ réalisent en fait les \og meilleures \fg{} approximations de $A$. Cela permet alors de majorer $\mu_n(A_J| e)_{k - g(A_J,e)}$ et ainsi de conclure la preuve de la propriété~\ref{7prop_int}.

\begin{lem}\label{7lem_meilleur_espaces}
Soit $\varepsilon >0$ et $C$ un sous-espace rationnel de dimension $e$ tel que :
\begin{align}\label{7hyp_lem_meilleure_approx}
 \psi_{k-g}(A_J,C) \leq H(C)^{- \frac{1}{\sum\limits_{j = f+1}^k \frac{1}{K_{j,v_{j} }} } - \varepsilon}. 
\end{align}
Alors, si $H(C)$ est assez grand en fonction de $\varepsilon$, il existe $N \in (\Nx)^{k}$ tel que 
 \begin{align*}
 C = C^J_{N}.
 \end{align*}
\end{lem}

\begin{preuve}
Dans cette preuve, on pose $K = \frac{1}{\sum\limits_{j = f+1}^k \frac{1}{K_{j,v_{j} }} } $. \\
Pour $i \in \llbracket 1, k \rrbracket$, on pose $N_i$ l'entier vérifiant :
 \begin{align}\label{7choix_N_i}
 \theta^{{\alpha_{i,{N_i+v_i -1 }}} } \leq H(C)^{ K+ \frac{\varepsilon}{2} - 1} < \theta^{{\alpha_{i,{N_i+v_i}}} }. 
\end{align}
On va montrer que $N = (N_1, \ldots, N_k)$ convient. Soit $Z_1, \ldots, Z_e$ une $\Zbase$ de $C\cap \Z^n$. On pose pour $N$ un entier et $i \in \llbracket 1, k \rrbracket$ :
\begin{align*}
 D_{N,i} = \| X_{N, i} \wedge Z_1 \wedge \ldots \wedge Z_e \|.
\end{align*}
On va montrer que pour tout $i \in \llbracket 1,k \rrbracket$ : 
\begin{align*}
 \forall \ell \in \llbracket 0, v_i -1 \rrbracket, \quad D_{N_{i}+\ell, i } < 1 
 \end{align*}
 et le lemme~\ref{2lem_X_in_B} permet de conclure.

 \bigskip 
 On fixe $i \in \llbracket 1, k \rrbracket $ et $\ell \in \llbracket 0, v_i -1 \rrbracket $. Le lemme~\ref{lem_phi_dim1} donne 
 \begin{align*}
 D_{N_{i}+\ell, i } &= \omega(X_{N_i+\ell, i}, C) \| X_{N_i+\ell,i} \| H(C).
 \end{align*}
 En utilisant l'inégalité triangulaire (lemme~\ref{2lem_ineg_triang}) on a 
 \begin{align*}
 \omega(X_{N_i+\ell, i}, C) \leq \omega(X_{N_i+\ell, i}, Y_i ) + \omega(Y_i, C). 
 \end{align*}
 Or la propriété~\ref{7prop_BNV_chap4} donne $$c_{\ref{7cons_norme_XNi_minor}}\theta^{\alpha_{i,N_i+\ell}} \leq 
 \| X_{N_i+\ell,i} \| \leq c_{\ref{7cons_norme_XNi_major}}\theta^{\alpha_{i,N_i+\ell}}$$ et $$
 c_{\ref{7cons_angle_Yi_XNi_minor}} \theta^{-{\alpha_{i,N_i+\ell+1}}} \leq \omega(X_{N_i+\ell,i},Y_i) \leq c_{\ref{7cons_angle_Yi_XNi_major}} \theta^{-{\alpha_{i,N_i+\ell+1}}}$$ avec des constantes indépendantes des $N_i$. 
 \\D'autre part, comme $Y_i \in A_J$ et $\dim(A_J) = k$ on a, par le lemme~\ref{lem_inclusion_croissance}:
 \begin{align*}
 \psi_{k-g}(A_J,C) = \omega_{k}(A_J, C) \geq \omega(Y_i, C). 
 \end{align*}
 On a donc :
 \begin{align}
 D_{N_{i}+\ell, i } &\leq c_{\ref{7cons_maj_DNell}} H(C) \theta^{\alpha_{i,N_i+\ell}} \left(\theta^{-\alpha_{i,N_i+\ell+1}} + H(C)^{-K - \varepsilon}\right) \nonumber \\
 &= c_{\ref{7cons_maj_DNell}} \left(\theta^{-\alpha_{i,N_i+\ell+1} + \alpha_{i,N_i+\ell}} H(C) + \theta^{\alpha_{i,N_i+\ell}}H(C)^{-K - \varepsilon +1 }\right) \label{7maj_DNiell_final}
 \end{align}
 avec $\cons \label{7cons_maj_DNell} $ indépendante de $N$, en utilisant l'hypothèse (\ref{7hyp_lem_meilleure_approx}) sur $C$.
 \\Le choix de $N_i$ en $(\ref{7choix_N_i})$ donne 
 \begin{align} \label{7inega_premier_terme}
 \theta^{\alpha_{i,N_i+\ell}}H(C)^{-K - \varepsilon +1 } \leq \theta^{\alpha_{i,N_i+v_i -1 }} \leq H(C)^{K + \frac{\varepsilon}{2} +1 }H(C)^{-K - \varepsilon -1 } = H(C) ^{-\frac{\varepsilon}{2}}.
 \end{align}
On s'intéresse maintenant au terme $\theta^{-\alpha_{i,N_i+\ell+1} + \alpha_{i,N_i+\ell}} H(C)$ de $(\ref{7maj_DNiell_final})$. Toujours par $(\ref{7choix_N_i})$, on a 
 \begin{align}
 \theta^{-\alpha_{i,N_i+\ell+1} + \alpha_{i,N_i+\ell}}H(C) &\leq \left(H(C)^{\frac{K + \frac{\varepsilon}{2} -1}{\alpha_{i,N_i +v_i}}}\right)^{-\alpha_{i,N_i+\ell-1} + \alpha_{i,N_i+\ell}} H(C) \nonumber \\
 &= H(C)^{ \frac{(K + \frac{\varepsilon}{2} -1)(-\alpha_{i,N_i+\ell+1} + \alpha_{i,N_i+\ell}) + \alpha_{i,N_i +v_i}}{\alpha_{i,N_i +v_i}}} \label{7inega_deuxieme_terme}.
 \end{align}
 On s'intéresse alors au numérateur de cet exposant ; en mettant de coté le terme en $\varepsilon$, on trouve :
 \begin{align*}
 (K -1)(-\alpha_{i,N_i+\ell+1} + \alpha_{i,N_i+\ell}) + \alpha_{i,N_i +v_i} = \alpha_{i,N_i + \ell} \left((K -1)(-\beta_{i,N_i+\ell+1} + 1) + \frac{\alpha_{i,N_i +v_i}}{\alpha_{i,N_i + \ell} }\right). 
 \end{align*}
 On va montrer que ce terme est négatif, en effet :
 \begin{align*}
 (K -1)(-\beta_{i,N_i+\ell+1} + 1) + \frac{\alpha_{i,N_i +v_i}}{\alpha_{i,N_i + \ell} } &= (K - 1)(-\beta_{i,N_i+\ell+1} + 1) + \beta_{i,N_i+\ell+1} \ldots \beta_{i,N_i+v_i} \\
 &= -\beta_{i,N_i+\ell+1} (K - 1 - \beta_{i,N_i+\ell+2} \ldots \beta_{i,N_i+v_i}) + K -1.
 \end{align*}
Comme $\ell \in \llbracket 0, v_i-1 \rrbracket$ on a $\beta_{i,N_i+\ell+2} \ldots \beta_{i,N_i+v_i} \leq \beta_{i,N_i+2} \ldots \beta_{i,N_i+v_i} \leq K_{i,v_i-1} $ car les $v_i-1$ facteurs du produit$\beta_{i,N_i+2} \ldots \beta_{i,N_i+v_i}$ sont consécutifs. 
On a alors 
 \begin{align}\label{7inegalite_presque_finale}
 (K -1)(-\beta_{i,N_i+\ell+1} + 1) + \frac{\alpha_{i,N_i +v_i}}{\alpha_{i,N_i + \ell} } &\leq -\beta_{i,N_i+\ell+1} (K - 1 - K_{i,v_i -1}) + K -1.
 \end{align}
 Or la propriété~\ref{7lem_min_KKi} donne 
 \begin{align*}
 K - 1 - K_{i,v_i-1} \geq \frac{K-1}{\min\limits_{ \ell \in \llbracket 1, m+1 \rrbracket} (\beta_{i,\ell})}.
 \end{align*}
 L'inégalité $(\ref{7inegalite_presque_finale})$ devient alors 
 \begin{align*}
 (K -1)(-\beta_{i,N_i+\ell+1} + 1) + \frac{\alpha_{i,N_i +v_i}}{\alpha_{i,N_i + \ell} } 
 &\leq -\beta_{i,N_i+\ell+1} \frac{K-1}{\min\limits_{ \ell \in \llbracket 1, m+1 \rrbracket} (\beta_{i,\ell})} + K-1 \\
 &\leq (K-1) \left(1 - \frac{\beta_{i,N_i+\ell+1} }{\min\limits_{ \ell \in \llbracket 1, m+1 \rrbracket} (\beta_{i,\ell})} \right) \\
 &\leq 0.
 \end{align*}
 \\ On reprend alors $(\ref{7inega_deuxieme_terme})$ et on a 
 \begin{align*}
 \theta^{-\alpha_{i,N_i+\ell+1} + \alpha_{i,N_i+\ell}}H(C) &\leq H(C)^{\frac{\varepsilon}{2}(\frac{-\alpha_{i,N_i + \ell +1} + \alpha_{i,N_i + \ell }}{\alpha_{i,N_i + v_i}})}. 
 \end{align*}
 Or on a $\beta_{i,j} \geq 2$ pour tout $j$, et $v_i \leq m+1 $ par le lemme~\ref{7lem_vj=m+1} donc en utilisant $(\ref{7_2bhypot_m=beta_m+1})$ on a $$\frac{\alpha_{i,N_i + \ell +1} - \alpha_{i,N_i + \ell }}{\alpha_{i,N_i + v_i}} = \frac{\beta_{i,N_i+ \ell +1} -1 }{\beta_{i,N_i + \ell +1} \ldots \beta_{i,N_i +v_i}} \geq \frac{1}{K_{i,v_i}} \geq \frac{1}{K_{i,m+1}} \geq \frac{1}{K_{1, m+1}}.$$ En posant $\cons \label{7cons_espsK} = \frac{1}{2K_{1, m+1}}$ on a alors 
 \begin{align}\label{7inega_second_terme}
 \theta^{-\alpha_{i,N_i+\ell+1} + \alpha_{i,N_i+\ell}}H(C) &\leq H(C)^{-c_{\ref{7cons_espsK}}\varepsilon}.
 \end{align}
 \bigskip \\
 En reprenant l'inégalité $(\ref{7maj_DNiell_final})$ et les estimations $(\ref{7inega_premier_terme})$ et $(\ref{7inega_second_terme})$ on a alors 
 \begin{align*}
 D_{N_{i}+\ell, i } &\leq c_{\ref{7cons_maj_DNell}}(H(C)^{-\frac{\varepsilon}{2}} + H(C)^{-c_{\ref{7cons_espsK}}\varepsilon})
 \end{align*}
 et donc si $H(C)$ est assez grand (en fonction de $\varepsilon, c_{\ref{7cons_maj_DNell}} $ et $c_{\ref{7cons_espsK}}$) on a pour tout $i \in \llbracket 1, k \rrbracket $
 \begin{align*}
 \forall \ell \in \llbracket 0, v_i-1 \rrbracket, \quad \| X_{N_i +\ell, i} \wedge Z_1 \wedge \ldots \wedge Z_e \| = D_{N_{i}+\ell, i } < 1.
 \end{align*}
 Le lemme~\ref{2lem_X_in_B} donne alors $X_{i,N_i+\ell} \in C$ pour tous $i \in \llbracket 1, k \rrbracket $ et $\ell \in \llbracket 0, v_i-1 \rrbracket$. 
 Alors $C^J_{N} \subset C$ et par égalité des dimensions $C = C^J_{N} $.

\end{preuve}

\begin{cor}\label{7cor_maj_exposant}
 On a $$ \mu_n(A_J| e)_{k - g(A_J,e)} \leq \frac{1}{\sum\limits_{i = f+1}^k \frac{1}{K_{i,v_{i} }} }.$$ 
\end{cor}

\begin{preuve}
 On suppose le contraire par l'absurde : soit $\varepsilon > 0$ tel que
 \begin{align*}
 \mu_n(A_J| e)_{k - g(A_J,e)} \geq \frac{1}{\sum\limits_{j= f+1}^k \frac{1}{K_{j,v_{j} }} } + 2\varepsilon.
 \end{align*}
 Il existe alors une infinité d'espaces $C$ rationnels de dimension $e$ tels que 
 \begin{align}\label{7maj_dunepart}
 \psi_{k-g}(A_J,C) \leq H(C)^{- \frac{1}{\sum\limits_{j = f+1}^k \frac{1}{K_{j,v_{j} }} } - \varepsilon}. 
 \end{align}
D'après le lemme~\ref{7lem_meilleur_espaces}, il existe $N \in (\Nx)^k$ tel que $C = C^J_{N}$ si $H(C)$ est assez grand.
\bigskip \\
D'autre part le lemme~\ref{7lem_prox_AJ_CN_H(CN)} donne :
\begin{align}\label{7min_dautrepart}
 \psi_{k-g}(A_J, C^J_{N}) \geq c_{\ref{7cons_prox_CN_AJ}} H(C^J_{N}) ^{ \frac{ - \min\limits_{j = f+1}^k \alpha_{j,M_j+1}} {\sum\limits_{j = f+1 }^k \alpha_{j,N_j} } }
\end{align}
 avec $ \cons \label{7cons_prox_CN_AJ} >0 $ indépendante de $N $. 
 Les inégalités $(\ref{7maj_dunepart}) $ et $ (\ref{7min_dautrepart})$ regroupées donnent alors
 \begin{align}\label{7ineg_comp_exposant}
 \frac{\min\limits_{j = f+1}^k \alpha_{j,N_j+v_j}} {\sum\limits_{j = f+1 }^k \alpha_{j,N_j} } \geq  \frac{1}{\sum\limits_{j = f+1}^k \frac{1}{K_{j,v_{j} }} } + \frac{\varepsilon}{2} 
 \end{align}
 en faisant tendre la hauteur vers $+ \infty$. 
\\Soit $j \in \llbracket f+ 1, k \rrbracket. $ On a :
\begin{align*}
 \frac{\min\limits_{i=f+1}^k {\alpha_{i,N_i +v_i} }}{K_{j,v_j} } \leq \frac{\alpha_{j,N_j+v_j}}{K_{j,v_j}} \leq \frac{ \alpha_{j,N_j} \beta_{j,N_j + 1} \ldots \beta_{j,N_j+v_j}}{K_{j,v_j}} \leq \alpha_{j,N_j}
 \end{align*}
car $K_{j,v_j}$ est le maximum des produits de $v_j$ termes consécutifs parmi les $\beta_{j,\ell}$. On somme maintenant sur les $j $ et on a 
\begin{align*}
  (\min\limits_{j=f+1}^k {\alpha_{j,N_j +v_j} } )\sum\limits_{j = f+1}^k \frac{1}{K_{j,v_j}} \leq \sum\limits_{j = f+1}^k \alpha_{j,N_j}. 
\end{align*}
Alors $\frac{\min\limits_{j = f+1}^k \alpha_{j,N_j+v_j}} {\sum\limits_{j = f+1 }^k \alpha_{j,N_j} } \leq  \frac{1}{\sum\limits_{j = f+1}^k \frac{1}{K_{j,v_{j} }} }$ ce qui contredit $(\ref{7ineg_comp_exposant})$ puisque $\varepsilon >0$, et termine la preuve du corollaire.

\end{preuve}

Les corollaires \ref{7cor_min_expos} et \ref{7cor_maj_exposant} concluent la preuve de la propriété~\ref{7prop_int}.


\section{Calcul final des exposants}\label{7section_calcul_expo}

On effectue la preuve du théorème~\ref{7theo_principale} en utilisant la propriété~\ref{7prop_int}, et le théorème~\ref{theo_somme_sev} pour en déduire les exposants $\mu_n(A|e)_{k -g(A,e)}$.

\bigskip 

Pour $j \in \llbracket 1,d \rrbracket$ on considère l'espace 
\begin{align*}
 R_j = \{0\}^{(j-1)(m+1)} \times \R^{m+1} \times \{0 \}^{(d-j)(m+1)} \subset \R^n.
\end{align*}
Les sous-espaces $R_j$ sont rationnels et en somme directe. On a donc $\bigoplus\limits_{j=1}^d R_j \subset \R^n$ et on remarque que pour tout $j \in \llbracket 1,d \rrbracket$ 
\begin{align*}
 \Vect(Y_j) \subset R_j.
\end{align*}
On applique le théorème~\ref{theo_somme_sev} avec les $R_j$ et les droites~$A_j = \Vect(Y_j)$. 

 \bigskip 
 Soit $e \in \llbracket 1, n-1 \rrbracket$ et $k \in \llbracket 1 + g(A,e), \min(d,e) \rrbracket$ tels que $e < k(m+1)$. Pour tout sous-ensemble $J \subset \llbracket 1,d \rrbracket$ de cardinal $k $, par le lemme~\ref{7lem_AJ_irr} 
 \begin{align*}
 A_J = \Vect_{j \in J}(Y_j) \text{ est } (e, k-g(A_J,e))\tir\text{irrationnel.}
 \end{align*} 
La première partie du théorème~\ref{theo_somme_sev} appliquée avec $k' = k - g(A,e)$ donne alors 
 $$A \in \II_n(d,e)_{k-g(A,e)}.$$ 
 On calcule maintenant l'exposant. Le même  théorème~donne 
 \begin{align*}
 \mu_n(A|e)_{k - g(A,e)} = \mu_n(A|e)_{k'} = \max\limits_{J \in \PP(k' + g(A,e),d) } \mu_n(A_J|e)_{k' + g(A,e) - g(A_J,e)}
 \end{align*}
 avec $\PP(\ell, d)$ l'ensemble des parties de cardinal $\ell$ de $ \llbracket 1,d \rrbracket$, et donc 
 \begin{align}\label{7calcul_exposants_presque_fin}
 \mu_n(A|e)_{k - g(A,e)} = \max\limits_{J \in \PP(k,d) } \mu_n(A_J|e)_{k - g(A_J,e)} = \max\limits_{J \in \PP(k,d) } \frac{1}{\sum\limits_{q = 1+ f(e,mk)}^{k} \frac{1}{K_{j_q,v_q }} }
 \end{align}
 par la propriété~\ref{7prop_int} appliquée avec $J = \{ j_1 < \ldots< j_{k} \}$.

Il reste alors à montrer le lemme suivant pour terminer la preuve du théorème~\ref{7theo_principale}.

\begin{lem}
 On a 
 $$\max\limits_{J = \{ j_1 < \ldots< j_{k} \} } \frac{1}{\sum\limits_{q = 1+ f(e,mk)}^{k} \frac{1}{K_{j_q,v_q }} } = \frac{1}{\sum\limits_{q = 1 + f(e,mk)}^{k} \frac{1}{K_{q + d- k +1,v_q }} }. $$
\end{lem}

\begin{preuve}
 On pose $J$ l'ensemble $\{d-k+1, \ldots, d\}$ qui correspond donc aux quantités $j_q = q + d- k $ pour $q \in \llbracket 1, k \rrbracket$. On a donc 
 \begin{align*}
 \max\limits_{J = \{ j_1 < \ldots< j_{k} \} } \frac{1}{\sum\limits_{q = 1+ f(e,mk)}^{k} \frac{1}{K_{j_q,v_q }} } \geq\frac{1}{\sum\limits_{q = 1 + f(e,mk)}^{k} \frac{1}{K_{q + d- k,v_q }} }.
 \end{align*}
De plus, les $K_{i,v}$ sont croissants en $i$ par $(\ref{7hypothèse_prop_princi3})$ donc en particulier 
\begin{align*}
\forall 1 \leq j_1 < \ldots < j_k \leq d, \quad \forall q \in \llbracket 1,k \rrbracket, \quad \forall v \in \llbracket 1, m \rrbracket, \quad K_{j_q, v } \leq K_{q +d-k, v}.
\end{align*}
En particulier 
\begin{align*}
 \max\limits_{J = \{ j_1 < \ldots< j_{k} \} } \frac{1}{\sum\limits_{q = 1+ f(e,mk)}^{k} \frac{1}{K_{j_q,v_q }} } \leq \frac{1}{\sum\limits_{q = 1 + f(e,mk)}^{k} \frac{1}{K_{q + d- k,v_q }} }
\end{align*}
 ce qui termine la preuve du lemme.

\end{preuve}

Soient $e \in \llbracket 1, n-1 \rrbracket$ et $k \in \llbracket 1 + g(A,e), \min(d,e) \rrbracket$ tels que $e < k(m+1)$, on déduit de ce lemme et de $(\ref{7calcul_exposants_presque_fin})$ :
\begin{align*}
 \mu_n(A|e)_{k - g(A,e)} = \frac{1}{\sum\limits_{q = 1 + f(e,mk)}^{k} \frac{1}{K_{q + d- k,v_q }} }.
\end{align*}
Cela termine la preuve du  théorème~\ref{7theo_principale}.

%% file: Chapitres/8_Independance_exposants.tex
\chapter{Indépendance algébrique d'exposants diophantiens}\label{chap9}

Les théorèmes démontrés dans les chapitres précédents permettent d'obtenir des résultats sur certains spectres joints. On s'intéresse ici aux propriétés d'indépendance algébrique des exposants composant ces spectres.

On fixe $n \in \Nx$ et $d \in \llbracket 1, n -1 \rrbracket$. On étudie la famille de fonctions
\begin{align*}
 (\mu_n(\cdot|e)_{k-g(d,e,n)})_{(e,k) \in U}
\end{align*}
pour $ U$ un sous ensemble de $$ V_{d,n} = \lbrace (e,k) \mid e \in \llbracket 1, n-1\rrbracket, k \in \llbracket 1 +g(d,e,n), \min(d,e) \rrbracket \rbrace.$$
Dans ce chapitre, on raisonne toujours à $d$ et $n$ fixés. Les fonctions que l'on étudie, notées $\mu_n(\cdot|e)_{k-g(d,e,n)}$, sont donc définies sur $\II_n(d, e)_{k-g(d,e,n)}$ l'ensemble des sous-espaces de $\R^n$ de dimension $d$ qui sont $(e,k-g(d,e,n))\tir$irrationnels. 

\bigskip
On pose pour $U \subset V_{d,n}$ l'ensemble $\II_U = \bigcap\limits_{(e,k)\in U} \II_n(d, e)_{k-g(d,e,n)}$ et l'application \begin{align*}
 M_U : \left| \begin{array}{ccc}
 \II_U &\longrightarrow & \R^{U} \\
 A & \longmapsto & (\mu_n(A|e)_{k-g(d,e,n)})_{(e,k) \in U}.
 \end{array}\right. 
\end{align*}

\bigskip 

On conjecture que la famille de fonctions $\left(\mu_n(\cdot|e)_{k-g(d,e,n)}\right)_{(e,k) \in V_{d,n}}$ est algébriquement indépendante sur $\R$. En direction de ce résultat on obtient le  théorème~suivant.

\begin{theo}\label{9theo_general}
Le degré de transcendance sur $\R$ de la famille $\left(\mu_n(\cdot|e)_{k-g(d,e,n)}\right)_{(e,k) \in V_{d,n}}$ est au moins égal à $n-d$.
\end{theo}

Cet énoncé signifie qu'il existe une sous-famille $U \subset V_{d,n}$ de cardinal $n-d$ telle que la famillle $\left(\mu_n(\cdot|e)_{k-g(d,e,n)}\right)_{(e,k) \in U}$ soit algébriquement indépendante sur $\R$. 
Dans la suite on démontre en effet des théorèmes plus précis sur de telles sous-familles.

\bigskip
Le  théorème~\ref{9theo_general} est directement une conséquence du théorème~\ref{9prop_premier_angle} énoncé ci-après. Celui-ci étend le résultat du corollaire~\ref{1cor_ind_alg_exposant_dim1} qui donne seulement le cas particulier $d = 1$ de ce théorème. 

\begin{theo}\label{9prop_premier_angle}
L'image de $(\mu_n(\cdot|1)_1, \ldots, \mu_n(\cdot|n-d)_1)$ contient un ouvert non vide de $\R^{n-d}$ et la famille de fonctions $(\mu_n(\cdot|e)_{1} )_{e \in \llbracket 1, n-d \rrbracket}$ est algébriquement indépendante sur $\R$.
\end{theo}

\bigskip
On étudie aussi des familles d'exposants correspondant au dernier angle ; c'est l'objet des deux théorèmes suivants. 

Le  théorème~\ref{9lem_sousfam_pour_theo2} est assez naturel, il a pour sujet la famille des exposants correspondant à $e + d \leq n$ et au dernier angle. 
\begin{theo}\label{9lem_sousfam_pour_theo2}
On suppose que $d$ divise $n$. 
\\Alors l'image de $(\mu_n(\cdot|1)_{\min(d,1)}, \ldots, \mu_n(\cdot|n-d)_{\min(d,n-d)})$ contient un ouvert non vide de $\R^{n-d}$ et la famille $(\mu_n(\cdot|e)_{\min(d,e)} )_{e \in \llbracket 1, n-d \rrbracket}$ est algébriquement indépendante sur $\R$. 
 
\end{theo}

\bigskip
Le  théorème~\ref{9lem_sousfam_pour_theo3} a pour sujet la famille des exposants correspondant au dernier angle dans le cas où  $j = d$ ; autrement dit $e \in \llbracket d, n-1 \rrbracket$. En particulier cette famille contient certains exposants vérifiant $g(d,e,n) = 0$ et d'autres $g(d,e,n) > 0$.

\begin{theo}\label{9lem_sousfam_pour_theo3}
On suppose que $d$ divise $n$. 
\\Alors l'image de $(\mu_n(\cdot|d)_{d}, \ldots, \mu_n(\cdot|n-1)_{d})$ contient un ouvert non vide de $\R^{n-d}$ et la famille $(\mu_n(\cdot|e)_{d} )_{e \in \llbracket d, n-1 \rrbracket}$ est algébriquement indépendante sur $\R$. 
\end{theo}

\bigskip

\begin{req}\label{9req_ouvert_implique_alg_ind}
Soit $U \subset V_{d,n}$. Si l'image de $M_U$ contient contient un ouvert non vide de $\R^{U}$, alors  la famille $(\mu_n(\cdot|e)_{k-g(d,e,n)})_{(e,k) \in U} $ est algébriquement indépendante sur $\R$. En effet on aurait sinon une relation polynomiale reliant ces fonctions et l'image de $M_U$ serait contenue dans une hypersurface algébrique de $\R^U$ qui est d'intérieur vide.  

Dans la suite les preuves s'attachent alors uniquement à montrer que les images des fonctions considérées contiennent un ouvert non vide.

\end{req}

\bigskip
Les théorèmes et propriété de ce chapitre se déduisent tous des constructions réalisées dans les chapitres \ref{chap6} et \ref{chap7}.
\\Le  théorème~\ref{9prop_premier_angle} est prouvé dans la section~\ref{9section_premier_angle} : on utilise le  théorème~\ref{6theo_principal} pour construire des espaces $A$ dont les exposants sont dans un ouvert non vide de $\R^{n-d}$. 
\\ Dans la section~\ref{9section_conditions_suffisantes} de ce chapitre, on énonce et montre deux propriétés donnant des conditions suffisantes sur $U$ pour qu'il existe un ouvert non vide dans l'image de $M_U$. Leurs preuves utilisent le  théorème~\ref{7theo_principale} Ces propriétés sont utilisées dans la section~\ref{9section_dernier_angle_theo} pour montrer les théorèmes \ref{9lem_sousfam_pour_theo2} et \ref{9lem_sousfam_pour_theo3}. Enfin on montre la propriété~\ref{9lem_sousfam_pour_theo} dans la section~\ref{9section_autre_resultat}, on y expose une nouvelle famille d'exposants algébriquement indépendants sur $\R$ qui fait notamment apparaître des exposants correspondants à tous les angles $(k \in \llbracket 1,d \rrbracket$).

\newpage
\section{Résultat pour le premier angle}\label{9section_premier_angle}
 Dans cette section, on fait la preuve du théorème~\ref{9prop_premier_angle} en utilisant les résultats du chapitre \ref{chap6}.

\begin{preuve}[ (théorème~\ref{9prop_premier_angle})]
 Le théorème~\ref{6theo_principal} donne un espace $A \in \II_n(d,n-d)_1$ vérifiant
 \begin{align*}
 \forall e \in \llbracket 1, n-d \rrbracket, \quad \mu_n(A|e)_1 = \gamma_e
 \end{align*}
 avec$(\gamma_1, \ldots, \gamma_{n-d}) \in \R^{n-d}$ vérifiant $\gamma_1 \geq C_d $ ainsi que 
 \begin{align*}
 \forall i \in \llbracket 2, n-d \rrbracket,& \quad \gamma_i \geq C_d \gamma_{i-1}, \\
 \forall (i,j) \in \llbracket 1, n-d-1 \rrbracket^2,& \quad i+j \leq n-1 \Longrightarrow \gamma_{i+j} \leq \gamma_i \gamma_j. 
 \end{align*}
En particulier l'image de $(\mu_n(\cdot|1)_1, \ldots, \mu_n(\cdot|n-d)_1)$ contient un l'ouvert non vide $O$ constitué de l'ensemble des $(\gamma_1, \ldots, \gamma_{n-d}) \in \R^{n-d}$ vérifiant $\gamma_1 > C_d $ ainsi que 
 \begin{align*}
 \forall i \in \llbracket 1, n-d-1 \rrbracket,& \quad C_d \gamma_{i} < \gamma_{i+1} < \min\limits_{j \in \llbracket 1, i \rrbracket} \gamma_{i+1-j}\gamma_{j}. 
 \end{align*}
 L'image de $(\mu_n(\cdot|1)_1, \ldots, \mu_n(\cdot|n-d)_1)$ contient donc un ouvert non vide et en particulier la famille de fonctions $$\mu_n(\cdot|1)_1, \ldots, \mu_n(\cdot|n-d)_1 $$ est algébriquement indépendante sur $\R$ par la remarque~\ref{9req_ouvert_implique_alg_ind}.
 
\end{preuve}


\section{Conditions sufisantes}\label{9section_conditions_suffisantes}

On développe ici deux lemmes techniques qui donnent des conditions suffisantes sur $U \subset V_{d,n}$ pour que $M_U(\II_U)$ contienne un ouvert non vide et que la famille d'exposants associée soit algébriquement indépendante sur $\R$. 

\bigskip

Dans cette section, on suppose que $n =d (m+1)$ avec $m \in \Nx$. On rappelle la notation $\PP(\llbracket 1,d \rrbracket \times \llbracket 1, m \rrbracket)$ pour l'ensemble des parties de $ \llbracket 1,d \rrbracket \times \llbracket 1, m \rrbracket$, et on pose 
\begin{align*}
 \chi : \left| \begin{array}{ccc}
 V_{d,n} &\longrightarrow & \PP(\llbracket 1,d \rrbracket \times \llbracket 1, m \rrbracket )\\
 (e,k) & \longmapsto & \llbracket 1 + f + d - k, u+d-k \rrbracket \times \llbracket 1, v+1 \rrbracket \cup \llbracket u+ 1 + d - k, d \rrbracket \times \llbracket 1, v \rrbracket 
 \end{array}\right.
\end{align*}
où $v$ et $u$ sont respectivement le quotient le reste de la division euclidienne de $e$ par $k$ et $f = f(e,mk)= \max(0, e-mk)$. 

De plus, on rappelle que dans le chapitre~\ref{chap7}, on a associé à $(e,k) \in V_{d,n}$ des quantités $v_1(e,k), \ldots, v_k(e,k)$ définies par :
\begin{align}\label{9eq_def_vj}
 v_q(e,k) &= \left\{
 \begin{array}{ll}
 v + 1 &\text{ si } q \in \llbracket 1, u \rrbracket \\
 v &\text{ si } q \in \llbracket u+1, k \rrbracket.
 \end{array}
 \right. 
\end{align}
Un couple $(q, \ell) \in \llbracket 1 ,d \rrbracket \times \llbracket 1,m \rrbracket$ appartient alors à $\chi(e,k)$ si et seulement si, on a 
\begin{align}\label{9apparaisant_dans_omega}
 q \geq 1 +f +d -k \text{ et } \ell \leq v_{q+k-d}(e,k).
\end{align}

\bigskip

\begin{lem}\label{9lem_technique}
 Soit $U \subset V_{d,n}$ avec $\#U \leq dm$ vérifiant pour tout $(e,k) \in U$
 \begin{align}
 e < k(m+1).\label{9hypothese_sur_U2}
 \end{align}
 On suppose de plus qu'il existe une relation d'ordre $<_U$ sur $U $ telle que 
 \begin{align}\label{9hypothese_sur_U} 
 \forall (e,k) \in U, \quad \chi(e,k) \smallsetminus \bigcup\limits_{(e ',k') <_U (e,k)} \chi(e',k') \neq \emptyset. 
 \end{align}
 Alors l'image de l'application $M_U$ contient un ouvert non vide de $\R^{U}$. 
 En particulier la famille $(\mu_n(\cdot|e)_{k-g(d,e,n)})_{(e,k) \in U}$ est algébriquement indépendante sur $\R$.
\end{lem}

Le lemme~\ref{9lem_technique} est en fait un cas particulier du lemme suivant. 

\begin{lem}\label{9lem_technique_faible}
 Soit $U \subset V_{d,n}$ avec $\#U \leq dm$ vérifiant pour tout $(e,k) \in U$
 \begin{align*}
 e < k(m+1).
 \end{align*}
 Pour $(e,k) \in U$ et $\beta = (\beta_{i,\ell})_{i \in \llbracket 1,d \rrbracket, \ell \in \llbracket 1, m \rrbracket} \in (\R_+^{*})^{dm}$, on pose $$\Omega_{(e,k)}(\beta) = \sum\limits_{q = 1 + f(e, mk)}^k \beta_{q+d-k,1} \ldots \beta_{q+d-k,v_q(e,k)} = \sum\limits_{q = 1 + f(e, mk) +d - k}^d \beta_{q,1} \ldots \beta_{q,v_{q+k-d}(e,k)}, $$
 ce qui définit une application $\Omega : (\R_+^{*})^{dm} \to \R^U, \beta \mapsto (\Omega_{(e,k)}(\beta))_{(e,k)\in U} $. \\
On suppose que la matrice jacobienne de $\Omega$ en $\beta$, $ J_{\Omega}(\beta) \in \MM_{\#U, dm}(\R)$ est de rang $\#U$ pour tout $\beta \in (\R_+^{*})^{dm}$. 
\\ Alors l'image de l'application $M_U$ contient un ouvert non vide de $\R^{U}$. 
 En particulier la famille $(\mu_n(\cdot|e)_{k-g(d,e,n)})_{(e,k) \in U}$ est algébriquement indépendante sur $\R$. 
\end{lem}

\begin{preuve}[ (lemme~\ref{9lem_technique_faible})]
D'après la remarque~\ref{9req_ouvert_implique_alg_ind}, il suffit de montrer que l'image de $M_U$ contient un ouvert non vide de $\R^{\#U}$. On pose alorsl'application suivante 
\begin{align*}
 \Omega' : \left| \begin{array}{ccc}
 (\R_+^{*})^{dm} &\longrightarrow & \R^{U}\\
 (\beta_{i, \ell})_{i \in \llbracket 1, d \rrbracket, \ell \in \llbracket 1, m \rrbracket} &\mapsto & \left(\frac{1}{\sum\limits_{q = 1 + f(e, mk)}^k (\beta_{q+d-k,1} \ldots \beta_{q+d-k,v_q(e,k)})^{-1}} \right)_{(e,k) \in U}.
 \end{array}\right. 
\end{align*}

On va appliquer le théorème~\ref{7theo_principale} pour construire des espaces $A$ de dimension $d$ dont la famille d'exposants $(\mu_n(A|e)_{k } )_{(e,k) \in U}$ est un point donné de l'image de $\Omega'$. 

Pour cela, on pose $H' \subset (\R_+^{*})^{dm}$ l'ensemble $(\beta_{i, \ell})_{i \in \llbracket 1, d \rrbracket, \ell \in \llbracket 1, m \rrbracket}$ vérifiant les hypothèses 
\begin{align}
 \forall i \in \llbracket 1,d \rrbracket,& \quad \beta_{i,1} > \ldots > \beta_{i, m}, \label{9hyp_croissance_beta} \\
 \min\limits_{\ell \in \llbracket 1, m \rrbracket}(\beta_{1,\ell}) > {(3d)^{\frac{c_{\ref{7cons_petite_hyp_theoc2}}}{c_{\ref{7cons_petite_hyp_theoc2}}-1}}} &\text{ et } \min\limits_{\ell \in \llbracket 1, m \rrbracket}(\beta_{1,\ell})^{c_{\ref{7cons_petite_hyp_theoc2c1}}} > \max\limits_{\ell \in \llbracket 1, m \rrbracket}(\beta_{1,\ell})^{c_{\ref{7cons_petite_hyp_theoc2}}}, \nonumber
\end{align}
ainsi que pour tout $i \in \llbracket 1,d-1 \rrbracket$ : 
\begin{align*}
 \min\limits_{\ell \in \llbracket 1, m \rrbracket}(\beta_{i,\ell})^{c_{\ref{7cons_petite_hyp_theoc2c1}}} > \max\limits_{\ell \in \llbracket 1, m \rrbracket}(\beta_{i+1,\ell})
 \\
 \text{ et } \min\limits_{\ell \in \llbracket 1, m \rrbracket}(\beta_{i+1,\ell}) > \max\limits_{\ell \in \llbracket 1, m \rrbracket}(\beta_{i,\ell})^{c_{\ref{7cons_petite_hyp_theoc2}}}
\end{align*}
avec $c_{\ref{7cons_petite_hyp_theoc2c1}}, c_{\ref{7cons_petite_hyp_theoc2}}$ les constantes définies dans le  théorème~\ref{7theo_principale}. L'ensemble $H'$ ainsi défini est un ouvert non vide de $(\R_+^{*})^{dm}$ car les inégalités ci-dessus sont strictes.

\bigskip
Pour $\beta \in H'$, le théorème~\ref{7theo_principale} donne un espace $A$ de dimension $d$ dans $\R^{n}$ tel que pour tous $e \in \llbracket 1, n-1\rrbracket $ et $k \in \llbracket 1 + g(d,e, n), \min(d,e) \rrbracket$ vérifiant $e < k (m+1) $ on a 
\begin{align*}
 A \in \II_n(d,e)_{k-g(d,e,n)} \text{ et } \mu_n(A|e)_{k-g(d,e,n)}
 &= \frac{1}{\sum\limits_{q = 1 + f(e, mk)}^k \frac{1}{K_{q+d-k,v_q(e,k)} }} 
\end{align*}
avec ici, du fait de l'hypothèse $(\ref{9hyp_croissance_beta})$
\begin{align}\label{9eq_def_Kv}
 \forall i \in \llbracket 1, d \rrbracket, \forall v \in \llbracket 1, m \rrbracket, \quad K_{i,v} = \beta_{i,1} \ldots \beta_{i,v}.
\end{align}

Ce résultat donne alors $\Omega'(H') \subset M_U(\II_U)$. Pour prouver le lemme~\ref{9lem_technique_faible} on va montrer que $\Omega'(H')$ contient un ouvert non vide. Quitte à composer $\Omega'$ au départ et à l'arrivée par l'application qui inverse chaque coordonnée il suffit de montrer que l'application de classe $\mathcal{C}^1$ sur $H = \lbrace (\beta_{i, \ell})_{i \in \llbracket 1, d \rrbracket, \ell \in \llbracket 1, m \rrbracket}, (\frac{1}{\beta_{i, \ell}})_{i \in \llbracket 1, d \rrbracket, \ell \in \llbracket 1, m \rrbracket} \in H' \rbrace$ définie par 
\begin{align*}
 \Omega : \left| \begin{array}{ccc}
 H &\longrightarrow & \R^{U}\\
 (\beta_{i, \ell})_{i \in \llbracket 1, d \rrbracket, \ell \in \llbracket 1, m \rrbracket} &\mapsto & \left(\sum\limits_{q = 1 + f(e, mk)}^k \beta_{q+d-k,1} \ldots \beta_{q+d-k,v_q(e,k)} \right)_{(e,k) \in U}
 \end{array}\right.
\end{align*}
a un ouvert non vide dans son image.
\\Comme par hypothèse $J_{\Omega}(\beta) $ est de rang $\#U$ pour tout $\beta \in H \subset (\R_+^{*})^{dm}$, le  théorème~de l'image ouverte assure alors que $\Omega(H)$ contient un ouvert non vide ce qui conclut la preuve du lemme.

\end{preuve}

\begin{preuve}[ (lemme~\ref{9lem_technique})]
Chaque $\Omega_{(e,k)}$, défini dans les hypothèses du lemme~\ref{9lem_technique_faible}, est un polynôme en les $\beta_{i, \ell}$. On remarque en utilisant $(\ref{9apparaisant_dans_omega})$ que $\chi(e,k)$ est en fait l'ensemble des indices des $\beta_{i,\ell}$ \og apparaissant\fg{} dans $\Omega_{(e,k)}(\beta)$. Cela donne en particulier 
\begin{align}\label{9rel_polynome_nul_ounon}
 \frac{\partial \Omega_{(e,k)}}{\partial \beta_{q,\ell}} = \left\{ \begin{array}{cl}
 \prod\limits_{\underset{p \neq \ell}{1 \leq p \leq v_{q +k -d}(e,k)}} \beta_{q,p} & \text{ si } (q,\ell) \in \chi(e,k)\\
 0 & \text{ sinon } 
 \end{array} \right.
\end{align}
D'après l'hypothèse $(\ref{9hypothese_sur_U})$, on peut poser $U =\lbrace \theta_1, \ldots, \theta_{\#U} \} $ de sorte que ces quantités respectent l'ordre $<_U$ : 
\begin{align*}
 \forall i < j, \quad \theta_i <_U \theta_j
\end{align*}
et que, pour tout $ j \in \llbracket 2, \#U \rrbracket$, il existe $(q_j,\ell_j) \in \chi(\theta_j) \smallsetminus \bigcup\limits_{i= 1}^{j-1} \chi(\theta_i)$. 
\\On a en particulier, en utilisant $(\ref{9rel_polynome_nul_ounon})$, pour tout $\beta \in (\R_+^{*})^{dm}$
\begin{align}\label{9valeur_pour_trianguler}
 \forall i < j, \quad \frac{\partial \Omega_{\theta_i}}{\partial \beta_{q_j,\ell_j}} =0 \text{ et } \frac{\partial \Omega_{\theta_j}}{\partial \beta_{q_j,\ell_j}} = \prod\limits_{\underset{p \neq \ell_j}{1 \leq p \leq v_{q_j +k -d}(\theta_j)}} \beta_{q_j,p} \neq 0. 
\end{align}
On note $J_{\Omega}(\beta) = \left(\frac{\partial \Omega_{\theta_i}}{\partial \beta_{j,\ell}}\right)_{i \in \llbracket 1, \#U \rrbracket, (j,\ell) \in \llbracket 1, d \rrbracket \times \llbracket 1,m \rrbracket} \in \MM_{\#U, dm}(\R)$, la matrice jacobienne de $\Omega$ en $\beta \in (\R_+^{*})^{dm}$. Ici l'indice $i$ correspond aux lignes de la matrice et $(j,\ell)$ aux colonnes, en ordonnant $\llbracket 1, d \rrbracket \times \llbracket 1,m \rrbracket$ selon l'ordre lexicographique usuel. 
\\On va montrer que $J_{\Omega}(\beta) $ est de rang maximal $\#U$, pour ensuite appliquer le lemme~\ref{9lem_technique_faible}.
On extrait de $J_{\Omega}(\beta) $ (quitte à intervertir des colonnes) la matrice $$G_U = \left(\frac{\partial \Omega_{\theta_i}}{\partial \beta_{q_j,\ell_j}}\right)_{i \in \llbracket 1, \#U \rrbracket, j \in \llbracket 1, \#U \rrbracket} \in \MM_{\#U}(\R).$$
D'après $(\ref{9valeur_pour_trianguler})$, on a 
\begin{align*}
 G_U = \begin{pmatrix}
 \frac{\partial \Omega_{\theta_1}}{\partial \beta_{q_1,\ell_1}} & 0 & 0 & \cdots & 0 \\ 
 \frac{\partial \Omega_{\theta_2}}{\partial \beta_{q_1,\ell_1}} & \frac{\partial \Omega_{\theta_2}}{\partial \beta_{q_2,\ell_2}} & 0 & \cdots & 0 \\
 \vdots & & \ddots & \ddots & \vdots \\
 \frac{\partial \Omega_{\theta_{\#U -1 }}}{\partial \beta_{q_1,\ell_1}} & \frac{\partial \Omega_{\theta_{\#U -1 }}}{\partial \beta_{q_2,\ell_2}} & \cdots &\frac{\partial \Omega_{\theta_{\#U -1 }}}{\partial \beta_{q_{\#U -1},\ell_{\#U -1 }}} & 0 \\
 \frac{\partial \Omega_{\theta_{\#U }}}{\partial \beta_{q_1,\ell_1}} & \frac{\partial \Omega_{\theta_{\#U }}}{\partial \beta_{q_2,\ell_2}} & & \cdots &\frac{\partial \Omega_{\theta_{\#U }}}{\partial \beta_{q_{\#U },\ell_{\#U }}} 
 \end{pmatrix}.
\end{align*}
La matrice $G_U$ est donc triangulaire inférieure et ses coefficients diagonaux sont non nuls d'après $(\ref{9valeur_pour_trianguler})$, elle est donc inversible.

On a donc montré que $J_{\Omega}(\beta)$ est de rang $\#U$ pour tout $\beta \in H$. Le lemme~\ref{9lem_technique_faible} donne alors que l'image de l'application $M_U$ contient un ouvert non vide de $\R^{U}$ ce qui conclut la preuve du lemme~\ref{9lem_technique}.

\end{preuve}

\section{Résultats pour le dernier angle}\label{9section_dernier_angle_theo}

Dans cette section, on prouve les théorèmes \ref{9lem_sousfam_pour_theo2} et \ref{9lem_sousfam_pour_theo3} en utilisant les résultats de la section précédente. On suppose dans cette section, que $d $ divise $n$ et on écrit $n = d(m+1)$ avec $m \in \Nx$.

\begin{preuve}[ (théorème~\ref{9lem_sousfam_pour_theo2})]
On pose
\begin{align*}
 U = \{ (e,\min(d,e)) \mid e \in \llbracket 1, n-d \rrbracket \}.
\end{align*}
On va montrer que $U$ vérifie l'hypothèse $(\ref{9hypothese_sur_U})$ du lemme~\ref{9lem_technique}.
\\On remarque tout d'abord que pour tout $(e,k) \in U$, on a $g(d,e,n) = 0$ et $f(e,mk) = 0$ car $e \in \llbracket 1, dm \rrbracket$ et $dm = n-d$. 
 \\ On pose l'ordre suivant sur $U$ :
 \begin{align*}
 (e',k') <_U (e,k) \Longleftrightarrow e' < e.
 \end{align*}
 Cela définit bien un ordre car les premières coordonnées des couples de $U$ sont deux à deux distinctes.
 \\ Soit $(e,k) \in U$. Si $e \leq d$ alors $k =e $ et 
 \begin{align*}
 \chi(e,k) = \llbracket 1+d-e, d \rrbracket \times \{ 1 \}.
 \end{align*}
 De même pour tout $(e',k') <_U (e,k)$, $\chi(e',k') = \llbracket 1+d-e', d \rrbracket \times \{ 1 \} $ et donc comme $ 1+ d -e' > 1+ d -e $, on a $$(1+d-e, 1) \in \chi(e,k) \smallsetminus \bigcup\limits_{(e',k') <_U (e,k)} \chi(e',k').$$
 
 \bigskip
 Sinon $e > d$ alors $k =d $ et on pose $e = d v + u$ la division euclidienne de $e $ par $d $. On a alors 
 \begin{align*}
 \chi(e,k) = \llbracket 1, u \rrbracket \times \llbracket 1, v+1 \rrbracket \cup \llbracket u+ 1, d \rrbracket \times \llbracket 1, v \rrbracket. 
 \end{align*}
 \\ Soit $(e',k') <_U (e,k)$. On pose $e' = k'v' + u'$ la division euclidienne de $e' $ par $k'$. On distingue les cas $u = 0$ et $u \neq 0$. 
 \\ \textbullet \, \underline{Si $u =0 $.} On va montrer que $(d, v) \notin \chi(e',k') $. 
 
 Comme $e > d $ et $u =0$ on a nécessairement $v > 1$. Si $k' =e' $ alors $v' = 1 $ et $u' = 0 $ donc $$\chi(e',k') = \llbracket 1+d-e', d \rrbracket \times \{ 1 \} $$ et $(d,v) \notin \chi(e',k')$ car $v >1 $.
 \\Sinon $k' = d $ et alors $1 \leq v' < v $ car $e' < e$. On a 
 $$ \chi(e',k') =\llbracket 1, u' \rrbracket \times \llbracket 1, v'+1 \rrbracket \cup \llbracket u'+ 1, d \rrbracket \times \llbracket 1, v' \rrbracket $$
 Alors $(d,v) \notin \chi(e',k') $ car $u' \leq d-1$ et $v' < v$.
 \\ \textbullet \, \underline{Si $u > 0 $.} On va montrer que $(u, v+1) \notin \chi(e',k') $. 
 \\Si $k' =e' $ alors $v' = 1$, $u'= 0 $ et $$\chi(e',k') = \llbracket 1+d-e', d \rrbracket \times \{ 1 \}. $$
 Comme $v+1 > 1$ on a $(u, v+1) \notin \chi(e',k') $ dans ce cas là.
 \\ Sinon $k' = d$ et alors $v' \leq v$ puisque $e' <e$. On a 
 $$ \chi(e',k') =\llbracket 1, u' \rrbracket \times \llbracket 1, v'+1 \rrbracket \cup \llbracket u'+ 1, d \rrbracket \times \llbracket 1, v' \rrbracket $$
 et donc 
 \begin{align*}
(u,v+1) \in \chi(e',k') \Longleftrightarrow [ v' = v \text{ et } u' \geq u ].
 \end{align*}
Comme $e' <e$, cette propriété n'est pas respectée et donc $(u,v+1) \notin \chi(e',k')$.

 \bigskip
 On a donc montré que pour tout $(e,k) \in U$, il existe $(q,\ell)$ tel que 
 \begin{align*}
 (q,\ell) \in \chi(e,k) \smallsetminus \bigcup\limits_{(e',k') <_U (e,k)}\chi(e',k').
 \end{align*}

Les hypothèses du lemme $\ref{9lem_technique}$ sont donc vérifiées : $M_U(\II_U)$ contient un ouvert non vide et la famille $(\mu_n(\cdot|e)_{k } )_{(e,k) \in U}$ est algébriquement indépendante sur $\R$.

\end{preuve}

\bigskip
Le théorème~\ref{9lem_sousfam_pour_theo3} se prouve en utilisant le lemme~\ref{9lem_technique_faible} car l'ensemble $U$ associé à la famille d'exposants ne vérifient pas l'hypothèse $(\ref{9hypothese_sur_U})$ du lemme~\ref{9lem_technique}.

\begin{preuve}[ (théorème~\ref{9lem_sousfam_pour_theo3})]
On pose 
\begin{align*}
 U = \llbracket d, n-1 \rrbracket \times \{ d \}.
\end{align*}
On va appliquer le lemme~\ref{9lem_technique_faible}. On va montrer l'inversibilité pour tout $\beta \in (\R_+^{*})^{dm} $ de la matrice jacobienne $ J_{\Omega}(\beta)\in \MM_{dm, dm}(\R)$ où 
 $$\Omega_{(e,d)}(\beta) =\sum\limits_{q = 1 + f(e, md) }^d \beta_{q,1} \ldots \beta_{q,v_{q}(e,d)} $$ pour $(e,d) \in U$ avec les $v_1(e,d), \ldots, v_d(e,d) $ définis en $(\ref{9eq_def_vj})$ et $f(e,md) = \max(0, e- md)$.
 \\ On a alors 
 \begin{align*}
 J_{\Omega}(\beta) = \begin{pmatrix}
 D_1 & \cdots & D_d
 \end{pmatrix} 
 \end{align*}
 avec pour $i \in \llbracket 1,d \rrbracket$ :
 \begin{align*}
 D_i = \left(\frac{\partial \Omega_{(e,d)}}{\partial \beta_{i,\ell}}\right)_{e \in \llbracket d, n-1 \rrbracket, \ell \in \llbracket 1,m \rrbracket} \in \MM_{dm,m}.
 \end{align*}
De par la forme de $\Omega_{(e,d)}$, on a pour $(i,\ell) \in \llbracket 1, d \rrbracket \times \llbracket 1, m \rrbracket$ 
\begin{align*}
 \frac{\partial \Omega_{(e,d)}}{\partial \beta_{i,\ell}} = \left\{ \begin{array}{cl}
 \prod\limits_{\underset{p \neq \ell}{1 \leq p \leq v_{i }(e,d)}} \beta_{i,p} & \text{ si } \ell \leq v_{i}(e,d)\\
 0 & \text{ sinon. } 
 \end{array} \right. 
\end{align*}
On remarque en particulier que pour $\ell \in \llbracket 1 , v_{i}(e,d) -1 \rrbracket$ on a 
\begin{align*}
\frac{\beta_{i, \ell+1}}{\beta_{i, \ell}} \frac{\partial \Omega_{(e,d)}}{\partial \beta_{i,\ell+1}} = \frac{\partial \Omega_{(e,d)}}{\partial \beta_{i,\ell}}.
\end{align*}
Notons $D_{i,1}, \ldots,D_{i,m}$ les colonnes de $D_i$. On fait successivement les opérations élementaires suivantes entre les colonnes de $J_\Omega(\beta)$ :
\begin{align*}
 \forall i \in \llbracket 1, d \rrbracket, \quad \forall \ell \in \llbracket 1, m-1 \rrbracket, \quad D_{i, \ell} \leftarrow \left(D_{i, \ell} - \frac{\beta_{i, \ell+1}}{\beta_{i, \ell}}D_{i, \ell+1} \right), 
\end{align*}
et on note $G$ la matrice obtenue. 
\\On a donc $\det(J_\Omega(\beta)) = \det(G)$ avec $ G =\left(g_{e,(i,\ell)} \right)_{ e \in \llbracket d, n-1 \rrbracket, (i,\ell) \in \llbracket 1,d \rrbracket \times \llbracket 1,m \rrbracket}$ et 
\begin{align}\label{9def_matrice_G}
 g_{e,(i,\ell)} \left\{ \begin{array}{ll}
 \neq 0 & \text{ si } v_i(e,d) = \ell \\
 = 0 & \text{ sinon}
 \end{array}\right. .
\end{align}
On étudie alors $G$ et pour cela on réordonne les colonnes de $G$, en posant sur $\llbracket 1,d \rrbracket \times \llbracket 1,m \rrbracket$ l'ordre suivant :
\begin{align*}
 (i, \ell) < (j,\ell') \Longleftrightarrow [ \ell < \ell' \text{ ou } (\ell = \ell' \text{ et } i <j)]
\end{align*}
qui est l'ordre lexicographique sur $\llbracket 1,d \rrbracket \times \llbracket 1,m \rrbracket$ avec comparaison prioritaire du second élement du couple. On note $G'$ la matrice ainsi obtenue et on a alors $\det(G) = \pm \det(G')$.

D'après l'ordre choisi pour les colonnes, les coefficients diagonaux de $G'$ sont ceux de la forme $g_{e,(u+1,v)} $ avec $e = dv + u$ la division euclidienne de $e$ par $d$ pour $e \in \llbracket d ,n -1 \rrbracket$. 
On va montrer que la matrice $G'$ est triangulaire supérieure à coefficients diagonaux non nuls. 

Soit $e \in \llbracket d ,n -1 \rrbracket$. On pose $e = dv + u$ la division euclidienne de $e$ par $d$.
\\D'après la définition des $v_i(e,d)$ en $(\ref{9eq_def_vj})$ on a $v_{u+1}(e,d) = v$. La construction de $G$ en $(\ref{9def_matrice_G})$ donne alors 
$$g_{e,(u+1,v)} \neq 0 $$
Soit maintenant $e' > e $. On va montrer que $g_{e',(u+1,v)} =0 $. On pose $e' = dv' + u'$ la division euclidienne de $e'$ par $d$ et on a alors $v' \geq v$. 
\\ \textbullet \, \underline{Si $v' > v $} alors comme $v_{u+1}(e',d) = v' \text{ ou } v' +1 $ on a $v_{u+1}(e',d) \neq v$ et donc 
$$g_{e',(u+1,v)} =0 $$
d'après $(\ref{9def_matrice_G})$.
\\ \textbullet \, \underline{Si $v' = v $} alors $u' > u$ et $u +1 \in \llbracket 1 , u' \rrbracket$. En particulier d'après la définition des $v_i(e',d)$ en $(\ref{9eq_def_vj})$, on a $v_{u+1}(e',d) = v' +1 \neq v$ et donc 
$$g_{e',(u+1,v)} =0 $$
d'après $(\ref{9def_matrice_G})$.
\\On a donc montré que tous les coefficients diagonaux de $G'$ sont non nuls et que $G'$ est triangulaire supérieure.
On déduit que $G$ est inversible car $\det(G) = \pm \det(G')$ et alors $J_\Omega(\beta)$ l'est aussi aussi. 
 \\Le lemme~\ref{9lem_technique_faible} permet alors de conclure que $M_U(\II_U)$ contient un ouvert non vide et que la famille $(\mu_n(\cdot|e)_{d} )_{e \in \llbracket dm, d(m+1)-1 \rrbracket}$ est algébriquement indépendante sur $\R$. 
 
\end{preuve}

\section{Famille d'exposants avec tous les angles}\label{9section_autre_resultat}
Dans cette section, on exhibe une dernière famille d'exposants algébriquement indépendants sur $\Q$. Celle-ci fait apparaître des exposants correspondant à tous les angles $ (k \in \llbracket 1, d \rrbracket$) et avec $g(d,e,n) = 0$. Précisement, on a des exposants au $k\tir$ième angle pour tout $k \in \llbracket 1, d \rrbracket$, mais on approxime seulement par des sous-espaces rationnels de dimension multiple de $k$. On remarque que comme dans les théorèmes \ref{9prop_premier_angle}, \ref{9lem_sousfam_pour_theo2} et \ref{9lem_sousfam_pour_theo3} on obtient l'indépendance algébrique d'une famille de $n-d$ exposants.

\begin{prop}\label{9lem_sousfam_pour_theo}
On suppose que $n = d(m+1)$ avec $m \in \Nx $. On pose 
\begin{align}\label{9def_Udans_lem_pour_theo}
 U =\{ (e,k) \mid k \in \llbracket 1, d \rrbracket, e \in \llbracket k , km \rrbracket, k |e \}.
\end{align}
Alors $M_U(\II_U)$ contient un ouvert non vide et la famille $(\mu_n(\cdot|e)_{k } )_{(e,k) \in U}$ est algébriquement indépendante sur $\R$. 
 \end{prop}

\begin{preuve}
On remarque tout d'abord que pour tout $(e,k) \in U$, on a $g(d,e,n) = 0$ et $f(e,mk) = 0$ car $e \in \llbracket 1, dm \rrbracket$ et $dm = n-d$. 
\\
Il suffit de vérifier que l'ensemble $U$ vérifie les hypothèses du lemme $\ref{9lem_technique}$ pour conclure. 
\\ Tout d'abord, si $(e,k) \in U$ alors $e \leq km < k(m+1)$. 

\bigskip
On introduit l'ordre lexicographique avec comparaison prioritaire du second élement du couple sur $U$ :
\begin{align*}
 \quad(e ',k') <_U (e,k) \Longleftrightarrow [k' < k \text{ ou } (k' = k \text{ et } e' < e)].
\end{align*}
Soit $(e,k) \in U$ ; alors $e = k \ell $ avec $\ell \in \llbracket 1, m\rrbracket$. On a donc 
\begin{align*}
 \chi(e,k) = \llbracket 1 + d - k, d \rrbracket \times \llbracket 1, \ell \rrbracket.
\end{align*}
On va montrer que pour tout $(e',k') <_U (e,k)$, on a $ (1+d-k, \ell) \notin \chi(e',k') $.
\\Soit $(e',k') <_U (e,k)$. 
\\ \textbullet \, \underline{Si $k' < k $} alors $ \chi(e',k') \subset \llbracket 1 +d -k', d \rrbracket \times \llbracket 1, m\rrbracket$.
\\En particulier, comme $ 1+ d -k \notin \llbracket 1 +d -k', d \rrbracket$ on a $ (1+d-k, \ell) \notin \chi(e',k') $.
\\ \textbullet \, \underline{Si $k' = k$} alors $ e' = k \ell'$ avec $\ell' < \ell $ et donc 
\begin{align*}
 \chi(e',k') = \llbracket 1 + d - k, d \rrbracket \times \llbracket 1, \ell' \rrbracket.
\end{align*}
Comme $\ell \notin \llbracket 1, \ell' \rrbracket$ on a $ (1+d-k, \ell) \notin \chi(e',k') $.

\bigskip
Les hypothèses du lemme $\ref{9lem_technique}$ sont donc vérifiées : $M_U(\II_U)$ contient un ouvert non vide et la famille $(\mu_n(\cdot|e)_{k } )_{(e,k) \in U}$ est algébriquement indépendante sur $\R$.

\end{preuve}

%% file: Chapitres/7_Dernier_angle_autre_exemple.tex
\chapter{Construction d'espaces avec exposants prescrits au dernier angle}\label{chap8}

Dans ce chapitre, on construit une nouvelle famille d'espaces $A$ dont on est capable de calculer les exposants correspondant au dernier angle, c'est-à-dire $\mu_n(A|e)_{\min(\dim(A), e)}$. On se limite ici au cas $g(d,e,n) =0$ c'est-à-dire $d +e \leq n$. 

\bigskip

On va montrer le  théorème~suivant.

\begin{theo}\label{8theo_dernier_angle}
 Soit $(d,q)\in (\Nx)^2$. On pose $n = (q+1)d$.
 \\ Pour tout $\alpha \geq 3d(d+4)$, il existe un sous-espace vectoriel $A$ de $\R^n$ de dimension $d$ tel que 
 \begin{align*}
 \forall e \in \llbracket d, n-d \rrbracket, \quad &\mu_n(A| e)_d = \frac{\alpha^{q_e +1 } }{r_e + (d-r_e)\alpha} \\
 \forall e \in \llbracket 1, d -1\rrbracket, \quad &\mu_n(A| e)_e = \frac{\alpha }{r_e } = \frac{\alpha}{e}
\end{align*}
où $q_e$ et $ r_e$ sont le quotient et le reste de la division euclidienne de $e $ par $d$.
\end{theo}

\begin{req}
 Le théorème~\ref{7theo_principale} donnait déjà un espace de dimension $d$ avec exposants diophantiens prescrits pour le dernier angle avec $e \in \llbracket 1, n-d \rrbracket$. 
 \\Ici on permet d'avoir d'autres restrictions sur les valeurs prises par ces exposants. On obtient en effet le même résultat que le corollaire~\ref{7cor_exposantsprescrits} avec $\alpha_1 = \ldots = \alpha_{d} = \alpha$ et on s'affranchit alors de l'hypothèse $$\forall i \in \llbracket 1,d-1 \rrbracket, \quad \alpha_{i}^{c_{\ref{7cons_petite_hyp_theoc2}}} < \alpha_{i+1} < \alpha_{i}^{c_{\ref{7cons_petite_hyp_theoc2c1}}}.$$
 De plus les preuves développées dans ce chapitre diffèrent de celles des autres chapitres et s'appuient notamment sur la géométrie des nombres.
\end{req}

 La construction de l'espace $A$ est réalisée dans la section~\ref{8premier_paragraphe_construction}, elle est similaire à celles des chapitres précédents. 
Pour démontrer le  théorème~\ref{8theo_dernier_angle}, on fait une disjonction de cas selon que $e <d $ ou $e \geq d$. 
\\Dans les deux cas, la minoration de l'exposant $\mu_n(A|e)_{\min(d,e)}$ est montrée en exhibant une famille d'espaces rationnels approchant bien $A$ (lemmes~\ref{8lem_psid_ACNe} et \ref{8lem_Angle_psie_DN}) de façon similaire aux chapitres précédents.
\\Pour $e \geq d$, on montre que les \og meilleurs \fg{} espaces $C$ approchant $A$ contiennent un certain espace rationnel $B_{N+1,q_e}$ (lemme~\ref{8lem_inclusion_BN_C}). On peut alors minorer la hauteur de l'espace $C$ (lemme~\ref{8minoration_somme_espace}) et on en conclut que $\mu_n(A|e)_{d}$ ne peut pas être trop grand grâce au lemme~\ref{8lem_minoration_exposant}.
\\Pour $e < d$, on montre que les \og meilleurs \fg{} espaces $C$ approchant $A$ intersectent non trivialement un certain espace rationnel $D_{N,d}$ (lemme~\ref{8lem_intersection_DN_C}). On aboutit à la minoration de $\mu_n(A|e)_{e}$ dans le lemme~\ref{8lem_minoration_exposant_e<d} en minorant $\psi_{1}(C \cap D_{N,d}, A) $ et donc a fortiori $\psi_e(C, A)$.

\bigskip
On définit la constante $C_d$ par
\begin{align*}
 C_d = 3d(d+4). 
\end{align*}

On énonce un lemme listant quelques inégalités plus ou moins triviales, vérifiées pour $\alpha \geq C_d$ et qui nous serviront pour prouver le théorème~\ref{8theo_dernier_angle}. Elles sont énoncées sous la forme utilisée dans les preuves de ce chapitre.

\begin{lem}\label{8lem_inegal_alpha}
Soit $\alpha \geq C_d$. Alors on a 
\begin{align}
 -\alpha^{2} + \alpha(2d+2) - d &\leq 0 \label{8lem_inegal_alpha1} \\
 - \frac{\alpha}{2} + d(d-1) +1 &\leq 0 \label{8lem_inegal_alpha2} \\
 -\alpha^2 +(1+2d)\alpha - d &\leq 0 \label{8lem_inegal_alpha3},
\end{align}
et pour tout $e \in \llbracket d, qd \rrbracket $ 
\begin{align}
 dr_e -(d-r_e) \alpha \leq 0 \label{8lem_inegal_alpha4} \\
 \frac{\alpha^{q_e}}{d -r_e + \frac{1}{2}} -\frac{\alpha^{q_e+1}}{r_e +(d -r_e)\alpha} + 1 \leq 0\label{8lem_inegal_alpha5}
\end{align}
en notant $q_e$ et $ r_e$ le quotient et le reste de la division euclidienne de $e $ par $d$.
\end{lem}

\begin{preuve}
 Les inégalités $(\ref{8lem_inegal_alpha1}), (\ref{8lem_inegal_alpha2}) $ et $(\ref{8lem_inegal_alpha3})$ sont vérifiées si 
 \begin{align}\label{8preuve_min_alpha}
 \alpha \geq \max\left(\frac{2d+2 + \sqrt{4d^2 + 4d + 4}}{2}, 2d(d-1)+2, \frac{1+2d + \sqrt{1+4d^2}}{2} \right).
 \end{align}
 On vérifie facilement que $C_d \geq \max\left(2d(d-1) +2, {2+\sqrt{3}}, \frac{3+\sqrt{5}}{2}\right)$. Cette constante est supérieure au maximum considéré en $(\ref{8preuve_min_alpha})$ ; les deux fractions y apparaissant valent en effet respectivement, $2+\sqrt{3}$ et $\frac{3+\sqrt{5}}{2}$ si $d = 1$, et le cas $d \geq 2 $ est trivial car alors $$C_d \geq 2d(d-1) +2 \geq \frac{2d+2 + \sqrt{4d^2 + 4d + 4}}{2} \geq \frac{1+2d + \sqrt{1+4d^2}}{2}.$$
L'inégalité $(\ref{8lem_inegal_alpha4})$ est vérifiée si 
\begin{align*}
 \alpha \geq \frac{dr_e}{d-r_e}.
\end{align*}
Or $\frac{dr_e}{d-r_e} \leq dr_e \leq d(d-1) \leq C_d$ ce qui permet de conclure.
\bigskip

Il reste alors à prouver l'inégalité $(\ref{8lem_inegal_alpha5})$. On remarque d'abord que 
\begin{align*}
 C_d &\geq (\sqrt{2} + 1)d(d+ \frac{3 }{2}) \\
 &\geq d-1 + d(d+\frac{1}{2}) + \sqrt{2}(d-1 +d(d+\frac{1}{2}))\\
 &\geq d-1 + d(d+\frac{1}{2}) + \sqrt{ 2(d-1 + d(d+\frac{1}{2}) )^2 } \\
 &\geq d-1 + d(d+\frac{1}{2}) + \sqrt{ (d-1 + d(d+\frac{1}{2}) )^2 + 2(d+\frac{1}{2})(d-1)}.
\end{align*}
En particulier, $C_d$ est supérieur à la plus grande racine du polynome 
\begin{align*}
 P_d(x) = \frac{x^2}{2} -\left(d-1 + (d + \frac{1}{2})d \right)x - (d + \frac{1}{2})(d-1) \geq 0.
\end{align*}
On en conclut que pour tout $x \geq C_d$, on a $P_d(x) \geq 0$.
L'inégalité $(\ref{8lem_inegal_alpha5})$ est équivalente à 
\begin{align*}
 &\alpha^{q_e}\left((d-r_e + \frac{1}{2}) \alpha - r_e - (d-r_e)\alpha \right) - (d-r_e + \frac{1}{2})(r_e + (d-r_e)\alpha) \geq 0 \\
 \Longleftrightarrow \quad &\alpha^{q_e}\left(\frac{\alpha}{2} - r_e \right) - (d-r_e + \frac{1}{2})(d-r_e)\alpha - (d-r_e + \frac{1}{2})r_e \geq 0.
\end{align*} 
Comme $e \in \llbracket d, qd\rrbracket$, on a $q_e \geq 1$. On cherche donc à montrer 
\begin{align*}
 \frac{\alpha^2}{2} -\left(r_e + (d-r_e + \frac{1}{2})(d-r_e) \right)\alpha - (d-r_e + \frac{1}{2})r_e \geq 0.
\end{align*}
Comme $0 \leq r_e \leq d-1$ et $(d-r_e + \frac{1}{2})r_e \leq (d + \frac{1}{2})(d-1)$ il suffit alors d'avoir : 
\begin{align*}
 \frac{\alpha^2}{2} -\left(d-1 + (d + \frac{1}{2})d \right)\alpha - (d + \frac{1}{2})(d-1) \geq 0.
\end{align*}
Or cette expression est $P_d(\alpha)$ et comme $\alpha \geq C_d$, on a $P_d(\alpha) \geq 0 $ ce qui termine la preuve.

 \end{preuve}

\section{Construction de l'espace \texorpdfstring{$A$}{}}\label{8premier_paragraphe_construction}
On rappelle que $n =(q+1)d$ et donc $qd = n-d$. 

On construit ici l'espace $A$ du théorème~\ref{8theo_dernier_angle} en utilisant les outils introduits dans le chapitre~\ref{chap3}.
Soit $\theta$ un nombre premier supérieur ou égal à $5$. 
\\Pour $j \in \llbracket 1, d \rrbracket$, on pose $\phi_j : \N \to \llbracket 0, qd-1 \rrbracket $ définie par :
\begin{align*}
 \phi_j(k) = k + (j-1)q \mod (qd)
\end{align*}
avec $x \mod (qd)$ le reste de la division euclidienne de $x$ par $qd$. \\
Dans la suite on note $\sigma_{i,j} = \sigma(\theta, u^{(i,j)}, (\alpha^k)_{k \in \N}) =  \sum\limits_{k = 0}^{+ \infty} \frac{u^{(i,j)}_k}{\theta^{ \floor{ \alpha^k} } }$ pour $i \in \llbracket 0, qd-1 \rrbracket$ et $j \in \llbracket 1, d \rrbracket$ avec des suites $u^{(i,j)}$ que l'on va construire .
D'après le lemme~\ref{3lem_sigma_alg_indep}, comme $\alpha >1$, il existe $qd$ suites $u^{(0,1)}, \ldots, u^{(qd-1, 1)}$ vérifiant : 
\begin{align*}
 \forall i \in \llbracket 0, qd-1 \rrbracket, \quad \forall k \in \N, \quad u^{(i,1)}_k \left\{ \begin{array}{ll}
 \in \{2, 3 \} &\text{ si } i = \phi_1(k)\\
 = 0 &\text{ sinon }
 \end{array} \right. 
\end{align*}
et telles que la la famille $(\sigma_{0, 1}, \ldots, \sigma_{qd-1,1})$ soit algébriquement indépendante sur $\Q$. \\
On construit alors des suites $u^{(0,2)},\ldots,u^{(qd-1,2)}, \ldots, u^{(0,d)}, \ldots, u^{(qd-1,d)}$ par récurrence.
\\ Soit $j \in \llbracket 1, d-1 \rrbracket $, on suppose que l'on a construit $u^{(0,1)},\ldots,u^{(qd-1,1)}, \ldots, u^{(0,j)}, \ldots, u^{(qd-1,j)}$ vérifiant 
\begin{align*}
 \forall i \in \llbracket 0, qd-1 \rrbracket, \quad \forall \ell \in \llbracket 1, j \rrbracket, \quad \forall k \in \N, \quad u^{(i,\ell)}_k \left\{ \begin{array}{ll}
 \in \{2, 3 \} &\text{ si } i = \phi_\ell(k)\\
 = 0 &\text{ sinon }
 \end{array} \right. 
\end{align*}
et telles que la famille $(\sigma_{0,1},\ldots,\sigma_{qd-1,1}, \ldots, \sigma_{0,j}, \ldots, \sigma_{qd-1,j})$ soit algébriquement indépendante sur $\Q$. Le lemme~\ref{3lem_sigma_alg_indep}, appliqué avec $\FF = \{ \sigma_{0,1},\ldots,\sigma_{qd-1,1}, \ldots, \sigma_{0,j}, \ldots, \sigma_{qd-1,j} \} $ donne alors $qd$ suites $u^{(0,j+1)}, \ldots, u^{(qd-1, j+1)}$ vérifiant : 
\begin{align*}
 \forall i \in \llbracket 0, qd-1 \rrbracket, \quad \forall k \in \N, \quad u^{(i,j+1)}_k \left\{ \begin{array}{ll}
 \in \{2, 3 \} &\text{ si } i = \phi_{j+1}(k)\\
 = 0 &\text{ sinon }
 \end{array} \right. 
\end{align*}

\bigskip 
On a donc construit $qd \times d =(n-d)d$ suites $u^{(0,1)},\ldots,u^{(qd-1,1)}, \ldots, u^{(0,d)}, \ldots, u^{(qd-1,d)}$ vérifiant $ \forall i \in \llbracket 0, qd-1 \rrbracket, \forall j \in \llbracket 1, d \rrbracket, \forall k \in \N $ : 
\begin{align}\label{8construc_suite_u}
u^{(i,j)}_k \left\{ \begin{array}{ll}
 \in \{2, 3 \} &\text{ si } i = k + (j-1)q \mod (qd)\\
 = 0 &\text{ sinon }
 \end{array} \right. 
\end{align}
et telles que la famille $(\sigma_{0,1},\ldots,\sigma_{qd-1,1}, \ldots, \sigma_{0,d}, \ldots, \sigma_{qd-1,d})$ soit algébriquement indépendante sur $\Q$.

\begin{req}\label{8req_def_uk}
 A $k$ et $j $ fixés, $i = k + (j-1)q \mod (qd) $ est l'unique entier $i \in \llbracket 0, qd-1 \rrbracket$ tel que $u_k^{(i,j)} \neq 0$. 
\end{req}

\begin{lem}\label{8lem_unique_kj_pour_i}
 Soient $i \in \llbracket 0, qd-1 \rrbracket$ et $N \in \N$.
 \\Il existe alors un unique couple $(k,j) \in \llbracket 0, q-1 \rrbracket \times \llbracket 1, d \rrbracket$ tel que $u_{N+k}^{(i,j)} \neq 0$.
\end{lem}

\begin{preuve}
 \textbullet \, \underline{Unicité :} Supposons qu'il existe $\ell_1, \ell_2 \in \llbracket 0, q-1 \rrbracket$ et $j_1, j_2 \in \llbracket 1,d \rrbracket$ tels que :
\begin{align*}
 u_{N + \ell_1}^{(i,j_1)} \neq 0 \text{ et } u_{N + \ell_2}^{(i,j_2)} \neq 0.
\end{align*}
Par définition des $u_k^{(i,j)}$ on a :
\begin{align*}
 N + \ell_1 +(j_1-1)q \equiv N + \ell_2 +(j_2-1)q \mod(qd).
\end{align*}
Par unicité de la division euclidienne par $q$, puisque $\ell_1, \ell_2 \in \llbracket 0, q-1 \rrbracket$, on a $\ell_1 = \ell_2$.
\\ On a donc $(j_1-1)q \equiv (j_2-1)q \mod(qd) $ et alors $(j_1 - j_2)q \equiv 0 \mod (qd)$. \\ Comme $j_1, j_2 \in \llbracket 1,d \rrbracket$, on a $j_1 = j_2 $.
 
\textbullet \, \underline{Existence :} On écrit les divisions euclidiennes de $i$ et de $N$ par $q $ : 
 \begin{align*}
 i = qu + v \text{ et } N =qu' + v'
 \end{align*}
 avec $v, v' \in \llbracket 0, q-1 \rrbracket$.
 \\Si $v \geq v' $, on pose $k = v -v' \in \llbracket 0, q-1 \rrbracket$ et $j = (u -u' \mod(d) ) +1 \in \llbracket 1, d \rrbracket$. 
 \\ On vérifie alors que $i = N+k + (j-1)q \mod(qd)$ : 
 \begin{align*}
 N+k + (j-1)q \mod(qd) &= qu' + v' + v -v' + q(u-u') \mod(qd) \\
 &= qu + v \mod(qd) = i.
 \end{align*}
 \\Si $v < v' $, on pose $k = v -v' +q \in \llbracket 0, q-1 \rrbracket$ et $j = (u -u' -1 \mod(d) ) +1 \in \llbracket 1, d \rrbracket$. 
 \\ On vérifie alors que $i = N+k + (j-1)q \mod(qd)$ : 
 \begin{align*}
 N+k + (j-1)q \mod(qd) &= qu' + v' + v -v' +q + q(u-u'-1) \mod(qd) \\
 &= qu + v \mod(qd) = i.
 \end{align*}
\end{preuve}

\bigskip
On pose maintenant pour $j \in \llbracket 1,d \rrbracket$, le vecteur de $\R^n$ :
\begin{align*}
 Y_j = \begin{pmatrix} 
 0 \\ \vdots\\
 0 \\ 
 1 \\
 0 \\
 \vdots \\
 0 \\
 \sigma_{0,j} \\
 \vdots \\
 \sigma_{qd-1,j}
 \end{pmatrix} \begin{array}{cc} \\ \\
 \\ \\ \leftarrow j\text{-ème ligne} \\ \\ \\ \\ \\ \\ \\ \\ \\ \end{array}.
\end{align*}
On définit alors l'espace $A$ du théorème~\ref{8theo_dernier_angle} par $A = \Vect(Y_1, \ldots, Y_d) $. En regardant les $d $ premières lignes des vecteurs $Y_j$, on a clairement $\dim(A) = d$.

\begin{lem}\label{8lem_espace_irrat}
 L'espace $A$ est $(e,1)\tir$irrationnel pour tout $e \in \llbracket 1, qd \rrbracket$.
\end{lem}

\begin{preuve}
On pose $M = \begin{pmatrix}
 Y_1 |& \cdots &|Y_d
\end{pmatrix} \in \MM_{n,d}(\R)$ la matrice dont les colonnes sont $Y_1, \ldots, Y_d$. 
\\Par construction on a $M= \begin{pmatrix} I_d \\ \Sigma \end{pmatrix}$ avec $\Sigma= (\sigma_{i,j} )_{i \in \llbracket 0, qd-1 \rrbracket, j \in \llbracket 1, d \rrbracket} $.
\\Comme les coefficients de $\Sigma $ forment une famille algébriquement indépendante sur $\Q$ et que $I_d \in \GL_{d}(\Q)$, le lemme~\ref{3lem_1_irrat} donne que $A$ est $(e,1)\tir$irrationnel pour tout $e \in \llbracket 1, n-1 \rrbracket$.

\end{preuve}

\section{L'espace \texorpdfstring{$B_{N,v}$}{} }
On reprend ici les idées des constructions effectuées dans le chapitre~\ref{chap5} pour construire les \og meilleurs \fg{} espaces approximant $A$ et étudier leurs propriétés. \\ 
Pour $ i \in \llbracket 0, qd-1 \rrbracket$, $ j \in \llbracket 1,d \rrbracket$ et $N \in \N^* $ on pose :
\begin{align*}
 \sigma_{i,j }^N =  \sum\limits_{k = 0}^{N} \frac{u_k^{(i,j)}}{\theta^{\floor{\alpha^k}}} \in \frac{1}{\theta^{\floor{\alpha^N}}}\Z
\end{align*}
 \\On pose maintenant pour $ j \in \llbracket 1, d \rrbracket$ le vecteur de $\Z^n$ : 
\begin{align*}
 X_N^j =\theta^{\floor{\alpha^N}} \begin{pmatrix}
 0 \\ \vdots\\
 0 \\ 
 1 \\
 0 \\
 \vdots \\
 0 \\
 \sigma_{0,j}^N \\
 \vdots \\
 \sigma_{qd-1,j}^N
 \end{pmatrix} \begin{array}{cc} \\ \\
 \\ \\ \leftarrow j\text{-ème ligne} \\ \\ \\ \\ \\ \\ \\ \\ \\ \end{array}.
\end{align*}
Dans la suite on note pour $N$ et $v$ deux entiers non nuls :
\begin{align}\label{8def_Bnv}
 B_{N,v} = \Vect(X_N^1, X_{N+1}^1,\ldots, X_{N+v-1}^1, X_N^2 \ldots, X_{N+v-1}^2, \ldots, X_N^d,\ldots, X_{N+v-1}^d)
\end{align} qui est un espace rationnel par définition. 
\\ Pour $j \in \llbracket 1,d \rrbracket$ on remarque que 
:
\begin{align}\label{8rec_Xn_Un}
 X_{N+1}^j = {\theta^{\floor{\alpha^{N+1}} - \floor{\alpha^{N}}}}X_N^j + U_{N+1}^j \text{ avec } U_{N+1}^j = \begin{pmatrix}
 0 \\
 \vdots\\
 0 \\ 
 u_{N+1}^{(0,j)}\\
 \vdots \\
 u_{N+1}^{(qd-1,j)}
 \end{pmatrix}.
\end{align}
On pose $V_N^j = \frac{U_N^j }{\| U_N^j\|} \in \Z^n$ car les vecteurs $U_N^j$ ont une unique coordonnée non nulle d'après la construction en (\ref{8construc_suite_u}) et la remarque~\ref{8req_def_uk}.
\\Les vecteurs $V_N^j$ sont donc des vecteurs de la base canonique de $\R^{n}$.
\\ On introduit aussi les vecteurs :
\begin{align}\label{8vecteur_Zn}
 Z_{N}^j = \frac{1}{\theta^{\floor{\alpha^N}}} X_N^j 
\end{align}
et on remarque que $Z_N^j \underset{N\to+\infty}{\longrightarrow} Y_j$ et plus précisement pour $j \in \llbracket 1,d \rrbracket$ :
\begin{align}\label{8major_Yj_XNj}
 \psi_1(\Vect(Y_j), \Vect(X_N^j) ) = \omega(Y_j, Z_N^j) \leq \frac{\|Y_j -Z_N^j \| }{\|Y_j \|} \leq c_{\ref{8cons_major_angle_Yj_XNj}}\theta^{- \alpha^{N+1}}
\end{align}
avec $\cons \label{8cons_major_angle_Yj_XNj}$ indépendante de $N$. On en déduit par ailleurs qu'il existe des constantes $\cons \label{8cons_minor_norme_XN}$ et $\cons \label{8cons_major_norme_XN}$ indépendantes de $N$ telles que pour tout $N \in \N$ :
\begin{align}\label{8norme_XN}
 c_{\ref{8cons_minor_norme_XN}} \theta^{\alpha^N} \leq \|X_N^j\| \leq c_{\ref{8cons_major_norme_XN}} \theta^{\alpha^N}.
\end{align}
\bigskip

\begin{lem}\label{8lem_dimBN_zbaseBN}
Soit $v \in \llbracket 1, q \rrbracket$. Alors l'espace $B_{N,v}$ est de dimension $dv$.
\\ De plus les vecteurs $(X_N^j)_{j \in \llbracket 1,d\rrbracket} \cup (V_k^j)_{j \in \llbracket 1,d\rrbracket, k \in \llbracket N+1, N+v-1\rrbracket}$ forment une $\Zbase$ de $B_{N,v} \cap \Z^n$.

\end{lem}

\begin{preuve}
Par récurrence sur $v$ et en utilisant (\ref{8rec_Xn_Un}) on a : 
\begin{align*}
 B_{N,v}= \Vect(X_N^1, V_{N+1}^1,\ldots, V_{N+v-1}^1, X_N^2, V_{N+1}^2,\ldots, V_{N+v-1}^2, \ldots, X_N^d, V_{N+1}^d,\ldots, V_{N+v-1}^d).
\end{align*}
La remarque~\ref{8req_def_uk} permet d'affirmer que les $V_k^j$ considérés ici sont tous différents. On rappelle de plus que ce sont des vecteurs de la base canonique. On en déduit en particulier que $(X_N^j)_{j \in \llbracket 1,d\rrbracket} \cup (V_k^j)_{j \in \llbracket 1,d\rrbracket, k \in \llbracket N+1, N+v-1\rrbracket}$ forment une famille libre.
 \bigskip \\
 On montre d'abord le lemme avec $v =q$. 
 \\Soit $(a_{j,k})_{{j \in \llbracket 1,d\rrbracket, k \in \llbracket 0, q-1\rrbracket}} \in \left[0,1\right]^{qd}$ tel que :
\begin{align}\label{8vecteur_U}
 U =  \sum\limits_{j = 1}^d a_{j,0} X_N^j + \sum\limits_{k=1}^{q-1} \sum\limits_{j = 1}^d a_{j,k} V_{N+k}^j \in \Z^n
\end{align}
En étudiant les $d$ premières coordonnées de $U$ on trouve que $a_{j,0} \theta^{\floor{\alpha^N}} \in \Z$ pour tout $j \in \llbracket 1,d \rrbracket$. 
\\On peut écrire $ a_{j,0} = \frac{x_j }{y_j}$ avec $y_j \big{|} \theta^{\floor{\alpha^N}}$ et $\pgcd(x_j, \theta) = 1 $ pour tout $j \in \llbracket 1,d \rrbracket$.
\\Soit $j_ {\max} \in \llbracket 1,d \rrbracket$ tel que $\max\limits_{j = 1}^d y_j = y_{j_ {\max}}$. 
\\ Considérons maintenant l'entier $ i = N + q(j_ {\max} -1) \mod(qd) \in \llbracket 0, qd-1 \rrbracket$. Par définition des $u_k^{(i,j)}$ en (\ref{8construc_suite_u}) on a :
\begin{align}\label{8i_max}
 u^{(i,j_ {\max})}_N \neq 0
\end{align}
et donc 
\begin{align*}
 \forall k \in \llbracket 1, q-1\rrbracket, \quad \forall j\in \llbracket 1,d \rrbracket, \quad u_{N+k}^{(i,j)} = 0 
\end{align*} d'après le lemme~\ref{8lem_unique_kj_pour_i}.
\\ Par définition des $U_{N+k} $ en (\ref{8rec_Xn_Un}) et des $V_{N+k} = \frac{U_{N+k}}{\| U_{N+k}\|}$, la $(d+i)-$ème coordonnée du vecteur 
 $ \sum\limits_{k=1}^{q-1} \sum\limits_{j = 1}^d a_{j,k} V_{N+k}^j$ est donc nulle.
\\ On a donc d'après (\ref{8vecteur_U}) :
\begin{align*}
 \sum\limits_{j = 1}^d a_{j,0} \theta^{\floor{\alpha^N}}\sigma_{i,j}^N = \sum\limits_{j = 1}^d \frac{x_j}{y_j}\theta^{\floor{\alpha^N}}\sigma_{i,j}^N \in \Z. 
\end{align*}
Comme les $y_j$ sont des puissances de $\theta$ et $\max\limits_{j = 1}^d y_j = y_{j_ {\max}}$ on a $\frac{y_{j_ {\max}}}{y_j} \in \Z$ et :
\begin{align}\label{8divisibilite}
 \sum\limits_{j = 1}^d {x_j}\frac{y_{j_ {\max}}}{y_j}\theta^{\floor{\alpha^N}}\sigma_{i,j}^N \in y_{j_ {\max}}\Z.
\end{align}
Maintenant on remarque que pour $j \neq j_ {\max}$, on a $i \not \equiv N +1 +q(j-1) [qd]$ et donc $u_N^{(i,j)} =0 $ dans ce cas.
Comme $\floor{ \alpha^{N-1}} < \floor{\alpha^N} $, on a alors pour tout $j \neq j_ {\max}$:
\begin{align*}
 \theta {\mid} \theta^{\floor{\alpha^N}}\sigma_{i,j}^N = \theta^{\floor{\alpha^N}}\sum\limits_{k = 0}^{N} \frac{u_k^{(i,j)}}{\theta^{\floor{\alpha^k}}} .
\end{align*}
Si $\theta \mid y_{j_ {\max}}$ on a alors en utilisant (\ref{8divisibilite}): 
\begin{align*}
 \theta \mid x_{j_ {\max}}\theta^{\floor{\alpha^N}}\sigma_{i,j_ {\max}}^N 
\end{align*}
Or comme $u_n^{(i,j_ {\max})} $ est non nul et est premier avec $\theta$ on a $\pgcd(\theta^{\floor{\alpha^N}}\sigma_{i,j_ {\max}}^N, \theta )= 1$. 
\\ On en déduit enfin que $\theta \mid x_{j_ {\max}}$ ce qui est contradictoire avec $\pgcd(x_{j_ {\max}}, \theta )= 1 $.
\\Tous les $y_j$ sont donc égaux à $1$ et alors tous les $a_{j,0}$ sont entiers.
\\ En reprenant (\ref{8vecteur_U}) on trouve :
\begin{align*}
  \sum\limits_{k=1}^{q-1} \sum\limits_{j = 1}^d a_{j,k} V_{N+k}^j \in \Z^n.
\end{align*}
Or tous les vecteurs $V_{N+k}^j$ sont distincts et proviennent de la base canonique, donc $a_{j,k} $ est entier pour tous $j \in \llbracket 1, d \rrbracket$ et $k \in \llbracket 0,q-1\rrbracket $. Cela montre donc que $$(X_N^j)_{j \in \llbracket 1,d\rrbracket} \cup (V_k^j)_{j \in \llbracket 1,d\rrbracket, k \in \llbracket N+1, N+q-1\rrbracket}$$ forme une $\Zbase$ de $B_{N,q} \cap \Z^n$ et en particulier que $\dim(B_{N,q}) = qd$.

\bigskip 
Soit maintenant $v \in \llbracket 1,q - 1 \rrbracket$. On a $$(X_N^j)_{j \in \llbracket 1,d\rrbracket} \cup (V_k^j)_{j \in \llbracket 1,d\rrbracket, k \in \llbracket N+1, N+v-1\rrbracket} \subset (X_N^j)_{j \in \llbracket 1,d\rrbracket} \cup (V_k^j)_{j \in \llbracket 1,d\rrbracket, k \in \llbracket N+1, N+q-1\rrbracket}.$$ La remarque~\ref{2req_sous_Zbase} donne alors que la famille $(X_N^j)_{j \in \llbracket 1,d\rrbracket} \cup (V_k^j)_{j \in \llbracket 1,d\rrbracket, k \in \llbracket N+1, N+v-1\rrbracket}$ forme une $\Zbase$ du $\Z\tir$module qu'elle engendre. Or ce $\Z\tir$module est $B_{N,v} \cap \Z^n$ et le lemme est prouvé.

\end{preuve}

\section{Calcul de l'exposant dans le cas \texorpdfstring{$e \geq d$}{ } }
Dans cette section, on considère $e \in \llbracket d, qd \rrbracket $. 
\\ On rappelle que l'on a $e = q_e d + r_e $ la division euclidienne de $e $ par $d$. On a donc en particulier $1 \leq q_e \leq q$.
\\ On calcule ici $\mu_n(A| e)_d$.
\subsection{Minoration de l'exposant}
On va introduire dans cette section une suite d'espaces rationnels de dimension $e$ approximant bien $A$ ce qui nous permet de minorer $\mu_n(A|e)_d$. 

Soit $N \in \N$. On définit l'espace $C_{N,e}$ par :
\begin{align}\label{8def_Cn}
 C_{N,e} &= \Vect(X_{N+1}^1, \ldots, X_{N+q_e}^1, \ldots, X_{N+1}^d, \ldots, X_{N+q_e}^d) \bigoplus \Vect(X_N^1, \ldots, X_N^{r_e}) \\
 &= B_{N+1,q_e} \bigoplus \Vect(X_N^1, \ldots, X_N^{r_e}) \nonumber
\end{align}
qui est un espace rationnel.
En utilisant (\ref{8rec_Xn_Un}) et en raisonnant par récurrence pour chaque $j \in \llbracket 1,d \rrbracket$ on a :
\begin{align}\label{8CN_rec}
 C_{N,e} = &\Vect(X_{N}^1, V_{N+1}^1, \ldots, V_{N+q_e}^1, \ldots, X_{N}^{r_e},V_{N+1}^{r_e}, \ldots, V_{N+q_e}^{r_e}) \nonumber\\ &\bigoplus \Vect(X_{N+1}^{r_e+1}, V_{N+2}^{r_e+1}, \ldots, V_{N+q_e}^{r_e+1}\ldots, X_{N+1}^{d},V_{N+2}^{d}, \ldots, V_{N+q_e}^{d}).
\end{align}

\begin{req}\label{8req_CN=BN}
 On remarque que dans le cas où $r_e =0 $, par la définition en $(\ref{8def_Cn})$ on a $C_{N,e} = B_{N+1,q_e}$ . 
\\ Dans tous les cas, on a $ B_{N+1,q_e} \subset C_{N,e} \subset B_{N,q_e+1}$.
\end{req}

\begin{lem}\label{8lem_hauteur_CN}
On a $\dim(C_{N,e}) = e $.
\\ De plus, il existe des constantes $\cons \label{8cons_haut_CN_minor} >0 $ et $\cons \label{8cons_haut_CN_major} >0 $ indépendantes de $N$ telles que 
\begin{align*}
 c_{\ref{8cons_haut_CN_minor}} \theta^{ r_e {\alpha^N} + (d-r_e) {\alpha^{N+1}}} \leq H(C_{N,e}) \leq c_{\ref{8cons_haut_CN_major}} \theta^{ r_e {\alpha^N} + (d-r_e) {\alpha^{N+1}}}.
\end{align*}
\end{lem}

\begin{req}\label{8HauteurBN}
Ce lemme permet en particulier d'avoir $$ c_{\ref{8cons_haut_CN_minor}} \theta^{d{\alpha^N}} \leq 
 H(B_{N,v}) \leq c_{\ref{8cons_haut_CN_major}} \theta^{d{\alpha^N}} $$ pour tout $v \in \llbracket 1, q \rrbracket$ et $N \in \N\smallsetminus \{0\}$. En effet $B_{N,v} = C_{N-1,dv}$.
\end{req}

\begin{preuve}
Si $q_e = q$, alors $C_{N,e} = B_{N+1,q}$ car $e =qd $ et $r_e = 0$ ; sinon $q_e < q $ et dans ce cas $C_{N,e} \subset B_{N,q}$. 
Dans chacun des cas, le lemme~\ref{8lem_dimBN_zbaseBN} donne que les vecteurs considérés en $(\ref{8CN_rec})$ proviennent d'une $\Zbase$ de $B_{N+1, q}\cap \Z^n$ ou $B_{N,q} \cap \Z^n$ respectivement. 

La relation (\ref{8CN_rec}) donne directement $\dim(C_{N,e}) = (q_e+1)r_e + (d-r_e)q_e = q_ed + r_e =e $ car toutes les familles $X_{N}^j, V_{N+1}^j, \ldots, V_{N+q_e}^j$ sont libres.

Par la remarque~\ref{2req_sous_Zbase}, ces vecteurs forment donc une $\Zbase $ de $C_{N,e} \cap \Z^n$.
\\ En reprenant la notation $ Z_{N}^j = \frac{1}{\theta^{\floor{\alpha^N}}} X_N^j $ on a :
\begin{align}\label{8Haut_produit}
 H(C_{N,e}) = &\theta^{r_e\floor{\alpha^N} +(d-r_e)\floor{\alpha^{N+1}}} \| H_N \| \leq \theta^{r_e{\alpha^N} +(d-r_e){\alpha^{N+1}}} \| H_N \| 
\end{align}
où $H_N$ est le produit extérieur des vecteurs 
\begin{align}\label{8vecteurs_dans_haut_produit}
 (Z_N^j)_{j \in \llbracket 1,r_e \rrbracket} \cup (Z_{N+1}^j)_{j \in \llbracket r_e +1, d \rrbracket} \cup (V_{N+k}^j)_{j \in \llbracket 1, r_e \rrbracket, k\in \llbracket 1,q_e\rrbracket} \cup (V_{N+k}^j)_{j \in \llbracket r_e +1, d \rrbracket, k\in \llbracket 2,q_e\rrbracket}.
\end{align}
On peut majorer cette norme par : 
\begin{align*}
 \|H_N\| \leq \| Z_N^1 \wedge \ldots \wedge Z_N^{r_e} \wedge Z_{N+1}^{r_e+1} \wedge \ldots \wedge Z_{N+1}^{d} \| 
\end{align*}
car les normes des vecteurs $V_k^j$ sont égales à $1$. 
\\Or la quantité $\| Z_N^1 \wedge \ldots \wedge Z_N^{r_e} \wedge Z_{N+1}^{r_e+1} \wedge \ldots \wedge Z_{N+1}^{d} \| $ converge vers $ \| Y_1 \wedge \ldots \wedge Y_d \|$ quand $N$ tend vers l'infini et est donc majorée indépendammant de $N$. Il existe alors $ c_{\ref{8cons_haut_CN_major}}> 0 $, indépendante de $N$ telle que :
\begin{align}\label{8major_HN}
 \|H_N \| \leq c_{\ref{8cons_haut_CN_major}}.
\end{align}

D'autre part, on pose la matrice $M$ dont les vecteurs colonnes sont les vecteurs de $(\ref{8vecteurs_dans_haut_produit})$. 
\\Alors, en reprenant la construction des vecteurs $V_N^j$ en $(\ref{8rec_Xn_Un})$, $M$ est la forme 
\begin{align*}
 M = \begin{pmatrix}
 \begin{matrix}
 I_d \\
 \Sigma_N
 \end{matrix} &
 \begin{matrix}
 0 \\
 V_N
 \end{matrix}
 \end{pmatrix}
\end{align*}
où $\Sigma_N$ est une matrice dont les coefficients sont des $\sigma_{i,j}^N$ ou des $\sigma_{i,j}^{N+1}$, et $V_N$ est une matrice de $\mathcal{M}_{qd, e -d }(\Z)$ de rang $e-d$ car ses colonnes sont $e-d$ vecteurs distincts de la base canonique. 
\bigskip \\
Soit $\Delta$ un mineur non nul de $V_N$ de taille $e-d$. 
 On peut alors extraire une matrice carrée $M'$ de taille $e$ de $M$ en sélectionnant les $d$ premières lignes et $e-d$ parmi les dernières, celles correspondant au mineur $\Delta$. Le déterminant de $M'$ est entier car produit de $\det(I_d)= 1$ et de $\Delta$ qui est un mineur d'une matrice de $M_{e-d}(\Z)$. On a donc $|\det(M') |\geq 1$.
 \\Or $\det(M')$ est un mineur de taille $e$ de $M$, en particulier on a $\|H_N\| \geq |\det(M')| \geq 1$. 
 On trouve alors en combinant avec (\ref{8Haut_produit}) et (\ref{8major_HN}) :
 \begin{align*}
 c_{\ref{8cons_haut_CN_minor}}\theta^{r_e{\alpha^N} +(d-r_e){\alpha^{N+1}}} \leq \theta^{r_e\floor{\alpha^N} +(d-r_e)\floor{\alpha^{N+1}}} \leq H(C_{N,e}) \leq c_{\ref{8cons_haut_CN_major}}\theta^{r_e{\alpha^N} +(d-r_e){\alpha^{N+1}}} 
\end{align*}
et en posant $ c_{\ref{8cons_haut_CN_minor}} = \theta^{-d}$.

\end{preuve}

On s'intéresse maintenant à l'angle $\psi_d(A,C_{N,e})$. Pour cela on étudie d'abord l'angle $\psi_1(\Vect(Y_1), C_{N,e})$.

\begin{lem}\label{8lem_minoration_Y1angleCNE}
 Il existe une constante $\cons \label{8cons_minor_angle_y1_CNe}$ indépendante de $N$ telle que :
 \begin{align*}
 \psi_1(\Vect(Y_1), C_{N,e}) \geq c_{\ref{8cons_minor_angle_y1_CNe}} \theta^{- \alpha^{N+q_e+1}}.
 \end{align*}
\end{lem}

\begin{preuve}
 Soit $X \in C_{N,e}$ tel que $\psi_1(\Vect(Y_1), C_{N,e}) = \omega(Y_1, X)$. \\On utilise la base de $C_{N,e}$ explicitée en $(\ref{8CN_rec})$ et on pose $(a_{k,j})$ une famille de réels tels que $X$ s'écrive
 \begin{align*}
 \sum\limits_{j = 1}^{r_e} a_{0, j} \theta^{-\floor{\alpha^N}} X_N^j+  \sum\limits_{j = r_e +1}^{d} a_{0, j} \theta^{-\floor{\alpha^{N+1}}}X_{N+1}^j +  \sum\limits_{j = 1}^{r_e} \sum\limits_{k = 1}^{q_e} a_{k,j}V_{N+k}^j +  \sum\limits_{j = r_e +1}^{d} \sum\limits_{k =1}^{q_e-1} a_{k,j}V_{N+k+1}^j.
 \end{align*}
Quitte à renormaliser $X$, on suppose que 
\begin{align}\label{8preuve_aik_carre1}
  \sum\limits_{j = 1}^{r_e} a_{0, j}^2+  \sum\limits_{j = r_e +1}^{d} a_{0, j} ^2 +  \sum\limits_{j = 1}^{r_e} \sum\limits_{k = 1}^{q_e} a_{k,j}^2 +  \sum\limits_{j = r_e +1}^{d} \sum\limits_{k =1}^{q_e-1} a_{k,j}^2 =1.
\end{align}
Les vecteurs $ \theta^{-\floor{\alpha^N}}X_{N}^j$ et $V_{N+k}^j$ sont de normes bornées par une constante indépendante de $N$ et on a $\omega(Y_1, X) = \frac{ \| Y_1 \wedge X \|}{ \| Y_1 \| \cdot \|X \|}$. 
Pour conlure, il suffit alors de montrer que si $N$ est assez grand 
\begin{align}\label{8minor_angle_y1_CNe_preuve}
 \| Y_1 \wedge X \| \geq c_{\ref{8cons_minor_angle_y1_CNe}} \theta^{-\floor{ \alpha^{N+q_e+1}}} \geq c_{\ref{8cons_minor_angle_y1_CNe}}\theta^{-{ \alpha^{N+q_e+1}}}.
\end{align}

On rappelle que $Y_1 = \begin{pmatrix}
 1 & 0 & \cdots & 0 & \sigma_{0,1} & \cdots & \sigma_{qd-1, 1}
\end{pmatrix}^\intercal $et on explicite $X$ :
\begin{align*}
 \begin{pmatrix}
 a_{0,1} \\
 \vdots \\ 
 a_{0,d} \\
  \sum\limits_{j = 1}^{r_e} a_{0, j} \sigma_{0,j}^N +  \sum\limits_{j = r_e +1}^{d} a_{0, j} \sigma_{0,j}^{N+1} +  \sum\limits_{j = 1}^{r_e} \sum\limits_{k = 1}^{q_e} a_{k,j}v_{N+k,0}^j +  \sum\limits_{j = r_e +1}^{d} \sum\limits_{k =1}^{q_e-1} a_{k,j}v_{N+k+1,0}^j
 \\
 \vdots \\
  \sum\limits_{j = 1}^{r_e} a_{0, j} \sigma_{qd-1,j}^N +  \sum\limits_{j = r_e +1}^{d} a_{0, j} \sigma_{qd-1,j}^{N+1} +  \sum\limits_{j = 1}^{r_e} \sum\limits_{k = 1}^{q_e} a_{k,j}v_{N+k,qd-1}^j +  \sum\limits_{j = r_e +1}^{d} \sum\limits_{k =1}^{q_e-1} a_{k,j}v_{N+k+1,qd-1}^j
 \end{pmatrix}
\end{align*}
en posant $V_{N+k}^j = \begin{pmatrix}
 0 & \cdots & 0 &v_{N+k,0}^j & \cdots & v_{N+k,qd-1}^j
\end{pmatrix}^\intercal$.
\\ On prouve alors $(\ref{8minor_angle_y1_CNe_preuve})$ en faisant une disjonction de cas. On note $\sigma \geq 1 $ un majorant des $\sigma_{i,j}$, en particulier pour tout $N \in \N$, $\sigma \geq \sigma_{i,j}^N$. 
\\ \textbullet \: \underline{Premier cas :} Si il existe $j \in \llbracket 2, d \rrbracket$ tel que 
\begin{align*}
 |a_{0,j}| \geq \frac{\theta^{-\floor{ \alpha^{N+q_e+1}}}}{(\sigma qd)^2}
\end{align*}
alors en minorant $\| X \wedge Y_1\| $ par le mineur de $(Y_1| X)$ correspondant à la première ligne et la $j$-ième ligne on trouve :
\begin{align*}
 \| X \wedge Y_1 \| \geq \left| \det \begin{pmatrix}
 1 & a_{0,1} \\
 0 & a_{0,j}
 \end{pmatrix} \right| \geq |a_{0,j}| \geq \frac{\theta^{-\floor{ \alpha^{N+q_e+1}}}}{(\sigma qd)^2}
\end{align*}
ce qui donne $(\ref{8minor_angle_y1_CNe_preuve})$.
\\ \textbullet \: \underline{Second cas :} On a sinon :
\begin{align*}
 \forall j \in \llbracket 2, d \rrbracket, \quad |a_{0,j}| < \frac{\theta^{-\floor{ \alpha^{N+q_e+1}}}}{(\sigma qd)^2}.
\end{align*}
D'après $(\ref{8preuve_aik_carre1})$ on a 
\begin{align*}
 a_{0, 1}^2 +  \sum\limits_{j = 1}^{r_e} \sum\limits_{k = 1}^{q_e} a_{k,j}^2 +  \sum\limits_{j = r_e +1}^{d} \sum\limits_{k =1}^{q_e-1} a_{k,j}^2 \geq 1 -(d-1) \left(\frac{\theta^{-\floor{ \alpha^{N+q_e+1}}}}{(\sigma qd)^2} \right)^2.
\end{align*}
En particulier, si $N$ est assez grand, il existe $(j', k')$ tel que $|a_{k',j'}| \geq \frac{1}{qd}$ avec 
\begin{align*}
 k' > 0 \text{ ou } (k'= 0 \text{ et } j' =1). 
\end{align*}
On suppose d'abord qu'on est dans le cas $k' >0$.
\\On pose alors $i = N+ k' + (j'-1)q \mod(qd)$. Par définition des $u_{N+k}^{(i,j)}$ (et donc des $v_{N+k,i}^j$) en $(\ref{8construc_suite_u})$ on a alors :
\begin{align*}
 v_{N+k',i}^{j'} = 1 \text{ et } \forall (j,k) \neq (j',k'), \quad v_{N+k,i}^j = 0.
\end{align*}
En minorant $\| X \wedge Y_1\| $ par le mineur de $(Y_1| X)$ correspondant à la première ligne et la $(i+1+d)$-ième ligne on trouve :
\begin{align*}
 \| X \wedge Y_1 \| &\geq \left| \det \begin{pmatrix}
 1 & a_{0,1} \\
 \sigma_{i,1} &  \sum\limits_{j = 1}^{r_e} a_{0, j} \sigma_{i,j}^N +  \sum\limits_{j = r_e +1}^{d} a_{0, j} \sigma_{i,j}^{N+1} + a_{k',j'}v_{N+k',i}^{j'} 
 \end{pmatrix} \right|\\
 &= \left| a_{0,1}(\sigma_{i,1}^N - \sigma_{i,1} ) +  \sum\limits_{j = 2}^{r_e} a_{0, j} \sigma_{i,j}^N +  \sum\limits_{j = r_e +1}^{d} a_{0, j} \sigma_{i,j}^{N+1} + a_{k',j'}v_{N+k',i}^{j'} \right|.
\end{align*}
Or $|\sigma_{i,1}^N - \sigma_{i,1}| =  \sum\limits_{k= N+1}^{+ \infty} \frac{u_{k}^{(i,1)}}{\theta^{\floor{\alpha^k}}} \leq 4 \theta^{-\floor{\alpha^{N+1}}}$ si $N$ est assez grand. On a donc 
\begin{align*}
 \| X \wedge Y_1 \| &\geq |a_{j',k'}v_{N+k',i}^{j'}| -| a_{0,1}(\sigma_{i,1}^N - \sigma_{i,1} )| -|  \sum\limits_{j = 2}^{r_e} a_{0, j} \sigma_{i,j}^N | - | \sum\limits_{j = r_e +1}^{d} a_{0, j} \sigma_{i,j}^{N+1} | \\
 &\geq \frac{1}{qd} - 4 \theta^{-\floor{\alpha^{N+1}}} - d \sigma \left(\frac{\theta^{-\floor{ \alpha^{N+q_e+1}}}}{(\sigma qd)^2}\right) \\
 &\geq \theta^{-\floor{ \alpha^{N+q_e+1}}}
\end{align*}
si $N$ est assez grand, ce qui donne $(\ref{8minor_angle_y1_CNe_preuve})$.
\\ On suppose maintenant que $k' = 0 $ et $j'= 1$ et donc $|a_{0,1}| \geq \frac{1}{qd}$. \\ On pose $ i = N+q_e +1 \mod(qd)$. Par définition des $u_{N+k}^{(i,j)}$ (et donc des $v_{N+k,i}^j$) en $(\ref{8construc_suite_u})$ on a alors :
\begin{align*}
 \forall k \in \llbracket 1, q_e \rrbracket, \quad \forall j \in \llbracket 2, d \rrbracket, \quad v_{N+k,i}^j = 0
\end{align*}
et $u_{N+1}^{(i,1)} = \ldots = u_{N+q_e}^{(i,1)} = 0$ et $u_{N+q_e+ 1}^{(i,1)} \in \{2,3\}$. 
\\En particulier $\sigma_{i,1}^N =  \sum\limits_{k=0}^{N} \frac{u_{k}^{(i,1)}}{\theta^{\floor{\alpha^k}}} = \sum\limits_{k=0}^{N+q_e} \frac{u_{k}^{(i,1)}}{\theta^{\floor{\alpha^k}}} = \sigma_{i,1}^{N+q_e}$ et $$ \sum\limits_{j = 1}^{r_e} \sum\limits_{k = 1}^{q_e} a_{k,j}v_{N+k,i}^j +  \sum\limits_{j = r_e +1}^{d} \sum\limits_{k =1}^{q_e-1} a_{k,j}v_{N+k+1,i}^j = 0.$$
En minorant $\| X \wedge Y_1\| $ par le mineur de $(X| Y_1)$ correspondant à la première ligne et la $(i+1+d)$-ième ligne on trouve :
\begin{align*}
 \| X \wedge Y_1 \| &\geq \left| \det\begin{pmatrix}
 1 & a_{0,1} \\
 \sigma_{i,1} &  \sum\limits_{j = 1}^{r_e} a_{0, j} \sigma_{i,j}^N +  \sum\limits_{j = r_e +1}^{d} a_{0, j} \sigma_{i,j}^{N+1} 
 \end{pmatrix} \right|\\
 &= \left| a_{0,1}(\sigma_{i,1}^N - \sigma_{i,1} ) +  \sum\limits_{j = 2}^{r_e} a_{0, j} \sigma_{i,j}^N +  \sum\limits_{j = r_e +1}^{d} a_{0, j} \sigma_{i,j}^{N+1} \right|.
\end{align*}
Or $|\sigma_{i,1}^N - \sigma_{i,1}| = |\sigma_{i,1}^{N+q_e} - \sigma_{i,1}| =  \sum\limits_{k= N+q_e+1}^{+ \infty} \frac{u_{k}^{(i,1)}}{\theta^{\floor{\alpha^k}}} \geq \frac{u_{N+q_e +1 }^{(i,1)}}{\theta^{\floor{\alpha^{N+q_e+1}}}} \geq 2 \theta^{-\floor{\alpha^{N+q_e+1}}}$. On a donc 
\begin{align*}
 \| X \wedge Y_1 \| &\geq |a_{0,1}(\sigma_{i,1}^N - \sigma_{i,1} )| - |  \sum\limits_{j = 2}^{r_e} a_{0, j} \sigma_{i,j}^N | - | \sum\limits_{j = r_e +1}^{d} a_{0, j} \sigma_{i,j}^{N+1} | \\
 &\geq \frac{2 \theta^{-\floor{\alpha^{N+q_e+1}}}}{qd}- d \sigma \left(\frac{\theta^{-\floor{ \alpha^{N+q_e+1}}}}{(\sigma qd)^2}\right) \\
 &\geq \frac{\theta^{-\floor{ \alpha^{N+q_e+1}}}}{qd}
\end{align*}
 ce qui prouve $(\ref{8minor_angle_y1_CNe_preuve})$.
\bigskip 
\\ On a donc $(\ref{8minor_angle_y1_CNe_preuve})$ dans tous les cas ce qui termine la preuve du lemme. 

\end{preuve}

\begin{lem}\label{8lem_psid_ACNe}
 Il existe des constantes $\cons \label{8cons_minor_angle_psid}$ et $ \cons \label{8cons_major_angle_psid}$ indépendantes de $N$ telles que 
 \begin{align*}
 c_{\ref{8cons_minor_angle_psid}} H(C_{N,e})^{\frac{- \alpha^{q_e +1}}{r_e + (d-r_e)\alpha}} \leq \psi_d(A, C_{N,e}) \leq c_{\ref{8cons_major_angle_psid}} H(C_{N,e})^{\frac{- \alpha^{q_e +1}}{r_e + (d-r_e)\alpha}}.
 \end{align*}
\end{lem}

\begin{preuve}
 On rappelle que pour $j \in \llbracket 1,d \rrbracket$ et $N \in \N$ 
 \begin{align*}
 \psi_1(\Vect(Y_j), \Vect(X_{N+q_e}^j)) \leq c_{\ref{8cons_major_angle_Yj_XNj}} \theta^{-\alpha^{N+q_e+1}}
 \end{align*}
 d'après $(\ref{8major_Yj_XNj})$.
 \\Par construction de $C_{N,e}$, on a $\Vect(X_{N+q_e}^1, \ldots, X_{N+q_e}^d) \subset C_{N,e}$ et donc $$ \psi_d(A,C_{N,e}) \leq \psi_d(A, \Vect(X_{N+q_e}^1, \ldots, X_{N+q_e}^d))$$ d'après le lemme~\ref{lem_inclusion_croissance}.
\\D'après la propriété~\ref{2prop_4.5Elio} on a :
 \begin{align*}
 \psi_d(A, \Vect(X_{N+q_e}^1, \ldots, X_{N+q_e}^d)) &\leq c_{\ref{8cons_maj_prop45_appliquee}} \sum\limits_{j=1}^{d} \psi_1(\Vect(Y_j), \Vect(X_{N+q_e}^j) ) 
 \end{align*}
 avec $\cons \label{8cons_maj_prop45_appliquee} $ dépendant de $Y_1, \ldots, Y_d $ et $n$.
On a donc 
\begin{align*}
 \psi_d(A, \Vect(X_{N+q_e}^1, \ldots, X_{N+q_e}^d)) &\leq c_{\ref{8cons_maj_prop45_appliquee}}c_{\ref{8cons_major_angle_Yj_XNj}} d \theta^{-{\alpha^{N+q_e +1 }} }\\
 &\leq c_{\ref{8cons_maj_prop45_appliquee}}c_{\ref{8cons_major_angle_Yj_XNj}} d c_{\ref{8cons_haut_CN_major}}^{\frac{{\alpha^{N+q_e +1 }}}{r_e{\alpha^N} +(d-r_e){\alpha^{N+1}} }} H(C_{N,e})^{\frac{-{\alpha^{N+q_e +1 }}}{r_e{\alpha^N} +(d-r_e){\alpha^{N+1}} }} 
 \\
 &= c_{\ref{8cons_major_angle_psid}} H(C_{N,e})^{\frac{-{\alpha^{q_e +1 }}}{r_e +(d-r_e)\alpha }} 
 \end{align*}
car $H(C_{N,e}) \leq c_{\ref{8cons_haut_CN_major}} \theta^{r_e{\alpha^N} +(d-r_e){\alpha^{N+1}}} $ d'après le lemme~\ref{8lem_hauteur_CN}, ce qui prouve la majoration du lemme avec $c_{\ref{8cons_major_angle_psid}} = c_{\ref{8cons_maj_prop45_appliquee}}c_{\ref{8cons_major_angle_Yj_XNj}} d c_{\ref{8cons_haut_CN_major}}^{\frac{{\alpha^{q_e +1 }}}{r_e +(d-r_e){\alpha} }} $.

\bigskip 
Pour montrer la minoration on utilise le fait que $\psi_d(A,C_{N,e}) \geq 
\psi_1(\Vect(Y_1), C_{N,e}) $ d'après le lemme~\ref{lem_inclusion_croissance}. Le lemme~\ref{8lem_minoration_Y1angleCNE} donne alors 
\begin{align*}
 \psi_d(A,C_{N,e}) \geq c_{\ref{8cons_minor_angle_y1_CNe}} \theta^{- \alpha^{N+q_e+1}}.
\end{align*}
Comme $H(C_{N,e}) \geq c_{\ref{8cons_haut_CN_minor}} \theta^{r_e{\alpha^N} +(d-r_e){\alpha^{N+1}}} $ d'après le lemme~\ref{8lem_hauteur_CN}, on a la minoration 
\begin{align*}
 \psi_d(A,C_{N,e}) \geq c_{\ref{8cons_minor_angle_y1_CNe}}c_{\ref{8cons_haut_CN_minor}}^{\frac{{\alpha^{N+q_e +1 }}}{r_e{\alpha^N} +(d-r_e){\alpha^{N+1}} }} H(C_{N,e})^{\frac{{-\alpha^{N+q_e +1 }}}{r_e{\alpha^N} +(d-r_e){\alpha^{N+1}} }} = c_{\ref{8cons_minor_angle_psid}}H(C_{N,e})^{\frac{-{\alpha^{q_e +1 }}}{r_e +(d-r_e)\alpha }} 
\end{align*}
avec $c_{\ref{8cons_minor_angle_psid}} = c_{\ref{8cons_minor_angle_y1_CNe}}c_{\ref{8cons_haut_CN_minor}}^{\frac{{\alpha^{q_e +1 }}}{r_e +(d-r_e){\alpha} }} $ ce qui prouve la minoration du lemme.

\end{preuve}

On a alors construit une construit une infinité d'espaces rationnels $C_{N,e}$ de dimension $e$ tels que 
\begin{align*}
 \psi_d(A,C_{N,e}) \leq c_{\ref{8cons_major_angle_psid}}H(C_{N,e})^{\frac{-{\alpha^{q_e +1 }}}{r_e +(d-r_e)\alpha }} 
\end{align*}
ce qui donne en particulier : 
\begin{align*}
 \mu_n(A|e)_d \geq \frac{{\alpha^{q_e +1 }}}{r_e +(d-r_e)\alpha }.
\end{align*}

\subsection{Majoration de l'exposant}
On montre dans cette section, que la minoration trouvée dans la partie précédente est optimale c'est-à-dire que l'on va majorer $\mu_n(A|e)_d$ par $\frac{{\alpha^{q_e +1 }}}{r_e +(d-r_e)\alpha }$. \\
On énonce un premier lemme technique qui nous sera utile dans la démonstration de la majoration. Ce lemme généralise en fait la minoration du lemme~\ref{8lem_hauteur_CN}.
\\ On rappelle que pour $N$ et $v$ deux entiers non nuls on a défini :
\begin{align*}
 B_{N,v} = \Vect(X_N^1, X_{N+1}^1,\ldots, X_{N+v-1}^1, X_N^2 \ldots, X_{N+v-1}^2, \ldots, X_N^d,\ldots, X_{N+v-1}^d).
\end{align*}

\begin{lem}\label{8minoration_somme_espace}
Soit $N \in \N$, $v \in \llbracket 1,q -1 \rrbracket $ et $r \in \llbracket 0, d-1 \rrbracket $.
\\Pour $W $ un sous-espace rationnel de $\Vect(X_N^1, \ldots, X_N^d)$ de dimension $r$, on a :
\begin{align*}
 H(B_{N+1,v}\oplus W) \geq c_{\ref{8cons_minoration_somme_espace}} \theta^{ r {\alpha^N} +(d-r) {\alpha^{N+1}}} 
\end{align*}
avec $\cons \label{8cons_minoration_somme_espace}$ indépendante de $N$ et de $W$.

 \end{lem}

 \begin{preuve}
 Soit $U_1, \ldots, U_r$ une $\Zbase$ de $W \cap \Z^n$. Comme $W \subset \Vect(X_N^1, \ldots, X_N^d)$ et que ces vecteurs forment une $\Zbase$ de $\Vect(X_N^1, \ldots, X_N^d) \cap \Z^n$ on peut écrire pour $i \in \llbracket 1, r \rrbracket$:
 \begin{align}\label{8baseW}
 U_i =  \sum\limits_{j = 1}^d a_{i,j}X_N^j
 \end{align}
 avec $a_{i,j} \in \Z$.
 \\ On remarque que $B_{N+1,v}\oplus W $ est de dimension $\dim(B_{N+1,v}) + r = vd + r $ d'après le lemme~\ref{8lem_dimBN_zbaseBN}.
\\ D'après le lemme~\ref{8lem_dimBN_zbaseBN}, les vecteurs $(X_{N+1}^j)_{j \in \llbracket 1,d\rrbracket} \cup (V_k^j)_{j \in \llbracket 1,d\rrbracket, k \in \llbracket N+2, N+v\rrbracket}$ forment une $\Zbase$ de $B_{N+1,v} \cap \Z^n$.
\\ En concaténant ces bases de $B_{N+1,v} $ et de $W$ on forme une base de l'espace vectoriel réel $B_{N+1,v}\oplus W$ (mais pas forcément du $\Ztir$module $(B_{N+1,v}\oplus W) \cap \Z^n$). De plus les vecteurs de cette base sont entiers, la formule énoncée dans la remarque~\ref{2req_maj_hauteur_base_entiere} donne la hauteur de $B_{N+1,v}\oplus W$ dans ce cas-là :
\begin{align}\label{8hauteur_somme}
 H(B_{N+1,v}\oplus W) = \frac{\left\| (\bigwedge\limits_{j =1}^d X_{N+1}^j) \wedge \left(\bigwedge\limits_{k = N+2}^{N+v} (V_k^1 \wedge \ldots \wedge V_k^d) \right) \wedge U_1 \wedge \ldots \wedge U_r \right\| }{N(I)}
\end{align}
où $I$ est l'idéal engendré par les coordonnées
grassmaniennes associées à cette base, qui, on le rappelle, sont les mineurs de taille $vd + r$ de la matrice associée aux vecteurs $$(X_{N+1}^j)_{j \in \llbracket 1,d\rrbracket} \cup (V_k^j)_{j \in \llbracket 1,d\rrbracket, k \in \llbracket N+2, N+v\rrbracket} \cup (U_i)_{i \in \llbracket 1,r \rrbracket}.$$
\\D'après (\ref{8baseW}) on a pour $i \in \llbracket 1,r \rrbracket$ : 
\begin{align*}
 U_i =  \sum\limits_{j = 1}^d a_{i,j}X_N^j =  \sum\limits_{j = 1}^d a_{i,j}\frac{X_{N+1}^j -U_{N+1}^j}{\theta^{\floor{\alpha^{N+1}} - \floor{\alpha^{N}}} }
\end{align*}
d'après la formule (\ref{8rec_Xn_Un}).
\\ Ainsi :
\begin{align}
 &\| (\bigwedge\limits_{j =1}^d X_{N+1}^j) \wedge \left(\bigwedge\limits_{k = N+2}^{N+v} (V_k^1 \wedge \ldots \wedge V_k^d) \right) \wedge U_1 \wedge \ldots \wedge U_r \| \nonumber \\&= \left(\frac{1}{\theta^{\floor{\alpha^{N+1}} - \floor{\alpha^{N}}}}\right)^r \| (\bigwedge\limits_{j =1}^d X_{N+1}^j) \wedge \left(\bigwedge\limits_{k = N+2}^{N+v} (V_k^1 \wedge \ldots \wedge V_k^d) \right) \wedge \sum\limits_{j = 1}^d a_{1,j}U_{N+1}^j\wedge \ldots \wedge \sum\limits_{j = 1}^d a_{r,j}U_{N+1}^j \| \nonumber 
 \\&= \theta^{ r \floor{\alpha^N} +(d-r) \floor{\alpha^{N+1}}} \| (\bigwedge\limits_{j =1}^d Z_{N+1}^j) \wedge \left(\bigwedge\limits_{k = N+2}^{N+v} (V_k^1 \wedge \ldots \wedge V_k^d) \right) \wedge \sum\limits_{j = 1}^d a_{1,j}U_{N+1}^j\wedge \ldots \wedge \sum\limits_{j = 1}^d a_{r,j}U_{N+1}^j \|\label{8minot_prod_ext}
\end{align}
avec $Z_{N+1}^j = \frac{1}{\theta^{\floor{\alpha^{N+1}} }}X_{N+1}^j$.
 \\ Or la norme du produit extérieur qui apparaît en (\ref{8minot_prod_ext}) peut être minorée par la valeur absolue de tout mineur de taille $dv + r $ de la matrice $M$ dont les colonnes sont les vecteurs $(Z_{N+1}^j)_{j \in \llbracket 1, d\rrbracket} \cup (V_k^j)_{j \in \llbracket 1,d\rrbracket, k \in \llbracket N+2, N+v\rrbracket} \cup (\sum\limits_{j = 1}^d a_{i,j}U_{N+1}^j)_{i \in \llbracket 1,r \rrbracket}$. On a
\begin{align*}
 M = \begin{pmatrix}
 \begin{matrix}
 I_d \\
 \Sigma_{N+1}
 \end{matrix} &
 \begin{matrix}
 A
 \end{matrix}
 \end{pmatrix} \in \MM_{n,dv+r}(\R)
\end{align*}
où $\Sigma_{N+1} = (\sigma_{i,j }^{N+1})_{ i \in \llbracket 0,qd-1 \rrbracket, j \in \llbracket 1,d\rrbracket } $ et $A \in \MM_{n,d(v-1)+r}(\R)$ dont les colonnes sont les vecteurs $(V_k^j)_{j \in \llbracket 1,d\rrbracket, k \in \llbracket N+2, N+v\rrbracket} \cup (\sum\limits_{j = 1}^d a_{i,j}U_{N+1}^j)_{i \in \llbracket 1,r \rrbracket}$. 
\\Or d'après la construction des vecteurs $V_k^j$ et $U_{N+1}^j$ en (\ref{8rec_Xn_Un}) et comme $a_{i,j} \in \Z$, il existe $A' \in \MM_{n-d, d(v-1) + r }(\Z)$ une matrice à coefficients entiers telle que : 
\begin{align*}
 M = \begin{pmatrix}
 \begin{matrix}
 I_d \\
 \Sigma_{N+1}
 \end{matrix} &
 \begin{matrix}
 A
 \end{matrix}
 \end{pmatrix} = \begin{pmatrix}
 \begin{matrix}
 I_d \\
 \Sigma_{N+1}
 \end{matrix} &
 \begin{matrix}
 0 \\
 A'
 \end{matrix}
 \end{pmatrix}.
\end{align*}
De plus la matrice $A'$ est de rang $d(v-1) + r$ car $\rang(M) = dv + r$ puisque les colonnes de $M$ sont des vecteurs d'une base. On peut donc extraire un mineur de $\begin{pmatrix}
 \begin{matrix}
 I_d \\
 \Sigma_{N+1}
 \end{matrix} &
 \begin{matrix}
 0 \\
 A'
 \end{matrix}
 \end{pmatrix}$ qui soit non nul et entier. En particulier il est minoré en valeur absolue par $1$ et en reprenant (\ref{8minot_prod_ext}) on a :
 \begin{align}\label{8derniere_mino_prod}
 \| (\bigwedge\limits_{j =1}^d X_{N+1}^j) \wedge \left(\bigwedge\limits_{k = N+2}^{N+v} (V_k^1 \wedge \ldots \wedge V_k^d) \right) \wedge U_1 \wedge \ldots \wedge U_r \| \geq \theta^{ r \floor{\alpha^N} +(d-r) \floor{\alpha^{N+1}}}.
 \end{align}
D'après $(\ref{8hauteur_somme})$, il reste alors à montrer que $N(I)$ est majorée par une constante dépendant seulement de $A$. 
\\ Pour cela on montre que la famille $(X_{N+1}^j)_{j \in \llbracket 1,d\rrbracket} \cup (V_k^j)_{j \in \llbracket 1,d\rrbracket, k \in \llbracket N+2, N+v\rrbracket} \cup (U_i)_{i \in \llbracket 1,r \rrbracket}$ forme \og presque \fg{} une $\Zbase$ de $(B_{N+1,v}\oplus W)\cap \Z^n$. 
\\ Soit $(u_j)_{j \in \llbracket 1,d\rrbracket} \cup (w_{j,k})_{j \in \llbracket 1,d\rrbracket, k \in \llbracket N+2, N+v\rrbracket} \cup (v_i)_{i \in \llbracket 1,r \rrbracket} \in [0,1]^{r + vd}$ tel que :
\begin{align}\label{8Xentier_preuve}
 X =  \sum\limits_{ j = 1}^d u_j X_{N+1}^j + \sum\limits_{\underset{j \in \llbracket 1,d\rrbracket }{k \in \llbracket N+2, N+v\rrbracket} } w_{j,k} V_k^j + \sum\limits_{ i= 1}^r v_i U_i \in \Z^n. 
\end{align}
On décompose $X$ dans la $\Zbase$ de $B_{N,v+1} \cap \Z^n$ que l'on a explicitée dans le lemme~\ref{8lem_dimBN_zbaseBN}.
\\On a $ U_i =  \sum\limits_{j = 1}^d a_{i,j}X_N^j $ d'après $(\ref{8baseW})$, et $X_{N+1}^j = \theta^{\floor{\alpha^{N+1}} - \floor{\alpha^{N}}}X_N^j + \| U_{N+1}^j\|V_{N+1}^j$ en reprenant la formule (\ref{8rec_Xn_Un}).
\\ Pour $k \in \llbracket N+2, N+v\rrbracket $, les $V_k^j$ n'apparaissent donc pas dans la décomposition des $U_i$ et $X_{N+1}^j$. Par définition d'une $\Zbase$ on a alors
\begin{align*}
 \forall j \in \llbracket 1,d\rrbracket, \quad \forall k \in \llbracket N+2, N+v\rrbracket, \quad w_{j,k} \in \Z .
\end{align*}
En particulier $ \sum\limits_{\underset{j \in \llbracket 1,d\rrbracket }{k \in \llbracket N+2, N+v\rrbracket} } w_{j,k} V_k^j \in \Z^n$ et la relation $(\ref{8Xentier_preuve})$ donne alors 
\begin{align*}
 \sum\limits_{ j = 1}^d u_j (\theta^{\floor{\alpha^{N+1}} - \floor{\alpha^{N}}}X_N^j + \| U_{N+1}^j\|V_{N+1}^j) + \sum\limits_{j = 1}^d \left(\sum\limits_{ i= 1}^r v_i a_{i,j} \right) X_N^j \in \Z^n.
\end{align*}
\\D'après le lemme~\ref{8lem_dimBN_zbaseBN}, les $V_{N+1}^j$ et $X_N^j$ forment une $\Zbase$ de $V \cap \Z^n$ où $V$ est le sous-espace vectoriel réel engendré par ces vecteurs. On a donc pour tout $j \in \llbracket 1, d\rrbracket$ :
\begin{align}
 u_j \| U_{N+1}^j\| \in \Z, \label{8inZ1} \\
 u_j \theta^{\floor{\alpha^{N+1}} - \floor{\alpha^{N}}} + \sum\limits_{ i= 1}^r v_i a_{i,j} \in \Z.\label{8inZ2}
\end{align}
Comme $\| U_{N+1}^j\| = 2 $ ou $3$, la relation $(\ref{8inZ1})$ donne $6u_j \in \Z $ pour tout $j \in \llbracket 1,d \rrbracket$. Par la deuxième relation $(\ref{8inZ2})$ on a alors $\sum\limits_{ i= 1}^r6 v_i a_{i,j} \in \Z$ pour tout $j \in \llbracket 1, d\rrbracket$. \\Enfin :
\begin{align*}
 \sum\limits_{i = 1 }^r 6 v_i U_i =  \sum\limits_{j= 1}^d \left(\sum\limits_{i = 1 }^r 6 v_i a_{i,j} \right) X_N^{j} \in \Z^n.
\end{align*}
Comme $U_1, \ldots, U_r$ est une $\Zbase$ de $W \cap \Z^n$ on trouve $ 6v_i \in \Z$ pour tout $i \in \llbracket 1, r\rrbracket$.
On a donc finalement :
\begin{align*}
 w_{j,k} \in \Z& \text{ pour } j \in \llbracket 1,d\rrbracket, k \in \llbracket N+2, N+v\rrbracket,\\
 6u_j \in \Z& \text{ pour } j \in \llbracket 1,d\rrbracket,\\
 6 v_i \in \Z& \text{ pour } i \in \llbracket 1,r\rrbracket. 
\end{align*}
En particulier cela donne $N(I) \leq 6^{d+r} \leq 6^{2d}$.
En combinant cela avec (\ref{8derniere_mino_prod}) dans (\ref{8hauteur_somme}) on trouve :
\begin{align*}
 H(B_{N+1,v} \oplus W) \geq 6^{-2d}\theta^{ r \floor{\alpha^N} +(d-r) \floor{\alpha^{N+1}}} \geq c_{\ref{8cons_minoration_somme_espace}} \theta^{ r {\alpha^N} +(d-r) {\alpha^{N+1}}}.
\end{align*}
ce qui est le résultat attendu avec $c_{\ref{8cons_minoration_somme_espace}} = 6^{-2d}\theta^{-d} $.

 \end{preuve}

 On peut maintenant minorer le $d\tir$ième angle que réalise $A$ avec tout espace rationnel de dimension $e$.

\begin{lem}\label{8lem_inclusion_BN_C}
 Soit $ \varepsilon >0$ et $C$ un espace rationnel de dimension $e$ tel :
 \begin{align*}
 \psi_d(A,C) \leq H(C) ^{-  \frac{{\alpha^{q_e +1 }}}{r_e +(d-r_e){\alpha}} - \varepsilon }.
 \end{align*}
 Supposons $H(C)$ assez grand, et notons $N \in \N$ 
 \begin{align}\label{8choixN}
 \theta^{\alpha^{N+q_e}} \leq H(C) ^{  \frac{\alpha^{q_e +1 }}{r_e +(d-r_e)\alpha}
 + \frac{\varepsilon}{2} -1 } < \theta^{{\alpha^{N+q_e +1}}}.
\end{align}
Alors $B_{N+1,q_e} \subset C$.
\end{lem}

\begin{preuve}
 Soit $N$ l'unique entier vérifiant $(\ref{8choixN})$. 
\\On pose $Z_1, \ldots, Z_{e} $ une $\Zbase$ de $C\cap \Z^n$.
\\
 Soit $j \in \llbracket 1,d \rrbracket $. On étudie 
 $D_{j,k} = \| X_{N+k}^j \wedge Z_1 \wedge \ldots \wedge Z_{e} \|$ pour $k \in \llbracket 1,q_e \rrbracket$.
 D'après le lemme~\ref{lem_phi_dim1}, on a :
 \begin{align*}
 D_{j,k} &= \psi_1(\Vect(X_{N+k}^j),C) \|X_{N+k}^j \| H(C) \\
 &\leq c_{\ref{8cons_major_norme_XN}} \theta^{\alpha^{N+k}}H(C) \left(\omega(X_{N+k}^j, Y_{j}) + \psi_1(\Vect(Y_j),C)\right)
\end{align*}
car $ \|X_{N+k}^j\| \leq c_{\ref{8cons_major_norme_XN}} \theta^{\alpha^{N+k}} $ d'après $(\ref{8norme_XN})$ et $\omega(X_{N+k}^j,C) \leq \omega(X_{N+k}^j, Y_{j}) + \psi_1(\Vect(Y_j),C)$ par l'inégalité triangulaire du lemme~\ref{2lem_ineg_triang}. \\
Or $\psi_1(Y_j,C) \leq \psi_d(A, C)$ d'après le lemme~\ref{lem_inclusion_croissance} car $e \geq d $ et $ \omega(X_{N+k}^j, Y_{j}) \leq c_{\ref{8cons_major_angle_Yj_XNj}} \theta^{-{\alpha^{N+k+1}}} $ d'après $(\ref{8major_Yj_XNj})$ et donc :
\begin{align*}
  D_{j,k} &\leq c_{\ref{8cons_maj_1Djk}}\theta^{{\alpha^{N+k}}} H(C) \left(\theta^{-{\alpha^{N+k+1}}} + H(C) ^{\frac{-\alpha^{q_e +1 }}{r_e +(d-r_e)\alpha} -   \varepsilon }\right)
\end{align*}
avec $\cons \label{8cons_maj_1Djk}$ indépendante de $N$.
\\ Or $\theta^{{\alpha^{N+k}}} \leq \theta^{{\alpha^{N+q_e}}} \leq H(C) ^{   \frac{\alpha^{q_e +1 }}{r_e +(d-r_e)\alpha}
 + \frac{\varepsilon}{2} -1 } $ par le choix de $N$ en $(\ref{8choixN})$ et donc :
\begin{align*}
 D_{j,k} &\leq c_{\ref{8cons_maj_1Djk}} \theta^{{\alpha^{N+k}} -{\alpha^{N+k+1}}} H(C) + c_{\ref{8cons_maj_1Djk}}H(C) ^{-\frac{ \varepsilon }{2}}.
\end{align*}
De plus $H(C) ^{   \frac{\alpha^{q_e +1 }}{r_e +(d-r_e)\alpha} + \frac{\varepsilon}{2} -1} \leq \theta^{{\alpha^{N+q_e + 1}}} $ d'après $(\ref{8choixN})$ et donc $$\theta \geq H(C)^{ \frac{ \frac{\alpha^{q_e +1 }}{r_e +(d-r_e)\alpha} + \frac{\varepsilon}{2} -1 }{{\alpha^{N+q_e + 1}}} }. $$ 
\\ On obtient alors la majoration : 
\begin{align}\label{8major_Djk}
 D_{j,k} &\leq c_{\ref{8cons_maj_1Djk}} H(C)^{1+ \frac{ ({\alpha^{N+k}} -{\alpha^{N+k+1}})(\frac{{\alpha^{q_e +1 }}}{r_e +(d-r_e)\alpha} + \frac{\varepsilon}{2} -1 )}{{\alpha^{N+q_e + 1}}}} + c_{\ref{8cons_maj_1Djk}}H(C) ^{-\frac{ \varepsilon }{2}}.
\end{align}
On étudie l'exposant et comme $ k \geq 1 $ :
\begin{align*}
 1+ \frac{ ({\alpha^{N+k}} -{\alpha^{N+k+1}})(\frac{{\alpha^{q_e +1 }}}{r_e +(d-r_e)\alpha} + \frac{\varepsilon}{2} -1 )}{{\alpha^{N+q_e + 1}}} &\leq 1+ \frac{ ({\alpha^{N+1}} -{\alpha^{N+2}})(\frac{{\alpha^{q_e +1 }}}{r_e +(d-r_e)\alpha} + \frac{\varepsilon}{2} -1 )}{{\alpha^{N+q_e + 1}}} \\
 &= 1+ \frac{ (1 -{\alpha})(\frac{{\alpha^{q_e +1 }}}{r_e +(d-r_e)\alpha} + \frac{\varepsilon}{2} -1 )}{{\alpha^{q_e }}}\\
 &\leq \frac{\alpha^{q_e } + (1 -{\alpha})(\frac{{\alpha^{q_e }}}{d} + \frac{\varepsilon}{2} -1 )}{{\alpha^{q_e }}}
\end{align*}
car $\frac{{\alpha^{q_e +1 }}}{r_e +(d-r_e)\alpha} \geq \frac{{\alpha^{q_e }}}{d} \geq 1 $ et $1-\alpha <-1$. De plus 
\begin{align*}
 \alpha^{q_e } + (1 -{\alpha})(\frac{{\alpha^{q_e }}}{d} -1 ) &\leq \frac{1}{d}\left(-\alpha^{q_e+1} + \alpha^{q_e}(d+1) + d(\alpha - 1) \right) \\
 &\leq \frac{1}{d}\left(-\alpha^{2} + \alpha(2d+2) - d \right) 
\end{align*}
car $e \geq d$ et donc $q_e \geq 1$. Par l'inégalité $(\ref{8lem_inegal_alpha1})$ du lemme~\ref{8lem_inegal_alpha} on a $-\alpha^{2} + \alpha(2d+2) - d \leq 0$. Alors 
\begin{align*}
 \frac{\alpha^{q_e } + (1 -{\alpha})(\frac{{\alpha^{q_e }}}{d} + \frac{\varepsilon}{2} -1 )}{{\alpha^{q_e }}} \leq \frac{(1-\alpha)\varepsilon}{2 \alpha^{q_e}} \leq \frac{-\varepsilon}{2 \alpha^{q}}
\end{align*}
car $1-\alpha <-1$ et $q_e \leq q$.
 \\On a donc en reprenant (\ref{8major_Djk}): 
\begin{align*}
 D_{j,k} &\leq c_{\ref{8cons_maj_1Djk}}H(C)^{- \frac{\varepsilon}{2\alpha^q}} + c_{\ref{8cons_maj_1Djk}}H(C) ^{-\frac{ \varepsilon }{2}} \leq 2c_{\ref{8cons_maj_1Djk}}H(C)^{- \frac{\varepsilon}{2\alpha^q}} .
\end{align*}
Pour $H(C)$ assez grand en fonction de $\varepsilon, \alpha, n $ et $c_{\ref{8cons_maj_1Djk}} $, on a $ D_{j,k} < 1$ pour tous $j \in \llbracket 1,d \rrbracket $ et $k \in \llbracket 1,q_e \rrbracket$. Le lemme~\ref{2lem_X_in_B} donne alors :
\begin{align*}
 \forall j \in \llbracket 1,d \rrbracket, \quad \forall k \in \llbracket 1,q_e \rrbracket, \quad X_{N+k}^j \in C.
\end{align*}
 Ces vecteurs engendrent $B_{N+1,q_e}$ donc $B_{N+1,q_e } \subset C$. 

 \end{preuve}

\begin{lem}\label{8lem_minoration_exposant}
 Soit $ \varepsilon >0$. Pour tous les espaces rationnels $C$ de dimension $e$ sauf un nombre fini on a :
 \begin{align*}
 \psi_d(A,C) > H(C) ^{- \frac{{\alpha^{q_e +1 }}}{r_e +(d-r_e){\alpha}} - \varepsilon }.
 \end{align*}
 
\end{lem}

\begin{preuve}
On suppose l'inverse par l'absurde. Alors pour un certain $\varepsilon>0$, il existe une infinité d'espaces rationnels $C$ de dimension $e$ vérifiant :
 \begin{align}\label{8hyp_absurde_preuve_lem14}
 \psi_d(A,C) \leq H(C) ^{  \frac{{-\alpha^{q_e +1 }}}{r_e +(d-r_e){\alpha}} - \varepsilon }.
 \end{align}
D'après le lemme~\ref{8lem_inclusion_BN_C}, si un tel $C$ a une hauteur suffisamment grande on a $B_{N+1, q_e} \subset C$ avec $N$ vérifiant $(\ref{8choixN})$.

On sépare la preuve en deux cas selon la valeur de $r_e$.

\bigskip

 \textbullet \: \underline{Premier cas $r_e =0$ :} Dans ce cas, on a $e = q_e d $ et $\frac{{\alpha^{q_e +1 }}}{r_e +(d-r_e){\alpha}} = \frac{\alpha^{q_e}}{d}$. \\On a déjà $B_{N+1,q_e} \subset C$; par égalité des dimensions et comme $r_e =0$ on a $$C = B_{N+1,q_e} = C_{N,e}$$ par la remarque~\ref{8req_CN=BN} et d'après le lemme~\ref{8lem_dimBN_zbaseBN}.
\\ On a donc $H(C) = H(C_{N,e})$ et par le lemme~\ref{8lem_psid_ACNe} on a :
\begin{align*}
 c_{\ref{8cons_minor_angle_psid}}H(C)^{\frac{-{\alpha^{q_e+1}}}{r_e+(d-r_e){\alpha}} }\leq \psi_d(A,C_{N,e}) = \psi_d(A,C) \leq H(C) ^{-\frac{{\alpha^{q_e +1 }}}{r_e +(d-r_e){\alpha}} -  \varepsilon }.
\end{align*}
 Cela donne alors $\frac{{\alpha^{q_e+1}}}{r_e +(d-r_e){\alpha}} \geq \frac{{\alpha^{q_e +1 }}}{r_e +(d-r_e){\alpha}} + \frac{\varepsilon }{2} $ car on peut faire tendre $H(C) $ vers $+ \infty$. On trouve alors
 $\frac{\varepsilon}{2} \leq 0$ ce qui est contradictoire.
 
\bigskip

 \textbullet \: \underline{Deuxième cas $r_e \neq 0$ :} Dans ce cas $e = q_e d + r_e $ et en particulier $1 \leq q_e \leq q-1$.\\
Montrons que $\Vect(X_N^1, \ldots, X_N^d) \cap C $ est de dimension supérieure ou égale à $r_e$ en montrant par récurrence sur $r \in \llbracket 0,r_e-1 \rrbracket$ qu'il existe $r+1$ vecteurs entiers linéairement indépendants dans $\Vect(X_N^1, \ldots, X_N^d) \cap C $.
\\Soit $r \in \llbracket 0, r_e -1 \rrbracket$. On suppose qu'il existe $U_1, \ldots, U_r$ vecteurs entiers linéairement indépendants dans $\Vect(X_N^1, \ldots, X_N^d) \cap C $ ; si $r = 0$ cette hypothèse est vide donc vraie. Montrons qu'il existe $U_{r+1} \in \Vect(X_N^1, \ldots, X_N^d) \cap C \cap \Z^n $ tel que $U_1, \ldots, U_{r+1}$ soient linéairement indépendants.
\bigskip

Comme $q_e \leq q-1$ et que la famille de vecteurs intervenant dans la définition de $B_{N+1, q_e}$ en $(\ref{8def_Bnv})$ est libre, on a $B_{N+1,q_e} \cap \Vect(X_N^1, \ldots, X_N^d) = \{0\}$, 
 et donc d'après le lemme~\ref{8lem_dimBN_zbaseBN}, $$\dim(B_{N+1,q_e} \oplus \Vect(U_1, \ldots, U_r)) = dq_e +r. $$
On note $G_r =B_{N+1,q_e} \oplus \Vect(U_1, \ldots, U_r) $ ainsi que $D_r = G_r^\perp \cap C$ l'orthogonal de $G_r$ dans $C$.
\\On a $\dim(D_r) = e - dq_e -r = r_e - r \geq 1$. 
\\Soit $\pi_r : C \to D_r $ la projection orthogonale sur $D_r$.
\\ On pose $\Delta_r = \pi_r(C \cap \Z^n)$. Alors $\Delta_r$ est un réseau euclidien de $D_r$ de déterminant :
\begin{align*}
 d(\Delta_r) = \frac{H(C)}{H(G_r)} ;
\end{align*}
on peut retrouver ce résultat dans la preuve du  théorème~2 de \cite{Schmidt}. \\ 
D'après le théorème~de Minkowski (théorème~\ref{2th_Minko}), il existe $X_r' \in \Delta_r \smallsetminus \{0 \} \subset D_r \cap \Q^n$ tel que :
\begin{align}\label{8normeXr}
 \| X_r' \| \leq c_{\ref{8cons_minko_1}} d(\Delta_r)^{\frac{1}{\dim(D_r)}} \leq c_{\ref{8cons_minko_1}} \left(\frac{H(C)}{H(G_r)}\right)^{\frac{1}{r_e-r}}
\end{align}
avec $\cons \label{8cons_minko_1}$ une constante ne dépendant que de $e$.
\\Comme $X_r' \in \Delta_r$, il existe $X_r \in C \cap \Z^n$ tel que $\pi_r(X_r) = X_r'$ ; on a $X_r \notin G_r$ donc $$X_r \notin B_{N+1,q_e}.$$
On pose $E_r$ le produit extérieur des vecteurs $(X_{N+k}^j)_{j \in \llbracket 1,d\rrbracket,k \in \llbracket 0, q_e \rrbracket} $ et de $X_r$. 
\\ On cherche à montrer que $E_r= 0$ et pour cela on étudie sa norme. 
On a :
\begin{align*}
 \| E_r \| &= \|X_r \wedge \bigwedge\limits_{k=0}^{q_e} (X_{N+k}^1 \wedge \ldots \wedge X_{N+k}^d) \| = \|\pi_r(X_r) \wedge \bigwedge\limits_{k=0}^{q_e} (X_{N+k}^1 \wedge \ldots \wedge X_{N+k}^d) \|
\end{align*}
car $X_r - \pi_r(X_r) \in G_r \subset B_{N,q_e+1}$.
Or $\| \pi_r(X_r) \| = \|X_r' \| \leq c_{\ref{8cons_minko_1}} \left(\frac{H(C)}{H(G_r)}\right)^{\frac{1}{r_e-r }}$ d'après (\ref{8normeXr}) et donc : 
\begin{align*}
 \| E_r\| &\leq \|\pi_r(X_r) \| \cdot \| \bigwedge\limits_{k=0}^{q_e} (X_{N+k}^1 \wedge \ldots \wedge X_{N+k}^d) \| \\
 &\leq c_{\ref{8cons_minko_1}} \left(\frac{H(C)}{H(G_r)}\right)^{\frac{1}{r-r_e}} \| (X_{N}^1 \wedge \ldots \wedge X_{N}^d) \wedge \bigwedge\limits_{k=1}^{q_e} (U_{N+k}^1 \wedge \ldots \wedge U_{N+k}^d) \|
\end{align*}
en utilisant la formule (\ref{8rec_Xn_Un}) sur les $X_N^j$ avec $U_{N+k}^j = \begin{pmatrix}
 0 &
 \cdots &
 0 &
 u_{N+k}^{(0,j)}& 
 \cdots & u_{N+k}^{(qd-1,j)}
 \end{pmatrix}^\intercal$.
\\ D'après la construction des $u_k^{(i,j)}$ en (\ref{8construc_suite_u}) on a $\| U_k^j \| \leq 3 $ pour tous $k $ et $j$. Cela donne donc :
\begin{align*}
 \| E_r \| &\leq c_{\ref{8cons_minko_1}} \left(\frac{H(C)}{H(G_r)}\right)^{\frac{1}{r_e-r}} \| X_{N}^1 \wedge \ldots \wedge X_{N}^d \| \prod\limits_{k=1}^{q_e} (\|U_{N+k}^1 \| \cdots \| U_{N+k}^d \|) \\&\leq 3^{dq}c_{\ref{8cons_minko_1}} \left(\frac{H(C)}{H(G_r)}\right)^{\frac{1}{r_e-r}} \| X_{N}^1 \wedge \ldots \wedge X_{N}^d \|.
\end{align*}
Or $\Vect(X_{N}^1, \ldots, X_{N}^d) = B_{N,1}$ d'après la construction en (\ref{8def_Bnv}). 
\\ D'après le lemme~\ref{8lem_dimBN_zbaseBN} $X_{N}^1, \ldots, X_{N}^d$ forment une $\Z-$base de $B_{N,1} \cap \Z^n$ et donc $H(B_{N,1}) = \| X_{N}^1 \wedge \ldots \wedge X_{N}^d \|$. 
\\ Or la remarque~\ref{8HauteurBN} indique que $H(B_{N,1}) \leq c_{\ref{8cons_haut_CN_major}} \theta^{d{\alpha^N}}$ et donc $\| X_{N}^1 \wedge \ldots \wedge X_{N}^d \| \leq c_{\ref{8cons_haut_CN_major}} \theta^{d{\alpha^N}}$.
\\ D'autre part, le lemme~\ref{8minoration_somme_espace} donne $H(G_r) = H(B_{N+1,q_e} \oplus \Vect(U_1, \ldots, U_r)) \geq
c_{\ref{8cons_minoration_somme_espace}} \theta^{ r {\alpha^N} +(d-r) {\alpha^{N+1}}} $. 
\bigskip
\\ On a donc :
 \begin{align*}
 \| E_r \| &\leq 3^{dq}c_{\ref{8cons_minko_1}}c_{\ref{8cons_minoration_somme_espace}}^{\frac{-1}{r_e-r}} c_{\ref{8cons_haut_CN_major}}\theta^{ d{\alpha^N} - \frac{r {\alpha^N} +(d-r) {\alpha^{N+1}}}{r_e-r} } H(C)^{\frac{1}{r_e-r} } \\
 &= 3^{dq}c_{\ref{8cons_minko_1}}c_{\ref{8cons_minoration_somme_espace}}^{\frac{-1}{r_e-r}} c_{\ref{8cons_haut_CN_major}}\theta^{ \frac{\alpha^N (d(r_e-r) - r - (d-r) {\alpha} )} {r_e-r} } H(C)^{\frac{1}{r_e-r} }.
 \end{align*}
 Or $ d(r_e-r) - r - (d-r) {\alpha} \leq dr_e -(d-r_e) \alpha \leq 0 $ d'après l'inégalité $(\ref{8lem_inegal_alpha4})$ du lemme~\ref{8lem_inegal_alpha}, et le choix de $N$ en (\ref{8choixN}) donne $$H(C)^{   \frac{{\alpha^{q_e }} }{\alpha^{N+q_e +1 }(d -r_e + \frac{1}{2})} } \leq H(C)^{  \frac{ \frac{\alpha^{q_e +1 }}{r_e +(d-r_e)\alpha } + \frac{\varepsilon}{2} -1 }{\alpha^{N+q_e +1 }} } \leq \theta$$
 car $\frac{\alpha^{q_e }}{d -r_e + \frac{1}{2} } \leq \frac{\alpha^{q_e +1 }}{r_e +(d-r_e)\alpha } -1 + \frac{\varepsilon}{2} $ d'après l'inégalité $(\ref{8lem_inegal_alpha5})$ de ce même lemme~\ref{8lem_inegal_alpha}. 
 L'inégalité devient alors :
 \begin{align*}
 \| E_r \| &\leq 3^{dq}c_{\ref{8cons_minko_1}}c_{\ref{8cons_minoration_somme_espace}}^{\frac{-1}{r_e-r}} c_{\ref{8cons_haut_CN_major}} H(C)^{ \frac{\alpha^N\alpha^{q_e } (d(r_e-r) - r - (d-r) {\alpha} )}{\alpha^{N+q_e +1} (r_e-r)(d -r_e + \frac{1}{2})} + \frac{1}{r_e-r} } \\
 &= 3^{dq}c_{\ref{8cons_minko_1}}c_{\ref{8cons_minoration_somme_espace}}^{\frac{-1}{r_e-r}} c_{\ref{8cons_haut_CN_major}} H(C)^{ \frac{ (d(r_e-r) - r - (d-r) {\alpha} ) + \alpha(d -r_e + \frac{1}{2}) }{\alpha(r_e-r)(d -r_e + \frac{1}{2})} } .
 \end{align*}
 On étudie l'exposant que l'on note $\delta =\frac{ (d(r_e-r) - r - (d-r) {\alpha} ) + \alpha(d -r_e + \frac{1}{2}) }{\alpha(r_e-r)(d -r_e + \frac{1}{2})} $ et on peut alors majorer 
 \begin{align*}
 \delta &= \frac{ 1}{\alpha(r_e-r)(d -r_e + \frac{1}{2})} (- \alpha(r_e-r - \frac{1}{2}) + d(r_e-r) - r ) \\
 &\leq \frac{ 1}{\alpha(r_e-r)(d -r_e + \frac{1}{2})} (- \frac{\alpha}{2} + d(d-1) )
 \end{align*}
 car $0 \leq r \leq r_e -1$ et $1 \leq r_e -r \leq d-1 $. Enfin d'après l'inégalité $(\ref{8lem_inegal_alpha2})$ du lemme~\ref{8lem_inegal_alpha}, on a $- \frac{\alpha}{2} + d(d-1) \leq -1$ et donc 
 \begin{align*}
 \delta &\leq \frac{ -1}{\alpha(r_e-r)(d -r_e + \frac{1}{2})} 
 \end{align*}
 et enfin $\| E_r \| \leq 3^{dq}c_{\ref{8cons_minko_1}}c_{\ref{8cons_minoration_somme_espace}}^{\frac{-1}{r_e-r}} c_{\ref{8cons_haut_CN_major}} H(C)^{ \frac{ -1}{\alpha(r_e-r)(d -r_e + \frac{1}{2}) }}$. En particulier si $H(C)$ est assez grand on a $\| E_r \| < 1$. 
 \\
Or les vecteurs considérés dans le produit extérieur $E_r = X_r \wedge \bigwedge\limits_{k=0}^{q_e} (X_{N+k}^1 \wedge \ldots \wedge X_{N+k}^d) $ sont entiers donc ce produit extérieur s'annule.
 \\Il existe alors $U_{r+1} \in (B_{N+1,q_e} \oplus \Vect(X_r)) \cap \Vect(X_N^1, \ldots, X_N^d) \smallsetminus \{ 0\} $. On rappelle que $X_r \notin B_{N+1,q_e}$ et comme les espaces considérés sont rationnels, on peut prendre $U_{r+1} \in \Z^n$. Comme $(B_{N+1,q_e} \oplus \Vect(X_r)) \subset C$ on a aussi $U_{r+1} \in C$.
\\ Montrons maintenant que les vecteurs $U_1, \ldots, U_{r+1}$ sont linéairement indépendants sur $\R$. Soit $ \sum\limits_{k =1}^{r+1} \lambda_k U_k =0 $ une relation de dépendance linéaire.
 \\On applique $\pi_r$ la projection orthogonale sur $D_r = (B_{N+1,q_e} \oplus \Vect(U_1, \ldots, U_r))^\perp \cap C $ et on trouve :
 \begin{align}\label{8indep}
 0= \pi_r(\sum\limits_{k =1}^{r+1} \lambda_k U_k) &= \sum\limits_{k =1}^{r+1} \lambda_k \pi_r(U_k) = \lambda_{r+1} \pi_r(U_{r+1}).
 \end{align}
 Or $U_{r+1} = Z + \mu X_r$ avec $Z \in B_{N+1,q_e} $ et $\mu \in \R$. 
 \\ De plus comme $U_{r+1} \in \Vect(X_N^1, \ldots, X_N^d) $ et que $\Vect(X_N^1, \ldots, X_N^d) \cap B_{N+1,q_e} = \{0\}$, on a nécessairement $\mu \neq 0$. \\
 On a donc $\pi_r(U_{r+1}) = \mu X_r' \neq 0$ ce qui donne $\lambda_{r+1} = 0 $ en utilisant (\ref{8indep}).
\\ Enfin on trouve $\lambda_1 = \ldots = \lambda_r = 0$ en utilisant l'hypothèse de récurrence donnant $U_1, \ldots, U_r$ linéairement indépendants. Ceci termine la récurrence. 

\bigskip
On a donc montré qu'il existe $W \subset \Vect(X_N^1, \ldots, X_N^d) \cap C$ un sous-espace rationnel de dimension $r_e$. 
\\ Comme $B_{N+1,q_e} \subset C$ et $\Vect(X_N^1, \ldots, X_N^d) \cap B_{N+1,q_e} = \{0\}$ on a donc $C = B_{N+1,q_e} \oplus W $ par égalité des dimensions. 
\\ Le lemme~\ref{8minoration_somme_espace} donne alors 
\begin{align}\label{8minor_hautC_rappel}
 H(C) \geq c_{\ref{8cons_minoration_somme_espace}} \theta^{ r_e {\alpha^N} +(d-r_e) {\alpha^{N+1}}}.
\end{align}
On remarque d'autre part que $C = B_{N+1,q_e} \oplus W \subset C_{N-1, (q_e+1)d}$ et d'après le lemme~\ref{8lem_minoration_Y1angleCNE} appliqué avec $N' = N-1$ et $e' = q_{e'}d$ où $q_{e'} = q_e +1$, on a 
\begin{align}\label{8minor_rappel_y1_C}
 \psi_1(\Vect(Y_1), C_{N-1,(q_e+1)d}) \geq c_{\ref{8cons_minor_angle_y1_CNe}} \theta^{- \alpha^{N'+ q_{e'}}} =
 c_{\ref{8cons_minor_angle_y1_CNe}} \theta^{- \alpha^{N+q_e+1}}
\end{align}
et donc, en utilisant le lemme~\ref{lem_inclusion_croissance} puisque $Y_1 \in A \smallsetminus \{0\} $ et $\dim(A) =d $ : 
\begin{align}\label{8minor_angle_C}
 \psi_d(A,C) \geq \psi_d(A,C_{N-1,(q_e+1)d}) \geq \psi_1(\Vect(Y_1),C_{N-1,(q_e+1)d}).
\end{align}
Les inégalités $(\ref{8minor_hautC_rappel})$, $(\ref{8minor_rappel_y1_C})$ et $(\ref{8minor_angle_C})$ donnent alors l'existence d'une constante $\cons \label{8cons_preuve_infinit_derniere_cons} >0 $ indépendante de $C$ telle que :
\begin{align*}
 \psi_d(A,C) \geq c_{\ref{8cons_preuve_infinit_derniere_cons}} H(C) ^{ \frac{-\alpha^{N+q_e+1}}{r_e {\alpha^N} +(d-r_e) {\alpha^{N+1}}}} = c_{\ref{8cons_preuve_infinit_derniere_cons}} H(C) ^{ \frac{-\alpha^{q_e+1}}{r_e +(d-r_e) {\alpha}}}.
\end{align*}
En reprenant l'hypothèse faite en $(\ref{8hyp_absurde_preuve_lem14})$ on a donc 
\begin{align*}
 c_{\ref{8cons_preuve_infinit_derniere_cons}} H(C) ^{ \frac{-\alpha^{q_e+1}}{r_e +(d-r_e) {\alpha}}} \leq H(C) ^{ \frac{-\alpha^{q_e+1}}{r_e +(d-r_e) {\alpha}} -  \varepsilon}.
\end{align*}
En particulier pour tout $C$ avec $H(C)$ assez grand on a $ c_{\ref{8cons_preuve_infinit_derniere_cons}} \leq H(C)^{-\varepsilon }$ et donc $$ c_{\ref{8cons_preuve_infinit_derniere_cons}} = 0$$
ce qui est contradictoire et termine la preuve du lemme.

\end{preuve}

Ce lemme~\ref{8lem_minoration_exposant} donne alors :
\begin{align*}
 \mu_n(A|C)_d \leq \frac{\alpha^{q_e +1 }}{r_e +(d-r_e){\alpha}}.
\end{align*}
Cela termine la preuve du théorème~\ref{8theo_dernier_angle} dans le cas $d \leq e$.

\section{Calcul de l'exposant dans le cas \texorpdfstring{$e < d$}{}}
Dans cette section, on montre $\mu_n(A|e)_e = \frac{\alpha}{e}$ pour $e \in \llbracket 1,d -1 \rrbracket $. On peut remarquer que dans le cas où $e = d $ cette égalité est aussi vérifiée ; en effet on a prouvé le théorème~\ref{8theo_dernier_angle} dans ce cas. Ce cas apparaissant aussi dans les preuves de section, on montrera ici de nouveau $\mu_n(A|d)_d = \frac{\alpha}{d}$. 
\\ Dans la suite on fixe $e \in \llbracket 1, d \rrbracket$.

\subsection{Minoration de l'exposant}
La minoration de l'exposant $\mu_n(A|e)_e$ reprend les mêmes idées que dans le cas $e \geq d$. \\
On introduit ici une suite d'espaces approchant très bien $A$.
\\ Pour $N \in \N$, on pose :
\begin{align}\label{def_DN}
 D_{N,e} = \Vect(X_N^1, \ldots, X_N^e)
\end{align}
qui est un espace rationnel de dimension $e$ car $X_N^i \in\Z^n$ pour tout $i$.

\begin{lem}\label{8lem_hauteur_DN} Les vecteurs $X_N^1, \ldots, X_N^{e}$ forment une $\Zbase$ de $D_{N,e} \cap \Z^n$ et 
\begin{align*}
 c_{\ref{8cons_minor_haut_DN}} \theta^{e {\alpha^N}} \leq H(D_{N,e}) \leq c_{\ref{8cons_major_haut_DN}} \theta^{e {\alpha^N}}
\end{align*}
\end{lem}
avec $\cons \label{8cons_minor_haut_DN} $ et $\cons \label{8cons_major_haut_DN}$ indépendantes de $N$.

\begin{preuve}
Par le lemme~\ref{8lem_dimBN_zbaseBN}, les vecteurs $X_N^1, \ldots, X_N^e $ sont des vecteurs d'une $\Zbase$ de $B_{N,v} \cap \Z^n$ introduit en (\ref{8def_Bnv}).
\\D'après la remarque~\ref{2req_sous_Zbase}, ils forment alors en particulier une $\Zbase$ de l'espace $W \cap \Z^n$ où $W$ est l'espace qu'ils engendrent, c'est-à-dire $D_{N,e}$.
\\ Par définition de la hauteur on a donc $H(D_{N,e}) = \| X_N^1 \wedge \ldots \wedge X_N^e\| $.
\\ On reprend la notation suivante pour $i \in \llbracket 1,e \rrbracket$ :
\begin{align}\label{vecteur_Zn}
 Z_{N}^j = \frac{1}{\theta^{\floor{\alpha^N}}} X_N^j 
\end{align}
et on a $Z_N^j \underset{N\to+\infty}{\longrightarrow} Y_j$.

En particulier on a 
\begin{align*}
 \| Z_N^1 \wedge \ldots \wedge Z_N^e\| \underset{N\to+\infty} {\longrightarrow} \| Y_1 \wedge \ldots \wedge Y_e \| \neq 0.
\end{align*}
On en déduit 
\begin{align*}
 \theta^{-e\floor{\alpha^N}} \| X_N^1 \wedge \ldots \wedge X_N^e\| = \| Z_N^1 \wedge \ldots \wedge Z_N^e\| \underset{N\to+\infty} {\longrightarrow} \| Y_1 \wedge \ldots \wedge Y_e \| \neq 0.
\end{align*}
Donc il existe $\cons \label{8cons_preuve_haut_DN1}$ et $\cons \label{8cons_preuve_haut_DN2}$ telles que pour tout $N \in \N $,
\begin{align*}
 c_{\ref{8cons_preuve_haut_DN1}} \theta^{e \floor{\alpha^N}} \leq H(D_{N,e}) \leq c_{\ref{8cons_preuve_haut_DN2}} \theta^{e \floor{\alpha^N}}
\end{align*}
ce qui permet de conclure puisque $\theta^{e {\alpha^N} -e } \leq \theta^{e \floor{\alpha^N}} \leq \theta^{e \alpha^N}$.

\end{preuve}

\begin{lem}\label{8lem_Angle_psie_DN}
Il existe une constante $\cons \label{8cons_major_psie_A_DNe}$ indépendante de $N$ telle que :
 \begin{align*}
 \psi_e(A,D_{N,e}) \leq c_{\ref{8cons_major_psie_A_DNe}} H(D_{N,e})^{{-\alpha}/{e}}.
 \end{align*}
 pour $N $ assez grand.
\end{lem}

\begin{preuve}
 On rappelle que pour $j \in \llbracket 1,d \rrbracket$ et $N \in \N$ 
 \begin{align*}
 \psi_1(\Vect(Y_j), \Vect(X_{N}^j)) \leq c_{\ref{8cons_major_angle_Yj_XNj}} \theta^{-\alpha^{N+1}}
 \end{align*}
 d'après $(\ref{8major_Yj_XNj})$.
\\ On a $\Vect(Y_1,\ldots, Y_e) \subset A $ donc $ \psi_e(A,D_{N,e}) \leq \psi_e(\Vect(Y_1,\ldots, Y_e), D_{N,e})$.
\\D'après la propriété~\ref{2prop_4.5Elio} on a :
 \begin{align*}
 \psi_e(\Vect(Y_1,\ldots, Y_e), D_{N,e}) &\leq c_{\ref{8cons_maj_prop45}} \sum\limits_{j=1}^{e} \psi_1(\Vect(Y_j), \Vect(X_{N}^j) ) 
 \end{align*}
 avec $\cons \label{8cons_maj_prop45} $ dépendant de $Y_1, \ldots, Y_e $ et $n$.
On a donc 
\begin{align*}
 \psi_e(A, D_{N,e}) &\leq c_{\ref{8cons_maj_prop45}} c_{\ref{8cons_major_angle_Yj_XNj}} e \theta^{-{\alpha^{N +1 }} }\\
 &\leq c_{\ref{8cons_maj_prop45}} c_{\ref{8cons_major_angle_Yj_XNj}}d c_{\ref{8cons_major_haut_DN}}^{\frac{{\alpha^{N+1 }}}{e{\alpha^N} }} H(D_{N,e})^{\frac{-{\alpha^{N +1 }}}{e{\alpha^N} }} 
 \\
 &= c_{\ref{8cons_major_psie_A_DNe}} H(D_{N,e})^{{-{\alpha}}/{e }} 
 \end{align*}
 avec $c_{\ref{8cons_major_psie_A_DNe}} = c_{\ref{8cons_maj_prop45}} c_{\ref{8cons_major_angle_Yj_XNj}}d c_{\ref{8cons_major_haut_DN}}^{\frac{{\alpha}}{e}} $ et car $H(D_{N,e}) \leq c_{\ref{8cons_major_haut_DN}} \theta^{e {\alpha^N}} $ d'après le lemme~\ref{8lem_hauteur_DN}.

\end{preuve}

On a donc construit une infinité d'espaces rationnels $D_{N,e}$ de dimension $e$ tels que 
\begin{align*}
 \psi_e(A,D_{N,e}) \leq c_{\ref{8cons_major_psie_A_DNe}} H(D_{N,e})^{{-\alpha}/{e}} ;
\end{align*}
cela donne en particulier : 
\begin{align*}
 \mu_n(A|e)_e \geq \frac{{\alpha}}{e }.
\end{align*}

\subsection{Majoration de l'exposant}
Le but de cette section est de majorer $\mu_n(A|e)_e$. 

\begin{lem}\label{8lem_intersection_DN_C}
 Soit $ \varepsilon >0$ et $C$ un espace rationnel de dimension $e$ tel :
 \begin{align*}
 \psi_e(A,C) \leq H(C) ^{-\frac{\alpha}{e} - \varepsilon }.
 \end{align*}
 Alors si $H(C)$ est assez grand, il existe $N \in \N$ et $Z \in \Z^n \smallsetminus \{0 \} $ tels que $Z \in C \cap D_{N,d}$ et 
 \begin{align*}
 \| Z \| \leq c_{\ref{8cons_enonce_lem_minko}}H(C)^{{1}/{e}}
 \end{align*}
 avec $\cons \label{8cons_enonce_lem_minko}$ indépendante de $Z$ et de $N$.
\end{lem}

\begin{preuve}
 Par le corollaire~\ref{cor_Mink2} du  théorème~de Minkowski il existe $Z \in C \cap \Z^n \smallsetminus \{0 \}$ tel que :
\begin{align}\label{8vecteur_Z}
 \| Z \| \leq c_{\ref{8cons_minko2}} H(C)^{{1}/{e}}
\end{align}
avec $\cons \label{8cons_minko2}$ une constante indépendante de $Z$.
\\ Il reste alors à montrer qu'il existe un certain $N$ tel que $Z \in D_{N,d} = \Vect(X_N^1, \ldots, X_N^d)$. 
\\ On pose $Z^A$ le projeté orthogonal de $Z$ sur $A$. On introduit $a_1, \ldots, a_d \in \R$ tels que :
\begin{align*}
 Z^A =  \sum\limits_{j = 1}^d a_j Y_j.
\end{align*}
On cherche à montrer qu'il existe $N$ tel que $\| Z \wedge X_N^1 \wedge \ldots \wedge X_N^d \|$ s'annule. Or $X_N^1, \ldots, X_N^d$ est une $\Zbase$ de $D_{N,d} \cap \Z^n$ d'après le lemme~\ref{8lem_hauteur_DN}. Pour $N \in \N $, on a donc d'après le lemme~\ref{lem_phi_dim1} 
\begin{align}
 \| Z \wedge X_N^1 \wedge \ldots \wedge X_N^d \| &= \psi_1(\Vect(Z), D_{N,d}) H(D_{N,d}) \|Z \| \nonumber\\
 & \leq \omega(Z, \sum\limits_{j = 1}^d a_j X_N^j) H(D_{N,d}) \|Z \| \nonumber\\
 & \leq \left(\omega(Z, Z^A) + \omega(Z^A, \sum\limits_{j = 1}^d a_j X_N^j) \right)H(D_{N,d}) \|Z \| \label{8angle3} .
\end{align}
En rappelant la notation $Z_N^j = \theta^{-\floor{\alpha^N}}X_N^j $, on a \begin{align*}
 \omega(Z^A, \sum\limits_{j = 1}^d a_j X_N^j) = \omega(Z^A, \sum\limits_{j = 1}^d a_j Z_N^j) \leq \frac{\left\|Z^A - \sum\limits_{j = 1}^d a_j Z_N^j \right\| }{ \| Z^A \| } 
 \leq \frac{ \sum\limits_{j = 1}^d |a_j| \left\|Y_j -Z_N^j \right\| }{ \left\| \sum\limits_{j = 1}^d a_j Y_j \right\| }.
\end{align*}
Or par construction des $Y_j$, on a $\left\| \sum\limits_{j = 1}^d a_j Y_j \right\| \geq \sqrt{ \sum\limits_{j=1}^d |a_j|^2} \geq \sum\limits_{j = 1}^d \frac{|a_j|}{d } $ ; en effet pour $i \in \llbracket 1,d \rrbracket$, la $i\tir$ème coordonnée de $Y_j$ est égale à $1$ si $i =j$ et à $0$ sinon. \\De plus, on a $\ \|Y_j-Z_N^j \| \leq \|Y_j \| c_{\ref{8cons_major_angle_Yj_XNj}}\theta^{-{\alpha^{N+1}}} $ par $(\ref{8major_Yj_XNj})$. On a donc 
\begin{align}\label{8angle_ZA_DN}
 \omega(Z^A, \sum\limits_{j=1}^d a_j X_N^j) \leq c_{\ref{8cons_angle_ZA_DN}}{\theta^{-{\alpha^{N+1}}}}
\end{align}
avec $\cons \label{8cons_angle_ZA_DN} = {c_{\ref{8cons_major_angle_Yj_XNj}}}{d} \max\limits_{j \in \llbracket 1,d \rrbracket} \|Y_j \| $. \\
D'autre part on a, en utilisant le lemme~\ref{lem_inclusion_croissance} puisque $Y_1 \in A \smallsetminus \{0\} $ et $\dim(C) =e $ :
\begin{align*}
 \omega(Z, Z^A) = \psi_1(\Vect(Z),A) \leq \psi_e(C,A).
\end{align*}
L'hypothèse du lemme donne alors: 
\begin{align}\label{8angle_ZA_C}
 \omega(Z, Z^A) \leq H(C) ^{ -\frac{\alpha}{e } - \varepsilon }.
\end{align}
On combine (\ref{8angle_ZA_DN}) et (\ref{8angle_ZA_C}) avec (\ref{8angle3}) et on a :
\begin{align*}
 \| Z \wedge X_N^1 \wedge \ldots \wedge X_N^d \| \leq \left(H(C) ^{-\frac{{\alpha}}{e } - \varepsilon } + c_{\ref{8cons_angle_ZA_DN}}{\theta^{-{\alpha^{N+1}}}}\right)H(D_{N,d}) \|Z \|.
\end{align*}
On sait de plus par le lemme~\ref{8lem_hauteur_DN} que $H(D_{N,d}) \leq c_{\ref{8cons_major_haut_DN}}\theta^{d{\alpha^N}}$ et par (\ref{8vecteur_Z}) que $\|Z\| \leq c_{\ref{8cons_minko2}}H(C)^{{1}/{e}}$. On a donc pour $N \in \N $ :
\begin{align*}
 \| Z \wedge X_N^1 \wedge \ldots \wedge X_N^d \| &\leq \left(H(C) ^{-\frac{{\alpha}}{e } - \varepsilon } + c_{\ref{8cons_angle_ZA_DN}}{\theta^{-{\alpha^{N+1}}}} \right) c_{\ref{8cons_major_haut_DN}}\theta^{d{\alpha^N}}c_{\ref{8cons_minko2}}H(C)^{{1}/{e}} \\
 &\leq c_{\ref{8cons_major_final_ext_preuve}} \left(H(C) ^{-\frac{{\alpha-1}}{e } - \varepsilon }\theta^{d{\alpha^N}} + \theta^{\alpha^N(d-\alpha) }H(C)^{{1}/{e}}\right)
\end{align*}
avec $\cons \label{8cons_major_final_ext_preuve} = c_{\ref{8cons_major_haut_DN}}c_{\ref{8cons_minko2}}\max(1,c_{\ref{8cons_angle_ZA_DN}})$ indépendante de $C$ et de $N$.
\\On choisit maintenant $N$. Soit $N$ l'entier tel que : 
\begin{align}\label{8defN}
 \theta^{d{\alpha^N}} \leq H(C) ^{\frac{{\alpha -1 }}{e } +\frac{\varepsilon}{2} } < \theta^{d{\alpha^{N+1}}}. 
\end{align}
Cela donne donc $\theta^{d{\alpha^N}} \leq H(C) ^{\frac{{\alpha -1 }}{e } +\frac{\varepsilon}{2} }$ et $\theta > H(C) ^{\frac{\frac{{\alpha -1 }}{e } +\frac{\varepsilon}{2} }{d{\alpha^{N+1}}} } $, d'où, puisque $\alpha > d$ :
\begin{align}
 \| Z \wedge X_N^1 \wedge \ldots \wedge X_N^d \| &\leq c_{\ref{8cons_major_final_ext_preuve}} \left(H(C)^{ -\frac{\varepsilon}{2}} + H(C)^{ \alpha^N(d-\alpha) \left(\frac{\frac{{\alpha -1 }}{e } +\frac{\varepsilon}{2} }{d{\alpha^{N+1}}}\right) + \frac{1}{e}} \right) \nonumber \\
 &\leq c_{\ref{8cons_major_final_ext_preuve}} \left(H(C)^{ -\frac{\varepsilon}{2}} + H(C)^{ \frac{(d-\alpha)(2\alpha - 2 + e\varepsilon) + 2d \alpha}{2ed{\alpha}} } \right) \label{8prod_exterieur_Z_2HC}.
\end{align}
On étudie alors l'exposant du deuxième terme :
\begin{align*}
 \frac{(d-\alpha)(2\alpha - 2 + e\varepsilon) + 2d \alpha}{2ed{\alpha}} &= \frac{1}{ed\alpha}\left(-\alpha^2 + (1 + 2d)\alpha - d \right) - \frac{\alpha -d }{2d\alpha} \varepsilon.
\end{align*}
Or d'après l'inégalité $(\ref{8lem_inegal_alpha3})$ du lemme~\ref{8lem_inegal_alpha}, on a $-\alpha^2 + (1 + 2d)\alpha - d \leq 0$ et donc 
\begin{align*}
 \frac{(d-\alpha)(2\alpha - 2 + e\varepsilon) + 2d \alpha}{2ed{\alpha}} \leq - c_{\ref{8cons_major_epsilon}}\varepsilon
\end{align*} 
avec $\cons \label{8cons_major_epsilon} = \frac{\alpha -d }{2d\alpha} > 0 $ car $\alpha > d$.
 \\ L'inégalité (\ref{8prod_exterieur_Z_2HC}) devient donc :
 \begin{align*}
 \| Z \wedge X_N^1 \wedge \ldots \wedge X_N^d \| &\leq c_{\ref{8cons_major_final_ext_preuve}} \left(H(C)^{ -\frac{\varepsilon}{2}} + H(C)^{- c_{\ref{8cons_major_epsilon}}\varepsilon} \right) 
 \end{align*}
avec $c_{\ref{8cons_major_final_ext_preuve}}$ et $c_{\ref{8cons_major_epsilon}}$ indépendante de $C$. Pour $H(C)$ assez grand on a donc $$\| Z \wedge X_N^1 \wedge \ldots \wedge X_N^d \| <1 $$ et le lemme~\ref{2lem_X_in_B} donne alors 
\begin{align*}
 Z \in \Vect(X_N^1, \ldots, X_N^d) = D_{N,d}
\end{align*}
pour $N$ vérifiant $(\ref{8defN})$, ce qui termine la preuve.

\end{preuve}

On peut maintenant montrer le résultat permettant de minorer l'exposant.
\begin{lem}\label{8lem_minoration_exposant_e<d}
 Soit $ \varepsilon >0$. Pour tous les espaces rationnels $C$ de dimension $e$ sauf un nombre fini on a :
 \begin{align*}
 \psi_e(A,C) > H(C) ^{-\frac{\alpha}{e} - \varepsilon }.
 \end{align*}
 
\end{lem}

\begin{preuve}
 On suppose le contraire par l'absurde. 
\\ Il existe alors une infinité d'espaces rationnels $C$ de dimension $e$ vérifiant : 
\begin{align}\label{8hyp_absurd}
 \psi_e(A,C) \leq H(C) ^{-\frac{\alpha}{e } - \varepsilon }.
\end{align}
Soit $C$ un tel espace. Le lemme $\ref{8lem_intersection_DN_C}$ donne alors $N \in \N$ et $Z \in \Z^n \smallsetminus \{0 \} $ tels que $Z \in C \cap D_{N,d}$ avec
 \begin{align*}
 \| Z \| \leq c_{\ref{8cons_enonce_lem_minko}}H(C)^{\frac{1}{e}}.
 \end{align*}
 \\Comme $Z \in D_{N,d} \cap \Z^n$ on peut écrire :
 \begin{align}\label{8ecriture_Z}
 Z =  \sum\limits_{j = 1}^d v_j X_N^j
 \end{align}
 avec $v_j \in \Z$ car les $X_N^j$ forment une $\Z-$base de $D_{N,d} \cap \Z^n$ d'après le lemme~\ref{8lem_hauteur_DN}.
 \\Soit, pour $ j\in \llbracket1,d\rrbracket$, $z_j = \theta^{\floor{\alpha^N}}v_j $, ce qui donne $Z =  \sum\limits_{j = 1}^d z_j Z_N^j $. 
 \\On souhaite maintenant minorer $\psi_1(Z,A)$. On rappelle la notation $Z^A$ pour le projeté orthogonal de $Z$ sur $A$. On introduit $a_1, \ldots, a_d \in \R$ tels que :
\begin{align*}
 Z^A =  \sum\limits_{j = 1}^d a_j Y_j.
\end{align*} On pose de plus : 
\begin{align*}
 \Delta = Z^A -  \sum\limits_{i = 1}^d z_j Y_j, 
\end{align*}
ainsi que 
\begin{align*}
 \omega = \| Z^A - Z \|.
\end{align*}
Comme les $Z_N^j$ et $Y_j$ ont leur $i\tir$ème coordonnée de $Y_j$ égale à $1$ si $i =j$ et à $0$ sinon, cela donne 
\begin{align*}
 Z^A-Z = \begin{pmatrix} a_1 -z_1 & \cdots & a_d -z_d& \star &\ldots & \star \end{pmatrix}^\intercal.
\end{align*}
\\ On a alors :
\begin{align*}
 \omega^2 \geq  \sum\limits_{j = 1}^d (a_j -z_j)^2, 
\end{align*}
et donc pour tout $j \in \llbracket 1,d\rrbracket$, $ | a_j - z_j | \leq \omega$. 
\\ Cela nous donne en particulier 
\begin{align}\label{8major_delta}
 \| \Delta \| = \left\|  \sum\limits_{j = 1}^d (a_j -z_j) Y_j \right\| \leq c_{\ref{8cons_maj_Delta_par_omega}}\omega 
\end{align}
 avec $\cons \label{8cons_maj_Delta_par_omega} = d\max\limits_{j \in \llbracket 1,d \rrbracket} \|Y_j\|$.
On peut alors minorer $\|Z\wedge Z^A \| $ : 
\begin{align}
\|Z\wedge Z^A \| &= \|Z \wedge ( \sum\limits_{j = 1 }^d z_j Y_j + Z^A-  \sum\limits_{j = 1 }^d z_j Y_j ) \| \nonumber \\ 
 &= \| Z \wedge  \sum\limits_{j = 1 }^d z_j Y_j + Z \wedge \Delta \| \nonumber \\
 &\geq \| Z \wedge  \sum\limits_{j = 1 }^d z_j Y_j \| - \| Z \wedge \Delta \|. \label{8minoration_1_angleZ}
\end{align}
D'un côté, pour tous $j_0 \in \llbracket 1,d \rrbracket$ et $i \in \llbracket 0, qd-1 \rrbracket$, en se rappelant des définitions des $Y_j$ et des $\sigma_j $ données dans la section $\ref{8premier_paragraphe_construction}$ et en considérant les $j_0\tir$ème et $(d +i)\tir$ème lignes de la matrice $\begin{pmatrix}
 Z & \left| \:  \sum\limits_{j= 1 }^d z_j Y_j \right.
\end{pmatrix}$ on a:
\begin{align*}
 \| Z \wedge  \sum\limits_{j= 1 }^d z_j Y_j \| &= \|  \sum\limits_{j = 1 }^d z_j Z_N^j \wedge  \sum\limits_{j = 1 }^d z_j Y_j \| \\
 &\geq \left| \det \begin{pmatrix}
 z_{j_0} & z_{j_0} \\
  \sum\limits_{j= 1 }^d z_j \sigma_{i,j}^N &  \sum\limits_{j = 1 }^d z_j \sigma_{i,j} 
 \end{pmatrix} \right|\\
 &= |z_{j_0}| \left|  \sum\limits_{j = 1 }^d z_j (\sigma_{i,j} - \sigma_{i,j}^N) \right| \\
 &= |z_{j_0}| \left|  \sum\limits_{j = 1 }^d z_j  \sum\limits_{ k = N+1 }^{+ \infty} \frac{u_k^{(i,j)}}{\theta^{\floor{\alpha ^k}}} \right|. 
\end{align*}
On choisit $j_0$ tel que $|z_{j_0}| = \max\limits_{j=1}^d{|z_j|} \neq 0$ ; un tel $j_0$ existe car $Z \neq 0$. 
\\En particulier on a $|z_{j_0}| = \theta^{\floor{\alpha^N}} |v_{j_0}| \geq \theta^{\floor{\alpha^N}}$ car $v_{j_0} \in \Z$ d'après (\ref{8ecriture_Z}).
\\ On choisit $i$ tel que $u_{N+1}^{(i,j_0)} \neq 0 $, c'est-à-dire $i = \phi_{j_0}(N+1)$ en reprenant les notations de la section~\ref{8premier_paragraphe_construction} . On a alors $u_{N+1}^{(i,j_0)} \geq 2 $ et $u_{N+1}^{(i,j)} = 0 $ pour tout $j \neq j_0 $ donc : 
\begin{align*}
 \| Z \wedge  \sum\limits_{j= 1 }^d z_j Y_j \| &\geq |z_{j_0}| \left|  \sum\limits_{j = 1 }^d z_j  \sum\limits_{ k = N+1 }^{+ \infty} \frac{u_k^{(i,j)}}{\theta^{\floor{\alpha ^k}}} \right| \\ 
 &\geq |z_{j_0}| \left| z_{j_0}\frac{u_{N+1}^{(i,j_0)}}{\theta^{\floor{\alpha^{N+1}}}} +  \sum\limits_{j = 1 }^d z_j  \sum\limits_{ k = N+2 }^{+ \infty} \frac{u_k^{(i,j)}}{\theta^{\floor{\alpha ^k}}} \right| \\
 &\geq |z_{j_0}| \left(|z_{j_0}|\frac{2}{\theta^{\floor{\alpha^{N+1}}}} - d|z_{j_0}|  \sum\limits_{ k = N+2 }^{+ \infty} \frac{u_k^{(i,j)}}{\theta^{\floor{\alpha ^k}}} \right) \\
 &\geq |z_{j_0}|^2 \frac{1}{\theta^{\floor{\alpha^{N+1}}}} 
\end{align*}
car si $N$ est assez grand on a $d  \sum\limits_{ k = N+2 }^{+ \infty} \frac{u_k^{(i,j)}}{\theta^{\floor{\alpha ^k}}} \leq \frac{1}{\theta^{\floor{\alpha^{N+1}}}} $.
\\On a donc :
\begin{align}\label{8minorer_prod_ZYj}
 \| Z \wedge  \sum\limits_{j= 1 }^d z_j Y_j \| &\geq \frac{\theta^{2\floor{\alpha^{N}}}}{\theta^{\floor{\alpha^{N+1}}}}.
\end{align}
D'autre part, d'après (\ref{8major_delta}), on a :
\begin{align}\label{8majorer_prod_Zomega}
 \| Z \wedge \Delta \| \leq \|Z\| \|\Delta \| \leq c_{\ref{8cons_maj_Delta_par_omega}} \omega \|Z \|.
\end{align} 
En combinant (\ref{8minoration_1_angleZ}) avec (\ref{8minorer_prod_ZYj}) et (\ref{8majorer_prod_Zomega}) on trouve :
\begin{align*}
 \|Z \wedge Z^A \| &\geq \frac{\theta^{2\floor{\alpha^{N}}}}{\theta^{\floor{\alpha^{N+1}}}} - c_{\ref{8cons_maj_Delta_par_omega}} \omega \|Z\|.
\end{align*}
D'après le lemme~\ref{angle_Xortho_sur_X}, $\omega = \| Z^A - Z \| = \| Z \| \omega(Z^A,Z)$ et alors comme $\|Z^A \|\leq \|Z \|$ :
\begin{align*}
 \omega &= \| Z \| \frac{\| Z^A \wedge Z \|}{\| Z \| \cdot \|Z^A \|} \\ 
 &\geq \ \frac{\| Z^A \wedge Z \|}{\| Z \| } \\
 &\geq \frac{\theta^{2\floor{\alpha^{N}}}}{\|Z\|\theta^{\floor{\alpha^{N+1}}}} - c_{\ref{8cons_maj_Delta_par_omega}} \omega .
\end{align*}
On a finalement :
\begin{align*}
 \omega \geq \frac{c_{\ref{8cons_minor_omega}}}{\|Z \|\theta^{\floor{\alpha ^{N+1}} -2 \floor{\alpha^N}}}
\end{align*}
avec $\cons \label{8cons_minor_omega} = (1 +c_{\ref{8cons_maj_Delta_par_omega}})^{-1} $.\\
On a enfin 
\begin{align*}
 \omega(Z, Z^A)=\frac{\omega}{\|Z\|} \geq \frac{c_{\ref{8cons_minor_omega}}}{\|Z \|^2\theta^{\floor{\alpha ^{N+1}} - 2 \floor{\alpha^N}}},
\end{align*}
avec $c_{\ref{8cons_minor_omega}} > 0$ indépendante de $N$ et de $Z$. \\
On rappelle que $Z \in C$, et donc $ \omega(Z, Z^A) = \psi_1(\Vect(Z),A) \leq \psi_e(C,A)$ et donc par (\ref{8hyp_absurd}) :
\begin{align}\label{8majorationderniereZ}
 \frac{c_{\ref{8cons_minor_omega}}}{\|Z \|^2\theta^{\floor{\alpha ^{N+1}} - 2 \floor{\alpha^N}}} \leq H(C)^{-\frac{\alpha}{e} -   \varepsilon}.
\end{align}
De plus on a $Z$ tel que $\|Z\| \leq c_{\ref{8cons_enonce_lem_minko}}H(C)^{{1}/{e}}$. L'inégalité (\ref{8majorationderniereZ}) devient donc :
\begin{align*}
 \frac{c_{\ref{8cons_minor_omega}}}{c_{\ref{8cons_enonce_lem_minko}}^{\alpha+ e\varepsilon }\|Z \|^2\theta^{\floor{\alpha ^{N+1}} - 2 \floor{\alpha^N}}} \leq \|Z\|^{-{\alpha} - e\varepsilon}. 
\end{align*}
En particulier on a 
\begin{align}\label{8derniere_contradiction}
 \frac{c_{\ref{8cons_derniere_contradiction}}}{\|Z \|^2\theta^{\floor{\alpha ^{N+1}} - 2 \floor{\alpha^N}}} \leq \|Z\|^{-{\alpha} - e\varepsilon} 
\end{align}
avec $ \cons \label{8cons_derniere_contradiction} = c_{\ref{8cons_minor_omega}}c_{\ref{8cons_enonce_lem_minko}}^{-\alpha-e\varepsilon }$.
 D'un autre côté, par construction des $Z_N^j$ on a 
 \begin{align*}
 \forall j \in \llbracket 1,d \rrbracket, \quad \|Z\| \geq |z_j|.
 \end{align*}
 En effet, pour $i \in \llbracket 1,d \rrbracket$ la $i\tir$ème coordonnée de $Z_N^j $ est égale à $1$ si $i =j$ et à $0$ sinon. En particulier on a $\|Z \| \geq |z_{j_0}| \geq \theta^{\floor{\alpha ^N}} $ car $z_{j_0} \neq 0$. On combine cela avec (\ref{8derniere_contradiction}) et on trouve : 
 \begin{align*}
 c_{\ref{8cons_derniere_contradiction}} &\leq \theta^{-\floor{\alpha ^N}(\alpha + e\varepsilon - 2) + \floor{\alpha ^{N+1}} - 2 \floor{\alpha^N}}
 \end{align*}
 avec, on le rappelle, $c_{\ref{8cons_derniere_contradiction}} >0 $ une constante indépendante de $N$. \\
 Or on a $$-\floor{\alpha ^N}(\alpha + e\varepsilon - 2) + \floor{\alpha ^{N+1}} - 2 \floor{\alpha^N} = \floor{\alpha ^{N+1}} -\floor{\alpha ^N}(\alpha + e\varepsilon) \underset{N \to + \infty}{\longrightarrow} - \infty.$$ 
 Si $H(C) $ est grand alors $N $ l'est aussi par $(\ref{8choixN})$. 
Donc pour $H(C)$ tendant vers $+ \infty$, on trouve $c_{\ref{8cons_derniere_contradiction}} = 0$ ce qui est contradictoire et termine la preuve du lemme.

\end{preuve}

Ce lemme~\ref{8lem_minoration_exposant_e<d} donne alors pour $e \in \llbracket 1,d-1 \rrbracket $ :
\begin{align*}
 \mu_n(A|C)_e \leq \frac{\alpha}{e}.
\end{align*}
On a donc prouvé le théorème~\ref{8theo_dernier_angle} dans le cas $ e< d$.